\documentclass[letterpaper]{article}
\usepackage{authblk}
\usepackage{url}
\usepackage[margin=1in,footskip=0.25in]{geometry}
\usepackage{amssymb}
\usepackage{graphicx}
\usepackage{amsmath}
\usepackage{algorithmic} 
\usepackage{algorithm}
\usepackage{float}
\usepackage{lipsum}
\usepackage{wrapfig}
\usepackage{color}
\usepackage{easylist}
\usepackage{lineno}
\usepackage{stmaryrd}
\usepackage{booktabs}
\usepackage{cleveref}
\usepackage{multirow}
\usepackage{hhline}
\usepackage[bbgreekl]{mathbbol}
\usepackage{subcaption}

\newcommand{\ba}{\begin{array}}
\newcommand{\ea}{\end{array}}
\newcommand{\be}{\begin{equation}}
\newcommand{\ee}{\end{equation}}

\newcommand{\Div}[1]{\nabla\cdot #1}
\newcommand{\DIV}[1]{\nabla_{\mathbf{X}} \cdot #1}

\newcommand{\etal}{\emph{et al}.}

\newcommand{\BB}{\mathbb{B}}
\newcommand{\CC}{\mathbb{C}}
\newcommand{\FF}{\mathbb{F}}

\newcommand{\PP}{\mathbb{P}}
\newcommand{\MM}{\mathbb{M}}
\newcommand{\FFb}{\overline{\FF}}

\newcommand{\ub}{\boldsymbol{u}}
\newcommand{\ubf}{\boldsymbol{u}^{\text{f}}}
\newcommand{\ubs}{\boldsymbol{u}^{\text{s}}}

\newcommand{\cb}{\boldsymbol{\chi}}
\newcommand{\xb}{\boldsymbol{x}}
\newcommand{\Xb}{\boldsymbol{X}}
\newcommand{\fb}{\boldsymbol{f}}
\newcommand{\Fb}{\boldsymbol{F}}
\newcommand{\Ub}{\boldsymbol{U}}
\newcommand{\Nb}{\boldsymbol{N}}
\newcommand{\nb}{\boldsymbol{n}}
\newcommand{\hb}{\boldsymbol{h}}
\newcommand{\ab}{\boldsymbol{a}}
\newcommand{\Ab}{\boldsymbol{A}}

\newcommand{\cauchy}{\bbsigma}
\newcommand{\cauchyf}{\bbsigma^{\text{f}}}
\newcommand{\cauchyv}{\bbsigma^{\text{v}}}
\newcommand{\cauchys}{\bbsigma^{\text{e}}}
\newcommand{\cauchysolid}{\bbsigma^{\text{s}}}
\newcommand{\ztensor}{\mathbb{0}}
\newcommand{\pstab}{\pi_{\text{stab}}}
\newcommand{\pf}{\pi^{\text{f}}}
\newcommand{\ps}{\pi^{\text{s}}}
\newcommand{\pphys}{p}
\newcommand{\peul}{\pi}
\newcommand{\mus}{\mu^{\text{s}}}
\newcommand{\muf}{\mu^{\text{f}}}
\newcommand{\rhos}{\rho^{\text{s}}}
\newcommand{\rhof}{\rho^{\text{f}}}

\newcommand{\ds}{\boldsymbol{d}^{\text{s}}}

\newcommand{\PPs}{\PP^{\text{e}}}

\newcommand{\soliddom}{\Omega^{\text{s}}_t}
\newcommand{\soliddomO}{\Omega^{\text{s}}_0}
\newcommand{\fluiddom}{\Omega^{\text{f}}_t}
\newcommand{\fluiddomO}{\Omega^{\text{f}}_0}
\newcommand{\fsinterface}{\Gamma^{\text{fs}}_t}

\newcommand{\subsolid}{\omega^{\text{s}}_t}
\newcommand{\subsolidO}{\omega^{\text{s}}_0}

\newcommand{\Wb}{\overline{W}}

\newcommand{\dev}{\text{dev}}
\newcommand{\DEV}{\text{DEV}}
\newcommand{\mfac}{M_\text{FAC}}
\newcommand{\tr}{\text{tr}}
\newcommand{\nus}{\nu_\text{stab}}
\newcommand{\kappas}{\kappa_\text{stab}}

\newcommand{\Gt}{G_\text{T}}
\newcommand{\Gl}{G_\text{L}}
\newcommand{\El}{E_\text{L}}

\title
{Stabilization approaches for the hyperelastic immersed boundary method for problems of large-deformation incompressible elasticity}

\author[1]{Ben Vadala-Roth} %\corref{cor1}}
\author[2]{Shashank Acharya}
\author[2]{Neelesh A. Patankar}
\author[1,4]{Simone Rossi} 
\author[3,4,5,6]{Boyce E.~Griffith} %\corref{cor2}}

\affil[1]{Department of Mathematics, University of North Carolina, Chapel Hill, NC, USA}
\affil[2]{Department of Mechanical Engineering, Northwestern University, Evanston, IL, USA}
\affil[3]{Departments of Mathematics, Applied Physical Sciences, and Biomedical Engineering, University of North Carolina, Chapel Hill, NC, USA}
\affil[4]{Carolina Center for Interdisciplinary Applied Mathematics, University of North Carolina, Chapel Hill, NC, USA}
\affil[5]{Computational Medicine Program, University of North Carolina, Chapel Hill, NC, USA}
\affil[6]{McAllister Heart Institute, University of North Carolina, Chapel Hill, NC, USA}
\affil[ ]{\texttt{bvadalaroth@gmail.com} and \texttt{boyceg@email.unc.edu}}
\begin{document}

\maketitle

\begin{abstract}
The immersed boundary method is a mathematical framework for modeling fluid-structure interaction. This formulation describes the momentum, viscosity, and incompressibility of the fluid-structure system in Eulerian form, and it uses Lagrangian coordinates to describe the structural deformations, stresses, and resultant forces. Integral transforms with Dirac delta function kernels connect the Eulerian and Lagrangian frames. The fluid and the structure are both typically treated as incompressible materials. Upon discretization, however, the incompressibility of the structure is only maintained approximately. To obtain an immersed method for incompressible hyperelastic structures that is robust under large structural deformations, we introduce a volumetric energy in the solid region that stabilizes the formulation and improves the accuracy of the numerical scheme. This formulation augments the discrete Lagrange multiplier for the incompressibility constraint, thereby improving the original method's accuracy. This volumetric energy is incorporated by decomposing the strain energy into isochoric and dilatational components, as in standard solid mechanics formulations of nearly incompressible elasticity. We study the performance of the stabilized method using several quasi-static solid mechanics benchmarks, a dynamic fluid-structure interaction benchmark, and a detailed three-dimensional model of esophageal transport. The accuracy achieved by the stabilized immersed formulation is comparable to that of a stabilized finite element method for incompressible elasticity using similar numbers of structural degrees of freedom.

\end{abstract}


%%%%%%%%%%%%%%%%%%%%%%%%%%%%%%%%%%%%%%%%%%%%%%%%%%%%%%%%%%%%%%%%%%%%%%%%%%%%%%%%%%%%%%%%%%%%%%%%%%%%%%%%%%%%%%%%%%%%%%%%%%%
\section{Introduction and overview}
\label{Intro}
\indent The immersed boundary (IB) method \cite{Peskin2002, Griffith2020} is a framework for modeling fluid-structure interaction (FSI) that was introduced by Peskin
to simulate blood flow through heart valves \cite{Peskin1972, Peskin1977}. This formulation describes the momentum, viscosity, and incompressibility of the fluid-solid system in Eulerian form, whereas Lagrangian variables describe the deformations, stresses, and resultant forces of the immersed structure.
Integral transforms with Dirac delta function kernels link the Eulerian and Lagrangian frames. When the IB equations are discretized, regularized delta functions are commonly used to construct coupling operators that connect the Eulerian grid and the Lagrangian mesh. This approach maintains a continuous velocity field across the fluid-solid interface while avoiding the need for body-fitted descriptions of the fluid and structure. The method was originally formulated to describe thin structures occupying zero volume within the fluid \cite{Peskin1972, Peskin1977}, and it was later extended to describe structures with finite volume \cite{Peskin2002}. In the work of Boffi \etal~\cite{Boffi2008}, the IB equations are systematically derived in the framework of large-deformation continuum mechanics. This formulation is particularly useful for applications that use nonlinear mechanics-based descriptions \cite{Holzapfel2009} of the immersed structure. Such models are commonly used, for instance, in biomedical applications, in which many physiologically realistic and experimentally validated models of soft tissue rely upon a continuum mechanics description. For this reason, in this work, we adopt the continuum formulation posed by Boffi \etal~\cite{Boffi2008}.\\
\indent A key feature of the immersed formulation for FSI introduced by Peskin, which is also used in the formulation considered herein, is that a common momentum equation is used for both the fluid and solid. This means that the incompressibility constraint is imposed in Eulerian form throughout the entire computational domain, automatically ensuring incompressibility in the solid region in the continuous IB formulation. This constraint is maintained via a Lagrange multiplier defined in the Eulerian frame. IB formulations have also been developed that treat compressible structures immersed in incompressible fluids \cite{Wang2009, Heltai2012}. As noted in various works on immersed methods for FSI \cite{Wang2009,Peskin1993,Griffith2012}, even if the solid is modeled as incompressible in the continuum equations, typical discretizations do not automatically satisfy this constraint. \\
\indent There are three main sources of error that lead to volume losses within the framework of immersed methods. One source of error is introduced by the coupling operators that link the Eulerian and Lagrangian descriptions of the structural velocity. Discrete coupling operators commonly used in practice generally do not preserve the divergence of the Eulerian velocity field when transferring it to the Lagrangian mesh, although specialized approaches have been developed \cite{Bao2017}. Another source of error is caused by the choice of approximation space used to represent the structural velocities. In commonly used IB methods, the Lagrangian velocity field is typically not represented in a way that allows for pointwise divergence-free discrete velocity fields (although see the work of Casquero \etal~\cite{Casquero2018}). Finally, time integration errors are another source of spurious changes in volume. Specifically, even if the discretized structural velocity is continuously divergence free, as in the IB formulation of Bao \etal~\cite{Bao2017}, time discretization errors can induce changes in the volume of the discretized structure. The approach introduced in this work addresses effects of volume errors originating from all these sources.\\
\indent The approach taken in this work to improve volume conservation is to reinforce the Eulerian incompressibility constraint ($\nabla \cdot \ub = 0$) by additional structural stresses that act to ensure that the Lagrangian incompressibility constraint ($J = 1$) is satisfied by the computed structural deformation. Specifically, the present work demonstrates that a simple stabilization method, active only in the solid region, can greatly improve the volume conservation properties of the hyperelastic immersed boundary method. The proposed method can be used with standard nodal Lagrangian finite element (FE) basis functions. To assess the performance of the proposed method, we employ standard quasi-static benchmark problems traditionally used for large-deformation incompressible elasticity. Our tests differ from standard elasticity benchmarks, however, in that the solid is immersed in a fluid and is treated dynamically. Comparisons to benchmark results are performed after allowing the fluid-structure system to reach equilibrium. This allows us to evaluate the accuracy of the structural component of our FSI calculations and enables direct comparisons to results obtained using FE methods specially designed for large-deformation incompressible elasticity. Additionally, we investigate the impact of the method on an FSI test that is driven by fluid forces, rather than a solid traction, and an FSI application involving a detailed three-dimensional model of esophageal transport \cite{Kou2017, Kou2017b}. \\
\indent The proposed stabilization method is rooted in the deviatoric-spherical decomposition of the Cauchy stress tensor. In nearly incompressible solid mechanics, a volumetric term is included in the stress to resist compressible deformations \cite{JavierBonetandRichardWood2008}. The magnitude of this additional term is related to a physical parameter, the bulk modulus $\kappa$, representing the resistance of the solid to compression. As $\kappa \rightarrow \infty$, the material can only experience incompressible deformations. We introduce a similar penalization method to stabilize the IB formulation. Specifically, our method penalizes changes in volume in the solid region by adding an additional isotropic stress to the structural stresses. For hyperelastic formulations, the pressure-like term can be derived from a volumetric term included in the strain energy as in standard nearly incompressible formulations of solid mechanics. In our IB formulation, the stabilization is controlled by a numerical parameter that we refer to as the \textit{numerical bulk modulus} $\kappas$. We emphasize that $\kappas$ \textit{does not} represent a physical parameter. Instead, $\kappas$ is a \textit{numerical} parameter that, in concert with the Eulerian Lagrange multiplier, reinforces the incompressibility constraint. 
%Therefore, our approach resembles an augmented Lagrangian method, in which the Lagrange multiplier is approximated by summing an initial approximation and a correction derived from a penalty function \cite{Glowinski1982, Simo1991}.
Furthermore, as detailed in Section \ref{Jumps}, this formulation satisfies both kinematic and dynamic conditions at fluid-solid interfaces, whether or not the volumetric stabilization terms are active in the solid stress.\\
\indent Numerical examples show that the proposed stabilization can improve simulations in which the elastic stress tensor is not deviatoric, but we also demonstrate that using a traceless Cauchy stress tensor generally results in further improvements in the volume conservation of the present formulation. We analyze two different strategies to achieve a traceless stress: the first is based on the Flory decomposition of the deformation gradient tensor \cite{Flory1961}; the second eliminates the volumetric contribution of the stress tensor using a deviatoric projection. Whereas the Flory decomposition is mainly used for hyperelastic materials, the deviatoric projection strategy is easily implemented for general elastic and hypoelastic materials. The Flory decomposition for hyperelastic materials is equivalent to a formulation that additively decouples the isochoric (volume-preserving) and dilatational (volume-changing) parts of the strain energy functional. The admissibility of such decompositions requires physical assumptions about the solid being studied: namely that uniform pressure only results in a change in size and does not result in changes in the shape \cite{Sansour2008}.\\
\indent It is well known that with many simple numerical methods for nearly incompressible elasticity, large values of the bulk modulus lead to volumetric locking or sub-optimal convergence rates in the computed displacement \cite{ThomasJRHughes2000}. In the proposed method, the penalization does not require $\kappas \rightarrow \infty$ because the penalty method acts alongside a Lagrange multiplier, and we demonstrate that locking can be avoided even with simple linear finite elements. Throughout this work, we determine the numerical bulk modulus from a parameter that we refer to as the \textit{numerical Poisson ratio} $\nus$ along with standard linearized elasticity relations. We therefore limit the values of $\nus$ to lie in the interval $\left[-1,\frac{1}{2}\right)$. In this work, $\nus = -1$ corresponds to no volumetric stabilization ($\kappas = 0$). We show that setting $\nus = -1$ ($\kappas = 0$) can result in large changes in the solid volume along with nonphysical solid deformations. We explore the effects of $\nus$ as a volumetric stabilization parameter used to reinforce the incompressibility of the structure. To achieve this goal, we empirically demonstrate that $\nus = 0.4$ is sufficient for several benchmark problems and across a range of grid spacings. We emphasize that $\nus$ is used to augment the Eulerian Lagrange multiplier field that imposes incompressibility within the structure, not to model a compressible material. Like $\kappas$, the numerical Poisson ratio $\nus$ is a numerical parameter and not a physical parameter of the model.\\ 
\indent The regularized delta function kernels used in the IB method to couple the Eulerian and Lagrangian frames can be interpreted as weighting functions. In this sense, the IB method resembles mesh free and particle based methods, including the element-free Galerkin method (EFG) \cite{Belytschko1994} and reproducing kernel particle methods (RKPM) \cite{Liu1995}, which use weighting functions to reconstruct continuum fields. In fact, neither the regularized delta functions typically used with the IB method nor the weighting functions typically used with the EFG and RKPM are interpolatory. Particle based methods may also be paired with traditional FE methods, yielding a combined method with the flexibility and smoothness of particle based methods but with the ability to more easily impose Dirichlet boundary conditions \cite{Wagner2001, Han2002, Zhang2002, Han2005}. As shown by Zhang \etal~\cite{Zhang2002}, they can be made computationally efficient and amenable to parallelization. In the EFG method, in particular, it is known that if the support of the weighting functions is small, volumetric locking may result. A simple way to alleviate volumetric locking for hexahedral elements, known as selective reduced integration (SRI) \cite{Malkus1978}, is to integrate the volumetric term associated with $\kappa$ using a quadrature rule with reduced order of accuracy. A procedure equivalent to SRI has been shown to correct this issue \cite{Dolbow1999}. We demonstrate herein that volumetric locking also occurs in our implementation of the IB method if we let $\nus \rightarrow \frac{1}{2}$. As already mentioned, however, the proposed method can circumvent issues with locking and obtain accurate solutions by using values of $\nus$ much smaller than $\frac{1}{2}$, even with low order elements. This is permissible because, as stated previously, the penalty parameter acts alongside the Eulerian Lagrange multiplier field to maintain incompressibility in the discrete case.\\
\indent The numerical method used here follows the one described by Griffith and Luo \cite{Griffith2017}. This method uses a finite difference scheme to approximate the Eulerian equations and an FE scheme for the Lagrangian equations, and it uses regularized delta functions in approximations to the integral transforms. In this method, an intermediate Lagrangian velocity field is projected onto FE basis functions to determine the velocity of the immersed structure. This FE and finite difference based method stands in contrast to an alternative numerical method, the immersed finite element method (IFEM) \cite{Zhang2004}, that uses the FE method to approximate both the Eulerian and Lagrangian equations. The IFEM can be directly applied to unstructured Eulerian grids, which is an advantage when working with complex computational domain geometries. This method achieves this through constructing delta functions using RKPM \cite{Liu1995, Liu2007}. For structured grids, a single kernel function may be used for the entire domain, which is also the case in the method considered here \cite{Griffith2017, Wang2004}. \\
\indent Another hybrid IB method combining finite difference and FE methods was proposed by Devendran and Peskin \cite{Devendran2012}. This method has only been formulated using a representation of the structure based on linear simplicial elements. Although different continuum material models may be selected, their implementation requires analytically calculating derivatives of the strain energy functional with respect to the coefficients of the FE representation of the displacement field. Other implementations include an extension using radial basis functions to represent the structure \cite{Shankar2015}, a particle based method to represent the structure \cite{Gil2010}, and the numerical methods of Boffi \etal~\cite{Boffi2008} and Roy \etal~\cite{Roy2015} that avoid the use of regularized delta functions, achieving regularization instead via the FE basis functions. A systematic study of the loss of incompressibility of the immersed solid was performed by Casquero \etal~\cite{Casquero2018} using an IB-type method with divergence-conforming B-splines. Through the use of these basis functions, their method achieves negligible changes of volume in the Eulerian frame and reduced incompressibility errors in the Lagrangian frame.\\
\indent The first two cases we employ in our numerical tests, the compressed block \cite{Reese1999} and Cook's membrane \cite{RDCook1974}, are two-dimensional problems that invoke the plane-strain assumption. Two additional quasi-static tests, an anisotropic extension to Cook's membrane \cite{Wriggers2016} and a torsion test \cite{Bonet2015}, are fully three-dimensional. We also consider a fully dynamic ``elastic band'' test, which is a two dimensional FSI test driven by fluid traction boundary conditions. As a more complex example, we also consider the impact of volumetric stabilization on a fully three-dimensional model of esophageal transport based upon the model of Kou \etal~\cite{Kou2017}. These tests all demonstrate that in the absence of stabilization, corresponding to $\kappas = 0$ (equivalently $\nus = -1$ for the cases considered here), leads to unphysical and inaccurate deformations and large volume conservation errors. Further, in the tests detailed herein, it is generally the case that using the Flory decomposition with a finite choice of $\kappas$ provides the best accuracy. An important result of these tests is that they clearly demonstrate that the accuracy of the structural response provided by the present stabilized IB formulation is comparable to that yielded by specialized methods for incompressible nonlinear elasticity.

%%%%%%%%%%%%%%%%%%%%%%
%%%%%%%%%%%%%%%%%%%%%%

\section{Continuous Formulation}
\subsection{Continuous Equations of Motion}
\label{Continuous Equations}
Let $\Omega = \fluiddom \cup \soliddom$ be the computational domain, in which $\fluiddom$ and $\soliddom$ are respectively the regions occupied by the fluid and the structure at time $t$. We describe the computational domain using Eulerian coordinates $\xb \in \Omega$. We describe the reference configuration of the structure using Lagrangian reference coordinates $\Xb \in \soliddomO$, in which $\soliddomO$ is the physical region occupied by the solid at time $t = 0$.  The mapping $\cb (\Xb,t): \soliddomO \mapsto \soliddom$ connects the reference configuration of the structure to its current configuration. The IB formulation for an immersed elastic structure employed in this work defines the Cauchy stress on the computational domain to be
\begin{equation}
  \cauchy (\xb,t) = \cauchyf(\xb,t) +  \begin{cases} \ztensor & \xb \in \fluiddom, \\ \cauchys(\xb,t) & \xb \in \soliddom. \end{cases} \label{cauchy-def}
\end{equation} 
Because we use a Lagrangian description of the structure, we use the first Piola-Kirchhoff stress $\PPs$ to describe the elastic response of the structure. Let $\FF = \frac{\partial \cb}{\partial \Xb} $ be the deformation gradient tensor, and let $J = \det(\FF)$. For the material models considered here, $\PPs$ is determined from a strain energy functional $\Psi(\FF)$ via $\PPs = \frac{\partial \Psi}{\partial \FF}$. The first Piola-Kirchhoff stress is related to the corresponding Cauchy stress by $\cauchys = \frac{1}{J}\PPs \FF^T$. We consider a Newtonian fluid with Cauchy stress given by $\cauchyf = - \peul\mathbb{I} + \mu\left(\nabla \ub + \nabla \ub ^T \right) $, in which $\peul$ is the Lagrange multiplier field and $\ub$ is the Eulerian velocity. We discuss the relationship between $\peul$ and the pressure $p$ in Section \ref{Volumetric Stabilization}.
The resulting IB form of the equations of motion, as derived by Boffi \etal~\cite{Boffi2008}, is:
\begin{align}
	\rho \frac{D \ub  }{Dt}(\xb, t) =& -\nabla \peul(\xb,t) + \mu \nabla ^2 \ub(\xb,t) + \fb( \xb ,t),  \label{ns} \\
	\nabla \cdot \ub( \xb ,t) =& \, 0, \label{divfree} \\
	\fb (\xb , t) =& \int_{\soliddomO} \nabla_{\Xb} \cdot \PPs (\Xb,t) \, \delta (\xb - \cb (\Xb,t)) \, d\Xb \label{elasticforce} \\
	& \mbox{} - \int_{\partial \soliddomO} \PPs (\Xb,t) \Nb(\Xb) \, \delta (\xb - \cb (\Xb,t)) \, d\Ab, \nonumber \\
	\frac{\partial \cb}{\partial t}(\Xb , t) =& \ \Ub(\Xb,t) = \int_{\Omega} \ub(\xb, t) \, \delta (\xb - \cb (\Xb,t)) \, d\xb. \label{noslip}
\end{align}
Here, $\ub(\xb,t)$ is the Eulerian velocity, $\Ub(\Xb,t)$ is the velocity of the structure, $\rho$ is the constant mass density, $\mu$ is the viscosity, $\fb(\xb,t)$ is the Eulerian form of the elastic body force of the immersed solid, and $\Nb(\Xb)$ is the outward unit normal along the solid boundary $\partial \soliddomO$ in the reference configuration. The Lagrange multiplier $\peul(\xb,t)$ is responsible for maintaining the incompressibility constraint, equation \eqref{divfree}, and within the fluid domain, $\peul$ is the physical pressure. The operators $\nabla ^2, \nabla \cdot \mbox{},$ and $\nabla$ are with respect to spatial coordinates, and $\frac{D}{Dt} = \frac{\partial}{\partial t} + \ub \cdot \nabla $ is the material time derivative. The differential operator $\nabla_{\Xb} \cdot \mbox{}$ is the Lagrangian divergence operator. Equations $(\ref{ns})$ and $(\ref{divfree})$ are the Navier-Stokes equations for an incompressible Newtonian fluid, augmented by elastic forces in the solid region. In the IB formulation, interactions between the Eulerian and Lagrangian variables occur via integral transforms with Dirac delta function kernels, equations ($\ref{elasticforce}$) and ($\ref{noslip}$). These relationships exploit the defining feature of the Dirac delta function as a linear functional. 
%To describe equation ($\ref{elasticforce}$), we say the force is \emph{prolonged}, or spread, from the solid region to the computational domain.For equation ($\ref{noslip}$), we say the solid \emph{interpolates} the Eulerian velocity. 
%Because $\ub$ is continuous space, which is implied by the $\mu\nabla^2\ub$ term, equation ($\ref{noslip} $) ensures that the fluid and solid move together along the fluid-solid interface. 
Equation \eqref{noslip} ensures that there is no slip along the solid boundary, and in partitioned methods for FSI, this no-slip condition is often referred to as the kinematic boundary condition.\\
\indent In this immersed formulation, the solid motion is exactly incompressible, which follows from equations ($\ref{divfree}$) and ($\ref{noslip}$). To demonstrate this, let $\subsolidO \subseteq \soliddomO$ be a subregion of the solid domain in the reference configuration, and let $\subsolid = \cb(\subsolidO,t)$ be the current configuration of this subregion at time $t$. The volume of this subregion in the current configuration is $V(t) = \int_{\subsolid } d\xb$, and its volume in the reference configuration is $V_0 = \int_{\subsolidO} d\Xb$. Because $\frac{\partial \cb}{\partial t} = \Ub$ and $\Div \ub = 0$, the Reynolds transport theorem %\cite{TedBelytschkoWingKamLiu2000}
implies:
\begin{equation}
  \frac{d}{dt} \int_{\subsolid} d\xb = \int_{\omega^s_t} \nabla \cdot \ub(\xb,t) \, d\xb = 0 \label{reynolds}, \\
\end{equation}
i.e. the volume of material region $\subsolidO$ does not change in time, and $V(t) = V_0$.
\noindent Because $\subsolidO \subseteq \soliddom$ is arbitrary, it also holds pointwise. It is also useful to recall that $V(t) = \int_{\subsolid } d\xb = \int_{\subsolidO} J(\Xb,t) \, d\Xb$. For $V(t) = V_0$ to hold for an arbitrary region, we must have $J(\Xb,t) \equiv 1$. Deviations from $J=1$, which can occur in discretizations of these equations, indicate local changes in volume. In the continuum equations, however, the present formulation consistently treats both the fluid and solid as exactly incompressible.\\
% \indent In practice, we solve a weak form of equation ($\ref{elasticforce}$) that is amenable to discretization via standard nodal FE methods. Rather than prolonging the divergence of the stress as in equation (\ref{elasticforce}), we instead prolong a force $\Fb$ that is weakly equivalent to $\DIV \PPs$. To determine $\Fb$, we introduce arbitrarily smooth test functions $\boldsymbol{V}(\Xb)$ and require
\indent The internal forces exerted by the structure $\soliddomO$ are defined by equation ($\ref{elasticforce}$). In practice, forces are evaluated using $\fb(\xb,t) = \int_{\soliddomO}\Fb(\Xb,t)\delta (\xb - \cb(\Xb,t)) \, d\Xb$, in which $\Fb(\Xb,t)$ is determined by a weak version of equation ($\ref{elasticforce}$). More precisely, we determine $\Fb(\Xb,t)$ by requiring:
\begin{equation}
      \int_{\soliddomO} \Fb(\Xb,t) \cdot \boldsymbol{V}(\Xb) \, d\Xb = \int_{\soliddomO} \left(\DIV \PPs (\Xb,t) \right) \cdot \boldsymbol{V}(\Xb) \, d\Xb - \int_{\partial \soliddomO} (\PPs(\Xb,t) \Nb(\Xb)) \cdot 
      \boldsymbol{V}(\Xb) \, d\Ab
\end{equation}
\noindent for all smooth $\boldsymbol{V}(\Xb)$. From the divergence theorem, we obtain
 \begin{equation}
	\int_{\soliddomO} \Fb(\Xb,t) \cdot \boldsymbol{V}(\Xb) \, d\Xb = -\int_{\soliddomO} \PPs(\Xb,t) : \nabla_{\Xb} \boldsymbol{V}(\Xb) \, d\Xb. \label{weak-force}
\end{equation}
\noindent This definition of the elastic forcing is amenable to discretization via standard nodal FE methods. In our numerical method, the force $\Fb(\Xb,t)$ is computed in a weak sense before being spread to the background grid, whereas $\fb(\xb,t)$ is used in the strong form of the equations of motion in the Eulerian frame (\ref{ns}). This method of calculating the internal force with respect to undeformed coordinates is effectively a total Lagrangian approach, as commonly used in nonlinear FE methods \cite{TedBelytschkoWingKamLiu2000}. The definition of the elastic body force $\fb(\xb,t)$ is called the \emph{unified weak formulation} by Griffith and Luo \cite{Griffith2017}.  \\

%%%%%%%%%%%%%%%%%%%%%%%%%%%
\subsection{FSI Sub-Problems and Boundary Conditions}
\label{Jumps}
A partitioned formulation of incompressible FSI, in Eulerian form, can be written alternatively as:
\begin{align}
\text{Fluid sub-problem}: \quad &\begin{cases}
	\rhof \frac{D \ubf  }{Dt}(\xb, t) &= -\nabla \pf(\xb,t) + \muf \nabla ^2 \ubf(\xb,t), \\
	\Div \ubf(\xb,t) &= 0, \, \\
\end{cases} \quad \xb \in \fluiddom \label{fsubproblem},\\
\text{Solid sub-problem}: \quad &\begin{cases}
	\rhos \frac{D \ubs  }{Dt}(\xb, t)  &= -\nabla \ps(\xb,t) + \mus \nabla ^2 \ubs(\xb,t) + \Div\cauchys(\xb,t), \\
	\Div \ubs(\xb,t) &= 0,\\
	\frac{D \ds  }{D t}(\xb, t) &= \ubs(\xb, t),
\end{cases} \quad \xb \in \soliddom, \label{ssubproblem}
\end{align}
in which $\ds(\xb,t)$ is the Eulerian displacement field of the solid. Superscripts of `f' and `s' indicate that a variable describes a fluid or solid quantity, respectively. Specifically, we are concerned with the case of $\muf = \mus$ and $\rhof = \rhos$. In the partitioned context, the quantities $\ub, \cauchy,$ and $\peul$ denote variables that have been extended to the entire domain $\Omega$, and when restricted to each sub-domain take on their respective fluid and solid counterparts.  As in the IB formulation, $\Omega = \fluiddom \cup \soliddom$. To complete the specification of the problem, we require continuity of velocity and traction at the interface $\fsinterface = \overline{\fluiddom} \cap\overline{\soliddom}$ and suitable initial conditions and boundary conditions along the remaining boundaries.\\
%As before, the domain $\Omega$ is decomposed into sub-domains $\fluiddom$ and $\soliddom$ such that  
%To relate this formulation to the immersed formulation considered herein, we would define $\cauchysolid$ via $\cauchysolid = \cauchyf + \cauchys$, that is we assume the solid to be a viscoelasic material. Note that $\cauchy$ restricted to $\fluiddom$ is $\cauchyf$ and restricted to $\soliddom$ is $\cauchysolid$.\\
\indent The continuity of velocity and traction conditions serve to couple the fluid and structure and are often referred to as boundary conditions for each sub-domain. As mentioned previously, the imposition of continuity of velocity $\ub$ across the fluid-structure interface is known as the kinematic boundary condition. Continuity of traction across the interface is known as the dynamic boundary condition. These conditions are equivalent to requiring that the jumps in these quantities are zero along the fluid-solid interface $\fsinterface$. Together, these two boundary conditions are:
\begin{align}
	\llbracket \ub(\xb,t) \rrbracket &= \boldsymbol{0}, \, \xb \in  \fsinterface \ \text{and} \label{kinematic-bc}  \\
	\llbracket \cauchy(\xb,t) \nb(\xb,t) \rrbracket &= \boldsymbol{0}, \, \xb \in  \fsinterface. \label{dynamic-bc}
\end{align}
The jump operator $\llbracket \cdot \rrbracket$ is defined for a function $g(\xb)$ on $\Omega$ as
\begin{equation}
	\llbracket g(\xb) \rrbracket = \lim_{\epsilon \downarrow 0} g(\xb + \epsilon \nb) - \lim_{\epsilon \downarrow 0} g(\xb - \epsilon \nb) = g_+(\xb) - g_-(\xb),
\end{equation}
in which $\nb$ is the outward normal at $\xb \in \fsinterface$ pointing into the fluid region $\fluiddom$. We use the notation $g_+$ and $g_-$ to denote the limits approaching $\fsinterface$ from the fluid region and from the solid region, respectively, and apply this notation to other quantities such as $\cauchyf$ and $\peul$.\\
\indent Sub-problems (\ref{fsubproblem}) -- (\ref{ssubproblem}) with interface conditions (\ref{kinematic-bc}) -- (\ref{dynamic-bc}) are equivalent to the IB equations (\ref{ns}) -- (\ref{noslip}). However, the two formulations are amenable to different discretization strategies. Previous work by Peskin and Printz \cite{Peskin1993} and Lai and Li \cite{Lai2001} establish the jump conditions within IB equations that imply the formulation satisfies conditions (\ref{kinematic-bc}) -- (\ref{dynamic-bc}). To briefly summarize, first note that in the immersed formulation, equation (\ref{kinematic-bc}) is maintained by definition of the Dirac delta function. Equation (\ref{dynamic-bc}) may be decomposed as $\llbracket \cauchy \nb \rrbracket = \llbracket \cauchyf\nb \rrbracket + \llbracket \cauchys \nb \rrbracket$ within the IB context. Noting that $\cauchys_+ = \ztensor$ by equation (\ref{cauchy-def}), we have the simplification that $\llbracket \cauchyf\nb \rrbracket = \cauchys_-\nb$ must hold for continuity of traction. The aforementioned work \cite{Peskin1993, Lai2001} demonstrate that the IB formulation yields jump conditions that, in the present context, imply the discontinuity in the fluid stress along $\fsinterface$ is precisely $\cauchys_-\nb$. In both formulations, the discontinuity in the Lagrange multiplier field $\peul$ across the interface $\fsinterface$ is specifically 
\begin{equation}
	\llbracket \peul \rrbracket = \nb^T \cauchys_- \nb, \label{pjump}
\end{equation}
and the viscous stress jump is 
\begin{equation}
	\mu \left\llbracket \mathbf{t}^T\left(\nabla\ub + \nabla\ub^T\right)\nb \right\rrbracket = \mathbf{t}^T\cauchys_- \nb, \label{njump}
\end{equation}
for any unit vector $\mathbf{t}$ tangent to $\fsinterface$. Using Nanson's relation, we may write the limiting value of the elastic traction in terms of $\PPs$
\begin{equation}
	\cauchys_- \nb = \frac{\PPs \Nb}{||\FF \mathbf{T}_1 \times \FF\mathbf{T}_2 ||},
\end{equation}
in which $\mathbf{T}_1$ and $\mathbf{T}_2$ are mutually orthogonal vectors that are tangent to $\fsinterface$ \cite{Gao2014}. In summary, the IB formulation given by equations (\ref{ns}) -- (\ref{noslip}) satisfies interface conditions (\ref{kinematic-bc}) -- (\ref{dynamic-bc}), and this is achieved through balancing the structural traction force with discontinuities in the viscous stress and Lagrange multiplier $\peul$.\\
\indent We remark that at steady-state, $\ub \equiv \boldsymbol{0}$, which implies the left hand side of each sub-problem in the partitioned formulation is zero. The stress in the solid region has the form $-\ps\mathbb{I} + \mus\left(\nabla \ubs + \left(\nabla \ubs\right)^T \right) + \cauchys$, so the equilibrium configuration defined zero velocity implies $- \nabla \ps + \Div \cauchys = \boldsymbol{0}$. This yields the same steady-state configuration as that of an uncoupled hyperelastic solid mechanics problem. This fact, which also holds true for the IB equations, forms the basis of our comparisons in Section \ref{benchmarks} to benchmark results computed using FE numerical methods for solid mechanics problems.
%%%%%%%%%%%%%%%%%%%%%%%%%%%

\subsection{Volumetric Stabilization}
\label{Volumetric Stabilization}
It is well known that for an artibrary second order tensor $\mathbb{T}$, there is a unique decomposition into deviatoric and isotropic parts, such that $\mathbb{T} = \dev[\mathbb{T}] - \varphi \mathbb{I}$. Here $\varphi = -\frac{\tr \left(\mathbb{T}\right)}{3}$ and 
\begin{equation}
	\dev[\boldsymbol{\mathbb{T}}] = \mathbb{T} - \frac{\tr \left(\mathbb{T} \right)}{3}\mathbb{I}. \label{dev}
\end{equation}
By construction, the deviatoric part will satisfy the property $\tr(\dev[\mathbb{T}]) = 0$. In continuum mechanics, the Cauchy stress may be similarly decomposed as
 \begin{equation}
	\cauchy = \dev[\cauchy] - \pphys\mathbb{I}, \label{cauchy-dev}
\end{equation}
in which $\pphys$ is the physical pressure and $\varphi = \pphys$. That is, the physical pressure $p$ is,
\begin{equation}
	\pphys = -\frac{\tr \left(\mathbb{\cauchy} \right)}{3},
\end{equation}
which means that $\pphys = \peul$ in the fluid region and $\pphys = \peul - \frac{\tr \left(\mathbb{\cauchys} \right)}{3}$ in the solid region in the immersed formulation of FSI. For incompressible motions, the pressure is defined by the incompressibility constraint. For compressible motions, the pressure encodes the volume change of the material. As described in Section \ref{Continuous Equations}, common IB formulations \cite{Griffith2017, Zhang2004} decompose the Cauchy stress via
\begin{equation}
	\cauchy = \cauchyv -\peul\mathbb{I} + \begin{cases} \ztensor & \xb \in \fluiddom, \\ \cauchys & \xb \in \soliddom, \end{cases} \label{cauchy-nodev}
\end{equation}
in which $\cauchyv = \mu\left(\nabla \ub + \nabla \ub ^T \right)$ is viscous stress and the elastic stress $\cauchys$ is not necessarily deviatoric. Note that $\cauchyv$ is already deviatoric because of equation (\ref{divfree}). As we will show in our simulations, using (\ref{cauchy-nodev}) in the discretized equations may lead to very poor numerical results, such as unphysical and sometimes extreme contractions of the immersed structure.\\
\indent As described in Section \ref{Intro}, loss of incompressibility in the solid region can be affected by the FSI coupling operators, choice of approximation space for the Lagrangian variables, and time integration errors. We wish to introduce a stabilization to the structural stress that corrects for the loss of incompressibility resulting from any of these sources of volume errors. We first take the deviatoric component of $\cauchys$ and then introduce a volumetric stabilization, such that the stress is
\begin{equation}
	\cauchy = \cauchyv - \peul\mathbb{I} + \begin{cases} \ztensor & \xb \in \fluiddom, \\ \dev[\cauchys] - \pstab\mathbb{I} & \xb \in \soliddom, \end{cases} \label{dev-stab}
\end{equation}
in which $\pstab\mathbb{I}$ is a stabilization term that acts like an additional pressure in the solid region. This will introduce an additional force in the structural region that, when spread to the background grid, will act to combat spurious compressible motions in the solid region. Considering equation (\ref{dev-stab}), we can see that the jump in the Eulerian Lagrange multiplier $\peul$ will now be $\llbracket \peul \rrbracket = \nb^T\dev\left[\cauchys_-\right]\nb - \pstab$, and the pressure in the solid region is $\pphys = \peul + \pstab$ for this model. Such volumetric stabilization can also be included if the deviatoric solid stress is not considered specifically,
\begin{equation}
	\cauchy = \cauchyv - \peul\mathbb{I} + \begin{cases} \ztensor & \xb \in \fluiddom, \\ \cauchys - \pstab\mathbb{I} & \xb \in \soliddom. \end{cases}\label{nodev-stab}
\end{equation}
\noindent Here, as well, we have a change in the jump for $\peul$. With formulation (\ref{nodev-stab}), the pressure jump is now $\llbracket \peul \rrbracket = \nb^T\cauchys_-\nb - \pstab$. The pressure in the solid region is $\pphys = \peul + \pstab - \frac{\tr \left(\mathbb{\cauchys} \right)}{3}$.\\
\indent Similar to treatments of nearly incompressible elasticity, we define $\pstab$ as a volumetric penalization term. More specifically, $\pstab$ is derived from a volumetric energy $U(J)$ that depends only on changes in the volume of the structure,
\begin{equation}
	\pstab = -\frac{\partial U(J)}{\partial J}.
\end{equation}
Thus, this formulation for stabilization parallels models of nearly incompressible elasticity \cite{JavierBonetandRichardWood2008}, except we have adopted the convention commonly used in fluid mechanics, in which the pressure has a negative sign in front of it.\\ 
%Furthermore, in the solid region $\pstab$ and $\peul$ effectively act together as an augmented Lagrangian approach \cite{Glowinski1982, Simo1991}. \\
\indent In nearly incompressible elasticity, restrictions are placed on $U(J)$ to achieve certain physically motivated properties. Specifically, a definition of $U(J)$ in which $U(1) = 0$ is physically meaningful because no additional energy is introduced if $\FF = \mathbb{I}$. Because the contribution of the volumetric energy to the first Piola-Kirchhoff stress is $ J U'(J) \FF^{-T}$, it is also common to require $U(J)$ to satisfy $U'(1) = 0$, which implies that no extra stress is introduced if $J = 1$. We also require that $\lim_{J \rightarrow \infty} U(J) = \infty$ and $\lim_{J \rightarrow 0} U(J) = \infty$, so that large dilatations and contractions are energetically unfavorable. Finally, we want to control the effect of $U(J)$ through a stabilization parameter, referred to as the numerical bulk modulus $\kappas$. A simple example of $U(J)$ that satisfies the above conditions is
\begin{equation}
	U(J) = \frac{\kappas}{2}(\ln J)^2. \label{vol-energy}
\end{equation}
\noindent To modulate the $\kappas$, we introduce the numerical Poisson ratio $\nus$. The two parameters are related via
\begin{equation}
\kappas = \frac{2G(1 + \nus)}{3(1 - 2\nus)} \label{kappa-nu},
\end{equation}
in which $G$ is the shear modulus. 
This relationship provides a mechanism for determining $\kappas$ that mimics the relationship between the physical Poisson ratio $\nu$ and the physical bulk modulus $\kappa$ in a compressible material model. Note that $\nus = -1$ yields $\pstab= 0$, retrieving the case with no stabilization. As discussed in Section \ref{Intro}, the numerical bulk modulus and the numerical Poisson ratio are not physical parameters of the model because the present formulation describes the immersed structure as incompressible in all cases.\\
\indent The previously mentioned models, (\ref{cauchy-nodev}) -- (\ref{nodev-stab}), are all capable of modeling incompressible materials in the continuous case. In fact, $\pstab = 0$ in the continuous IB formulation because $J \equiv 1$; any value of $\kappas \ge 0$ describes an incompressible solid in this formulation. It is interesting to study all these formulations, however, because we may lose discrete incompressibility of the solid even if we maintain a discretely divergence-free Eulerian velocity field. Therefore, (\ref{cauchy-nodev}) -- (\ref{nodev-stab}) may result in different structural deformations in the discretized equations.\\
\indent As mentioned previously in Section \ref{Jumps}, including the stabilization pressure $\pstab$ and using the $\dev[\cdot]$ operator will change the jumps in the Eulerian Lagrange multiplier. Continuity of traction, equation (\ref{dynamic-bc}), is maintained in all cases, however. This is because condition (\ref{dynamic-bc}) may be simplified to either $\llbracket \cauchyf \nb \rrbracket = \left(\dev\left[\cauchys_- \right] - \pstab\mathbb{I}\right)\nb$ or $\llbracket \cauchyf \nb \rrbracket = \left(\cauchys_- - \pstab\mathbb{I}\right)\nb$, depending on whether we use (\ref{dev-stab}) or (\ref{nodev-stab}), respectively, for the solid stress formulation. Likewise, it may be shown that the jump in the fluid traction is the structural force evaluated at the interface; see previous work by Lai and Li \cite{Lai2001}. Explicitly, the jump in fluid traction is either $\left(\dev\left[\cauchys_- \right] - \pstab\mathbb{I}\right)\nb$ or $\left(\cauchys_- - \pstab\mathbb{I}\right)\nb$, depending on our choice of solid formulation. Therefore, the dynamic boundary condition is satisfied for all formulations introduced in this section. Finally, these changes to the elastic part of the stress will not alter the jumps in the viscous stress because the changes involve isotropic tensors $\xi \mathbb{I}$ (for some scalar $\xi$) that disappear when multiplied by two mutually orthogonal vectors (e.g. $\mathbf{t}^T\xi \mathbb{I}\nb = 0$).

\subsubsection{Unmodified Model}
Model \eqref{cauchy-nodev} and \eqref{nodev-stab} can be derived by assuming the energy functional:
\begin{align}
	\Psi &= W(\FF), \label{bad_energy} \text{ and} \\
	\Psi &= W(\FF) + U(J) \label{stab_energy},
\end{align}
respectively. Note that $U(J)$ represents the stabilization energy proposed in our method. Model \eqref{cauchy-nodev} and therefore energy \eqref{bad_energy} lead to poor numerical results (see Section \ref{benchmarks}).

\subsubsection{Modified Model}
\indent It is possible to obtain models (\ref{cauchy-dev}) and (\ref{dev-stab}) in different ways. One way is through the Flory decomposition $\FFb= J^{-1/3}\FF$ \cite{Flory1961}. Note that $\det(\overline{\FF}) = 1$ by construction. We reformulate strain energy functionals (\ref{bad_energy}) and (\ref{stab_energy}), respectively, as:
\begin{align}
	\Psi &= W(\FFb), \label{mod_energy} \text{ and} \\
	\Psi &= W(\FFb) + U(J). \label{decoupled_energy}
\end{align}
The energies given by \eqref{mod_energy} and (\ref{decoupled_energy}) completely decouple energy associated to volume changing and volume preserving motions and achieves the desired split in the Cauchy stress. This decoupling is motivated by the physical assumption that a uniform pressure only produces changes in size but not changes in shape. Work by Sansour  \cite{Sansour2008} explores this physical assumption in depth, but we include only a brief explanation here. Specifically, with (\ref{decoupled_energy}), we obtain an additive split in the Cauchy stress into purely deviatoric and dilatational stresses. This means that the only contributions to the stabilizing pressure will come from $U(J)$.\\
\indent We show that using the model with the Flory decomposition has a similar effect as using the deviator operator. Let $\Wb$ denote $W(\FFb)$. The derivative of $\Wb$ is $\frac{\partial \Wb}{\partial \FF} = \frac{\partial \Wb}{\partial \overline{\FF}}:\frac{\partial \FFb}{\partial \FF}$, in which $\frac{\partial \FFb}{\partial \FF}$ is a fourth order tensor. Explicitly,
\begin{equation}
	\frac{\partial \FFb}{\partial \FF} = J^{-1/3}\left(\mathcal{I} - \frac{1}{3} \left( \FF^{-T} \otimes \FF \right)\right), \label{dFFbar-dFF}
\end{equation}
with $\mathcal{I}$ denoting the fourth order identity tensor. By contracting ($\ref{dFFbar-dFF}$) with $\frac{\partial \Wb}{\partial \FFb}$, we obtain
\begin{equation}
	\frac{\partial \Wb}{\partial \FFb}:\frac{\partial \FFb}{\partial \FF} = J^{-1/3}\left(\frac{\partial \Wb}{\partial \overline{\FF}} - \frac{1}{3} \left(\frac{\partial \Wb}{\partial \overline{\FF}} :\FF \right)\FF^{-T}\right). \label{dWbar-dFF}
\end{equation}
Pushing forward ($\ref{dWbar-dFF}$), it is evident that using the Flory decomposition yields an elastic stress whose corresponding Cauchy stress is traceless,
\begin{equation}
	\tr \left[ \frac{1}{J}\left( \frac{\partial \Wb}{\partial \FFb}:\frac{\partial \FFb}{\partial \FF} \right) \FF^T \right] = \tr \left[ J^{-4/3}\left(\frac{\partial \Wb}{\partial \overline{\FF}}\FF^T - \frac{1}{3} \tr \left(\frac{\partial \Wb}{\partial \overline{\FF}} \FF^T \right)\mathbb{I}\right) \right] = 0. \label{cauchy-traceless}
\end{equation}
Alternatively, note that applying the deviator operator (\ref{dev}) and subsequently the trace operator to $J^{-4/3}\frac{\partial \Wb}{\partial \overline{\FF}}\FF^T$ yields the same result as (\ref{cauchy-traceless}). These calculations are well established in literature on nonlinear solid mechanics and FE methods for nearly incompressible elasticity \cite{JavierBonetandRichardWood2008, TedBelytschkoWingKamLiu2000}.

\subsubsection{Deviatoric Projection}
\indent For material models in which an elastic energy cannot be defined, it is possible to define a deviatoric Cauchy stress using the deviatoric projection implicitly defined by equation (\ref{dev}). In our tests, deviatoric projections of hyperelastic models will be constructed by using the deviator operator for the first Piola-Kirchhoff stress,
\begin{equation}
	\DEV[\PPs] = \PPs - \frac{1}{3}\left(\PPs:\FF \right)\FF^{-T}. \label{DEV}
\end{equation}
Note that (\ref{DEV}) resembles (\ref{dWbar-dFF}) except for the $J^{-1/3}$ pre-factor. Also notice that pushing forward (\ref{DEV}) yields a traceless tensor,
\begin{equation}
	\DEV[\PPs]\FF^T = \PP \FF^T - \frac{1}{3}\tr \left(\PPs \FF^T \right)\mathbb{I},
\end{equation}
in which we have used the identity $ \tr\left(\PPs\FF \right) = \PPs:\FF^T$. Thus we have the relationship
\begin{equation}
	\DEV\left[\PPs\right]\FF^T = \dev\left[\PPs\FF^T\right]. \label{DEV-dev}
\end{equation}
The quantity on the right of equation (\ref{DEV-dev}) is equal to $\dev[J\cauchys]$. In our numerical computations, the $\DEV[\cdot]$ operator defined by (\ref{DEV}) is applied to the first Piola-Kirchhoff stress derived from equation (\ref{bad_energy}) and yields a Cauchy stress with the desired split. We study these models with and without the volumetric stabilization.

%%%%%%%%%%%%%%%%%%%%%

\subsection{Constitutive Laws}
For incompressible solids with isotropic material responses, we express $\Psi$ as a function of the first two tensor invariants of the right Cauchy-Green tensor, $\CC = \FF^T \FF$. These invariants are $I_1 = \tr(\CC)$ and $I_2 = \frac{1}{2}\left(I_1^2 - \tr(\CC^2) \right)$. This relationship between $\Psi$ and $\FF$ ensures for material frame indifference \cite{JavierBonetandRichardWood2008}. In incompressible cases, $J \equiv 1$, and there is no dependence on the third invariant, $I_3 = \det(\CC) = J^2$. In the compressible regime, of course, this is not the case. Often, the energy for compressible and nearly incompressible materials is written as a function of $\bar{I}_1= J^{-2/3}I_1$ and $\bar{I}_2 = J^{-4/3}I_2$, which are the invariants of $\overline{\CC} = \FFb^T\FFb$. Modifying the invariants in this way removes information about the volume change. Thus we refer to the invariants of $\overline{\CC}$ as the modified invariants.\\
\indent The invariant based models used in this work are of the form
\begin{align}
	\Psi &= W(I_1,I_2), \text{ and}\\
	\Psi &= W(I_1,I_2) + U(J)
\end{align}
for the deviatoric projection and unmodified models and 
\begin{align}
	\Psi &= W(\bar{I}_1,\bar{I}_2), \text{ and} \\
	\Psi &= W(\bar{I}_1,\bar{I}_2) + U(J)
\end{align}
for the modified models. We briefly describe specific constitutive models that do and do not have the desired deviatoric split.
\subsubsection{Neo-Hookean Models}
The neo-Hookean model is a simple hyperelastic model that depends only on the first invariant. Using the unmodified invariants, its energy and first Piola-Kirchhoff stress with stabilization are
\begin{align}
	\Psi &= \frac{G}{2} \left( I_1 - 3\right) + \frac{\kappas}{2}(\ln J)^2 \label{nh_energy} \text{ and}\\
	\PPs &=  G\FF + \kappas \ln(J) \FF^{-T}. \label{nh_stress}
\end{align}
\noindent The Young's modulus is commonly used to describe neo-Hookean materials. The Young's modulus $E$ is related to $G$ via $G = \frac{E}{2(1+\nu)}$. Here we use $\nu = \frac{1}{2}$ to relate $G$ and $E$ because we are modeling a material whose motions are incompressible. \\
\indent When using modified invariants, the energy and stress are
\begin{align}
	\Psi &= \frac{G}{2}\left(\bar{I}_1 - 3\right) + \frac{\kappas}{2}(\ln J)^2 \text{ and}\\
	\PPs &=  G J^{-2/3}\left(\FF - \frac{I_1}{3}\FF^{-T}\right) + \kappas \ln (J) \FF^{-T} \label{nh_stress_mod}.
\end{align}
Taking the deviatoric projection of (\ref{nh_stress}), we have
\begin{equation}
\noindent \PPs =  G\left(\FF - \frac{I_1}{3}\FF^{-T}\right) + \kappas \ln (J) \FF^{-T}. \label{nh_stress_dev}
\end{equation}\\
Note the similarities between equations (\ref{nh_stress_mod}) and (\ref{nh_stress_dev}); the difference is that equation (\ref{nh_stress_mod}) includes the $J^{-2/3}$ pre-factor.
\subsubsection{Mooney-Rivlin Models}
Mooney-Rivlin material models also include linear dependence on $I_2$ in the energy. The unmodified invariant case is given by
\begin{align}
	\Psi&= c_1 \left( I_1 - 3\right) + c_2(I_2 - 3) + \frac{\kappas}{2}(\ln J)^2 \label{mr_energy} \text{ and}\\
	\PPs &=  2c_1\FF + 2c_2(I_1\FF - \FF \CC) + \kappas \ln (J) \FF^{-T}, \label{mr_stress}
\end{align}
in which $c_1$ and $c_2$ are material constants. Using modified invariants yields
\begin{align}
	\Psi &= c_1\left(\bar{I}_1 - 3\right) + c_2(\bar{I}_2 - 3) + \frac{\kappas}{2}(\ln J)^2 \text{ and}\\
	\PPs &=  2c_1 J^{-2/3}\left(\FF - \frac{I_1}{3}\FF^{-T}\right) + 2c_2J^{-4/3}\left(I_1\FF - \FF \CC - \frac{2I_2}{3}\FF^{-T}\right) + \kappas \ln (J) \FF^{-T} \label{mr_stress_mod}.
\end{align}
Taking the deviatoric projection of (\ref{mr_stress}) yields
\begin{equation}
\PPs =  2c_1 \left(\FF - \frac{I_1}{3}\FF^{-T}\right) + 2c_2\left(I_1\FF - \FF \CC - \frac{2I_2}{3}\FF^{-T}\right) + \kappas \ln (J) \FF^{-T} \label{mr_stress_dev}.
\end{equation}
\indent As with the neo-Hookean models, we study material models described in this work with and without volumetric stabilization. For the Mooney-Rivlin material law, this will require being able to relate material constants to $\kappas$. For consistency between the small deformation (linear) and large deformation (nonlinear) regimes, we set $G = 2(c_1 + c_2)$ when calculating $\kappas$. This allows the use of the same formula, equation (\ref{kappa-nu}), that relates $\kappas$ and $\nus$ to a material quantity. \\
\subsubsection{Modified Standard Reinforcing Model}
\indent To examine the effects of anisotropy, we use the modified standard reinforcing model \cite{Murphy2013}. This model describes transversely isotropic materials with fibers given by a material vector $\Ab$ in the reference configuration and $\ab = \FF \Ab$ in the current configuration. The effect of the anisotropy appears through the anisotropic invariants $I_4$ and $I_5$:
\begin{align}
	I_4 &= \Ab^T\CC \Ab = \ab ^T \ab,  \text{ and} \\
 	I_5 &= \Ab^T \CC^2 \Ab =  \ab^T \BB \ab,
\end{align}
in which $\BB = \FF \FF^T$ is the left Cauchy-Green strain. Because $\ab$ is the stretched and rotated material vector, $I_4$ measures the stretch of the fiber, whereas $I_5$ encodes information related to the shear as well as the stretch \cite{Merodio2005}.  The modified standard reinforcing model is
\begin{align}
\Psi &= \frac{\Gt}{2}(I_1 - 3) + \frac{\Gt - \Gl}{2}(2I_4 - I_5 - 1) + \frac{\El + \Gt - 4\Gl}{8}(I_4 - 1)^2 + \frac{\kappas}{2}(\ln J)^2 \text{ and} \label{sr_energy}  \\
 \PPs &= \Gt \FF + (\Gt - \Gl)\left(2\FF \MM - \FF \MM\CC - \FF \CC \MM\right) +\frac{\El + \Gt - 4\Gl}{2}(I_4 - 1)\FF \MM + \kappas \ln (J) \FF^{-T}, \label{sr_stress}
\end{align}
\noindent in which $\MM = \Ab \otimes \Ab$. Here, $\Gt$ is the shear modulus of the material in the plane transverse to the fibers, and $\Gl$ is the shear modulus along the length of the fibers. To determine $\kappas$, $\Gt$ is used in equation ($\ref{kappa-nu}$) because this material model does not involve an isotropic shear modulus $G$. The parameter $\El$ is similar to a Young's modulus but in the direction of the fiber. 
When we modify $I_1$, we instead obtain:
%\begin{align}
%\Psi =& \frac{\Gt}{2}(\bar{I}_1 - 3) + \frac{\Gt - \Gl}{2}(2I_4 - I_5 - 1) + \frac{\El + \Gt - 4\Gl}{8}(I_4 - 1)^2 + \frac{\kappas}{2}(\ln J)^2,  \text{ and}\\
% \PPs =& \Gt J^{-2/3} \left(\FF - \frac{I_1}{3} \FF^{-T}\right) + (\Gt - \Gl)\left(2\FF \MM - \FF \MM\CC - \FF \CC \MM\right) \label{sr_stress_mod} \\
% &\mbox{} +\frac{\El + \Gt - 4\Gl}{2}(I_4 - 1)\FF \MM + \kappas \ln (J) \FF^{-T}. \nonumber
%\end{align}
\begin{align}
\Psi =& \frac{\Gt}{2}(\bar{I}_1 - 3) + \frac{\Gt - \Gl}{2}(2I_4 - I_5 - 1) + \frac{\El + \Gt - 4\Gl}{8}(I_4 - 1)^2 + \frac{\kappas}{2}(\ln J)^2,  \text{ and}\\
 \PPs =& \Gt J^{-2/3} \left(\FF - \frac{I_1}{3} \FF^{-T}\right) + (\Gt - \Gl)\left(2\FF \MM - \FF \MM\CC - \FF \CC \MM\right) +\frac{\El + \Gt - 4\Gl}{2}(I_4 - 1)\FF \MM + \kappas \ln (J) \FF^{-T}. \label{sr_stress_mod}
\end{align}
\noindent Likewise, for the modified standard reinforcing model, the deviatoric projection is:
\begin{align}
\PPs =& \Gt  \left(\FF - \frac{I_1}{3} \FF^{-T}\right) + (\Gt - \Gl)\left(2\FF \MM - \FF \MM\CC - \FF \CC \MM\right) +\frac{\El + \Gt - 4\Gl}{2}(I_4 - 1)\FF \MM +  \kappas \ln (J) \FF^{-T}.  \label{sr_stress_dev} 
\end{align}
\indent The anisotropic models considered for the benchmarks herein take the following forms:
\begin{align}
\Psi &= W(\bar{I}_1, \bar{I}_2, I_4, I_5) + U(J),\\
\Psi &= W(I_1, I_2, I_4, I_5) + U(J), \\
\Psi &= W(\bar{I}_1, \bar{I}_2, I_4, I_5), \text{or}\\
\Psi &= W(I_1, I_2, I_4, I_5), \label{bad_energy_aniso}
\end{align}
with the exception of the deviatoric projection of the standard reinforcing model, which does not arise from an energy functional. If we modify both $I_4$ and $I_5$ by using $\overline{\CC}$ and use the volumetric stabilization, we arrive at a Cauchy stress with an additive deviatoric-spherical split. In cases of uniform pressure, however, it is possible that a body will undergo a shape-changing deformation if the material is anisotropic. Thus, the volumetric split is not appropriate for the anisotropic part of the stress. \\
\indent We remark that in biomechanics literature, the standard reinforcing model is often defined as $\Psi =  c_1(I_1 - 3)  + c_4(I_4 - 1)^2$, without any dependence on $I_5$ \cite{Merodio2005}. It can be shown that omitting this anisotropic invariant implies that the linearized shear moduli in the direction of the fibers and perpendicular to the fibers must be the same. It can also be shown that the three modes of shear characteristic of transversely isotropic materials are not represented if $I_5$ is omitted \cite{Murphy2013}. The modified standard reinforcing model is arrived at by augmenting the standard reinforcing model in a way that allows for the consistency between the linear and finite regimes \cite{Merodio2005}. \\

%%%%%%%%%%%%%%%%%%%%%%
%%%%%%%%%%%%%%%%%%%%%%

\section{Numerical Methods}
\label{Numerical Methods}
This work focuses on improving the treatment of incompressible, nonlinearly elastic structures of the immersed boundary-finite element method of Griffith and Luo \cite{Griffith2017}. The benchmark problems that we consider here do not have analytic solutions, and so benchmark solutions are obtained using an FE method for large-deformation incompressible elasticity \cite{Chiumenti, Masud2013}. We briefly describe both numerical methods here. For FSI simulations, we use IBAMR \cite{IBAMR, Griffith2007}, which is an open-source adaptive and distributed-memory parallel implementation of the IB method. Specifically, we use the ``IBFE" module in IBAMR, which allows the use of volumetric structures. The quasi-static finite element benchmark solutions are computed using BeatIt \cite{Beatit}. Both IBAMR and BeatIt rely on the parallel C++ finite element library libMesh \cite{libmesh} and
on linear and nonlinear solver infrastructure provided by the PETSc library \cite{Petsc}. 

\subsection{Immersed Boundary-Finite Element Formulation}
 The Immersed Boundary-Finite Element (IBFE) method for FSI has been described in detail previously \cite{Griffith2017, Griffith2009}. Briefly, a staggered-grid finite difference method is used to discretize the Eulerian equations, and a nodal FE method is used to discretize the Lagrangian equations \cite{HowardElmanDavidSilvester2014}. In this scheme, the Eulerian velocity $\ub$ is approximated at the cell edges (faces in three spatial dimensions), and the Eulerian Lagrange multiplier field $\peul$ is approximated at cell centers. We use standard second order accurate finite differences to discretize the Eulerian incompressible Navier-Stokes equations (\ref{ns}) and (\ref{divfree}) \cite{Griffith2009, Leveque2007}. In our computations, we use the unified weak formulation of the hyperelastic IB method, which corresponds to determining the structural force $\Fb$ via equation (\ref{weak-force}) \cite{Griffith2017}. \\
\indent We discretize the structure $\soliddomO$ via a triangulation $\mathcal{T}_h = \cup_e K^e$, in which $K^e$ are isoparametric elements. On $\mathcal{T}_h$, we define Lagrange basis functions $\{\phi_\ell(\Xb)\}_{\ell=1}^m$, in which $m$ is the number of FE nodes in our mesh. In the computations performed in this study, these functions belong to the common FE spaces of \textbf{P1}, \textbf{P2}, \textbf{Q1}, and \textbf{Q2}, which in two spatial dimensions denote the spaces of linear, quadratic, bilinear, and biquadratic basis functions, respectively; see, for instance, Elman \etal~\cite{HowardElmanDavidSilvester2014}. In three dimensions, we only use \textbf{P1} and \textbf{Q1}, which are the spaces of linear and trilinear basis functions. We include low order spaces to demonstrate the robustness of the stabilization method and show that even low order approximations may also achieve favorable volume conservation.\\
\indent In all cases, the mapping $\cb$ is approximated by $\cb_h(\Xb,t) = \sum_{\ell=1}^m \cb_{\ell}(t)\phi_{\ell}(\Xb)$. The Lagrangian force is also approximated with the same FE basis via $\Fb_h(\Xb,t) = \sum_{\ell=1}^m \Fb_{\ell}(t) \phi_{\ell}(\Xb)$. Let $[\Fb]$ be the vector representation of $\Fb_h$, and let $ [\boldsymbol{B}]$ be the vector with entries $-\int_{\soliddomO} \PPs \cdot \nabla_{\Xb}\phi_l(\Xb) \, d\Xb$. Equation (\ref{weak-force}) may be written discretely as
\begin{equation}
	[\mathcal{M}][\Fb] = [\boldsymbol{B}], \label{discrete-force}
\end{equation}
in which $[\mathcal{M}]$ is the vector mass matrix consisting of scalar mass matrices as the diagonal blocks. The scalar mass matrix has entries $\int_{\soliddomO} \phi_l(\Xb) \phi_m(\Xb) \, d\Xb$. Gaussian quadrature rules are used to integrate the integral equations. For each case, we choose integration orders that would exactly integrate the bilinear form that arises in the weak form of static linear elasticity. Specifically, in this work, we do not employ selective reduced integration \cite{ThomasJRHughes2000}.\\
\indent To describe equation ($\ref{elasticforce}$), we say the force is \emph{prolonged}, or spread, from the solid region to the computational domain. For equation ($\ref{noslip}$), we say the solid \emph{interpolates} the Eulerian velocity. The force prolongation and velocity restriction operators of the coupling are constructed to be adjoints. This means that there is conservation of power as we map data between the Eulerian grid and Lagrangian mesh \cite{Peskin2002}. The discrete force prolongation operator $\mathcal{S}$ relates $\fb$ and $\Fb$ via $\fb = \mathcal{S}\left[\cb\right]\Fb$. In two spatial dimensions, the operator $\mathcal{S}\left[\cb\right]$ is implicitly defined by
\begin{align}
	(f_1)_{i-\frac{1}{2},j} &= \sum_{K^e \in \mathcal{T}_h}\sum_{Q=1}^{N^e} F_1(\Xb_Q^e,t)\delta_h(\xb_{i-\frac{1}{2},j} - \cb_h(\Xb_Q^e,t))w_Q^e, \label{discrete-force1}\\
	(f_2)_{i,j-\frac{1}{2}} &= \sum_{K^e \in \mathcal{T}_h}\sum_{Q=1}^{N^e} F_2(\Xb_Q^e,t)\delta_h(\xb_{i,j-\frac{1}{2}} - \cb_h(\Xb_Q^e,t))w_Q^e, \label{discrete-force2}
\end{align}
in which $\Xb^e_Q$ are quadrature points, $w^e_Q$ are quadrature weights, and $N^e$ is the number of quadrature points on element $e$. $\delta_h(\xb) = \prod_{i=1}^2\delta_h(x_i)$ denotes a two-dimensional regularized delta function corresponding to the four point kernel function used by Peskin \cite{Peskin2002}. We denote by $\xb_{i-\frac{1}{2},j}$ and $\xb_{i,j-\frac{1}{2}}$ the Eulerian grid point at $(i\Delta x, (j+\frac{1}{2})\Delta x)$ and $((i+\frac{1}{2})\Delta x, j\Delta x)$, respectively. $\Delta x$ is the Eulerian grid spacing. The first component of the Eulerian force density is $(f_1)_{i-\frac{1}{2},j}$, and it is evaluated at $\xb_{i-\frac{1}{2},j}$; $(f_2)_{i,j-\frac{1}{2}}$ is evaluated at $\xb_{i,j-\frac{1}{2}}$. The discrete velocity restriction operator $\mathcal{J}$ relates $\Ub$ and $\ub$ via $\Ub = \mathcal{J}\left[\cb\right]\ub$. The operator $\mathcal{J}\left[\cb\right]$ is implicitly defined by first determining an intermediate velocity field $\Ub^\text{IB}$ via
\begin{align}
	U_1^{\text{IB}}(\Xb,t) &= \sum_{i,j} (u_1)_{i-\frac{1}{2},j}\delta_h(\xb_{i-\frac{1}{2},j} - \cb_h(\Xb,t))\Delta x^2, \\
	U_2^{\text{IB}}(\Xb,t) &= \sum_{i,j} (u_2)_{i,j-\frac{1}{2}}\delta_h(\xb_{i,j-\frac{1}{2}} - \cb_h(\Xb,t))\Delta x^2,
\end{align}
in which $(u_1)_{i-\frac{1}{2},j}$ and $(u_2)_{i,j-\frac{1}{2}}$ are the components of the discrete Eulerian velocity field evaluated at $\xb_{i-\frac{1}{2},j}$ and $\xb_{i,j-\frac{1}{2}}$, respectively. $\Ub^{\text{IB}}$ is a Lagrangian velocity field that is in general not a linear combination of the basis functions $\{\phi(\Xb)\}_{\ell=1}^m$. To obtain a velocity field $\Ub_h$ that can be represented using the FE basis functions, $\Ub^{\text{IB}}$ is projected onto $\{\phi(\Xb)\}_{\ell=1}^m$ in an $L^2$ sense. \\ 
%Although $\mathcal{S}$ is defined by equations (\ref{discrete-force1}) -- (\ref{discrete-force2}), $\mathcal{J}$ is implicitly defined via a discrete conservation of power identity.
% Specifically,
%\begin{equation}
%	\left(\Fb, \mathcal{J} \ub\right)_{\Xb} = \left(\mathcal{S}\Fb, \ub\right)_{\xb}, \label{adjoints}
%\end{equation}
%in which $(,)_{\Xb}$ is the Lagrangian inner produce and $(,)_{\xb}$ is the Eulerian inner product.
%Let $[\ub]$ be the vector representation of the Eulerian velocity, and let $[\mathcal{S}]$ and $[\mathcal{J}]$ be the matrix representations of the coupling operators. The explicit version of equation (\ref{adjoints}) is
%\begin{equation}
%	[\Fb]^T[\mathcal{M}][\mathcal{J}][\ub] = \left([\mathcal{S}][\Fb]\right)^T[\ub]h^2.
%\end{equation}
%This relationship defines $\mathcal{J}$ to be 
%\begin{equation}
%	[\mathcal{J}] = [\mathcal{M}]^{-1}[\mathcal{S}]^Th^2.
%\end{equation}
%However, because we only need the action of $\mathcal{J}$ on $\ub$, we do not explicitly construct $\mathcal{J}$.\\
\indent We use the notation $\cb^n$ to denote our discrete approximation to the deformation at time $n\Delta t$, in which $\Delta t$ is the time-step size, and apply this convention to all quantities of interest; $\cb^{n+1/2}$ indicates an intermediate value used to calculate the solution at the next time-step. The operators $\nabla_h, \nabla_h \cdot$, and $\nabla_h^2$ denote second order finite difference operators as defined in previous works \cite{Griffith2009, Leveque2007}. The convective term $\mathcal{N}^{n+1/2} = \frac{3}{2}\ub^n \cdot \nabla_h \ub^n - \frac{1}{2} \ub^{n-1} \cdot \nabla_h \ub^{n-1}$ is computed using a piecewise parabolic method-type (PPM) spatial approximation \cite{Griffith2009} along with an Adams-Bashforth temporal discretization. The basic time stepping scheme is summarized in Algorithm \ref{ts-algorithm}. Notice that this scheme requires the solution of two mass matrix systems (to determine the Lagrangian structural force $\Fb$ and the Lagrangian structural velocity $\Ub$) along with the incompressible Stokes equations (to determine the Eulerian material velocity field $\ub$ and Lagrange multiplier field $\pi$). Notice also that it is not possible to determine the structural deformation without accounting for the effects of the force-spreading and velocity-interpolation operators. Although it is possible to formally eliminate the Eulerian velocity and Lagrange multiplier fields from the equations of motion \cite{Usabiaga2016, Kallemov2016}, e.g.~by forming an appropriate Schur complement system, doing so affects the algorithms used in solving the discrete equations but not the solutions themselves. Algorithm \ref{ts-algorithm} uses multiple previous values for $\ub$ and $\peul$ to compute subsequent approximations, so this scheme cannot be used for the initial time-step. Instead, a predictor-corrector method is used. For further details, we refer to the description of the method by Griffith and Luo \cite{Griffith2017}.\\
\indent All Dirichlet boundary conditions for the structure are imposed via a penalty method. Specifically, surface forces of the form $\boldsymbol{T} = \kappa_{\text{D}}(\cb_{\text{D}} - \cb)$ are applied to the boundary of the structure to approximate Dirichlet boundary conditions given by $\cb_{\text{D}}$. In this approach, the parameter $\kappa_{\text{D}}$ denotes a stiffness that penalizes deviations from the desired value. We use the scaling $\kappa_{\text{D}} \propto \frac{\Delta x}{(\Delta t)^2}$, so that the stiffness parameter increases as the mesh is refined. For all tests, we use a proportionality constant of $2.5$ to determine $\kappa_{\text{D}}$, yielding $\kappa_{\text{D}} = 2.5\frac{\Delta x}{(\Delta t)^2}$. \\
\begin{algorithm}
\caption{Time-Stepping Scheme for advancing from $\cb^n, \ub^n$, $\peul^n$ to $\cb^{n+1}, \ub^{n+1}$, $\peul^{n+1}$}
\label{ts-algorithm}
 Given $\cb^n, \ub^n$ \\
 Update $\cb^{n+1/2} = \cb^n + \frac{\Delta t}{2} \mathcal{J}\left[\cb^n\right] \ub^n$ \\
% Update $\left[\mathbf{B}\right]^{n+1/2} = \mathcal{A}\left(\cb^{n+1/2}\right)$ \\
 Evaluate $\Fb^{n+1/2}$ by using $\cb^{n+1/2}$ to solve equation (\ref{discrete-force})\\
 Update $\fb^{n+1/2} = \mathcal{S}\left[\cb^{n+1/2} \right]\Fb^{n+1/2}$ \\
 Solve for $\left(\ub^{n+1},\peul^{n+1/2}\right)$: $\begin{cases} \rho \left(\frac{\ub^{n+1} - \ub^n}{\Delta t} + \mathcal{N}^{n+1/2} \right) &= -\nabla_h \peul^{n+1/2} +\mu \nabla_h^2 \left(\frac{\ub^{n+1} + \ub^n}{2} \right) + \fb^{n+1/2}\\
 \nabla_h \cdot \ub^{n+1} &= 0 \end{cases}$ \\
 Update $\cb^{n+1} = \cb^n + \Delta t \mathcal{J}\left[ \cb^{n+1/2} \right] \left(\frac{\ub^{n+1} + \ub^n}{2}\right) $
\end{algorithm}

%%%%%%%%%%%%%%%%%%%%%%

\subsection{Structural Mechanics}
The results from the FSI calculations are compared to quasi-static fully incompressible elasticity finite element simulations. Specifically, we use a mixed displacement-pressure formulation to enforce the incompressibility constraint \cite{FrancoBrezziandMichaelFortin2012}. Piecewise linear polynomials are used to interpolate the displacement and pressure fields. This choice leads to an unstable numerical method, as the corresponding \textbf{P1}/\textbf{P1} elements (indicating a piecewise linear approximation for both the displacement and pressure) do not satisfy the Ladyzhenskaya-Brezzi-Babuska (LBB) condition \cite{HowardElmanDavidSilvester2014}. Work by Hughes, Franca, and Balestra \cite{Hughes1986} circumvented the LBB condition for the Stokes equations using a simple stabilization method suitable for equal order approximations and recovering the optimal order of convergence. Their method has been extended to linear and nonlinear elasticity \cite{Franca1988, Klaas1999}, and it has been reinterpreted in a variational multiscale framework \cite{Chiumenti, Masud2013}. We use this stabilized method \cite{Chiumenti, Masud2013} to solve the quasi-static incompressible nonlinear elasticity equations given by
\begin{align}
	\DIV \PPs(\Xb,t) &= 0, \label{solid-mech1} \\
	J(\Xb,t) &= 1, \quad \Xb \in \soliddomO. \label{solid-mech2}
\end{align}
Here $\PPs = \frac{\partial \Psi}{\partial \FF}$ and $\Psi = W - p(J-1)$. Modified invariants are used in this FE method by defining $W$ as either $W = W(\bar{I}_1,\bar{I}_2)$ or $W = W(\bar{I}_1,\bar{I}_2,I_4,I_5)$, implying that $\tr\left[\frac{\partial W}{\partial \FF}\FF^T \right] = 0$. Note that the variable $p$ serves as both the physical pressure and the Lagrange multiplier for the incompressibility constraint. We emphasize that this numerical scheme is \emph{not} used for the FSI computations, but is merely employed to generate structural deformations used for comparison against our FSI numerical methods. The FSI numerical methods, as detailed above, independently compute the deformations of the immersed solid.

%%%%%%%%%%%%%%%%%%%%%%
%%%%%%%%%%%%%%%%%%%%%%

\section{Benchmarks}
\label{benchmarks}
Most of our benchmark tests are derived from standard benchmark problems for incompressible elasticity drawn from the solid mechanics literature, except that here the solid bodies are embedded in an incompressible Newtonian fluid. As mentioned in Section \ref{Jumps}, the steady-state solutions of the FSI problems are the same as those from pure solid mechanics formulations. This follows from that the fact that at steady-state ($\ub = \mathbf{0}$), $\Div \cauchysolid = \mathbf{0}$ implies $\Div \cauchys - \nabla \peul = \mathbf{0}$, which has the same solution as the solid mechanics problem described by (\ref{solid-mech1}) -- (\ref{solid-mech2}).\\
\indent We report the computed displacements of a point of interest and the total volume conservation for each benchmark. Where relevant, we study the same points of interest as employed in previous studies. We also report the displacement results from the fully incompressible elasticity calculations, labeled as ``FE (\textbf{P1}/\textbf{P1})" in the displacement plots. The majority of the deformations presented in the figures do not show the Eulerian velocity field of the computational domain because the fluid flow field will be zero at steady state. When listing the ranges of the percent change in volume, we omit the coarsest discretizations of the structure. We also present the deformations of each benchmark for all relevant formulations of the energy functional (e.g.~equations (\ref{bad_energy}) -- (\ref{decoupled_energy}) for the isotropic material models). These visualizations depict the average values of $J$ for each element, and the extents of the color bar indicate cutoff values. The average of $J$ for element $e$ is calculated as:
\begin{equation}
	\text{Avg} \ J^e = \left( \sum_Q J(\Xb^e_Q) w^e_Q \right) \big/ \left( \sum_Q w^e_Q \right), \label{avgJ}
\end{equation}
recalling that $\Xb^e_Q$ and $w^e_Q$ are the quadrature points and quadrature weights associated with the element $e$. We use the same Gaussian quadrature rule here as for the approximation of the integrals in other parts of our method. We also report the stresses, pressures, and principal stretches for a representative selection of the benchmark tests.\\ 
\indent We use material models for the structure with modified isotropic invariants, unmodified isotropic invariants, and the deviatoric projection of the isotropic material response of the elastic stress. These models are studied with varying levels of volumetric stabilization that are tuned through different choices of $\nus$. Except where otherwise noted, we consider numerical Poisson ratios of $\nus = -1, 0, 0.4$, and $.49995$. We recall that $\nus = -1$ corresponds to the case of zero numerical bulk modulus and thus zero volumetric energy-based stabilization. We use $\nus = .49995$ to study the effect of volumetric locking (if any), and $\nus = 0$ and $\nus = .4$ are studied as intermediate values between the two extremes. Preliminary tests indicated that the displacement solutions begin to deteriorate (experience locking) slightly past $\nus = 0.4$ in many cases, especially on coarser computational grids. In fact, finer discretizations can withstand larger values of $\nus$ without suffering from locking, whereas coarser cases are more sensitive to values of $\nus$ close to $0.5$. This is consistent with the computational mechanics literature \cite{ThomasJRHughes2000}. \\
\indent Except where otherwise indicated, the computational domain is $\Omega = [0, L]^d$, in which $d = 2,3$ is the spatial dimension and $L$ is the domain length. Each structure is placed in the center of this computational domain. The Eulerian grid spacing is $\Delta x = \frac{L}{N}$, in which $N$ is the number of cells in one spatial dimension. The Lagrangian mesh width is $\Delta X$, and the mesh factor ratio $\mfac = \frac{\Delta X}{\Delta x}$ describes the relative grid spacing between the Eulerian and Lagrangian meshes. In our tests, the choice of $\mfac \leq 1$ is used for pressure driven cases and $\mfac \approx 2$ is used for shear driven cases. These choices of $\mfac$ were made based on preliminary tests. More specifically, we use $\mfac =1$ for the compression block test (Section \ref{Compression Test}); $\mfac = 2$ for the Cook's membrane, the anisotropic Cook's membrane, and the torsion test (Sections \ref{Cook's Membrane} -- \ref{Torsion}); and $\mfac = 0.5$ for the elastic band test (Section \ref{Elastic Band}).\\
\indent Most tests use zero velocity boundary conditions on the computational domain \cite{Griffith2009}. This allows for the fluid velocity to decay to zero. Additionally, for some tests viscous damping is used in the solid region to dampen oscillations and achieve near-critical damping of the structure's motion. This substantially decreases simulation time to reach steady state. Viscous damping is included to the solid region by adding a term $-\eta \Ub$ to the Lagrangian force $\Fb$, in which $\eta > 0$ is the damping coefficient. The simulations are run until a final time $T_{\text{f}}$, which is chosen such that the velocity is approximately zero. We use the polynomial $q(t) = -2\left(\frac{t}{T_{\text{l}} }\right)^3 + 3\left(\frac{t}{T_{\text{l}}}\right)^2$ to load the structure. At time $T_{\text{l}} = \alpha T_{\text{f}}$ the load is fully applied and sustained until the end of the simulation $T_{\text{f}}$. Here $\alpha \in (0,1)$. Between times $T_{\text{l}}$ and $T_{\text{f}}$, we let the structure relax to its resting configuration. Except where otherwise noted, the density is $\rho = 1.0 \ \frac{\text{g}}{\text{cm}^3}$, and the viscosity is $\mu = .01 \ \frac{\text{dyn} \cdot \text{s}}{\text{cm}^2}$, corresponding to water. In all tests, the structure and fluid are given the same density. This has no effect on the steady state cases but is a convenient choice that allows for efficient constant-coefficient linear solvers. \\

%%%%%%%%%%%%%%%%%%%%%%

\subsection{Compression Test}
\label{Compression Test}
\begin{figure}
\centering
\includegraphics[width=.7\linewidth, trim={80 0 70 0}, clip]{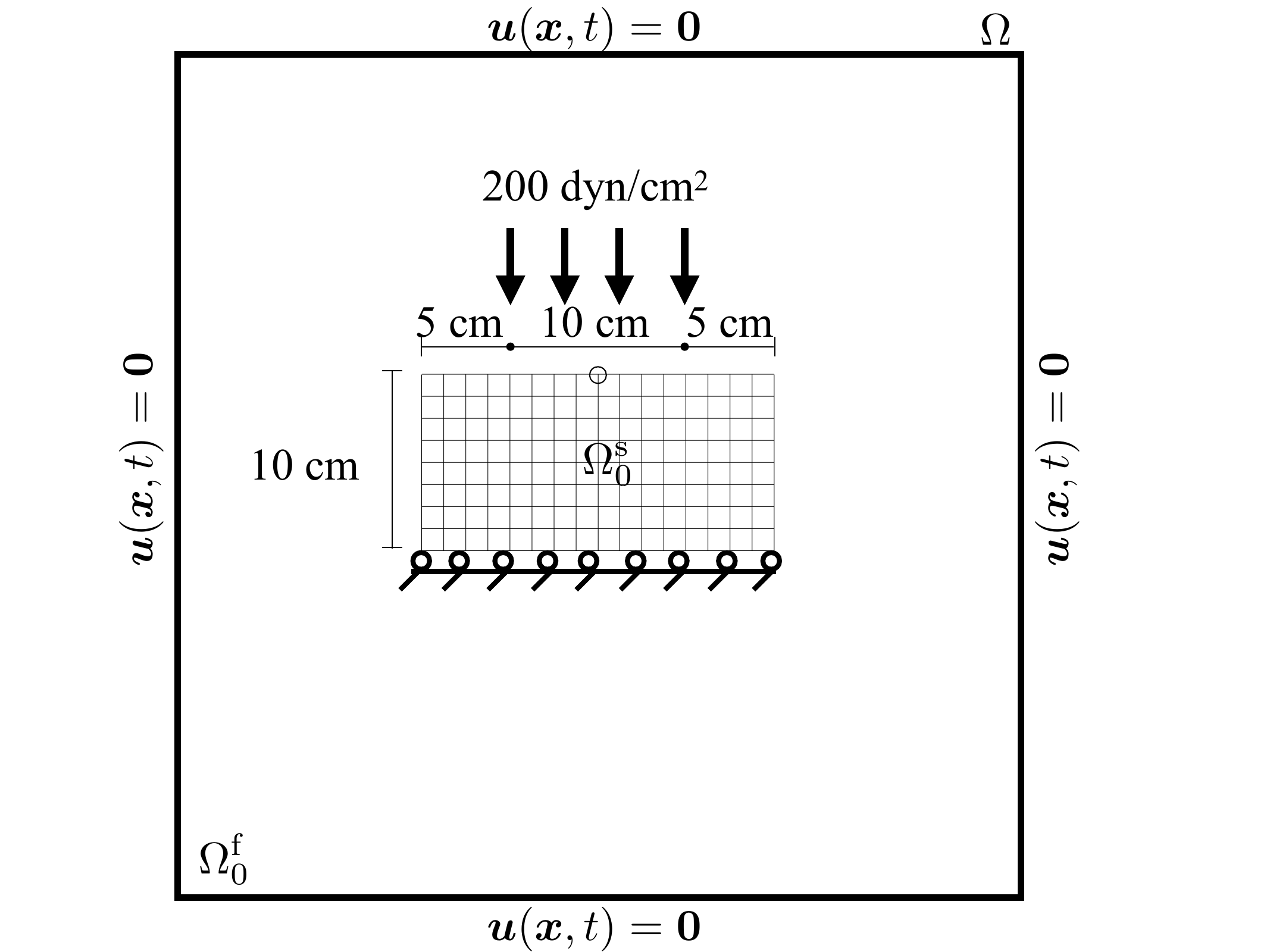}
\caption{Specifications of the compressed block benchmark (Section \ref{Compression Test}). The quantity of interest is the $y$-displacement as measured at the encircled point. The structure, shown here in its initial configuration and denoted by $\soliddomO$, is immersed in a fluid denoted by $\fluiddomO$. The entire computational domain is $\Omega = \fluiddom \cup \soliddom$. Zero fluid velocity is enforced on the boundary of $\Omega$.}
\label{comp_mesh}
\end{figure}

\begin{figure}
\begin{tabular}{l c c}
& \textbf{Modified Invariants} & \textbf{Unmodified Invariants} \\
\rotatebox{90}{\qquad\qquad \textbf{$\nus = .4$} }&
\subcaptionbox{\label{sfig:testa}} {\includegraphics[width=.45\linewidth, trim={200 200 70 400}, clip]{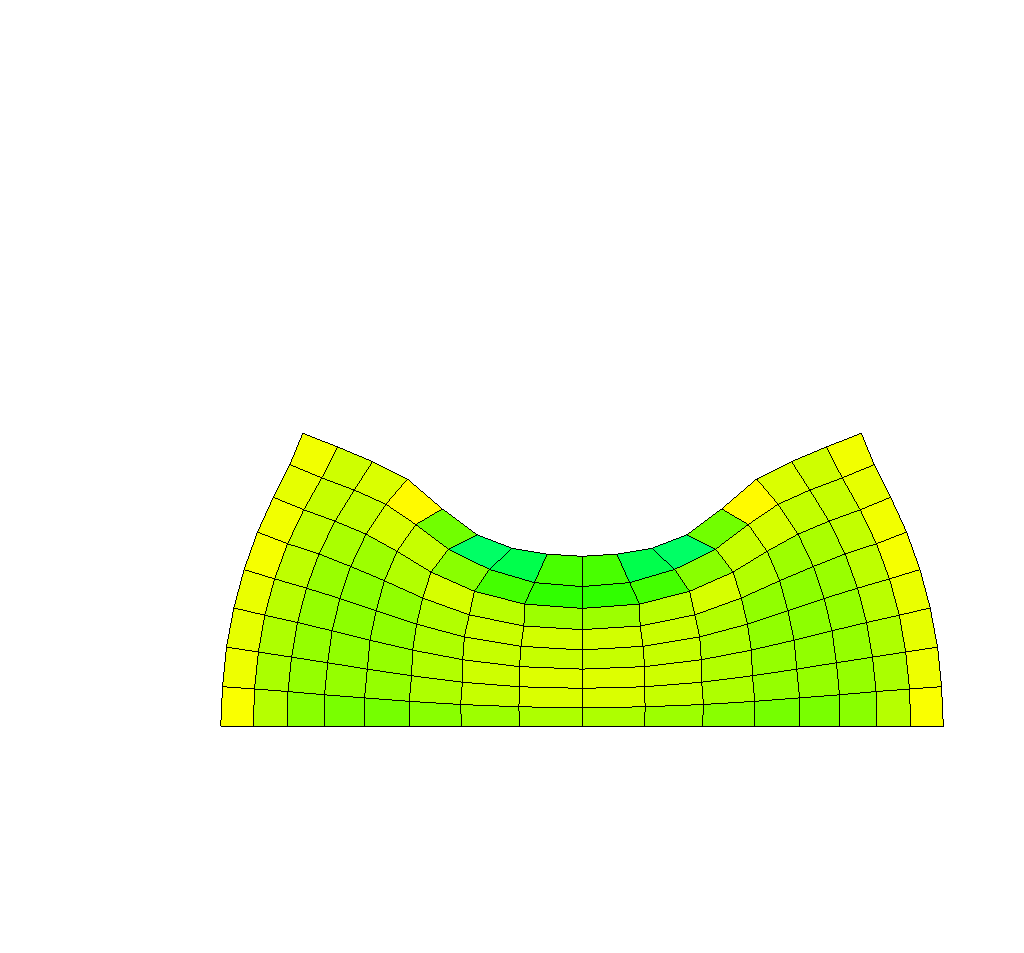}} &
\subcaptionbox{\label{sfig:testb}} {\includegraphics[width=.45\linewidth, trim={200 200 70 400}, clip]{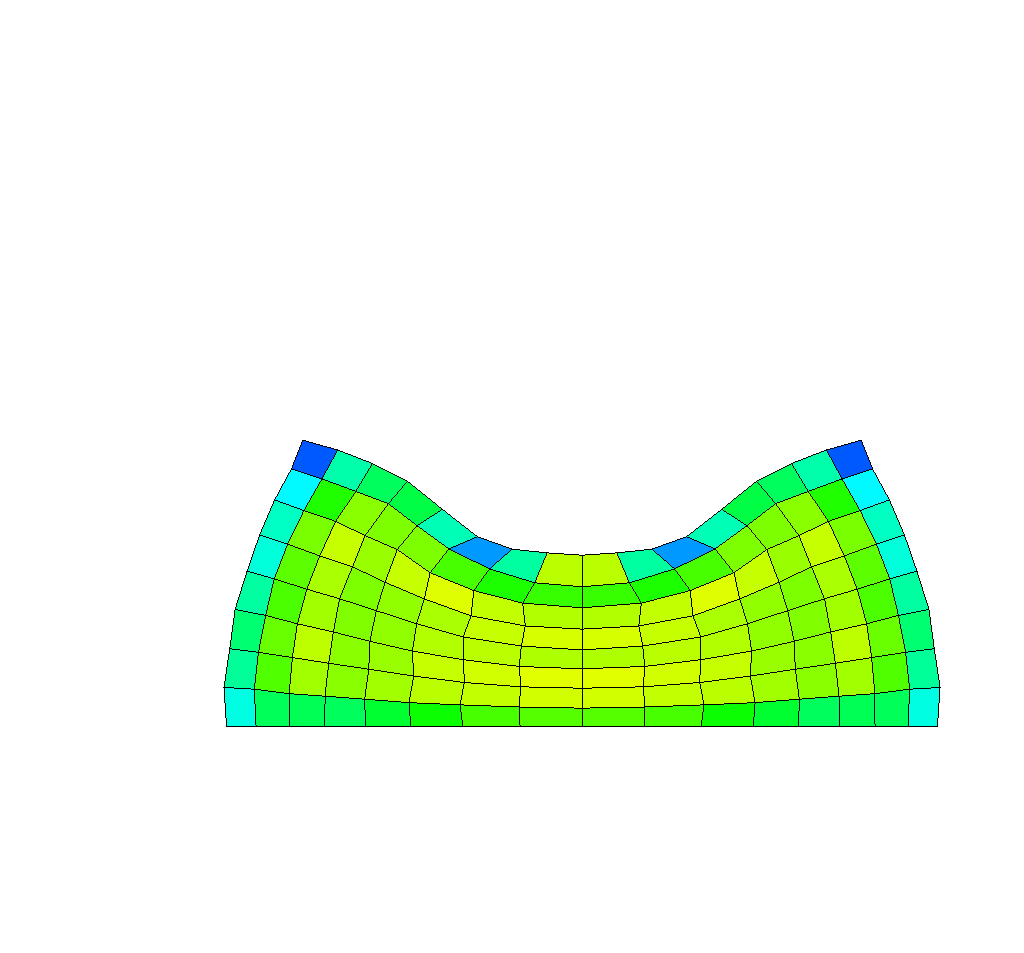}}\\

\rotatebox{90}{\qquad\qquad\textbf{$\nus = -1$}} &
\subcaptionbox{\label{sfig:testc}}{\includegraphics[width=.45\linewidth, trim={200 200 70 350}, clip]{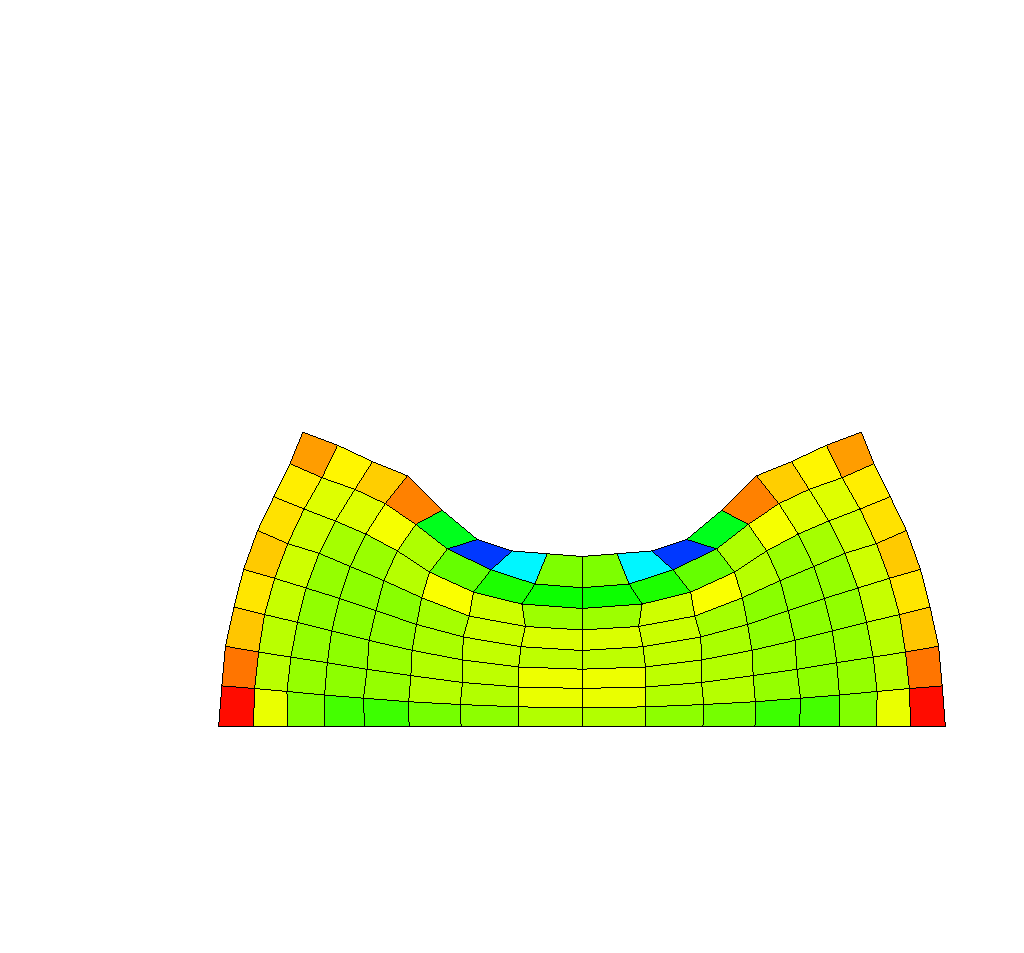}} &
\subcaptionbox{\label{sfig:testd}}{\includegraphics[width=.45\linewidth, trim={200 200 70 350}, clip]{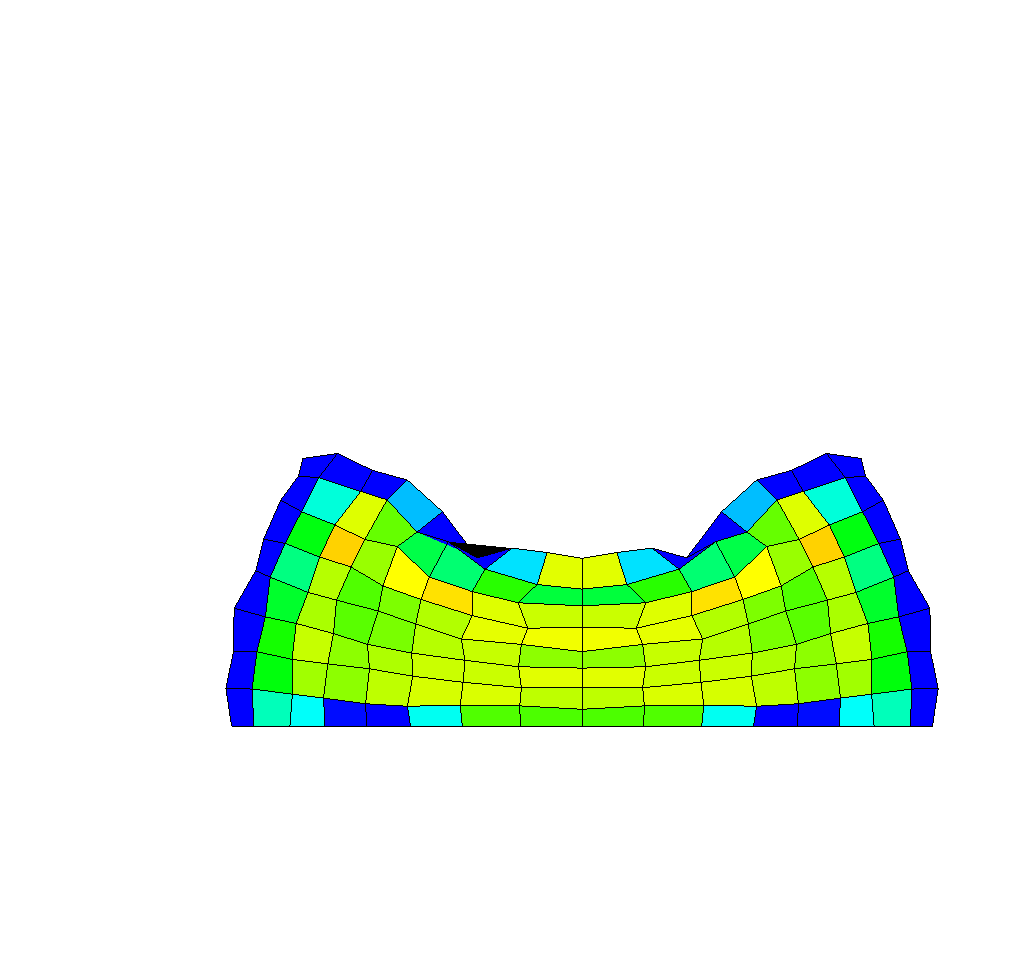}}  \\
\end{tabular}
%trim={left bottom right top}
\begin{centering}
Avg $J$ \\
\includegraphics[width=2.5in, trim={0 5in 0 5in}, clip]{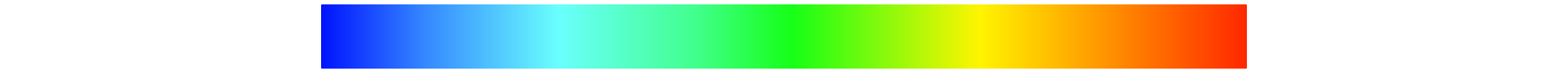}  \\
%0.40 \ \ \ \ \ \ \ \ \ \ \ \ \ \ \ \ 1.10
0.90 $\qquad\qquad\qquad\qquad$ 1.05

\end{centering}
\caption{Deformations of the compressed block benchmark (Section \ref{Compression Test}), along with mean values of $J$ within each element calculated using equation \eqref{avgJ}, using a neo-Hookean material model, equations (\ref{nh_energy}) -- (\ref{nh_stress_mod}), with $G = 80.194 \ \frac{\text{dyn}}{\text{cm}^2}$. The background Eulerian grid is not shown. Shown here are solid deformations computed using \textbf{Q1} elements and $m = 153$ solid degrees of freedom. The first row shows cases with $\nus = .4$, and the second row shows cases with $\nus = -1$ (here equivalent to $\kappas = 0$ and no volumetric-based stabilization). The first column depicts cases with modified invariants, and the second column depicts cases with unmodified invariants. Notice that the case with modified invariants with nonzero numerical bulk modulus has the smoothest deformations, whereas those of the case with unmodified invariants and zero numerical bulk modulus behave unphysically.}

\label{cb}
\end{figure}

\begin{figure}
\begin{tabular}{l c c c}
&$\cauchys_{00}$ & $\cauchys_{01}$& $\cauchys_{11}$ \\

\rotatebox{90}{$\quad\;\;$\textbf{m = 561} }&
\subcaptionbox{\label{sfig:testa}}{\includegraphics[width=.3\linewidth]{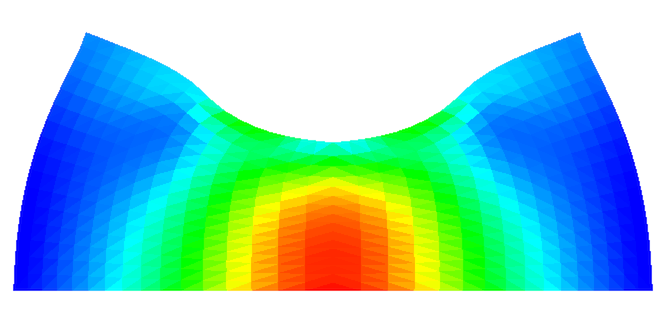}}&
\subcaptionbox{\label{sfig:testb}}{\includegraphics[width=.3\linewidth]{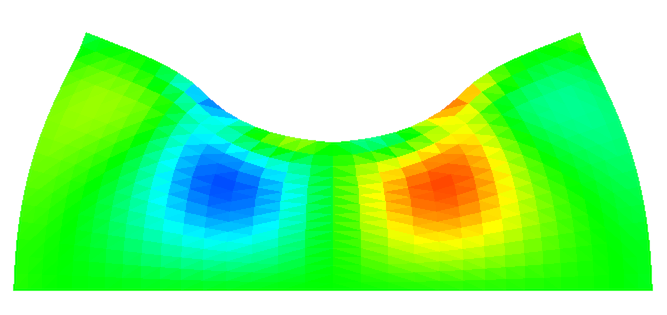}}&
\subcaptionbox{\label{sfig:testa}}{\includegraphics[width=.3\linewidth]{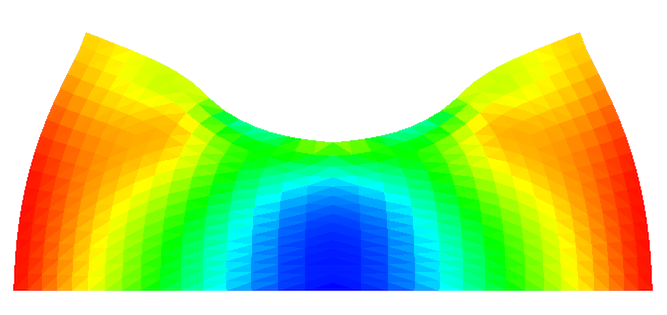}}\\

\rotatebox{90}{$\quad\;\;$ \textbf{m = 2145} }&
\subcaptionbox{\label{sfig:testa}}{\includegraphics[width=.3\linewidth]{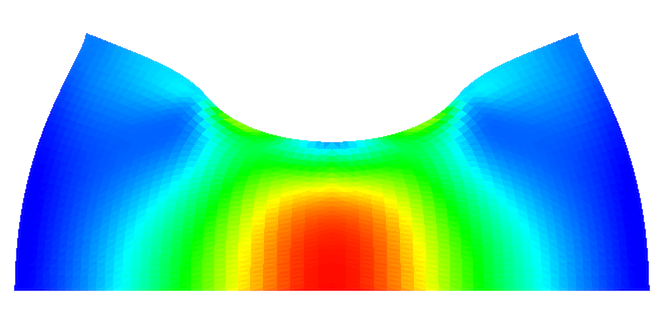}}&
\subcaptionbox{\label{sfig:testb}}{\includegraphics[width=.3\linewidth]{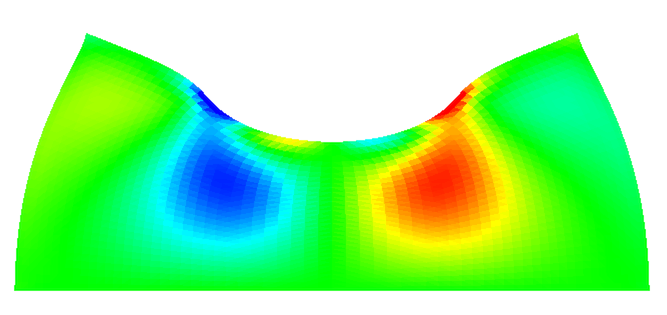}}&
\subcaptionbox{\label{sfig:testa}}{\includegraphics[width=.3\linewidth]{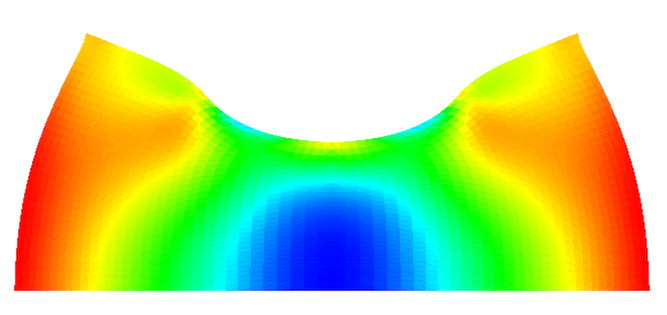}}\\

\rotatebox{90}{$\quad\;\;$ \textbf{m = 4753} }&
\subcaptionbox{\label{sfig:testa}}{\includegraphics[width=.3\linewidth]{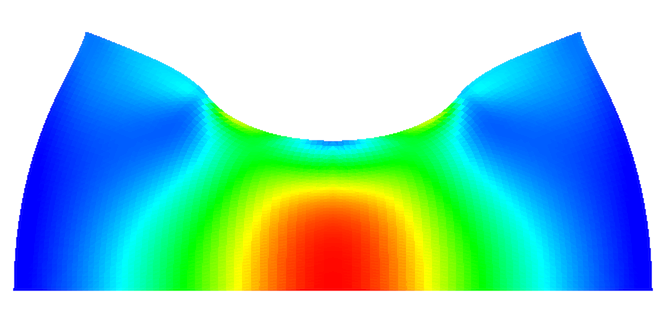}}&
\subcaptionbox{\label{sfig:testb}}{\includegraphics[width=.3\linewidth]{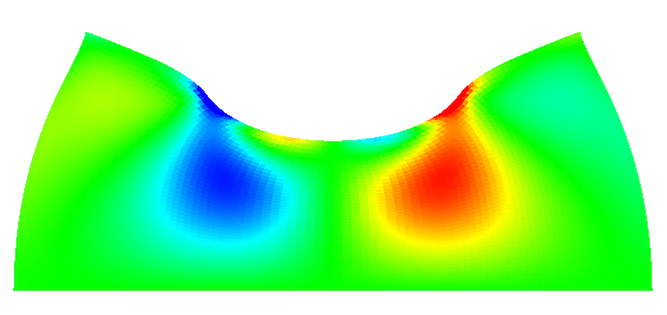}}&
\subcaptionbox{\label{sfig:testa}}{\includegraphics[width=.3\linewidth]{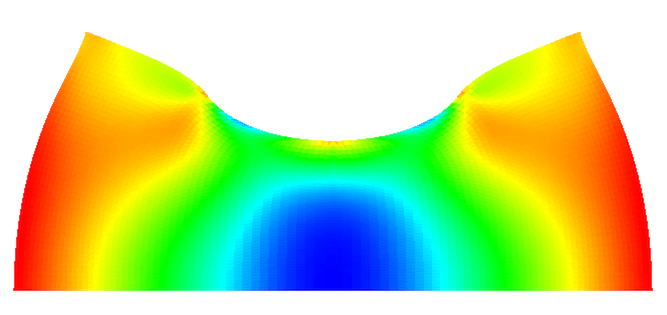}}\\

\\
\hline
\rotatebox{90}{$\quad\;\;$\textbf{FE (P1/P1)} }&
\subcaptionbox{\label{sfig:testa}}{\includegraphics[width=.3\linewidth]{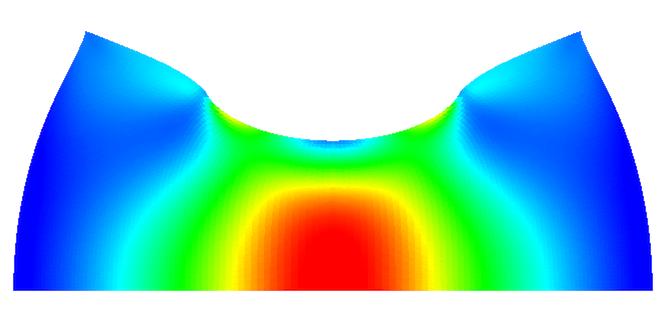}}&
\subcaptionbox{\label{sfig:testb}}{\includegraphics[width=.3\linewidth]{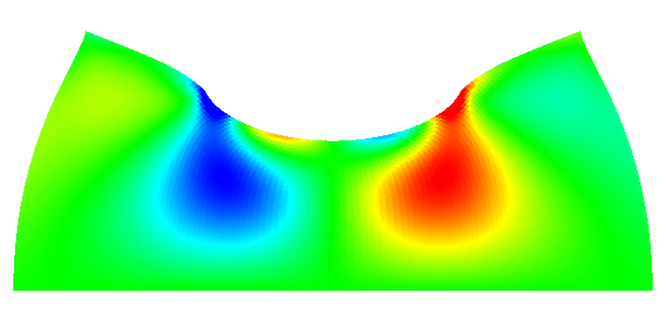}}&
\subcaptionbox{\label{sfig:testa}}{\includegraphics[width=.3\linewidth]{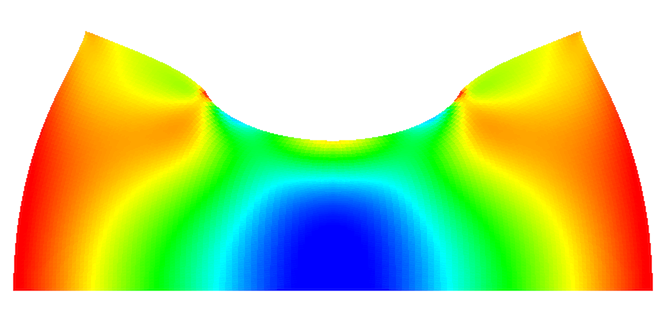}}\\

&
\includegraphics[width=.3\linewidth, trim={0 5in 0 5in}, clip]{color_bar.pdf}&
\includegraphics[width=.3\linewidth, trim={0 5in 0 5in}, clip]{color_bar.pdf}&
\includegraphics[width=.3\linewidth, trim={0 5in 0 5in}, clip]{color_bar.pdf}  \\
&-15 $\qquad\qquad\qquad$ 145& -45 $\qquad\qquad\qquad$ 45 & -100 $\qquad\qquad\qquad$ 20 \\
\end{tabular}
\caption{The three components of the deviatoric part of $\cauchys$ for the compressed block benchmark (Section \ref{Compression Test}). The IBFE method uses modified invariants and volumetric stabilization ($\nus = 0.4$). Each row is labeled with the solid degrees of freedom (DOF), and the bottom row depicts the FE (\textbf{P1/P1}) solution with $m = 8385$ solid DOF (higher than that of the highest resolution IBFE results presented in this figure). We use \textbf{P1} elements for each method. The results from the IBFE formulation are clearly converging to the higher-resolution FE solution.}
\label{cb_dev}
\end{figure}

\begin{figure}
\begin{tabular}{l c c c}
&\textbf{FE} \textbf{(P1}/\textbf{P1)} & \textbf{Modified Invariants \&} & \textbf{Unmodified Invariants \&} \\
&& \textbf{Stabilized} & \textbf{Unstabilized} \\
\rotatebox{90}{$\quad\;$ \textbf{Coarse} }&
\subcaptionbox{\label{sfig:testa}}{\includegraphics[width=.3\linewidth]{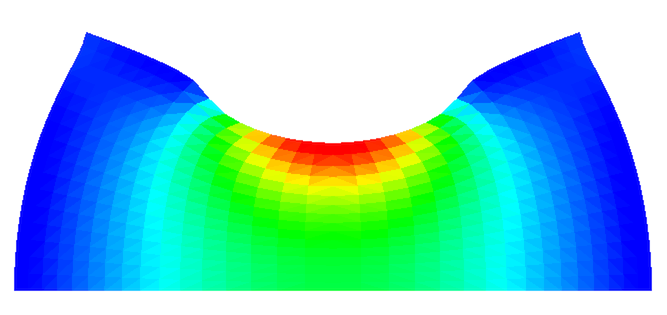}}&
\subcaptionbox{\label{sfig:testb}}{\includegraphics[width=.3\linewidth]{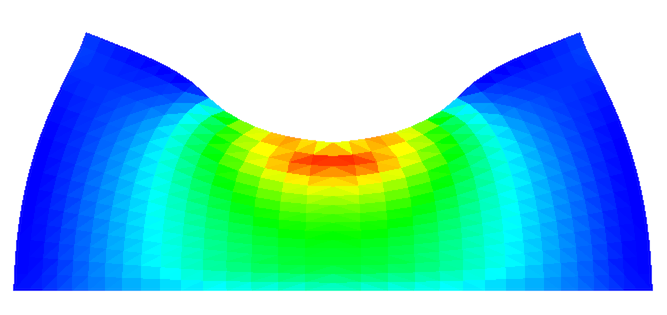}}&
\subcaptionbox{\label{sfig:testa}}{\includegraphics[width=.3\linewidth]{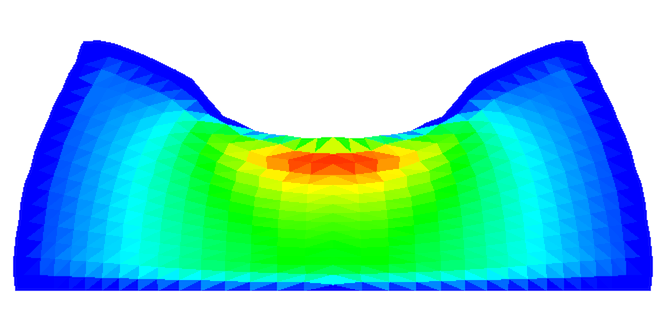}}\\

\rotatebox{90}{$\quad\;$ \textbf{Fine} }&
\subcaptionbox{\label{sfig:testa}}{\includegraphics[width=.3\linewidth]{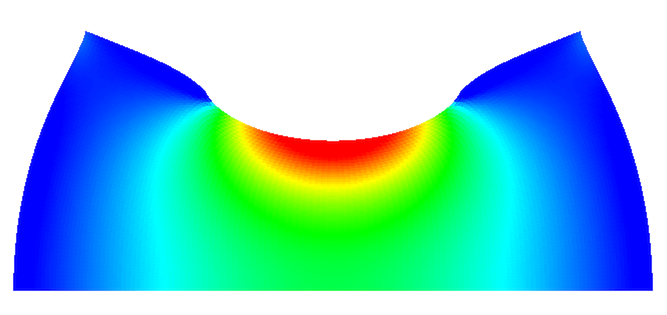}}&
\subcaptionbox{\label{sfig:testb}}{\includegraphics[width=.3\linewidth]{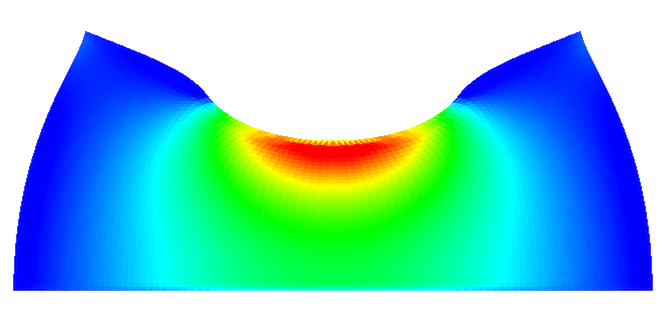}}&
\subcaptionbox{\label{sfig:testa}}{\includegraphics[width=.3\linewidth]{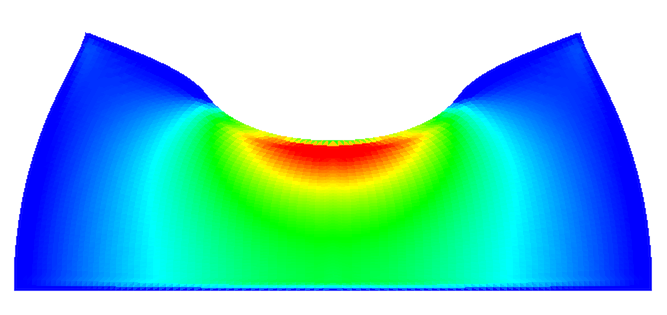}}\\
\end{tabular}
\vskip .5cm
\centering
$p \ \left(\frac{\text{dyn}}{\text{cm}^2}\right)$ \\
\includegraphics[width=2.5in, trim={0 5in 0 5in}, clip]{color_bar.pdf}  \\
-50 $\qquad\qquad\qquad\qquad$ 135
%\end{centering}
\caption{The pressure distributions of the compressed block benchmark (Section \ref{Compression Test}) for the FE (\textbf{P1}/\textbf{P1}) solution (left) and for IBFE solutions corresponding to the modified case with stabilization (middle) and the unmodified case without stabilization (right). We show a coarse and fine discretization and use \textbf{P1} elements for each method. Results from all formulations appear to be converging to the same solution. The IBFE method uses $m = 1089$ DOF and $m = 9409$ DOF for the coarse and fine cases, respectively. The FE method uses $m = 1089$ DOF and $m = 16,641$ DOF for the coarse and fine cases, respectively.}
\label{cb_p}
\end{figure}

\begin{figure}
$\qquad\qquad\qquad\;\;\;\;$ \textbf{P1} $\qquad\qquad\qquad\qquad\quad$  \textbf{Q1} $\qquad\qquad\qquad\qquad\;\;\;$  \textbf{P2} $\qquad\qquad\qquad\qquad\quad\;$ \textbf{Q2}\\
\rotatebox{90}{$\quad$ \textbf{$\nus = .49995$} }
   \rotatebox{90}{$\quad\;\;$ Disp. (cm) }
\includegraphics[width=.225\linewidth, trim={40 190 25 200}, clip]{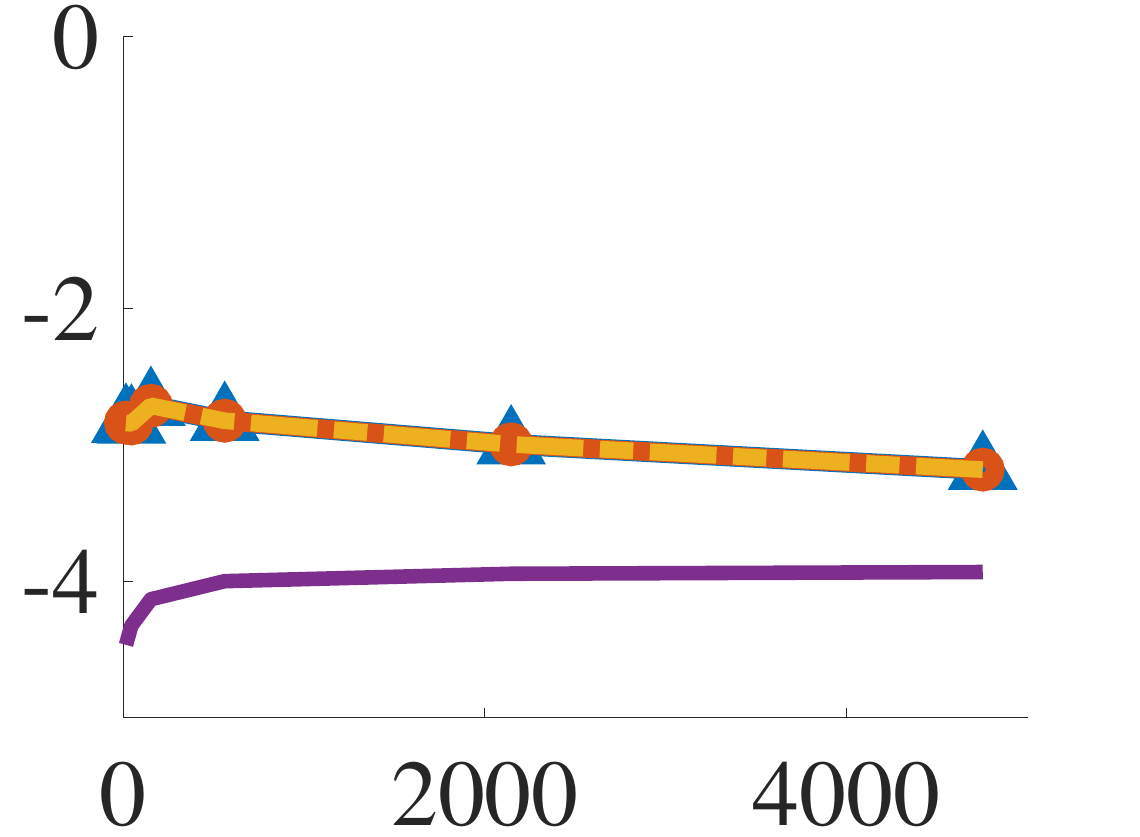} 
\includegraphics[width=.225\linewidth, trim={40 190 25 200}, clip]{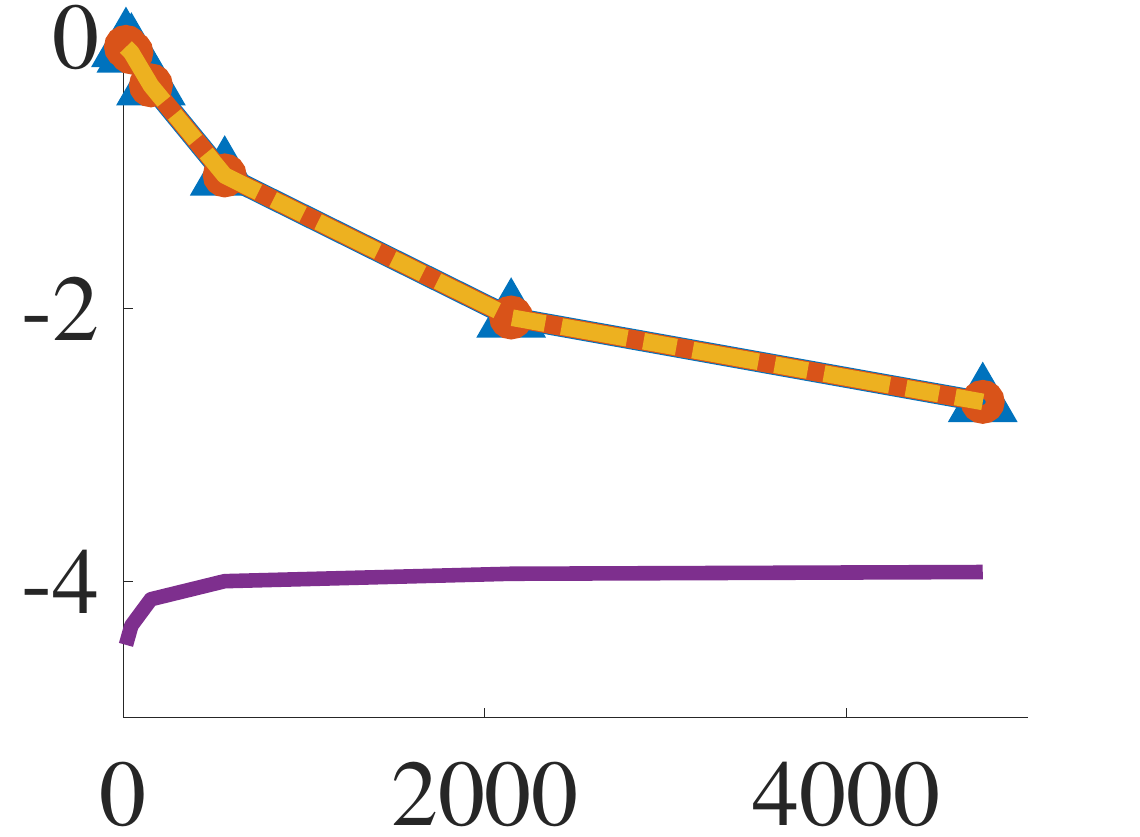}
\includegraphics[width=.225\linewidth, trim={40 190 25 200}, clip]{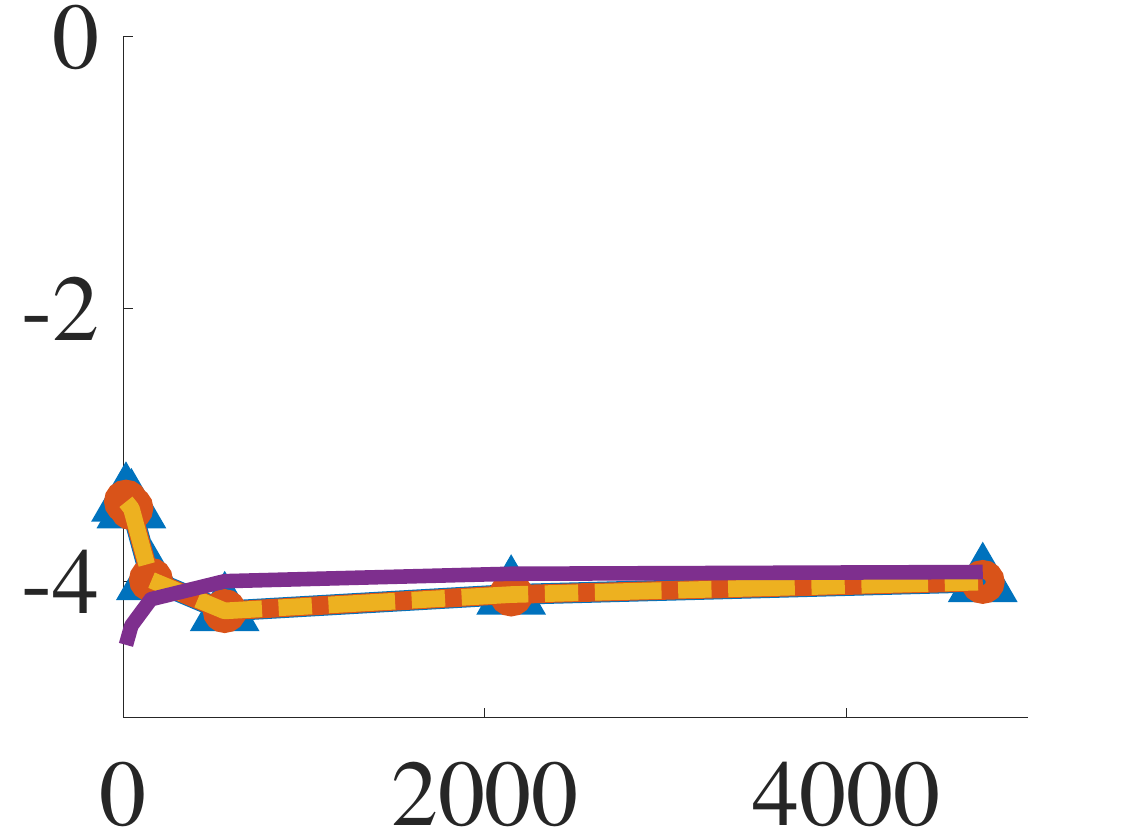}
\includegraphics[width=.225\linewidth, trim={40 190 25 200}, clip]{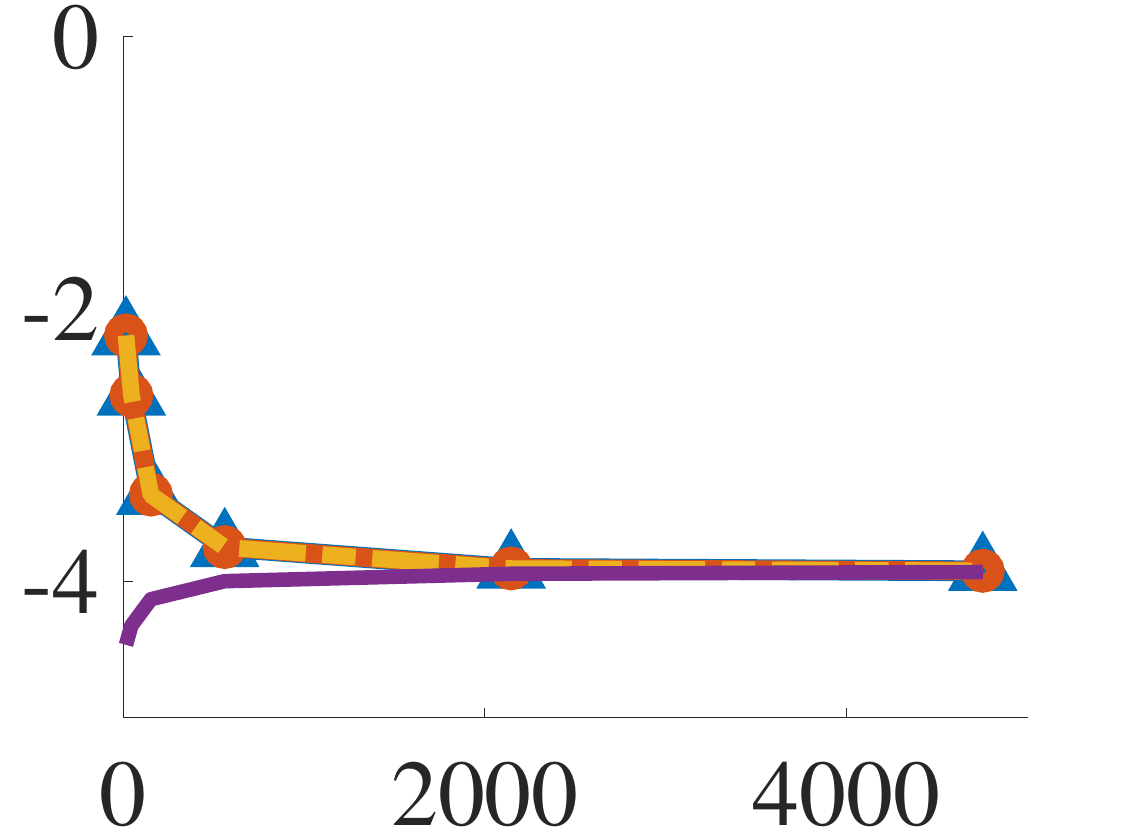} \\

\rotatebox{90}{$\qquad$ \textbf{$\nus = .4$} }
   \rotatebox{90}{$\quad\;\;$ Disp. (cm) }
\includegraphics[width=.225\linewidth, trim={40 190 25 200}, clip]{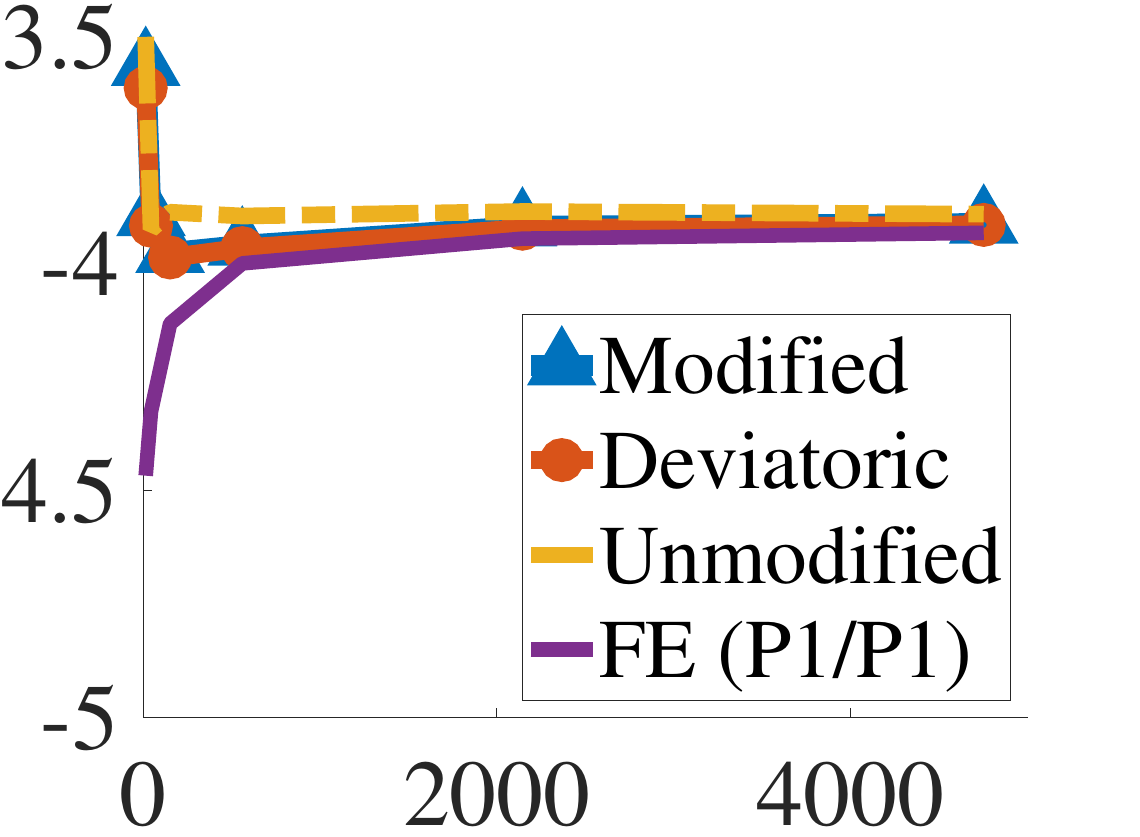} 
\includegraphics[width=.225\linewidth, trim={40 190 25 200}, clip]{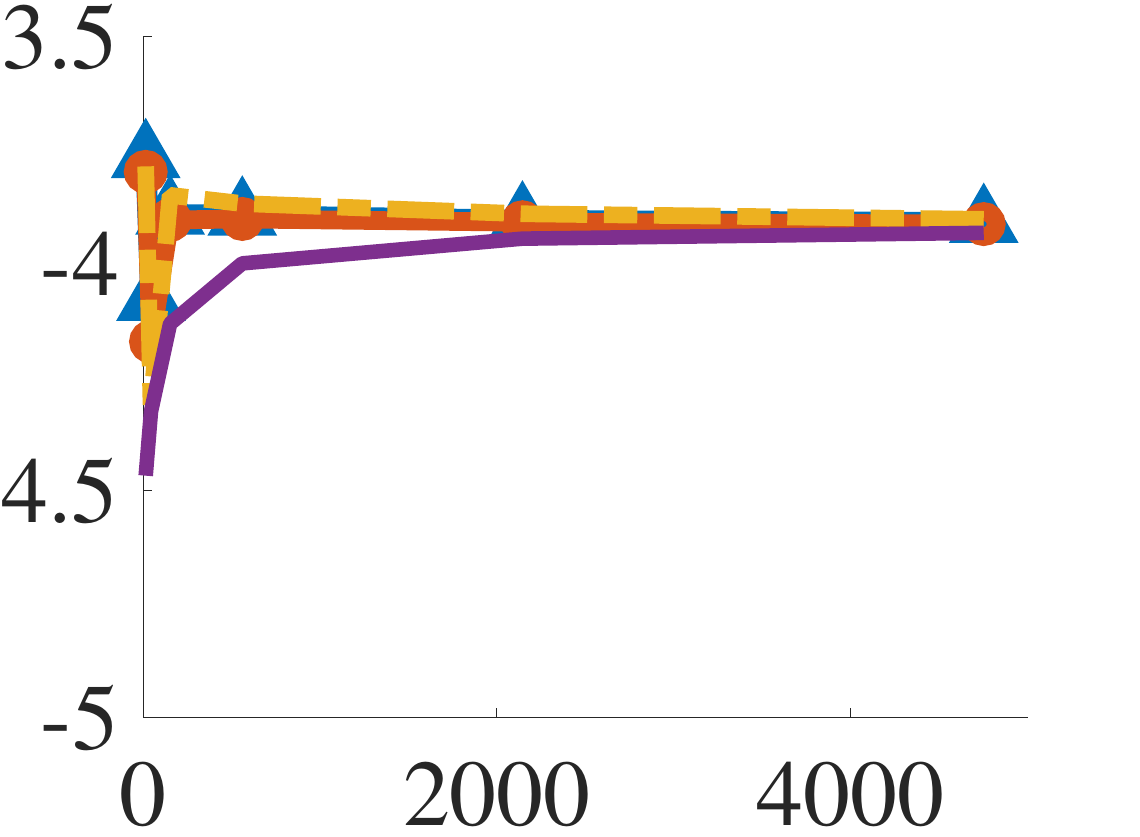}
\includegraphics[width=.225\linewidth, trim={40 190 25 200}, clip]{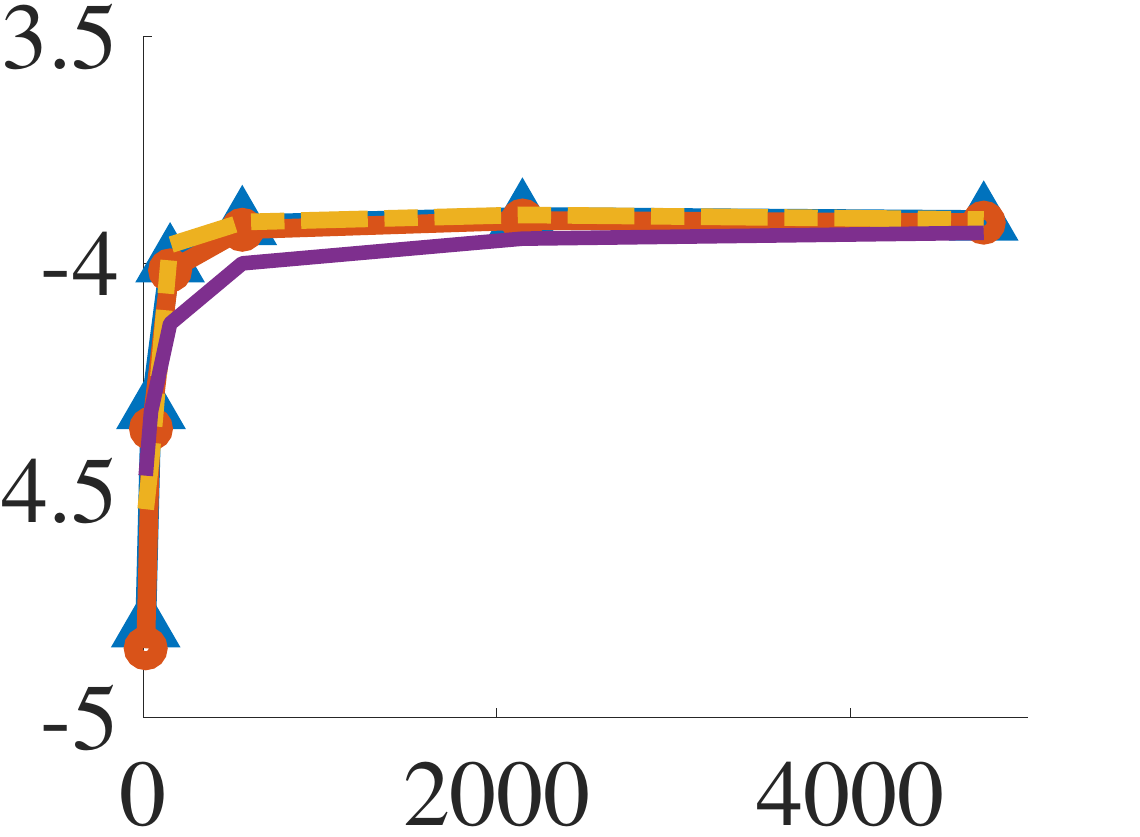}
\includegraphics[width=.225\linewidth, trim={40 190 25 200}, clip]{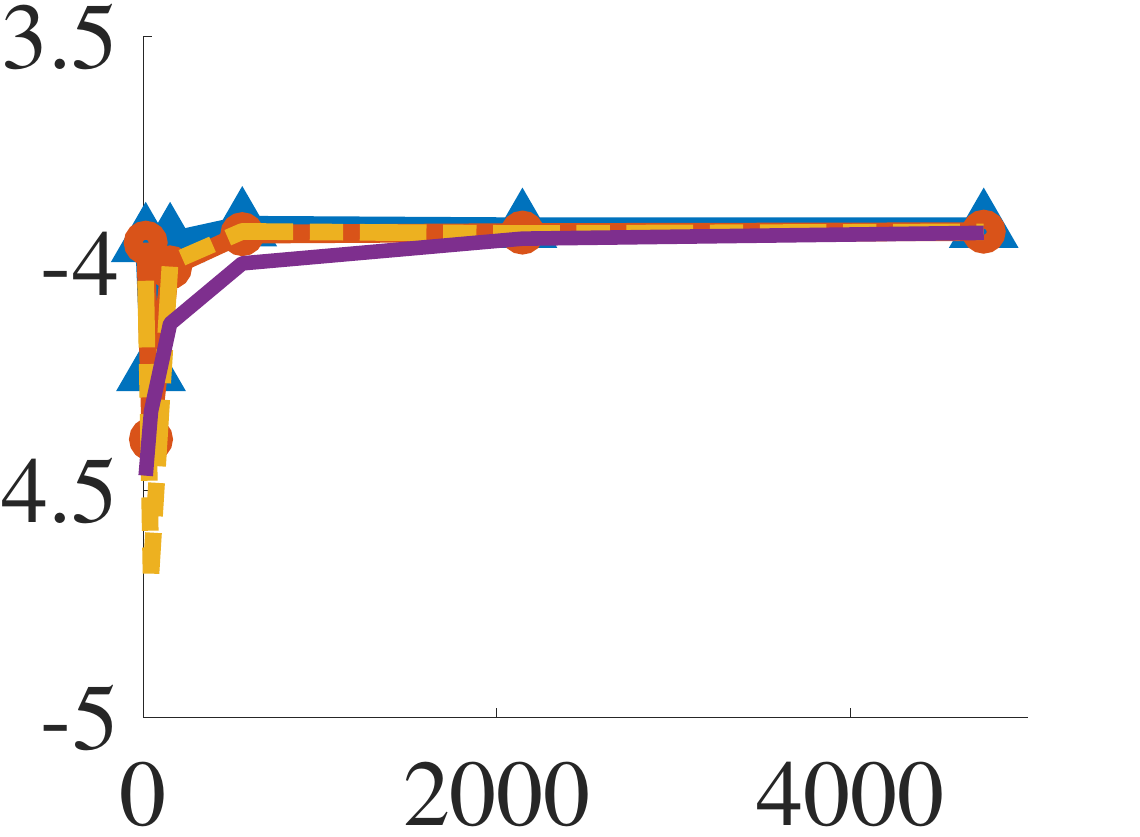}\\

\rotatebox{90}{$\qquad$ \textbf{$\nus = 0$} }
   \rotatebox{90}{$\quad\;\;$ Disp. (cm) }
\includegraphics[width=.225\linewidth, trim={40 190 25 200}, clip]{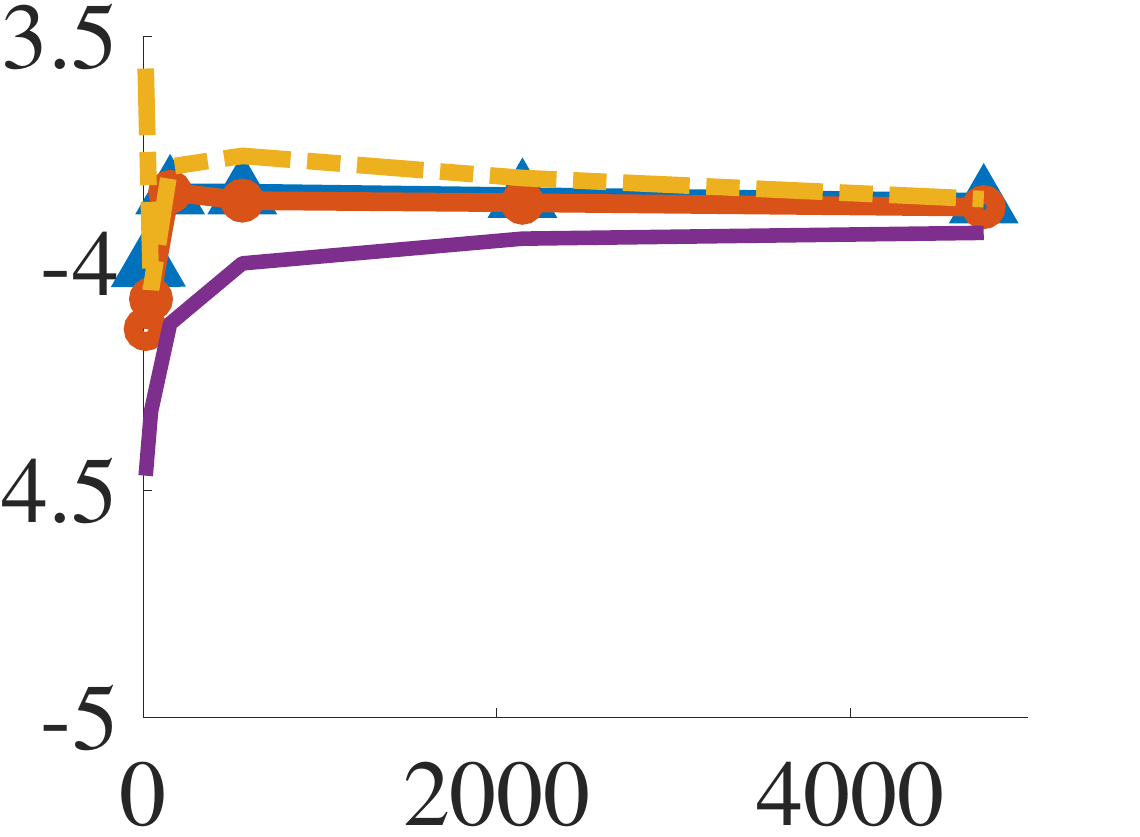}
\includegraphics[width=.225\linewidth, trim={40 190 25 200}, clip]{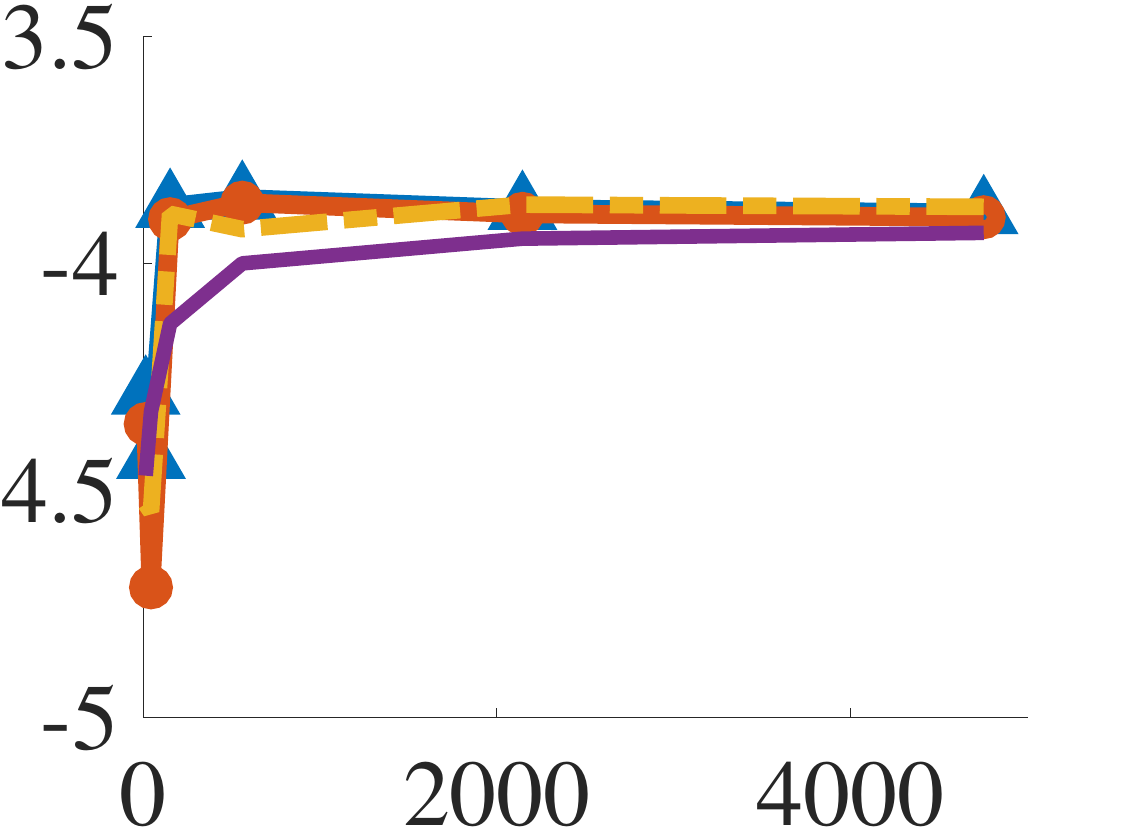}
\includegraphics[width=.225\linewidth, trim={40 190 25 200}, clip]{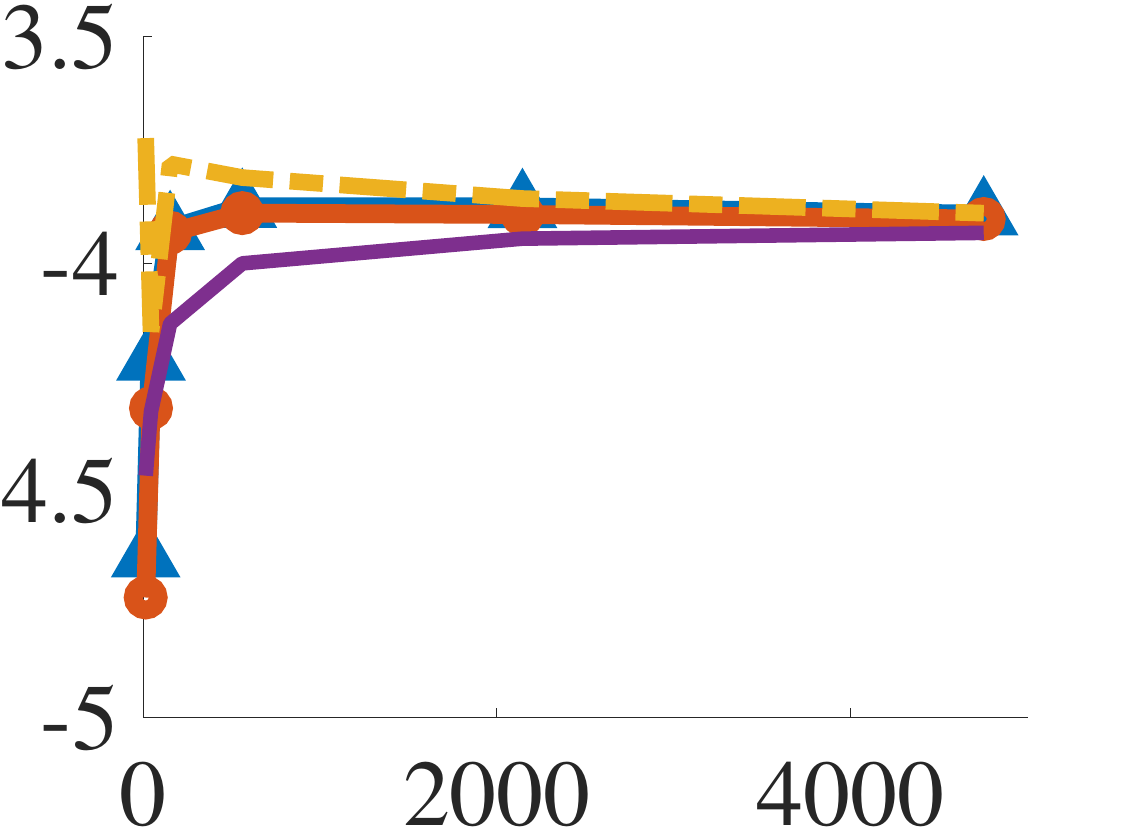}
\includegraphics[width=.225\linewidth, trim={40 190 25 200}, clip]{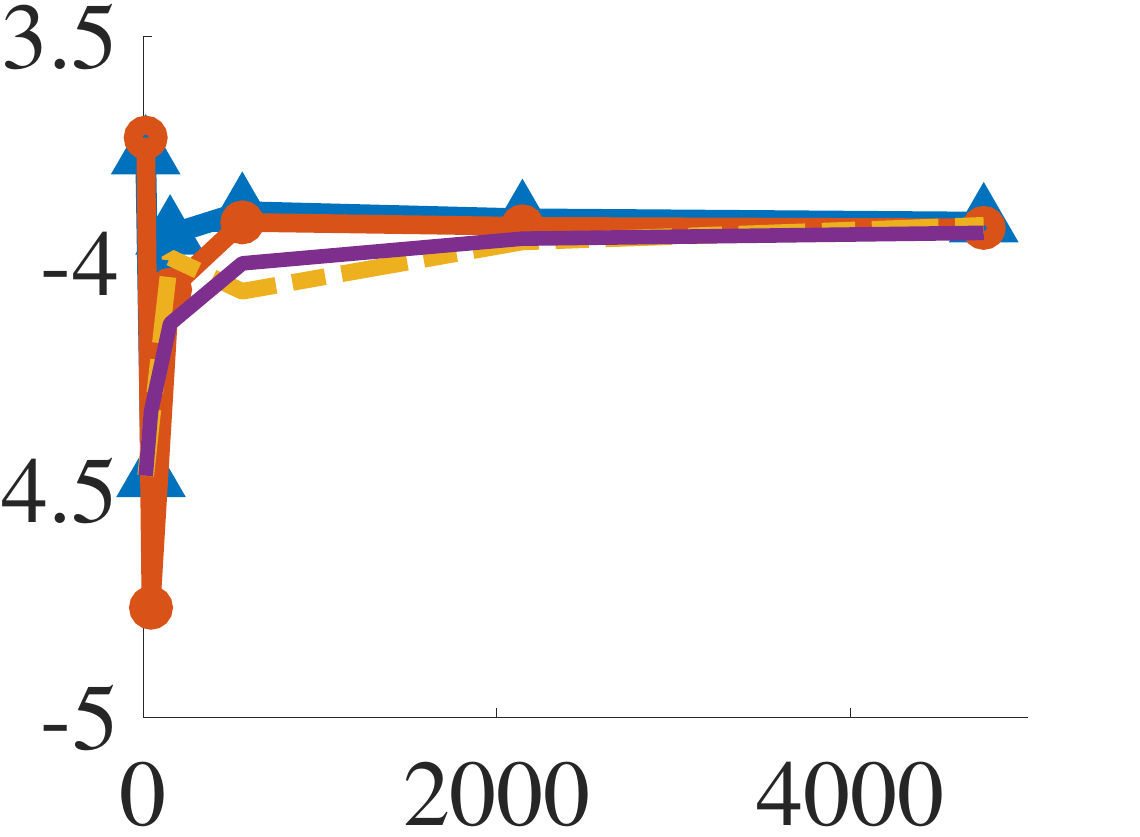}\\

\rotatebox{90}{$\quad\;\;\;$ \textbf{$\nus = -1$} }
   \rotatebox{90}{$\quad\;\;$ Disp. (cm) }
\includegraphics[width=.225\linewidth, trim={40 190 25 200}, clip]{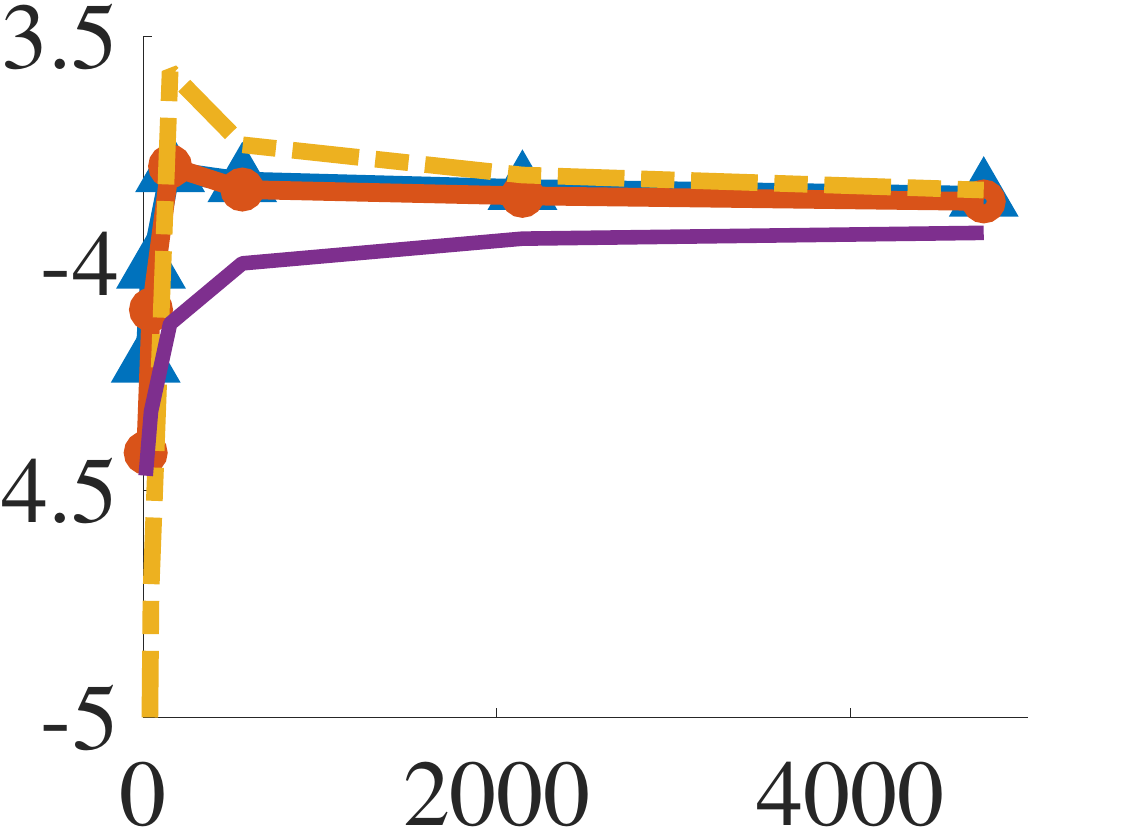}
\includegraphics[width=.225\linewidth, trim={40 190 25 200}, clip]{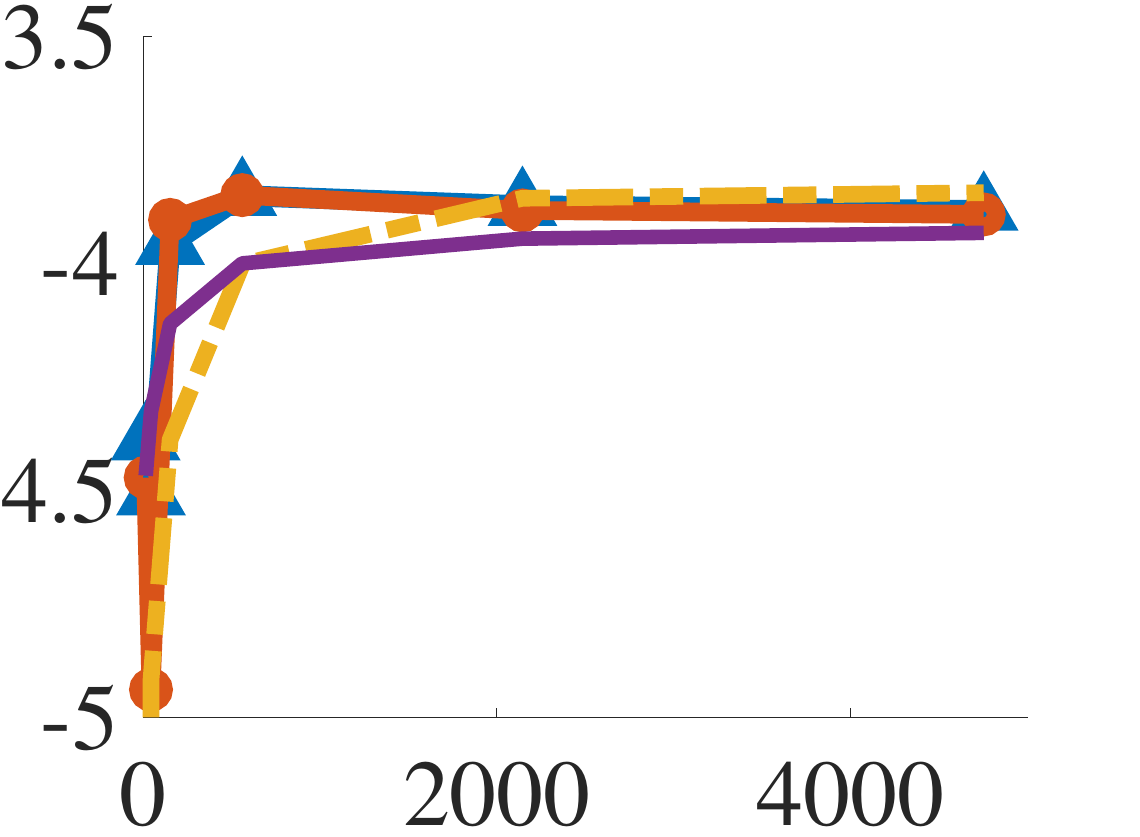}
\includegraphics[width=.225\linewidth, trim={40 190 25 200}, clip]{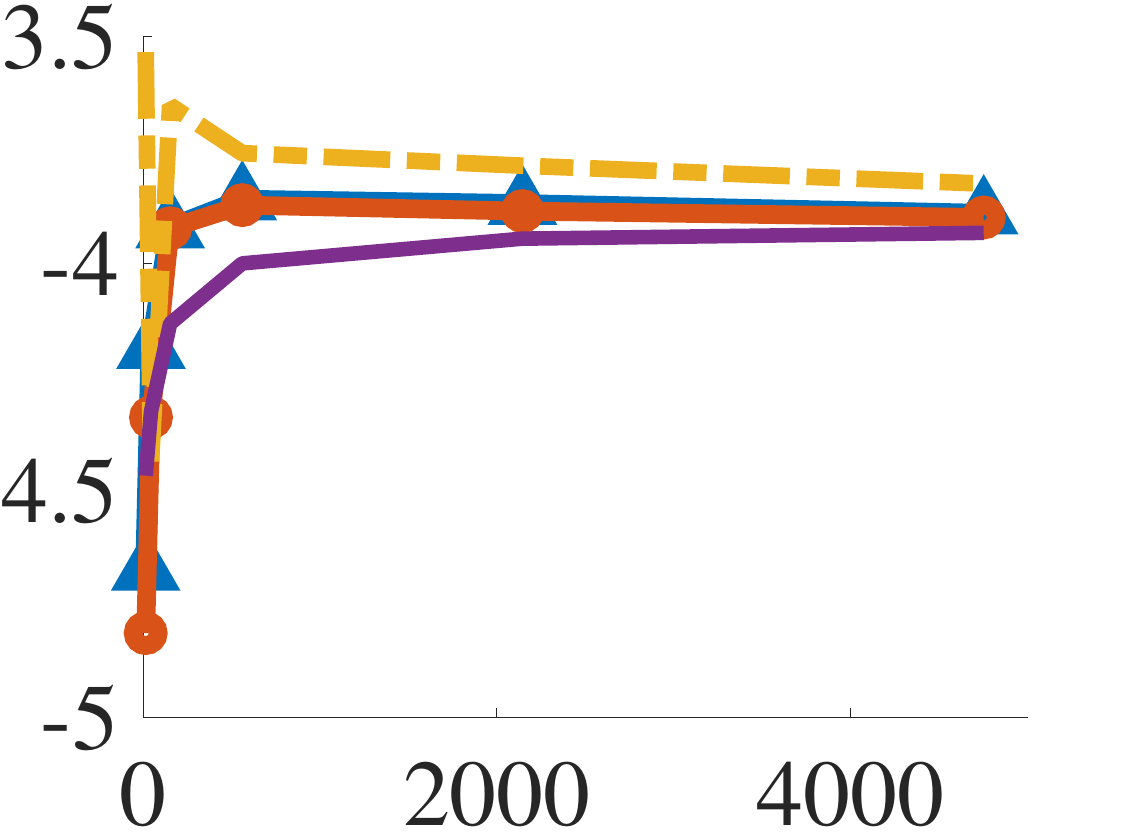}
\includegraphics[width=.225\linewidth, trim={40 190 25 200}, clip]{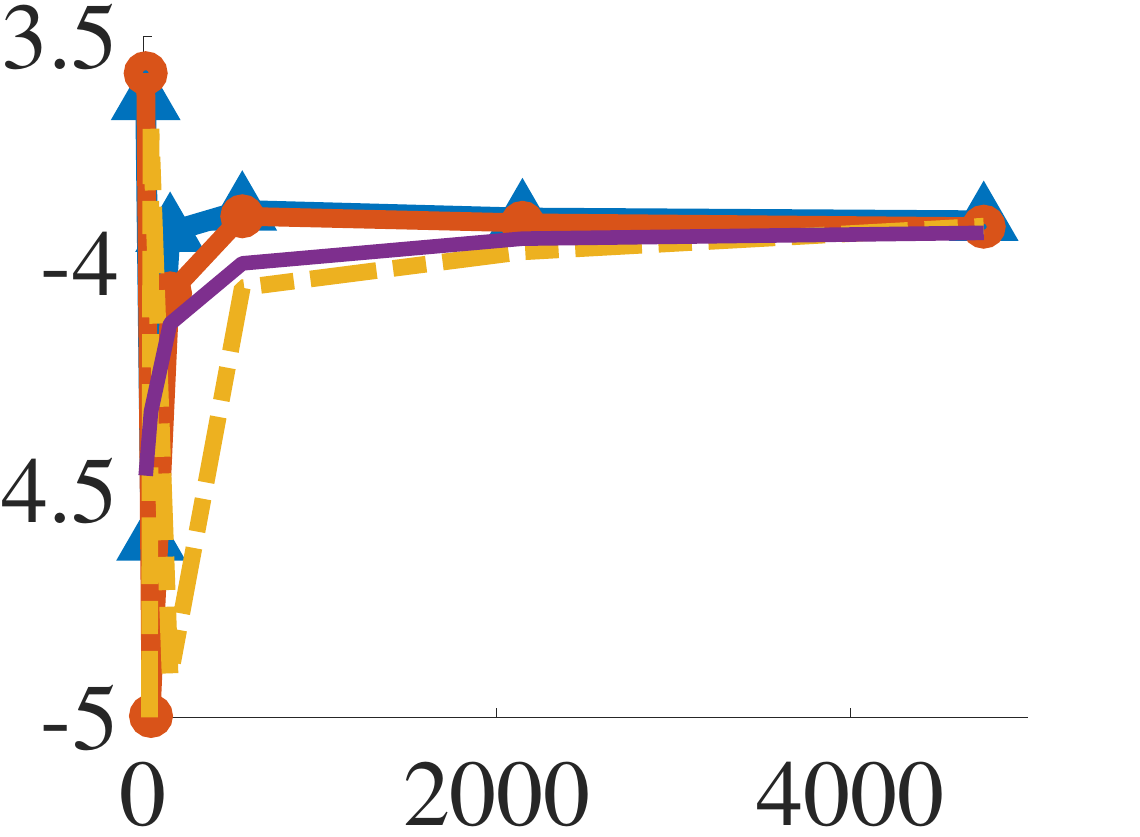}\\

$\qquad\qquad\quad$ \# Solid DOF $\qquad\qquad\quad\;$ \# Solid DOF $\qquad\qquad\quad$ \# Solid DOF $\qquad\qquad\quad\;$ \# Solid DOF
\caption{Displacement of the center point in Figure (\ref{comp_mesh}) for the compressed block benchmark (Section \ref{Compression Test}) for different choices of elements and numerical Poisson ratio. The solid degrees of freedom (DOF) range from $m = 15$ to $4753$. Notice that each row has the same extents. If a value of $\nus$ is close to $\frac{1}{2}$, low order elements produce volumetric locking, and higher order elements are needed for convergence at reasonable numbers of DOF.}
\label{comp_disp}
\end{figure}

\begin{figure}
\begin{tabular}{l c c}
& \textbf{$m=153$} & \textbf{$m=2145$} \\
\rotatebox{90}{$\qquad\qquad\qquad$ Disp. (cm) }&
\includegraphics[width=.45\linewidth, trim={50 190 25 200}, clip]{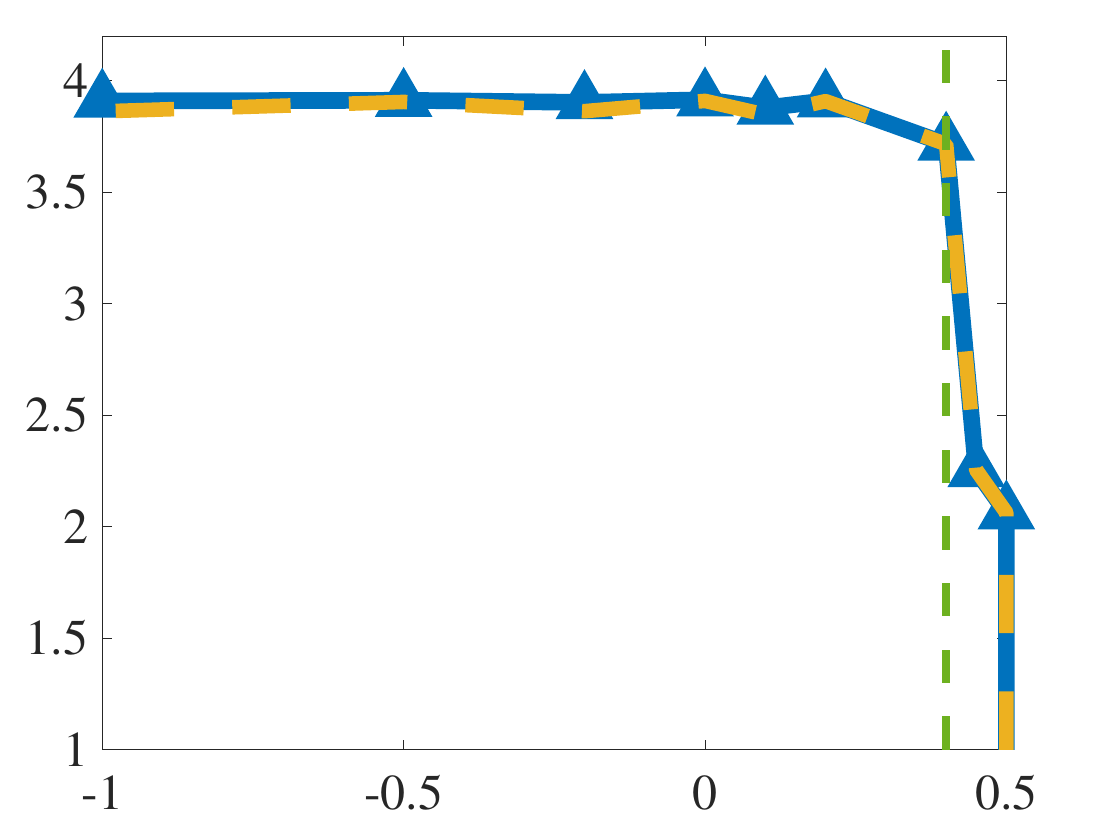} &
\includegraphics[width=.45\linewidth, trim={50 190 25 200}, clip]{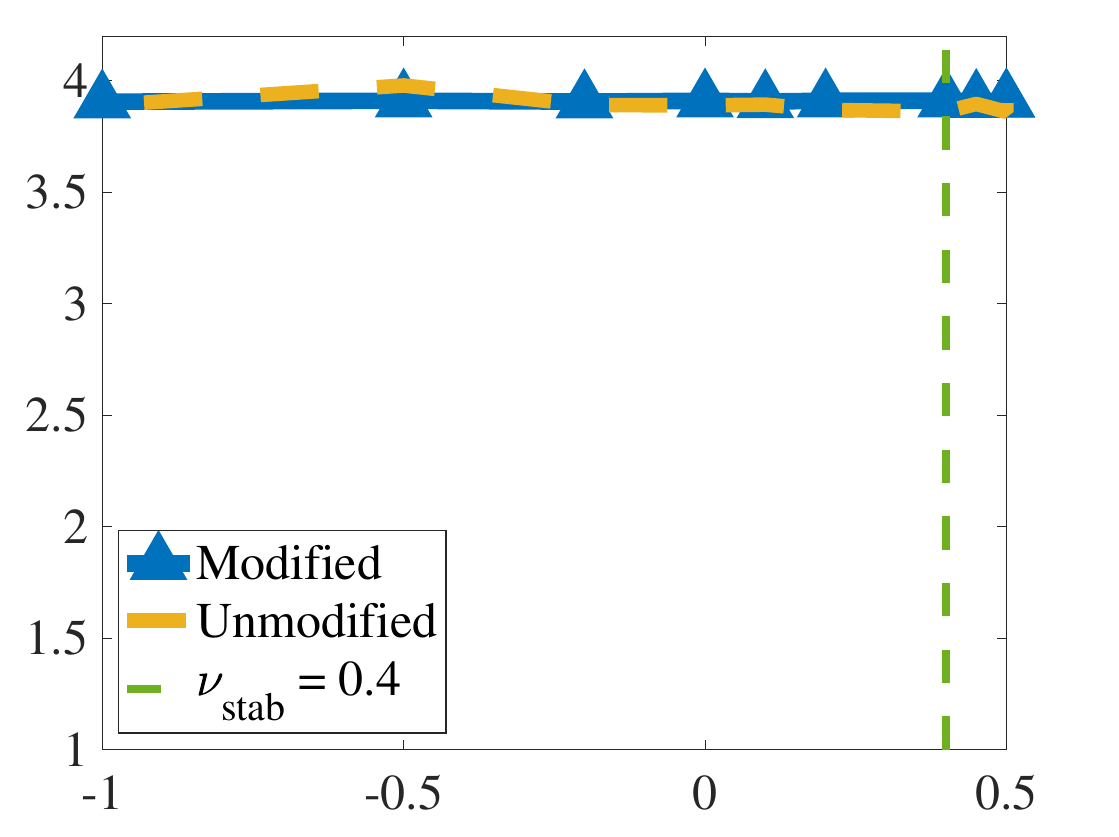}\\

& \textbf{$\nus$} & \textbf{$\nus$}
\end{tabular}
\caption{Displacement of the encircled point, shown in Figure (\ref{comp_mesh}), plotted against different values of $\nus$ for the compressed block benchmark (Section \ref{Compression Test}). Results indicate that for a discretization with $m = 2145$ DOF, the displacement is insensitive to larger values of $\nus$. For the coarser case of $m = 153$ DOF, volumetric locking begins to pollute the solution just past $\nus = 0.4$. Results are for \textbf{Q1} elements. }
%trim={left bottom right top}
\label{disp_v_nu}
\end{figure}

\begin{figure}
$\qquad\qquad\qquad\;\;\;\;$ \textbf{P1} $\qquad\qquad\qquad\qquad\quad$  \textbf{Q1} $\qquad\qquad\qquad\qquad\;\;\;$  \textbf{P2} $\qquad\qquad\qquad\qquad\quad\;$ \textbf{Q2}\\
\rotatebox{90}{$\quad$ \textbf{$\nus = .49995$} }
   \rotatebox{90}{$\quad$ Vol Change \% }
\includegraphics[width=.225\linewidth, trim={30 190 25 200}, clip]{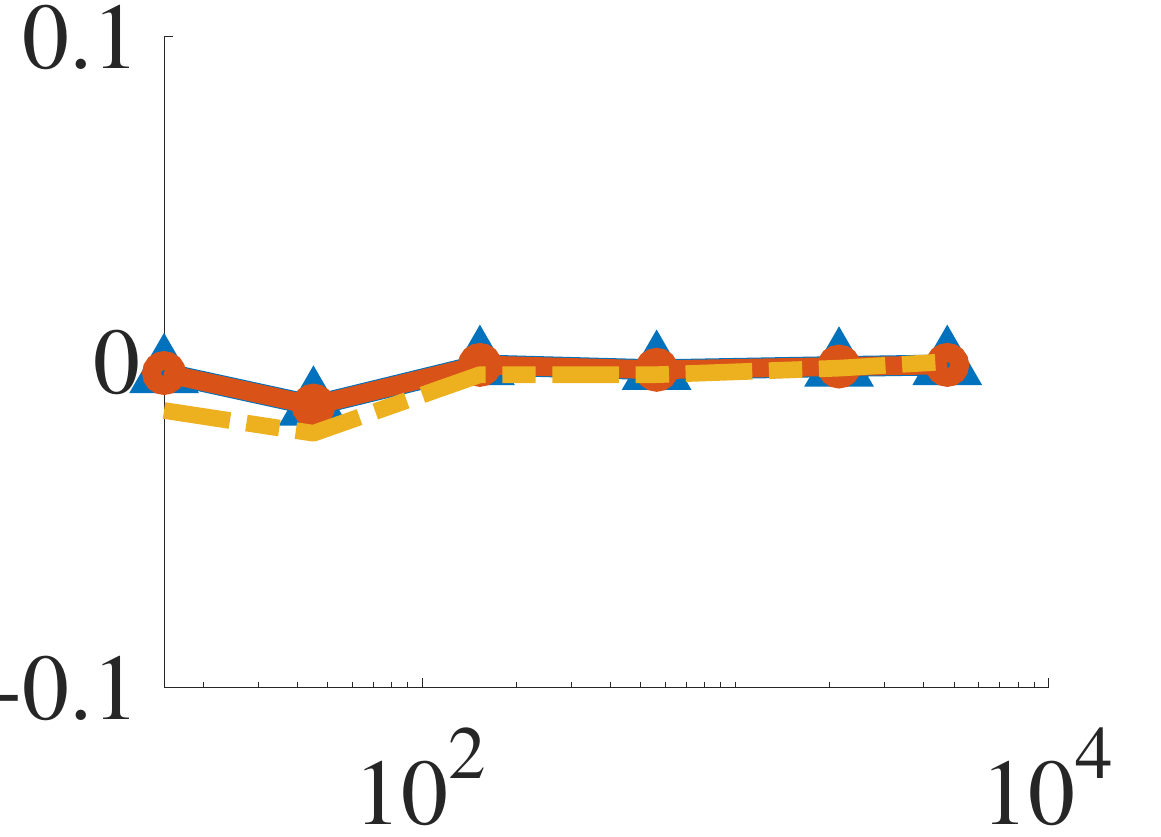} 
\includegraphics[width=.225\linewidth, trim={30 190 25 200}, clip]{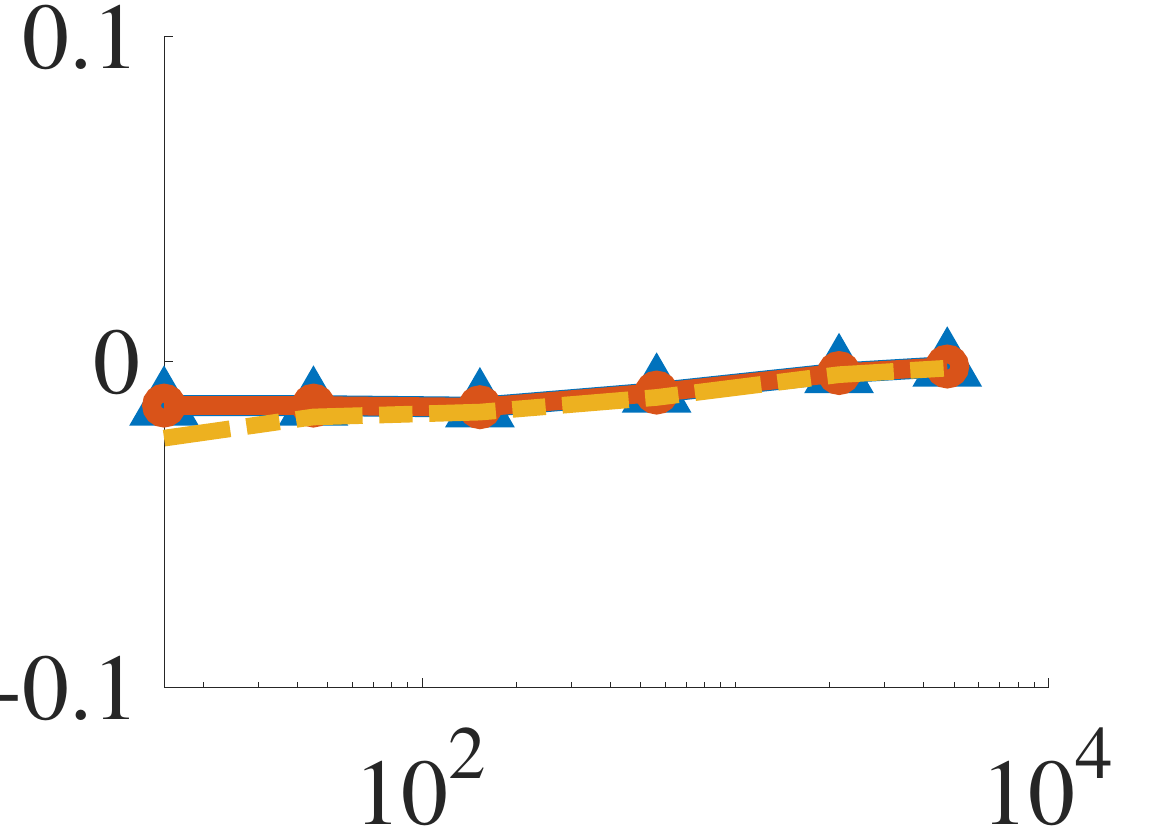} 
\includegraphics[width=.225\linewidth, trim={30 190 25 200}, clip]{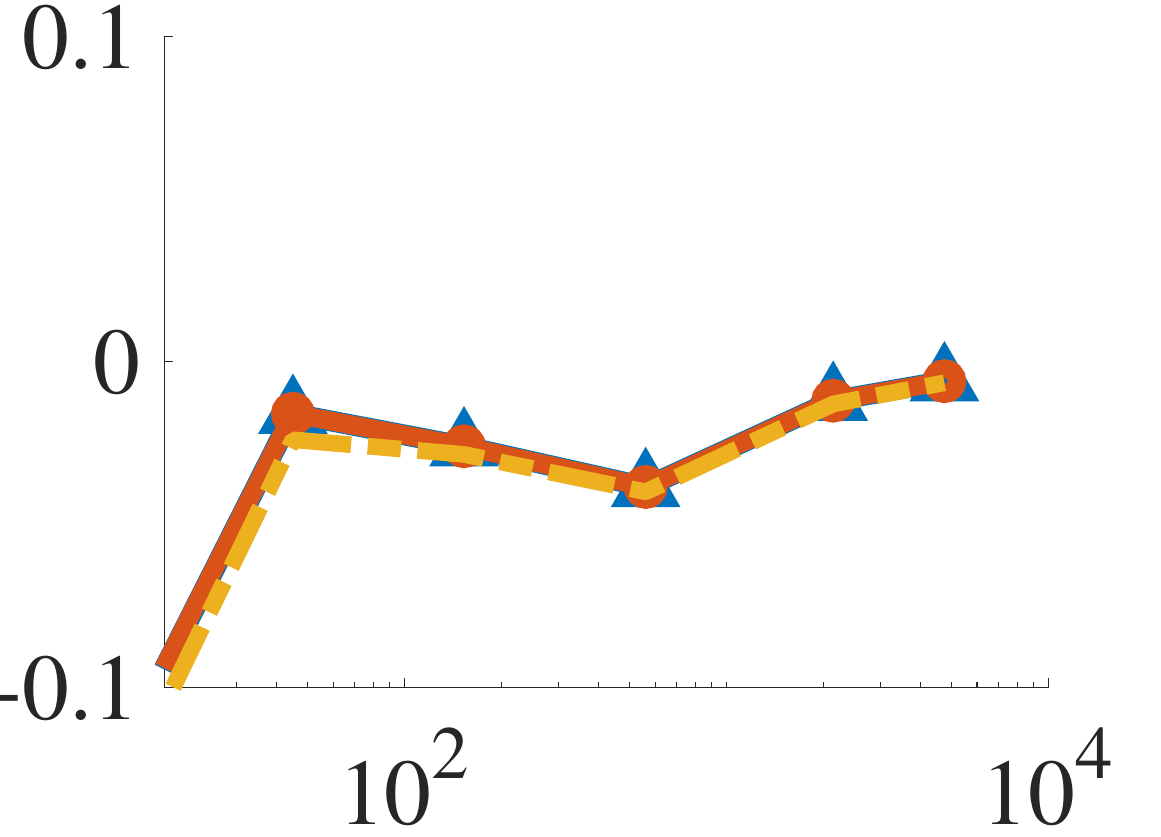} 
\includegraphics[width=.225\linewidth, trim={30 190 25 200}, clip]{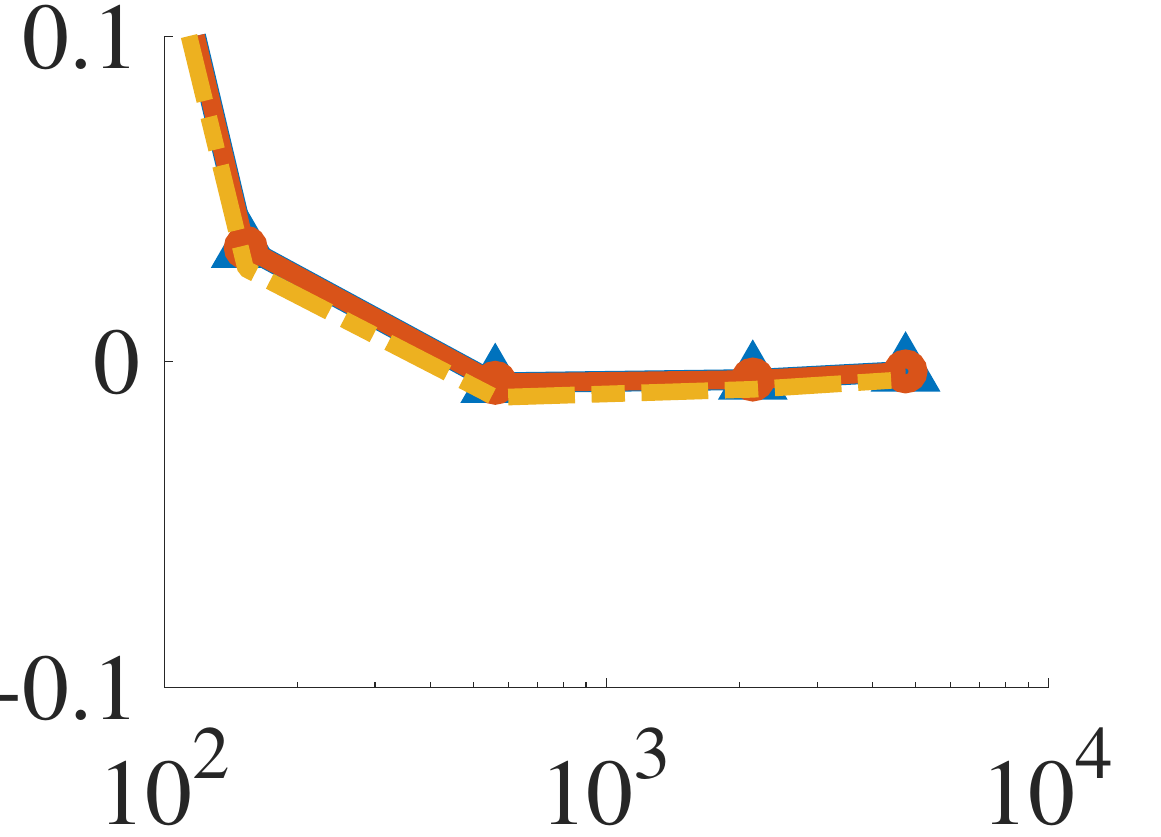} \\
\rotatebox{90}{$\qquad\;$ \textbf{$\nus = .4$} }
   \rotatebox{90}{$\quad$ Vol Change \% }
\includegraphics[width=.225\linewidth, trim={30 190 25 200}, clip]{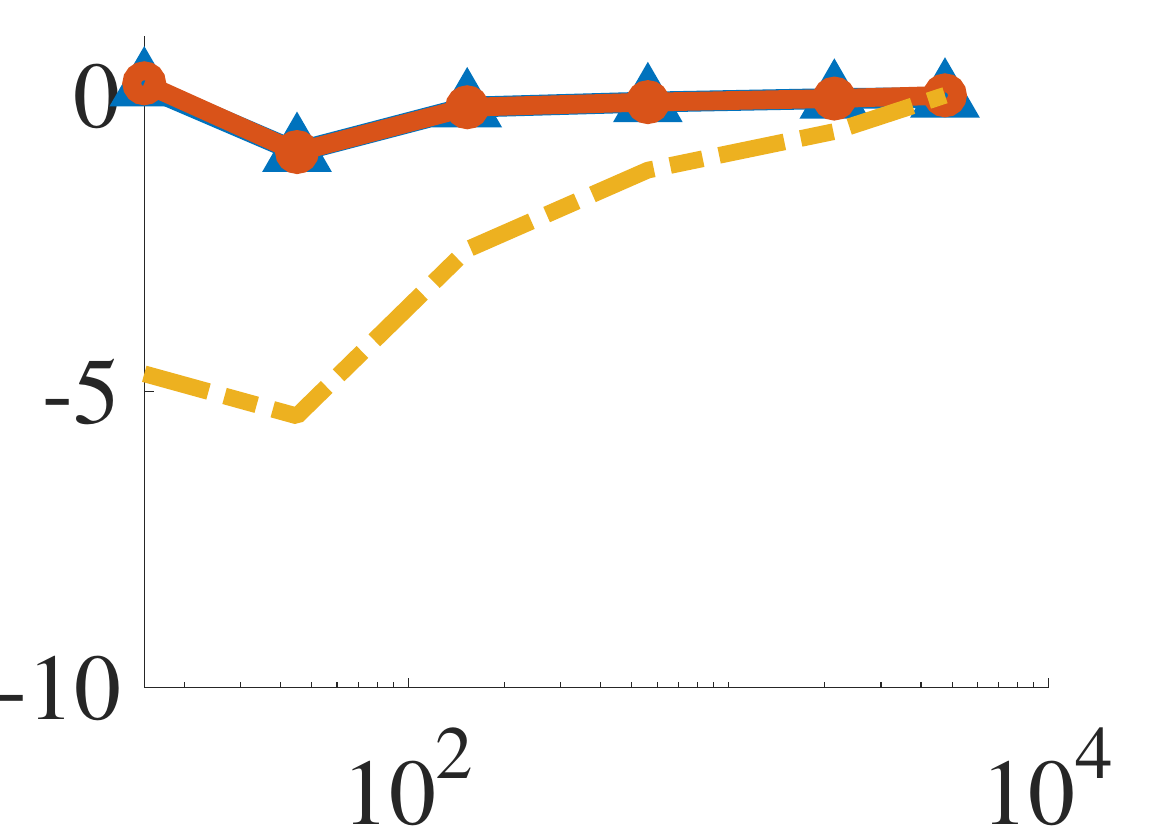} 
\includegraphics[width=.225\linewidth, trim={30 190 25 200}, clip]{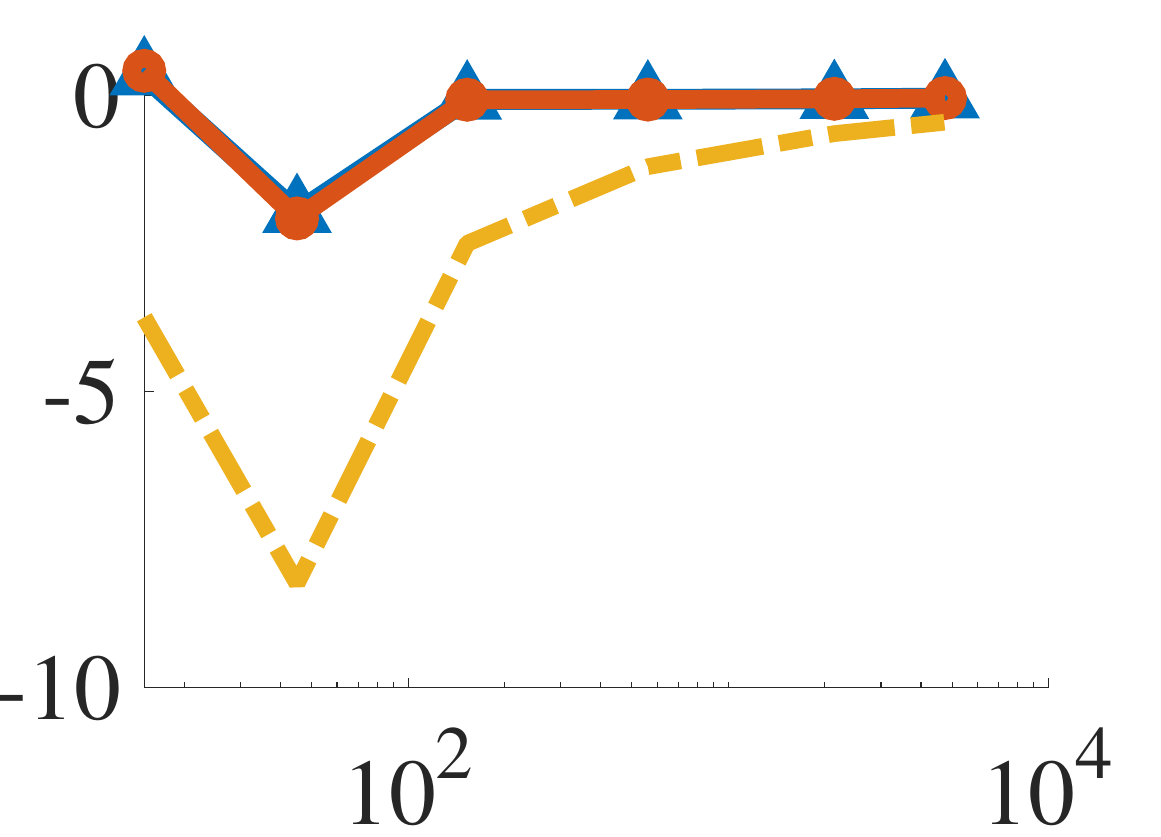}
\includegraphics[width=.225\linewidth, trim={30 190 25 200}, clip]{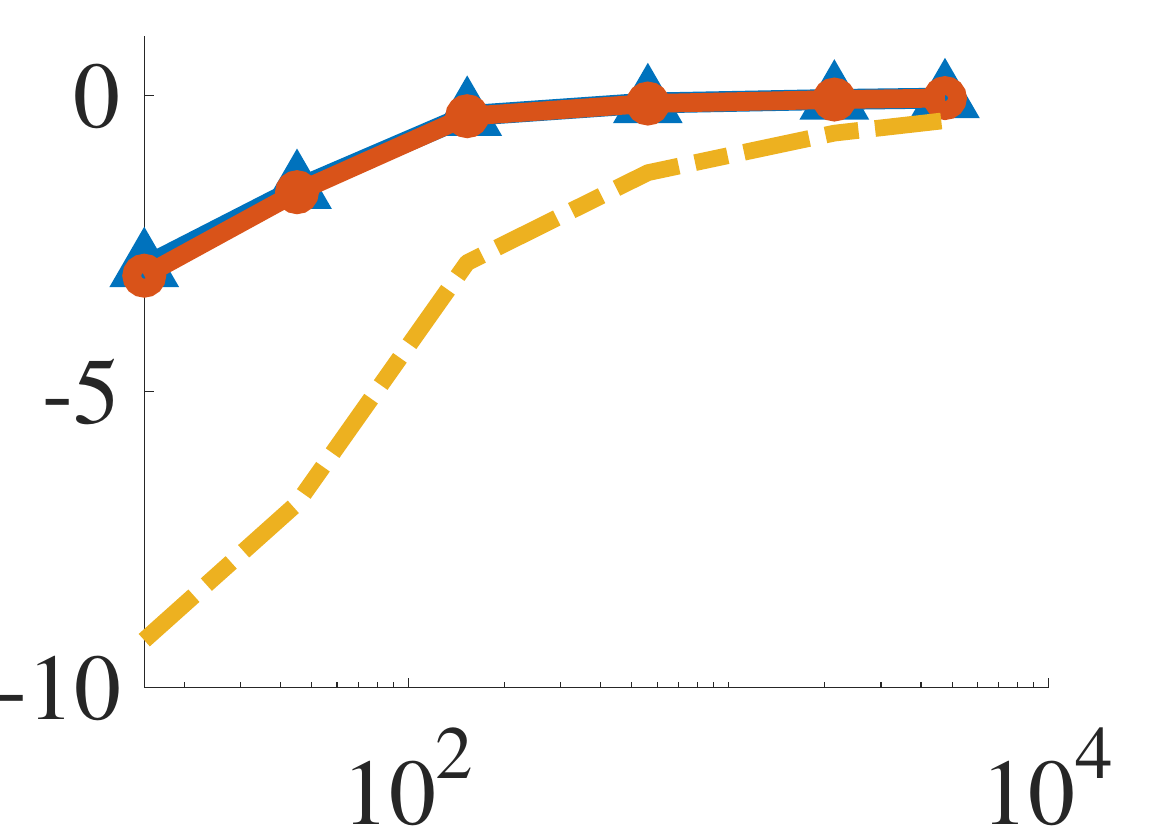} 
\includegraphics[width=.225\linewidth, trim={30 190 25 200}, clip]{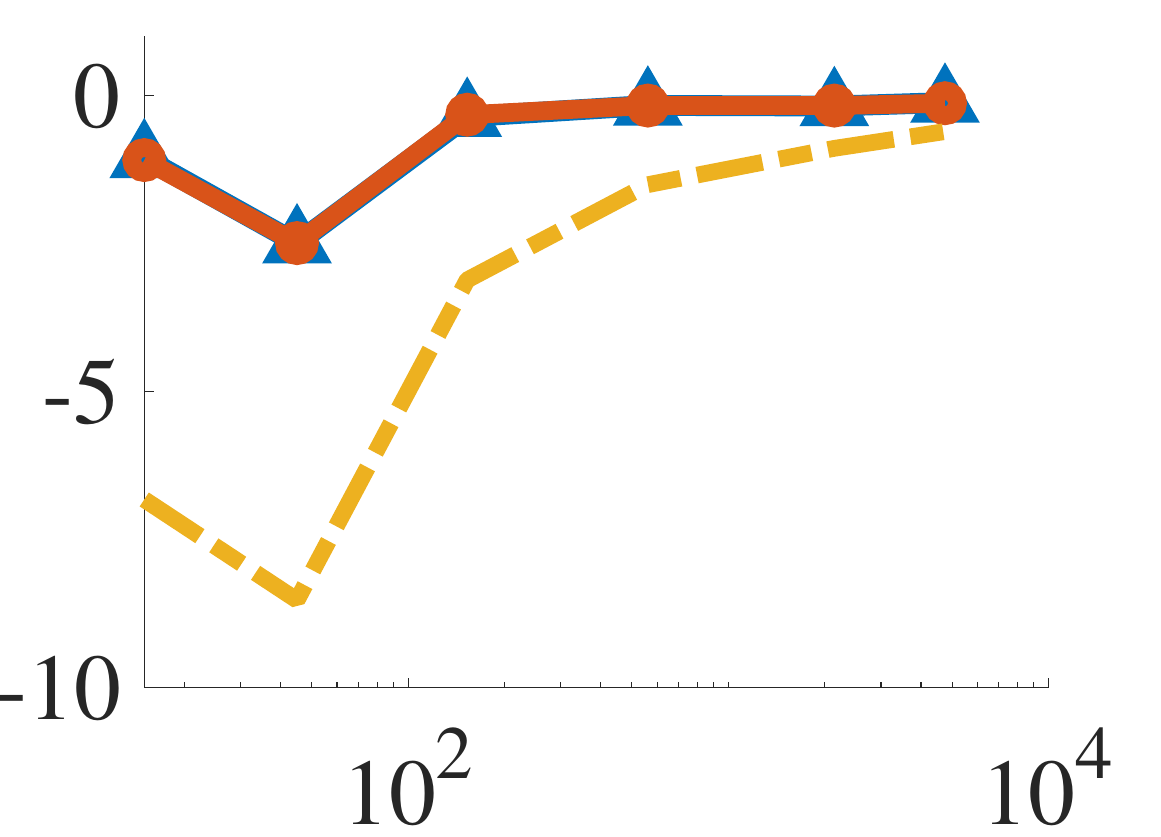}\\

\rotatebox{90}{$\qquad\;$ \textbf{$\nus = 0$} }
   \rotatebox{90}{$\quad$ Vol Change \% }
\includegraphics[width=.225\linewidth, trim={30 190 25 200}, clip]{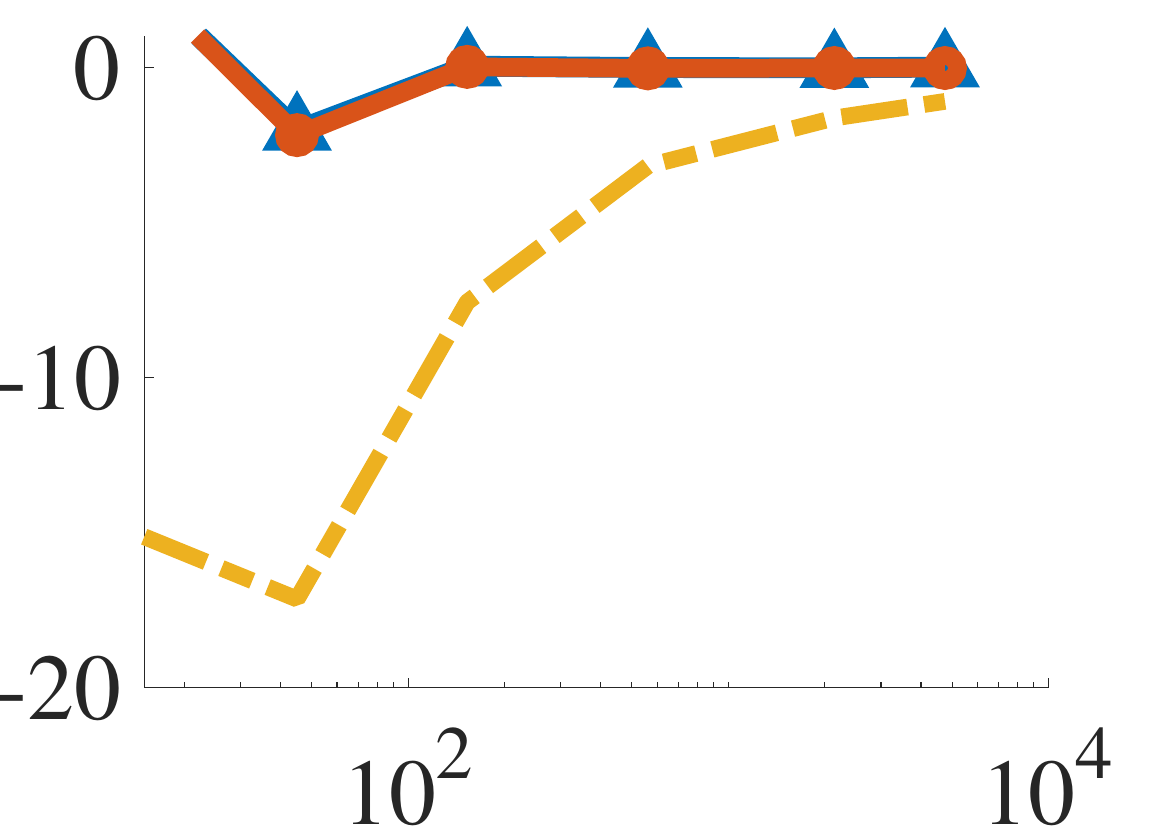} 
\includegraphics[width=.225\linewidth, trim={30 190 25 200}, clip]{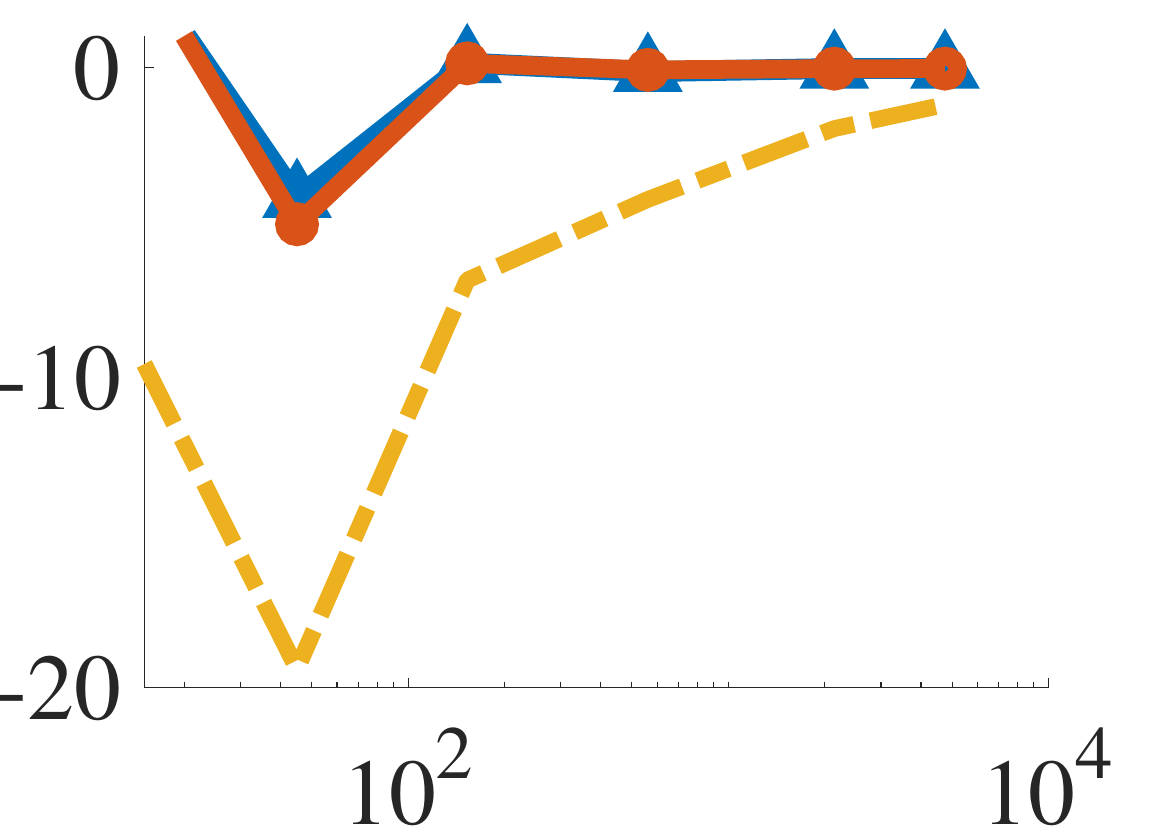} 
\includegraphics[width=.225\linewidth, trim={30 190 25 200}, clip]{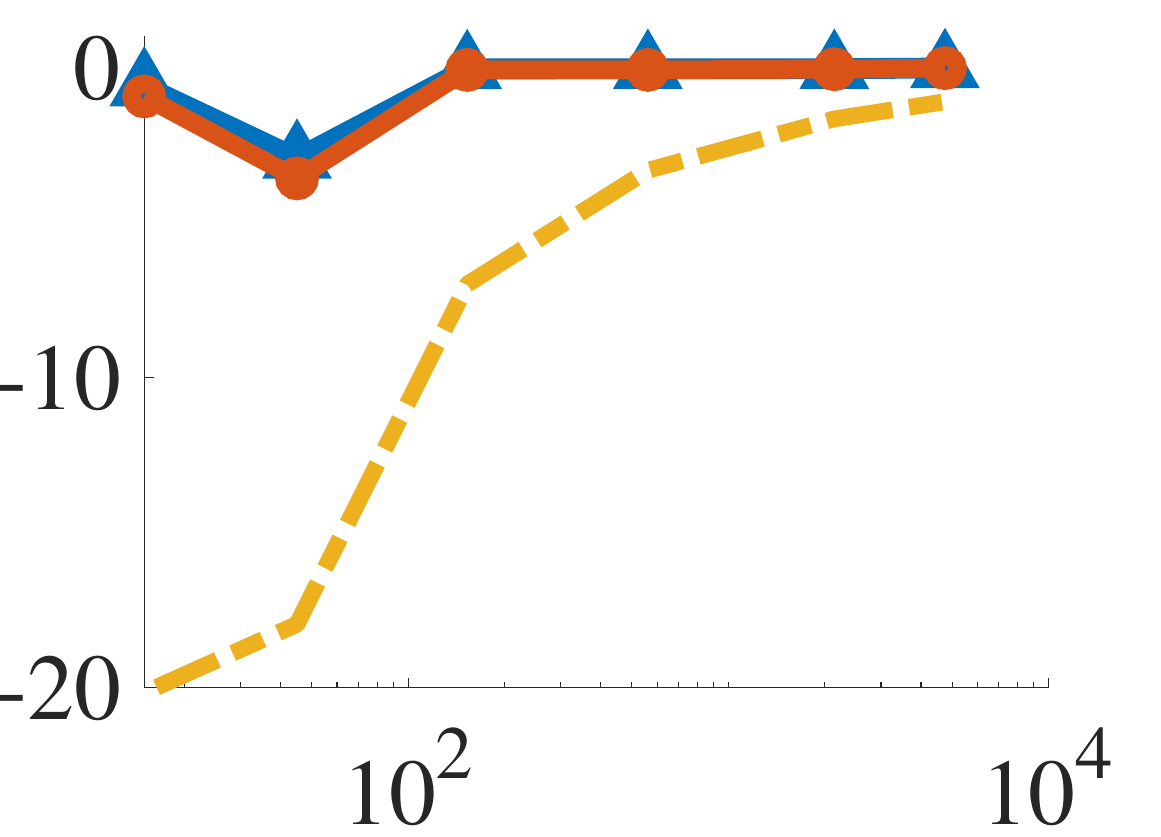} 
\includegraphics[width=.225\linewidth, trim={30 190 25 200}, clip]{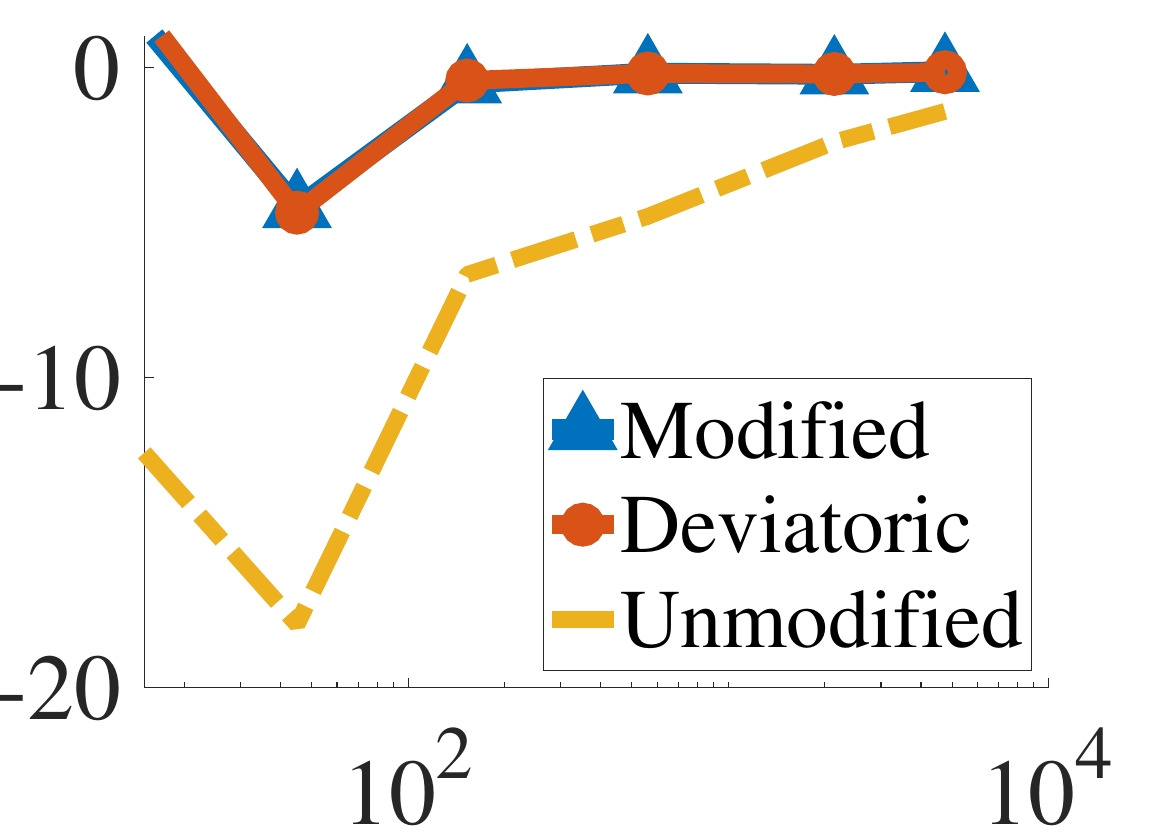}\\

\rotatebox{90}{$\qquad$ \textbf{$\nus = -1$} }
   \rotatebox{90}{$\quad$ Vol Change \% }
\includegraphics[width=.225\linewidth, trim={30 190 25 200}, clip]{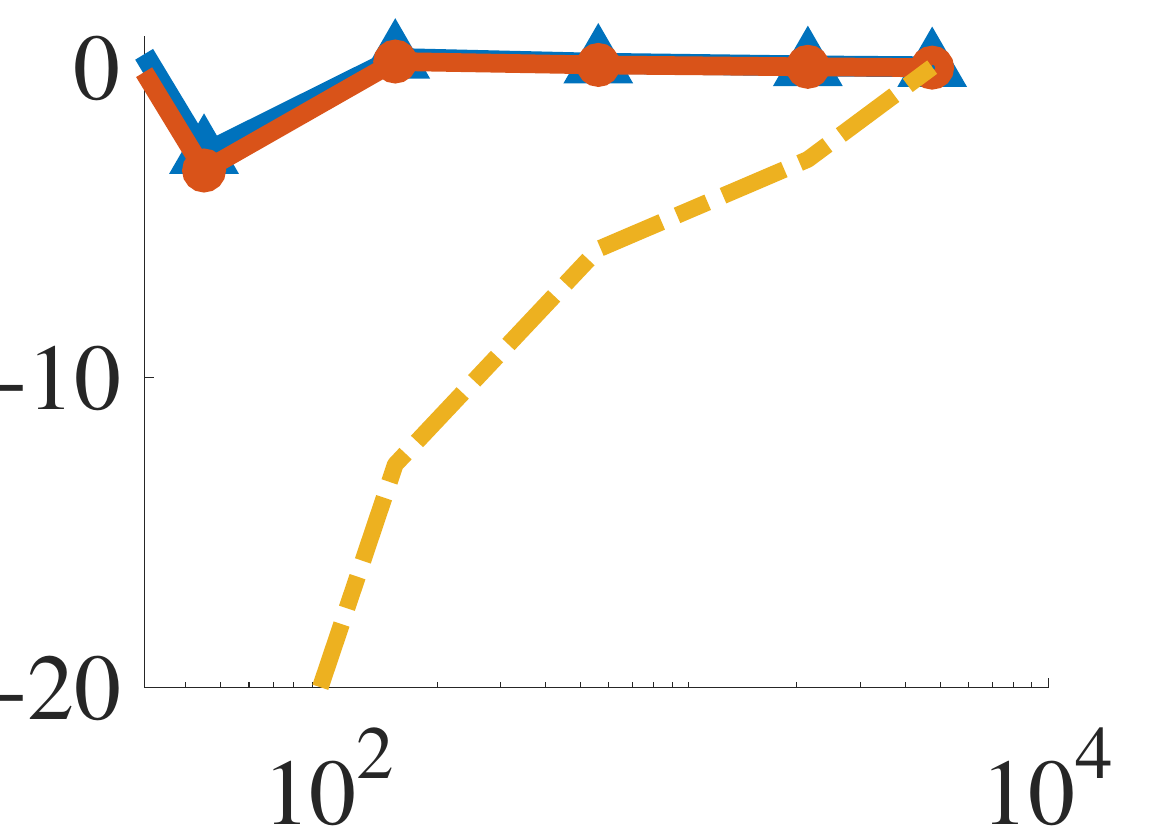} 
\includegraphics[width=.225\linewidth, trim={30 190 25 200}, clip]{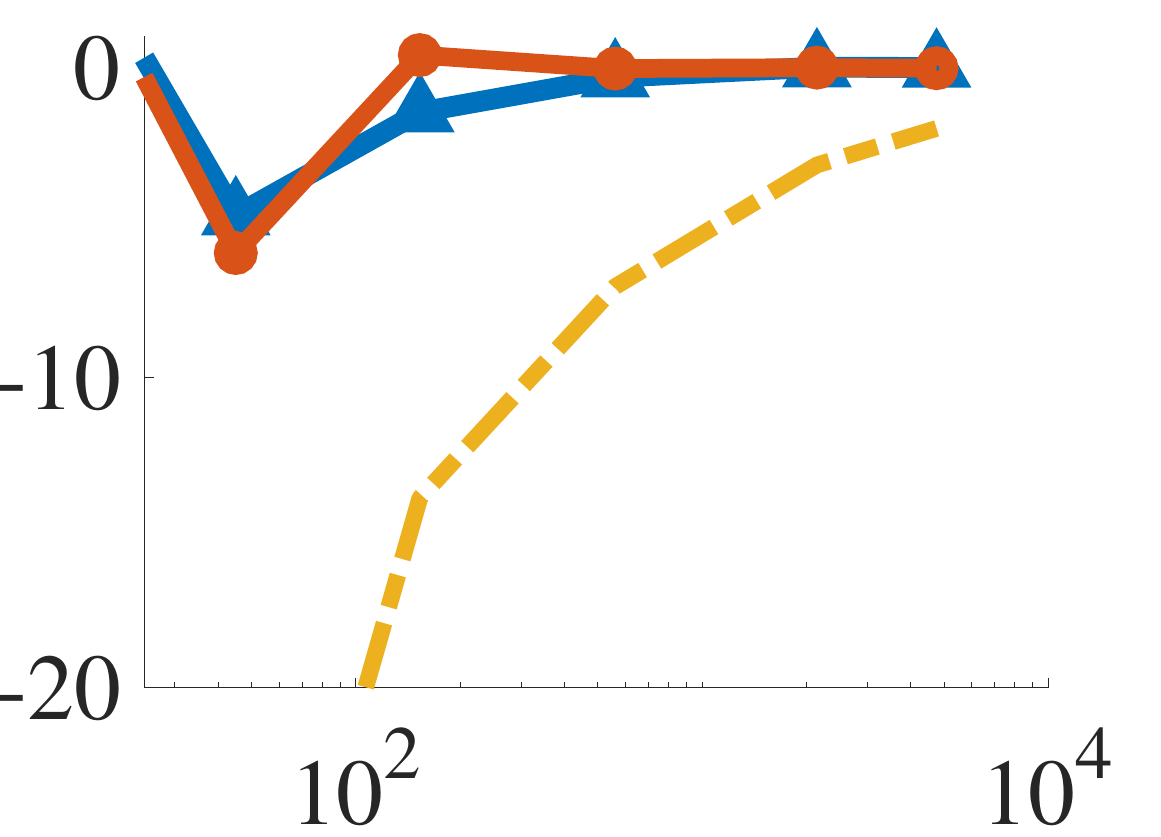} 
\includegraphics[width=.225\linewidth, trim={30 190 25 200}, clip]{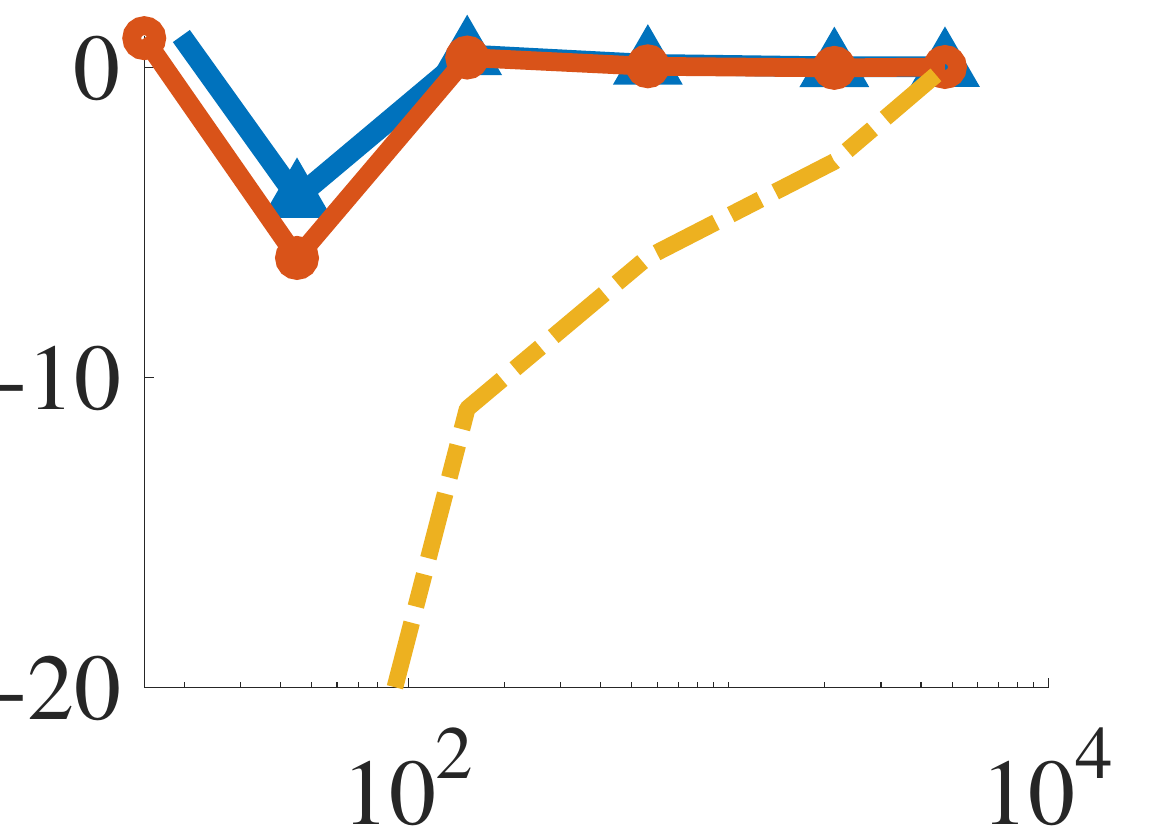} 
\includegraphics[width=.225\linewidth, trim={30 190 25 200}, clip]{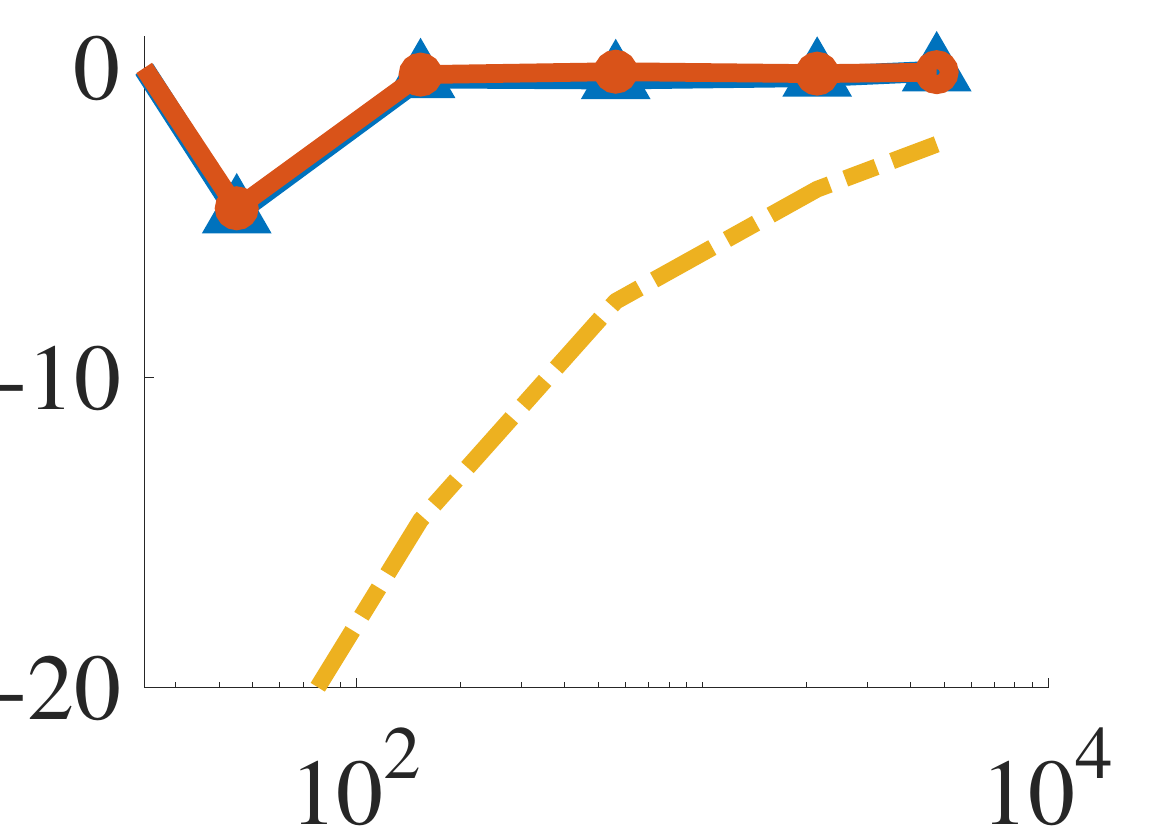}\\

$\qquad\qquad\quad$ \# Solid DOF $\qquad\qquad\quad\;$ \# Solid DOF $\qquad\qquad\quad$ \# Solid DOF $\qquad\qquad\quad\;$ \# Solid DOF
\caption{Percent change in total area for different numbers of solid degrees of freedom (DOF) for the compressed block benchmark (Section \ref{Compression Test}) after deformation. The DOF range from $m = 15$ to $4753$, and the $x$ axis is on a log scale. Omitting the coarsest discretizations ($m=15$), the largest deviations in total volume among all element types used are approximately $2.1\%$ for the modified case, $14\%$ for the unmodified case, and $2.4\%$ for the deviatoric case.}
\label{comp_area}
\end{figure}

\indent This test is a plane strain problem involving a rectangular block with a downward traction applied in the center of the top side of the mesh and zero vertical displacement applied on the bottom boundary; see Figure (\ref{comp_mesh}) for the loading configuration and dimensions of the structure. Zero horizontal displacement is also imposed along the top side. All other boundaries have zero traction applied. This test was used by Reese \etal~\cite{Reese1999} to test a stabilization technique for low order finite elements. A neo-Hookean model is used with shear modulus set to $G = 80.194 \ \frac{\text{dyn}}{\text{cm}^2}$, and damping is set to $\eta = 4.0097 \frac{\text{g}}{\text{s}}$ for this test.  The downward traction has magnitude $200 \ \frac{\text{dyn}}{\text{cm}^2}$. The computational domain is $\Omega = [0,L]^2$ with $L = 40 \ \text{cm}$. The numbers of solid degrees of freedom (DOF) range from $m = 15$ to $m = 4753$ for the FE (\textbf{P1}/\textbf{P1}) results and all the IBFE results. Specifically, we use a sequence of meshes that yields the same node locations for each element type (e.g. for $m = 15$, a \textbf{P1} mesh and a \textbf{Q2} mesh have FE nodes located in the same positions). \\
\indent The primary quantity of interest is the $y$-displacement at the center of the top face. Figure (\ref{cb}) shows deformations at time $T_{\text{f}} = 500$ s. The load time is $T_{\text{l}} = 100$ s. Figure (\ref{cb_dev}) shows the deviatoric stresses for both the FE method and the IBFE method with modified invariants and stabilization. The states of stress are clearly converging to distributions that are in excellent agreement with the FE results. Additionally, we report the pressure field of the IBFE method and the pressure field of the FE method in Figure (\ref{cb_p}). The pressure fields are qualitatively in agreement with the exception of the boundary, which is less accurate for the IBFE method. This is a result of the use of regularized delta functions in the Lagrangian-Eulerian coupling operators, which smooth discontinuities in $\peul$ that can occur at fluid-structure interfaces. Figure ($\ref{comp_disp}$) shows the behavior of the point of interest under refinement. Note the performance of the unmodified invariants in the final row of plots. The convergence behavior of the single recorded point is satisfactory although the overall deformations are unphysical in these cases; again, see Figure ($\ref{cb}$). Particularly noticeable in this benchmark is the effect of using modified invariants versus unmodified invariants while using a nonzero numerical bulk modulus. This is apparent in Figures (\ref{cb}a) and (\ref{cb}b), in which the deformations of the elements are smoothest in (\ref{cb}a) in the case in which modified invariants are used. As expected for values of $\nus$ close to $\frac{1}{2}$, volumetric locking plagues the lower order elements, resulting in poor convergence. Locking is avoided, however, for different values of $\nus$, corresponding to smaller numerical bulk moduli, even for low order elements. Figure (\ref{disp_v_nu}) depicts the displacement of this point as a function of $\nus$ for the compression test. Note that the locking behavior appears for values of $\nus$ larger than $\nus = 0.4$. \\
\indent Figure (\ref{comp_area}) reports the percent change in total volume. Formulations using modified invariants and deviatoric projection yield superior volume conservation. The percent change for all element types considered range between $.0004 \%$ and $2.1 \%$ for the modified invariants, between $.001 \%$ and $14 \%$ for the unmodified invariants, and between $.0005 \%$ and $2.4\%$ for the deviatoric projection. These ranges account for change in area in an absolute sense, whereas the plots display whether the change in area was a gain or loss.

%%%%%%%%%%%%%%%%%%%%%

\subsection{Cook's Membrane}
\label{Cook's Membrane}
\begin{figure}
\centering
\includegraphics[width=.7\linewidth, trim={60 0 70 0}, clip]{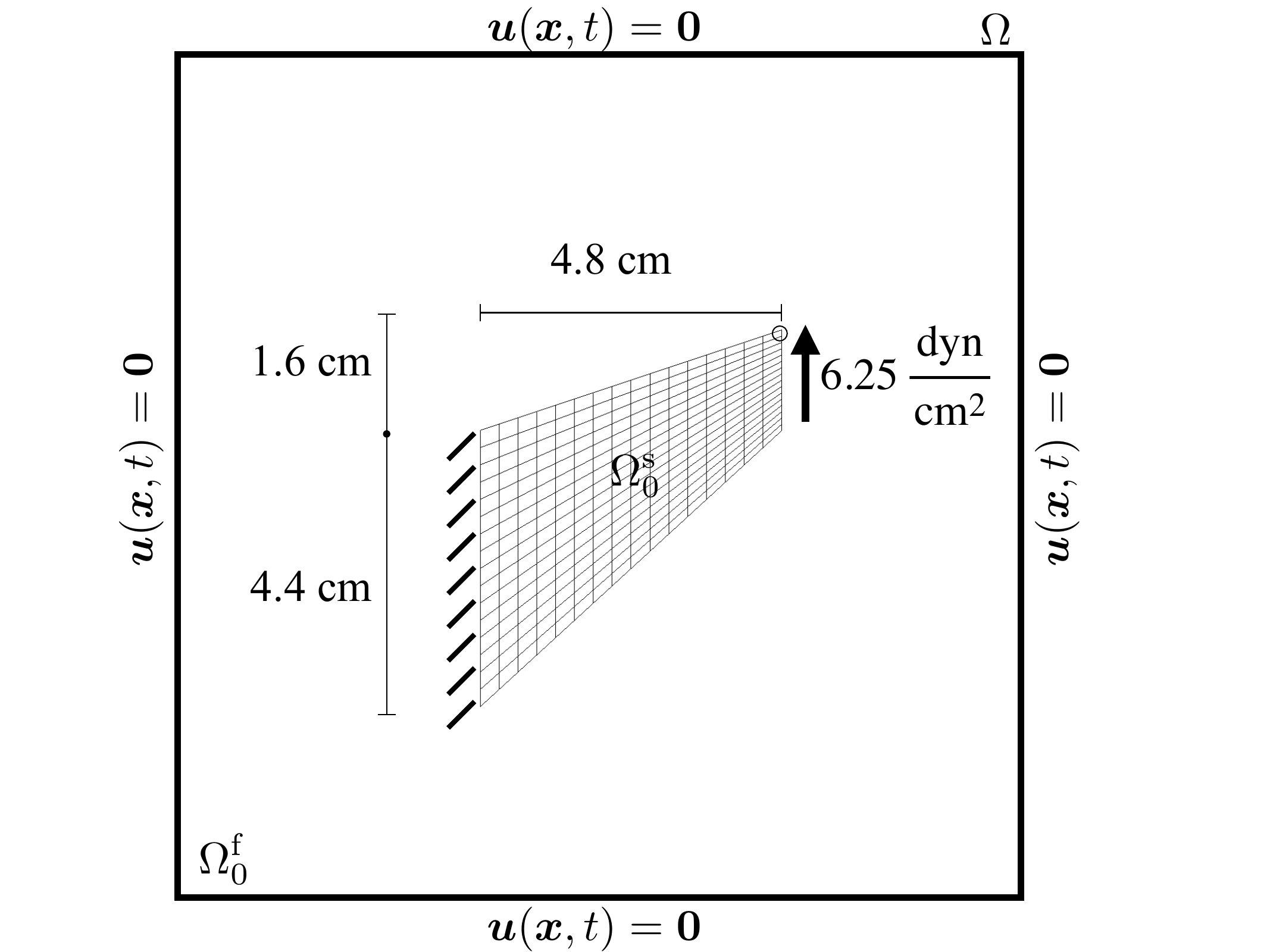}
\caption{Specifications of the Cook's membrane benchmark (Section \ref{Cook's Membrane}). The primary quantity of interest is the $y$-displacement as measured at the upper right hand corner, indicated by the circle. The structure, shown here in its initial configuration and denoted by $\soliddomO$, is immersed in a fluid denoted by $\fluiddomO$. The entire computational domain is $\Omega = \fluiddom \cup \soliddom$. Zero fluid velocity is enforced on the boundary of $\Omega$.}
\label{cooks}
\end{figure}

\begin{figure}
\begin{tabular}{l c c}
& \textbf{Modified Invariants} & \textbf{Unmodified Invariants} \\
\rotatebox{90}{\qquad\qquad\qquad\qquad \textbf{$\nus = .4$} }&
\subcaptionbox{\label{sfig:testa}} {\includegraphics[width=.45\linewidth, trim={50 140 0 20}, clip]{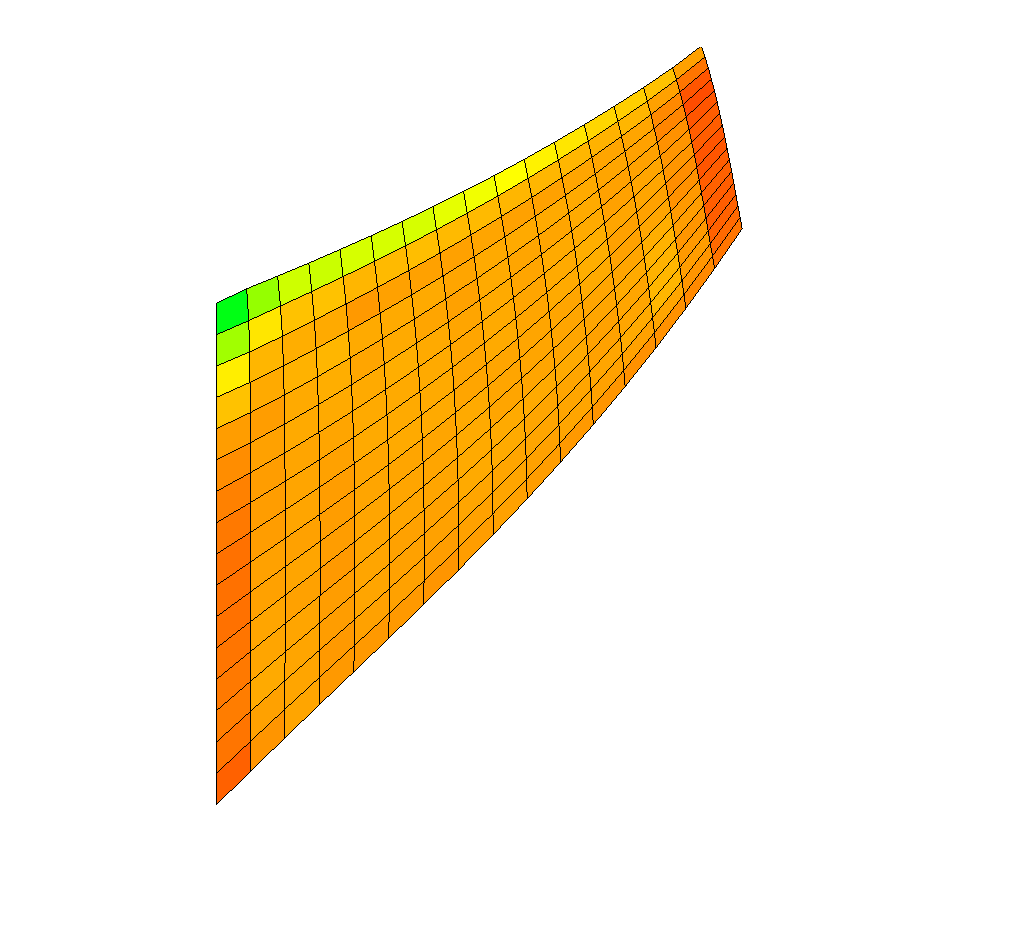}} &
\subcaptionbox{\label{sfig:testb}} {\includegraphics[width=.45\linewidth, trim={50 140 0 20}, clip]{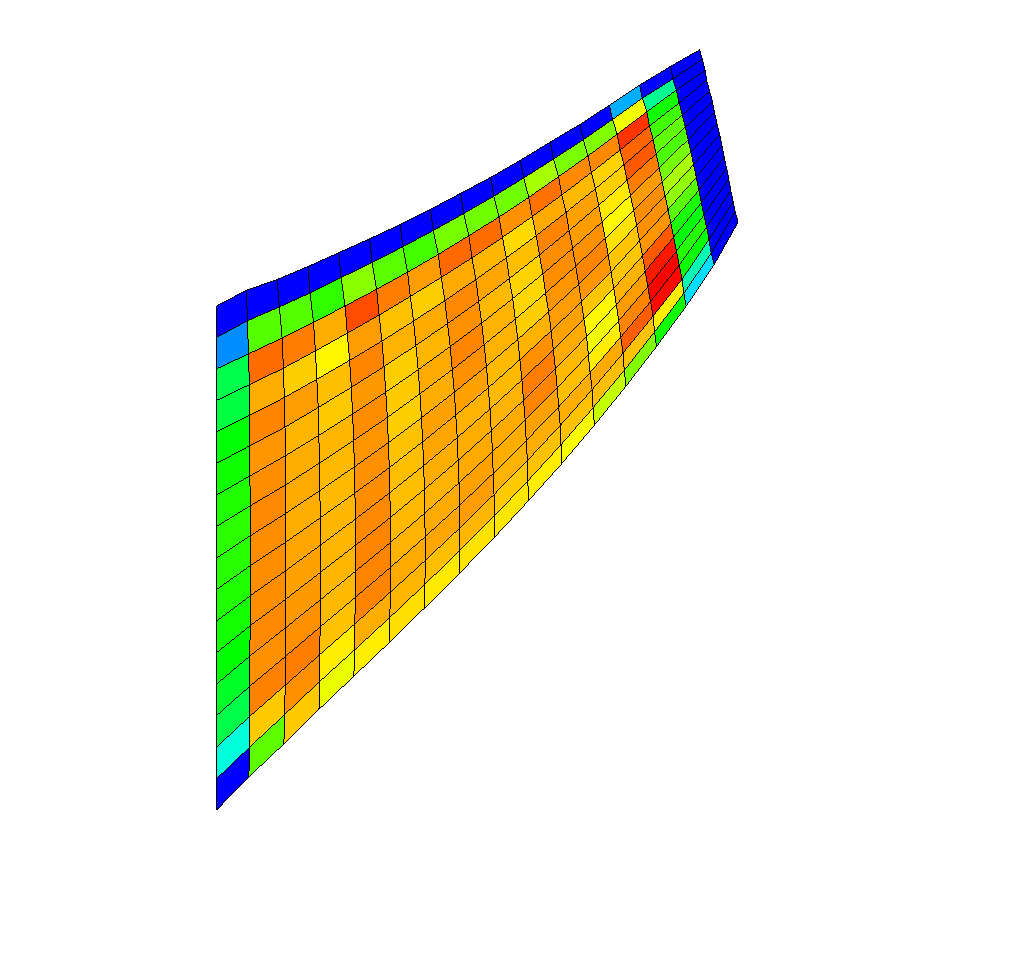}}\\

\rotatebox{90}{\qquad\qquad\qquad\qquad \textbf{$\nus = -1$}} &
\subcaptionbox{\label{sfig:testc}}{\includegraphics[width=.45\linewidth, trim={50 140 0 20}, clip]{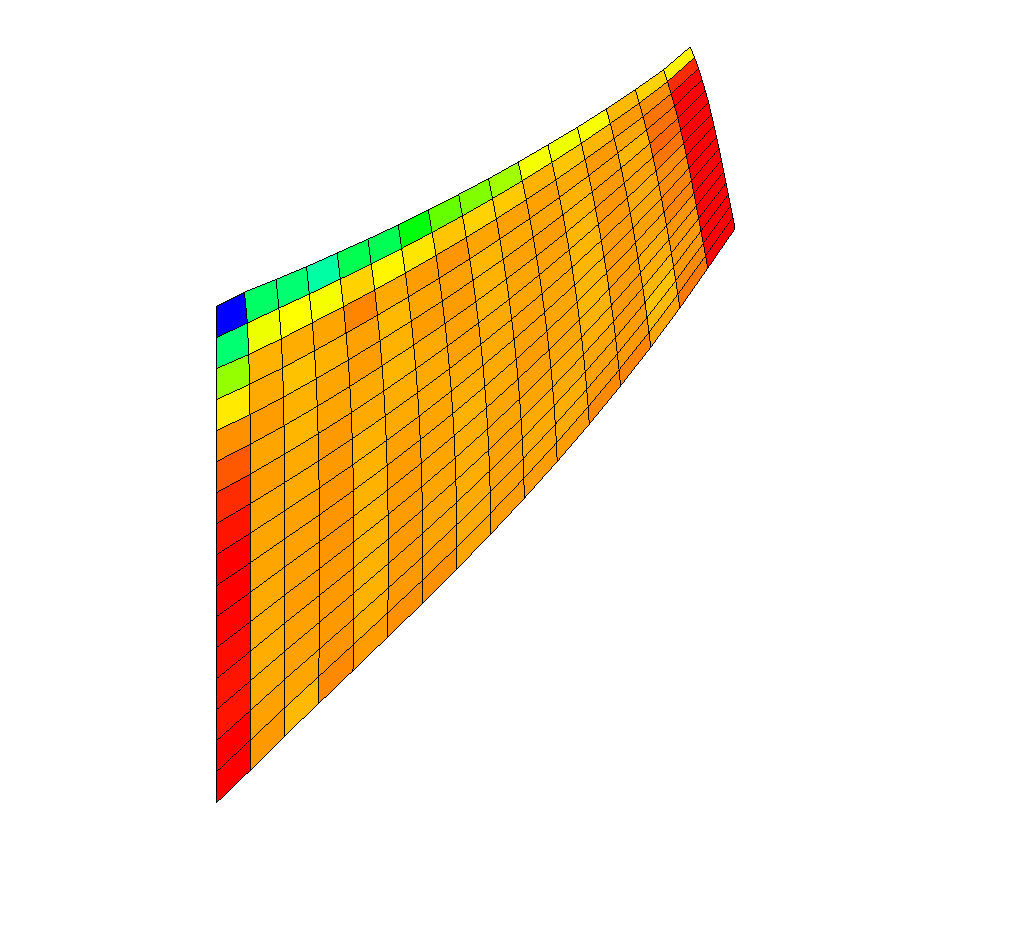}} &
\subcaptionbox{\label{sfig:testd}}{\includegraphics[width=.45\linewidth, trim={50 140 0 20}, clip]{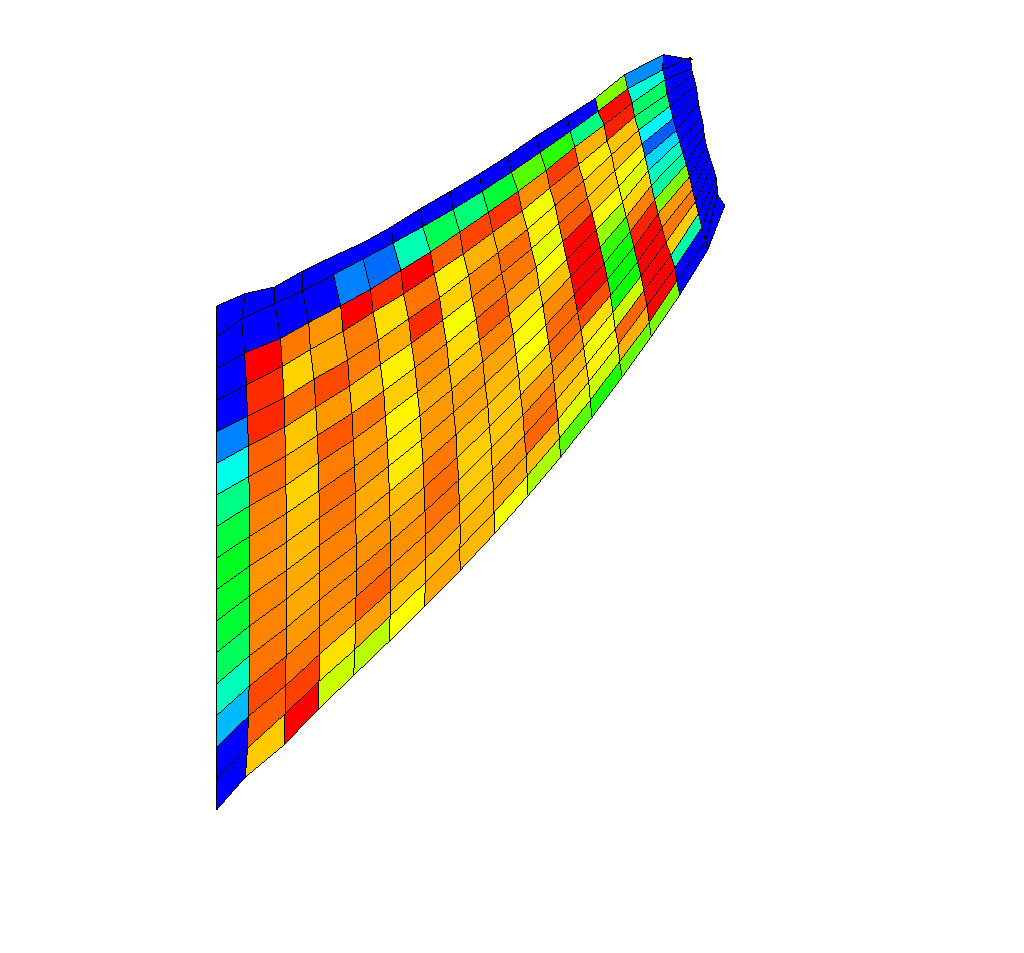}}  \\
\end{tabular}
%trim={left bottom right top}
\begin{centering}
Avg $J$ \\
\includegraphics[width=2.5in, trim={0 5in 0 5in}, clip]{color_bar.pdf}  \\
%0.40 \ \ \ \ \ \ \ \ \ \ \ \ \ \ \ \ 1.10
0.95$\qquad\qquad\qquad\qquad$ 1.01

\end{centering}
\caption{Deformations of the Cook's membrane benchmark (Section \ref{Cook's Membrane}), along with mean values of $J$ within each element calculated via equation \eqref{avgJ}, using a neo-Hookean material model, equations (\ref{nh_energy}) -- (\ref{nh_stress_mod}), with $G = 83.3333 \ \frac{\text{dyn}}{\text{cm}^2}$. The background Eulerian grid is not shown. Shown here are solid meshes with \textbf{Q1} elements and $m = 289$ DOF. The first row shows cases with $\nus = .4$, and the second row shows cases with $\nus = -1$ (here equivalent to $\kappas = 0$ and no volumetric-based stabilization). The first column depicts cases with modified invariants, and the second column depicts cases with unmodified invariants. Notice that the case with modified invariants with nonzero numerical bulk modulus have the smoothest deformations and provides the best volume conservation, whereas those of the case with unmodified invariants and zero numerical bulk modulus behave unphysically.}
\label{cm}
\end{figure}

\begin{figure}
\begin{tabular}{l c c c}
&$\cauchys_{00}$ & $\cauchys_{01}$& $\cauchys_{11}$ \\

\rotatebox{90}{$\qquad\;$ \textbf{m = 289} }&
\subcaptionbox{\label{sfig:testa}}{\includegraphics[width=.2\linewidth]{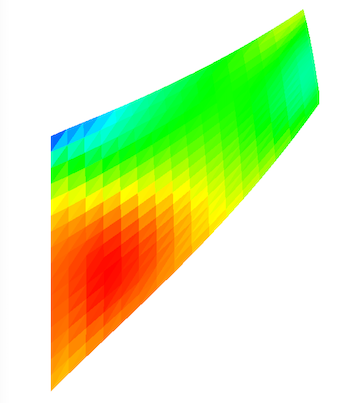}}&
\subcaptionbox{\label{sfig:testb}}{\includegraphics[width=.2\linewidth]{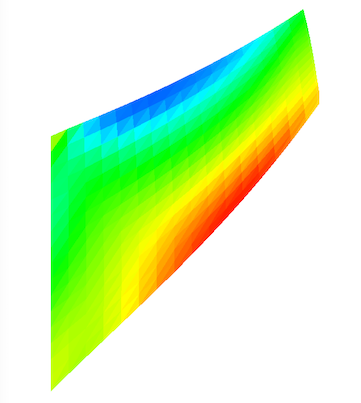}}&
\subcaptionbox{\label{sfig:testa}}{\includegraphics[width=.2\linewidth]{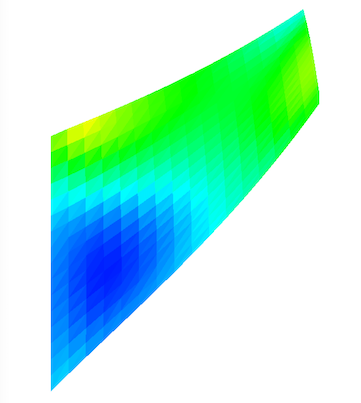}}\\

\rotatebox{90}{$\qquad\;$ \textbf{m = 1089} }&
\subcaptionbox{\label{sfig:testa}}{\includegraphics[width=.2\linewidth]{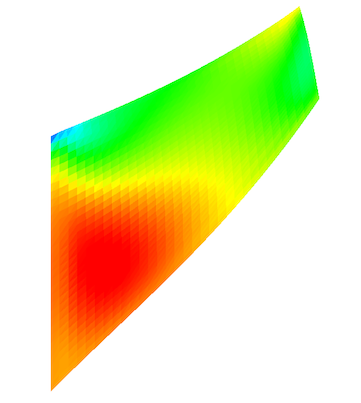}}&
\subcaptionbox{\label{sfig:testb}}{\includegraphics[width=.2\linewidth]{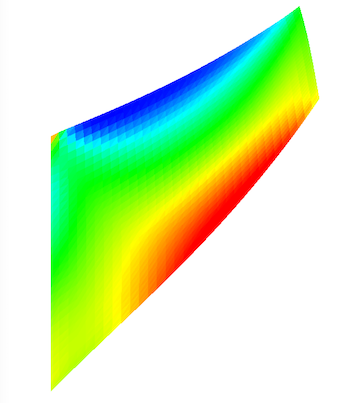}}&
\subcaptionbox{\label{sfig:testa}}{\includegraphics[width=.2\linewidth]{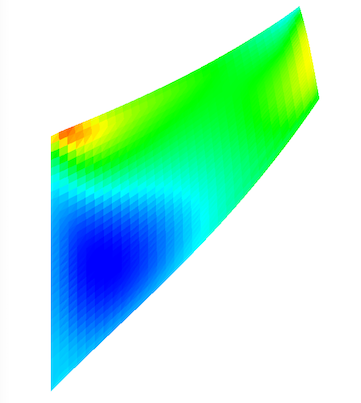}}\\

\rotatebox{90}{$\qquad\;$ \textbf{m = 4225} }&
\subcaptionbox{\label{sfig:testa}}{\includegraphics[width=.2\linewidth]{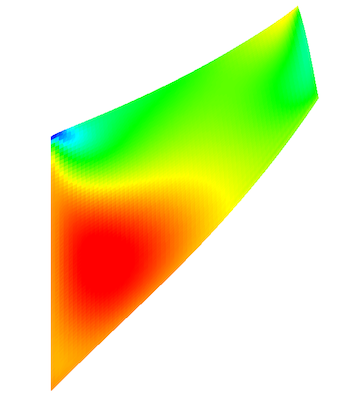}}&
\subcaptionbox{\label{sfig:testb}}{\includegraphics[width=.2\linewidth]{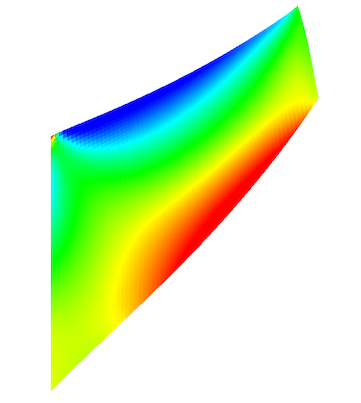}}&
\subcaptionbox{\label{sfig:testa}}{\includegraphics[width=.2\linewidth]{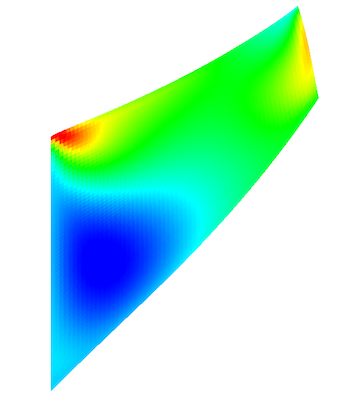}}\\

\\
\hline
\rotatebox{90}{$\qquad\;$ \textbf{FE (P1/P1)} }&
\subcaptionbox{\label{sfig:testa}}{\includegraphics[width=.2\linewidth]{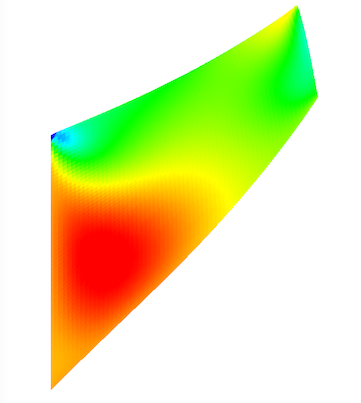}}&
\subcaptionbox{\label{sfig:testb}}{\includegraphics[width=.2\linewidth]{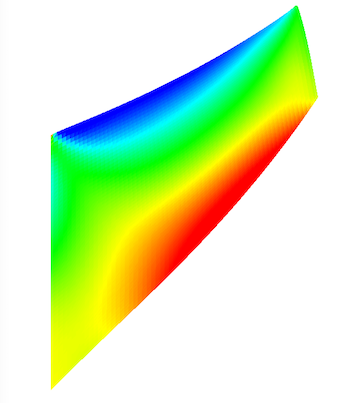}}&
\subcaptionbox{\label{sfig:testa}}{\includegraphics[width=.2\linewidth]{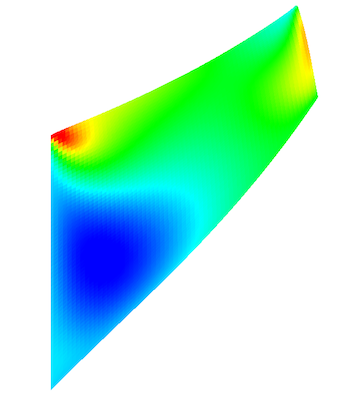}}\\

&
\includegraphics[width=.3\linewidth, trim={0 5in 0 5in}, clip]{color_bar.pdf}&
\includegraphics[width=.3\linewidth, trim={0 5in 0 5in}, clip]{color_bar.pdf}&
\includegraphics[width=.3\linewidth, trim={0 5in 0 5in}, clip]{color_bar.pdf}  \\
&-12 $\qquad\qquad\qquad$ 2.5& -6 $\qquad\qquad\qquad$ 10 & -2.5 $\qquad\qquad\qquad$ 10 \\
\end{tabular}
\caption{Three components of the deviatoric part of $\cauchys$ stress for the Cook's membrane benchmark (Section \ref{Cook's Membrane}). The IBFE method uses modified invariants and volumetric stabilization ($\nus = 0.4$). Each row is labeled with the DOF, and the bottom row depicts the FE (\textbf{P1/P1}) solution with $m = 4225$ solid DOF (equal to that of the highest resolution IBFE results presented in this figure). We use \textbf{P1} elements for each method. The results from the IBFE formulation are clearly converging to the high-resolution FE solution.}
\label{cm_dev}
\end{figure}

\begin{figure}
$\qquad\qquad\qquad\;\;\;\;$ \textbf{P1} $\qquad\qquad\qquad\qquad\quad$  \textbf{Q1} $\qquad\qquad\qquad\qquad\;\;\;$  \textbf{P2} $\qquad\qquad\qquad\qquad\quad\;$ \textbf{Q2}\\
\rotatebox{90}{$\quad\;$ \textbf{$\nus = .49995$} }
   \rotatebox{90}{$\qquad$ Disp. (cm) }
\includegraphics[width=.225\textwidth, trim={30 190 25 200}, clip]{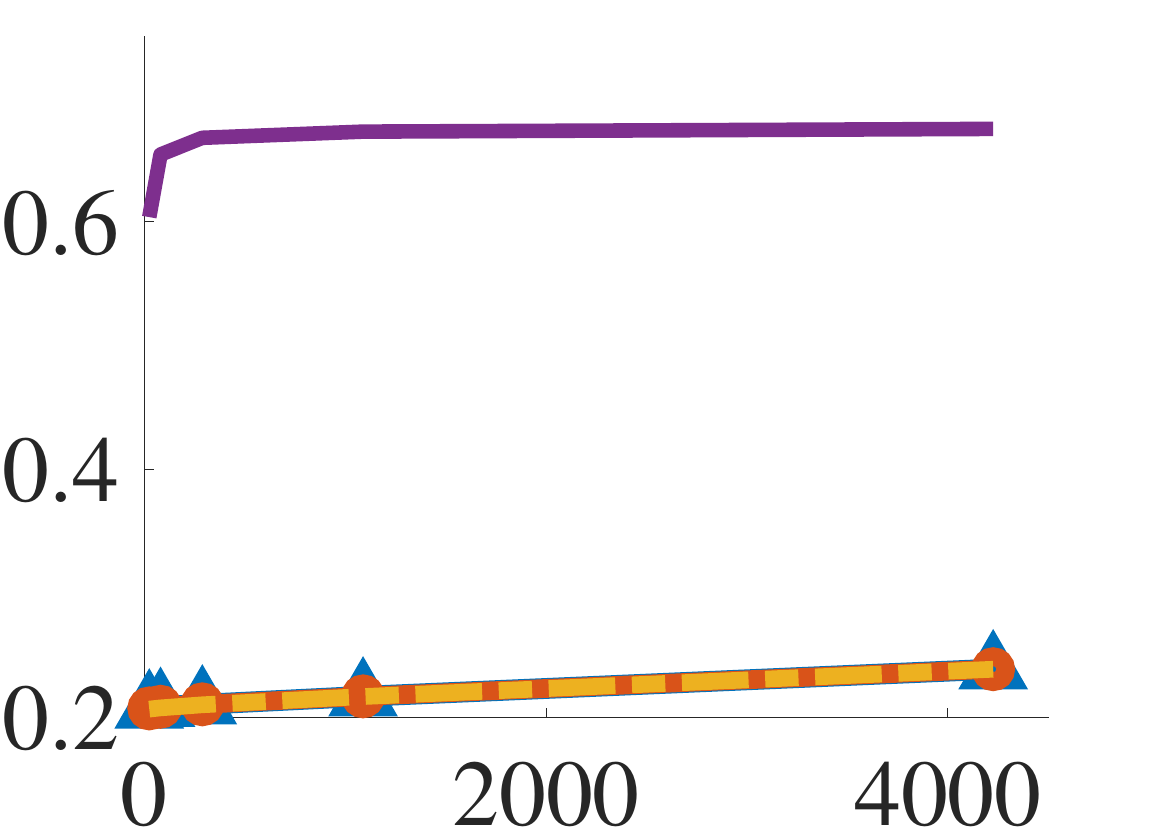}
\includegraphics[width=.225\textwidth, trim={30 190 25 200}, clip]{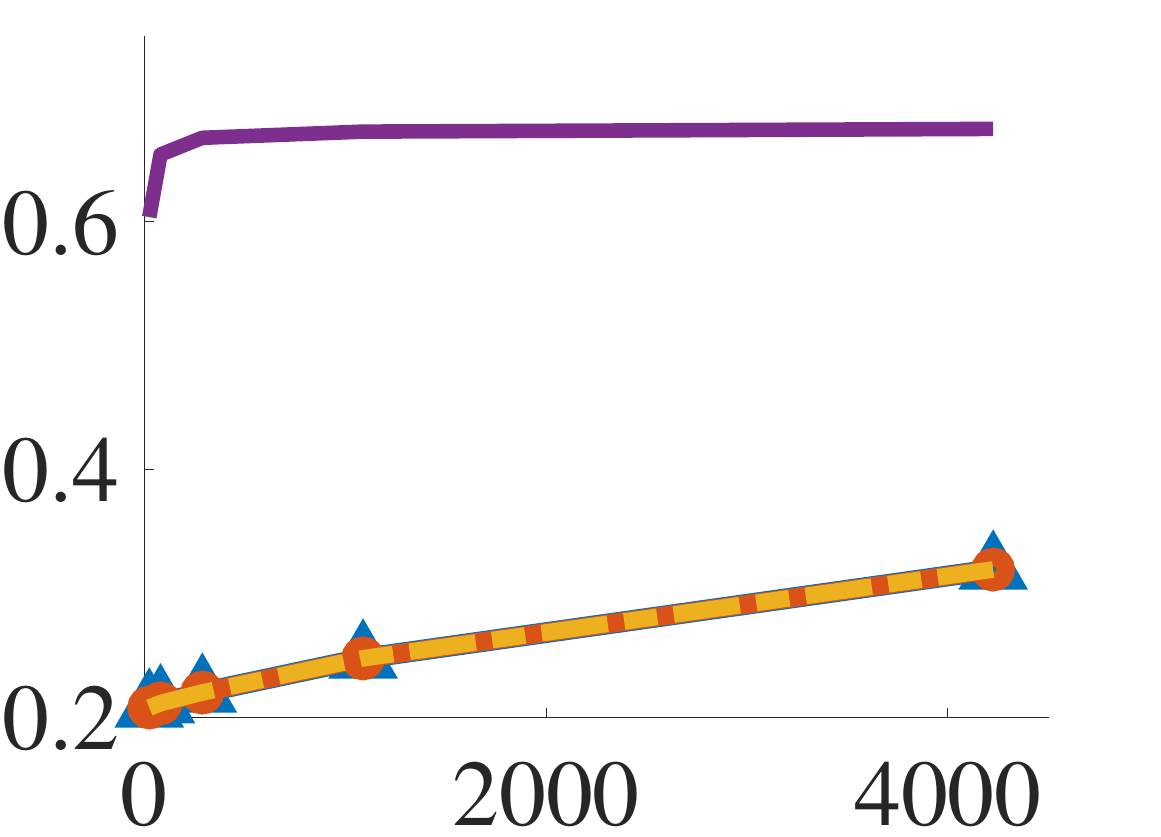}
\includegraphics[width=.225\textwidth, trim={30 190 25 200}, clip]{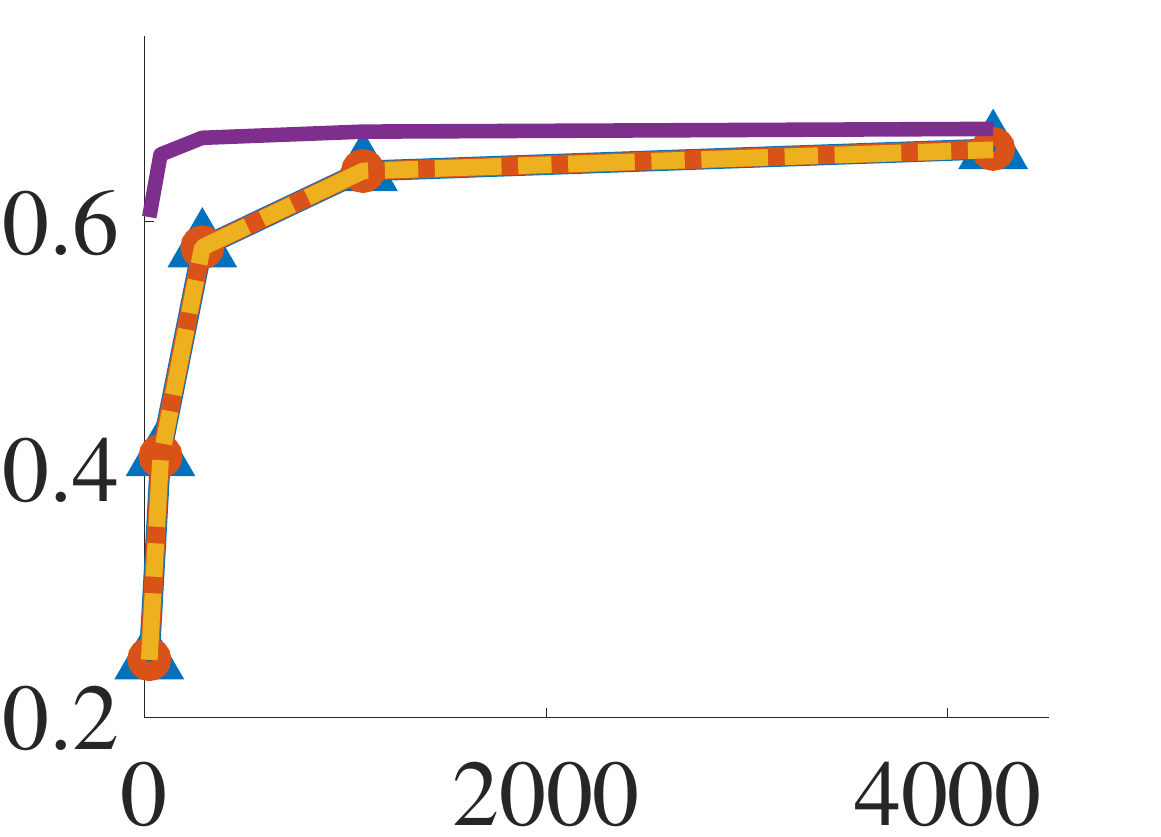} 
\includegraphics[width=.225\textwidth, trim={30 190 25 200}, clip]{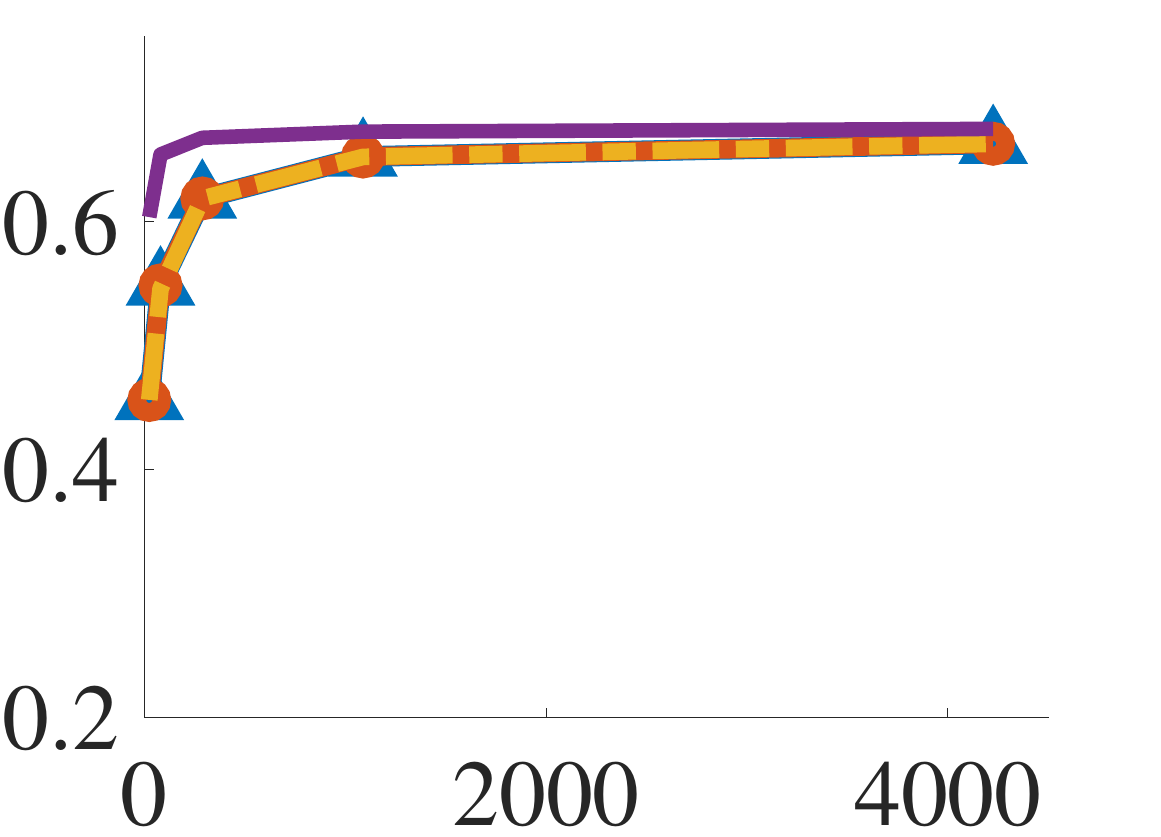} \\
\rotatebox{90}{$\qquad\;\;$ \textbf{$\nus = .4$} }
   \rotatebox{90}{$\qquad$ Disp. (cm) }
\includegraphics[width=.225\linewidth, trim={30 190 25 200}, clip]{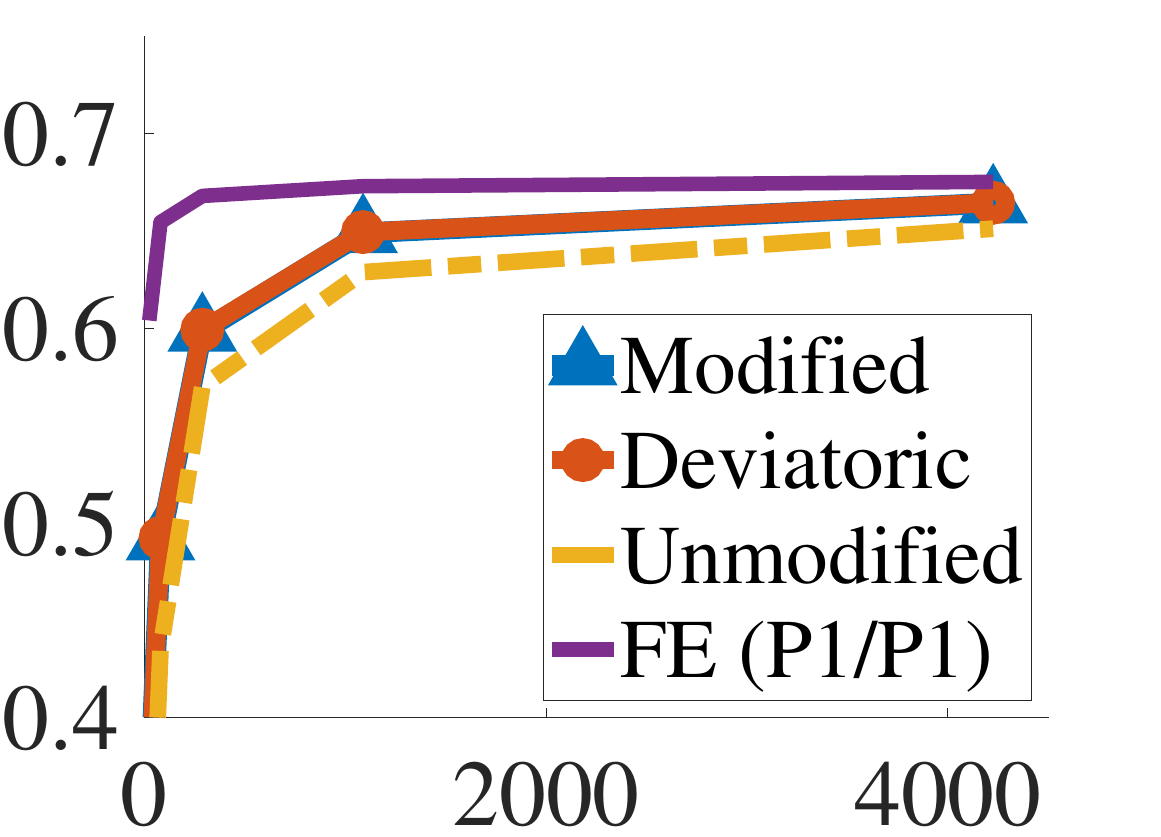}
\includegraphics[width=.225\linewidth, trim={30 190 25 200}, clip]{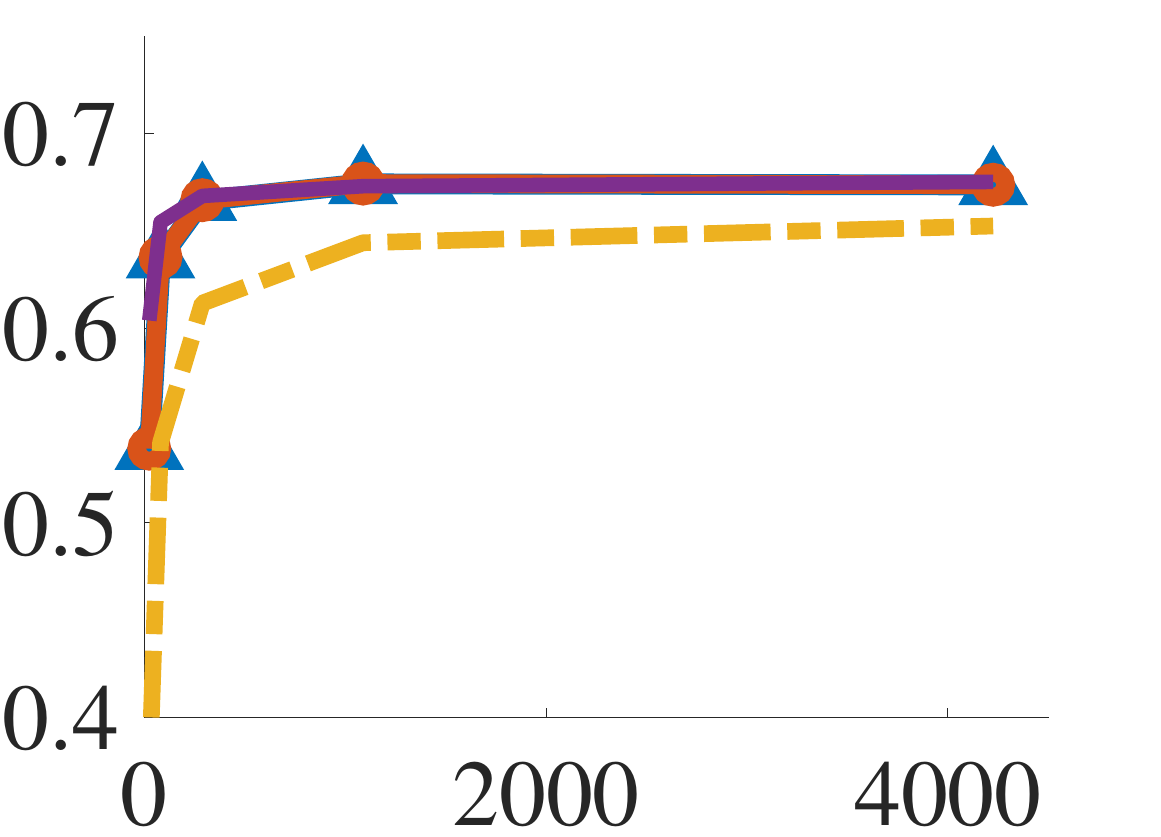}
\includegraphics[width=.225\linewidth, trim={30 190 25 200}, clip]{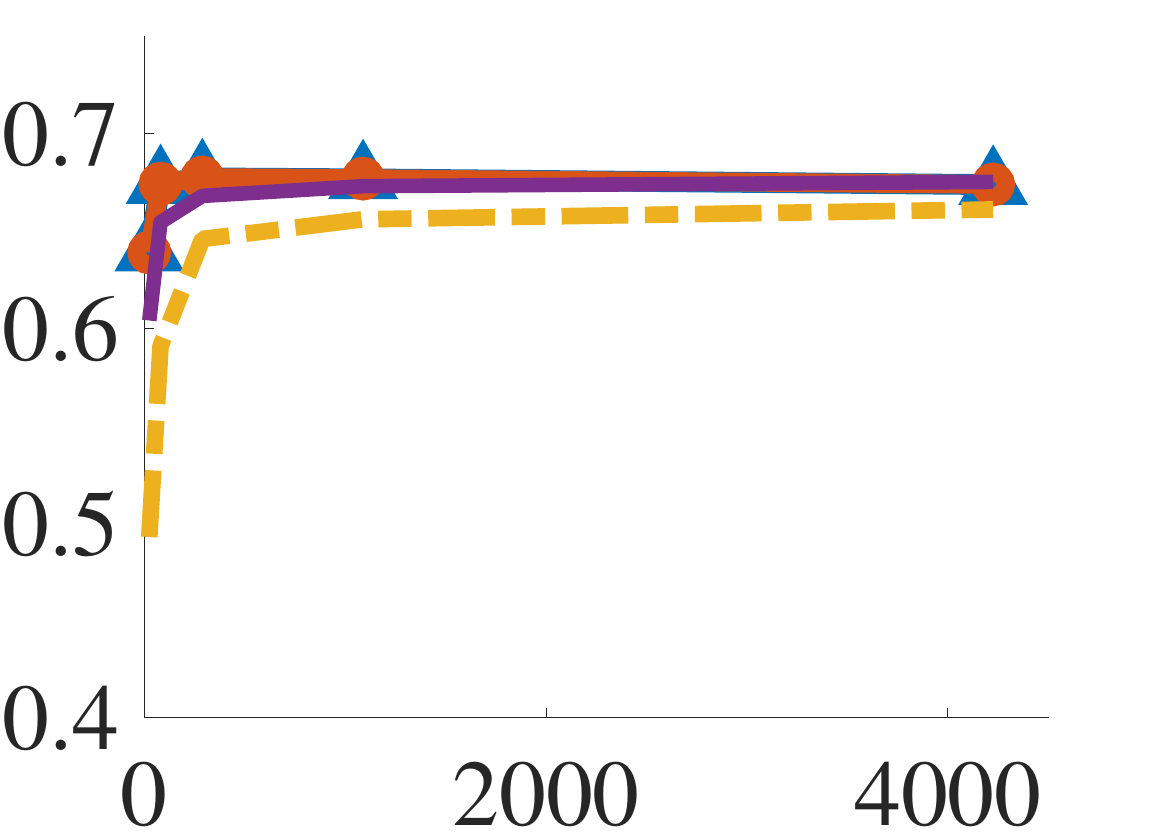}
\includegraphics[width=.225\linewidth, trim={30 190 25 200}, clip]{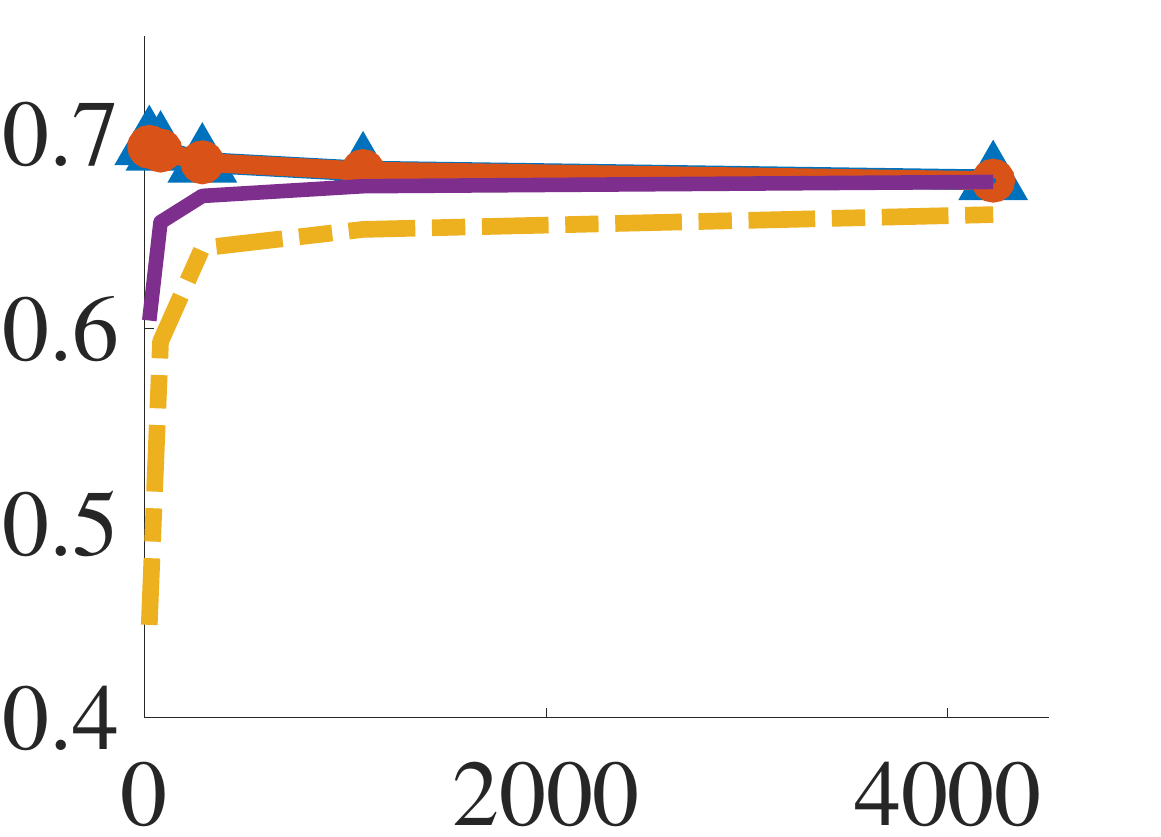}\\

\rotatebox{90}{$\qquad\;\;$ \textbf{$\nus = 0$} }
   \rotatebox{90}{$\qquad$ Disp. (cm) }
   \includegraphics[width=.225\linewidth, trim={30 190 25 200}, clip]{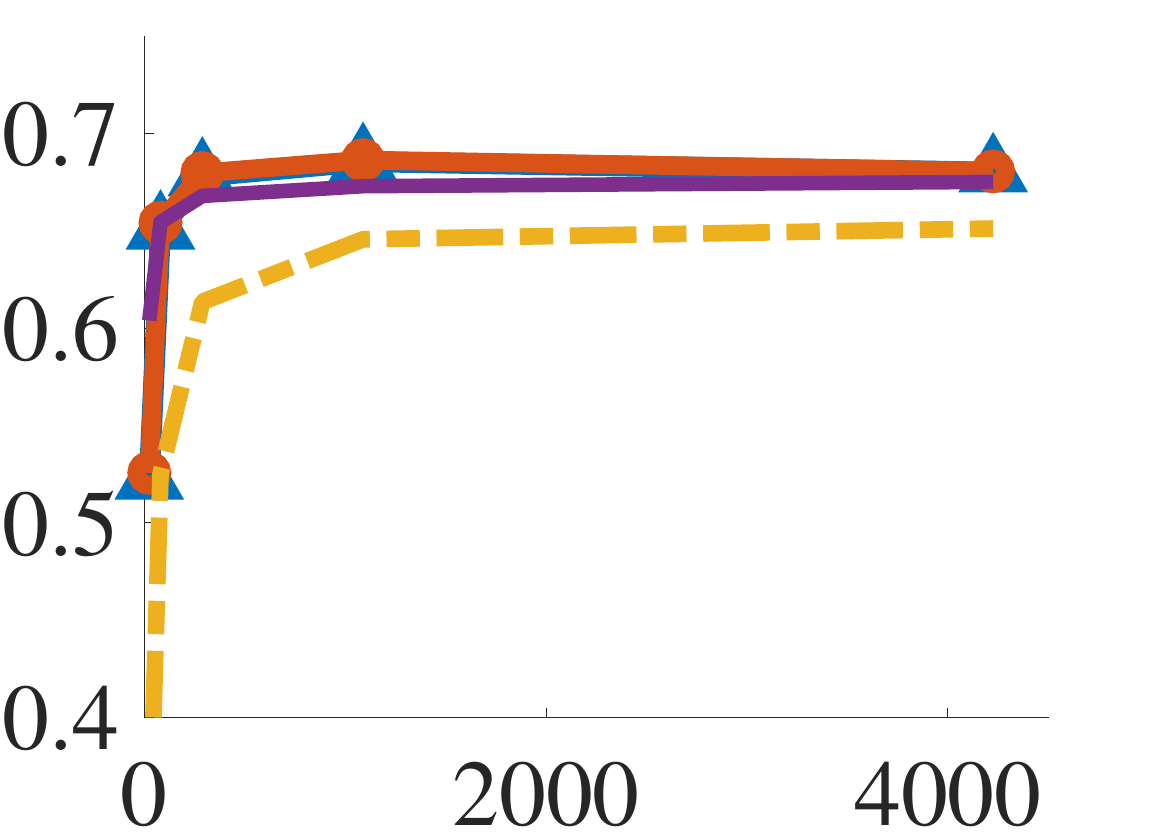}
\includegraphics[width=.225\linewidth, trim={30 190 25 200}, clip]{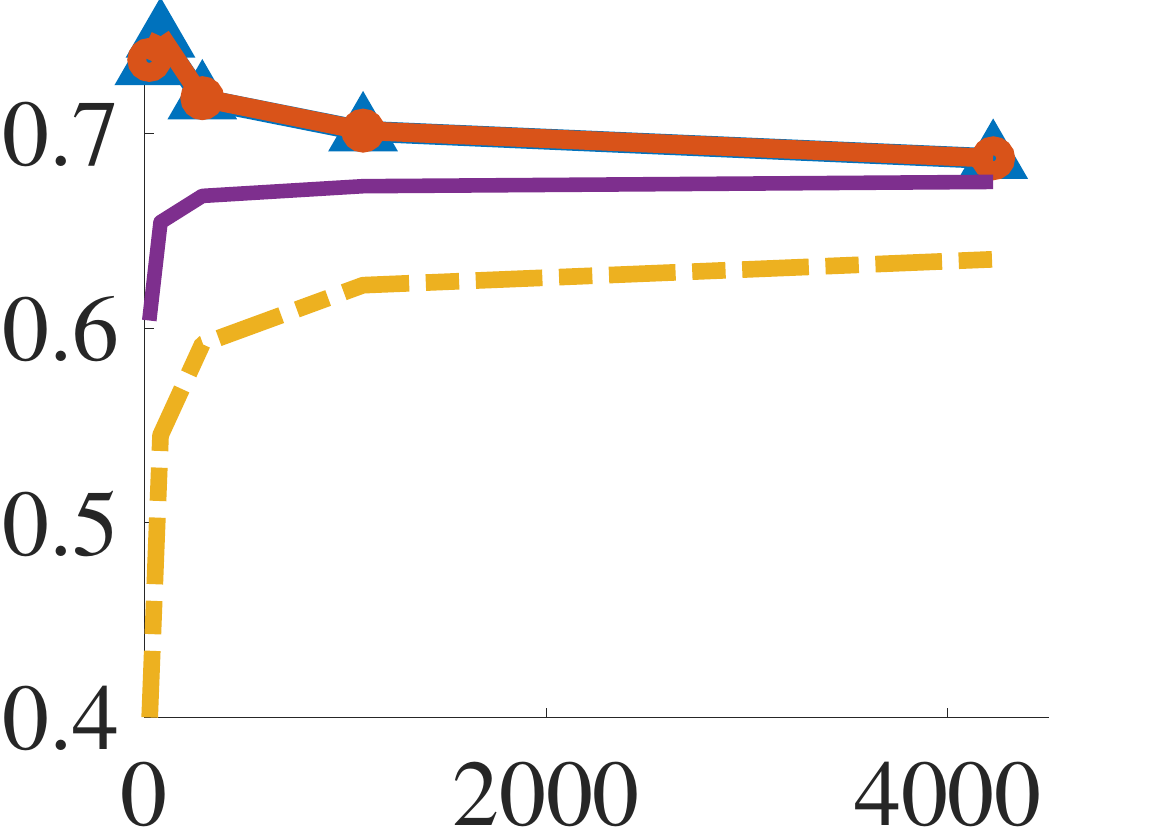}
\includegraphics[width=.225\linewidth, trim={30 190 25 200}, clip]{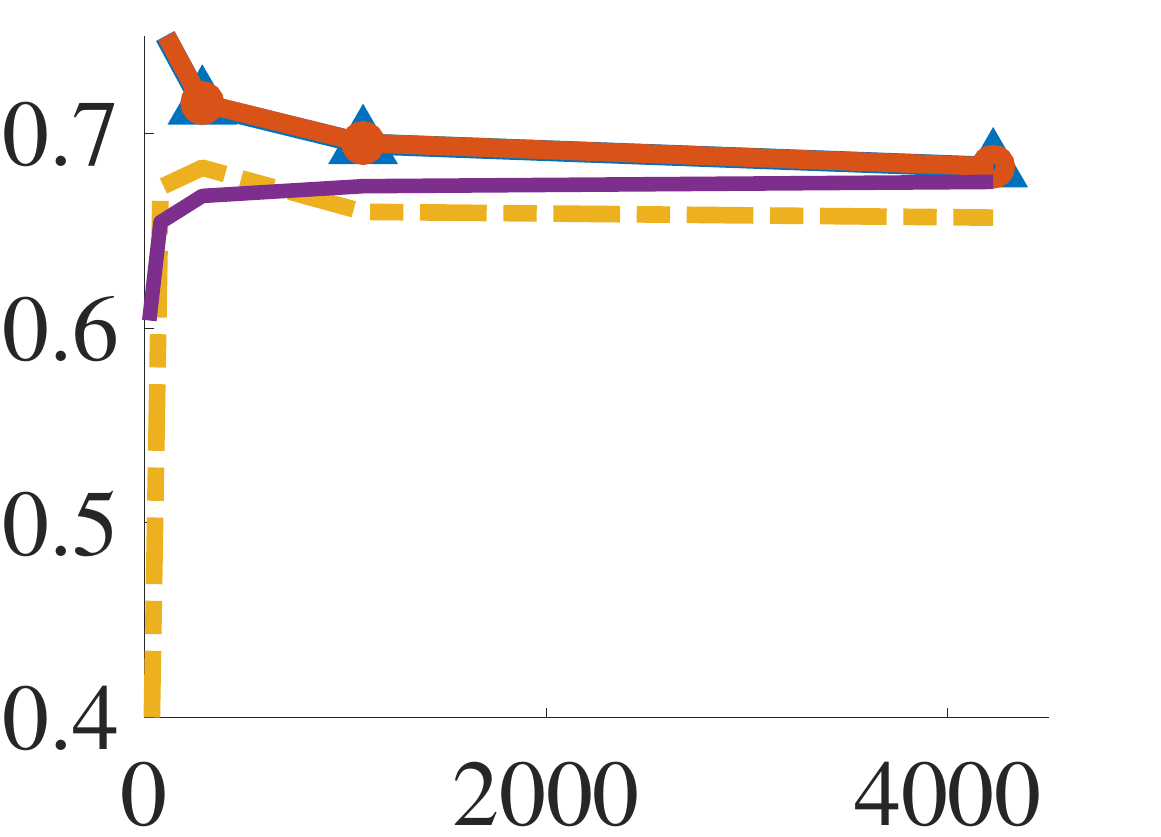}
\includegraphics[width=.225\linewidth, trim={30 190 25 200}, clip]{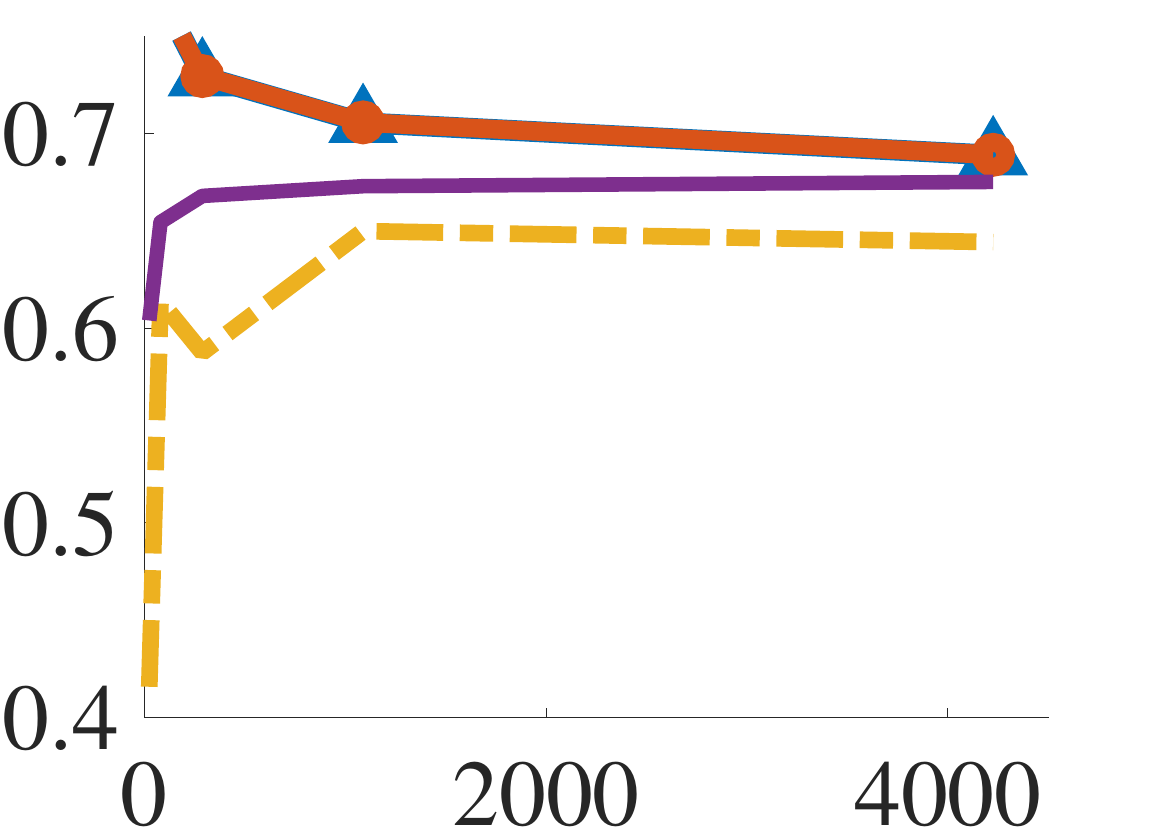}\\

\rotatebox{90}{$\qquad\;$ \textbf{$\nus = -1$} }
   \rotatebox{90}{$\qquad$ Disp. (cm) }
\includegraphics[width=.225\linewidth, trim={30 190 25 200}, clip]{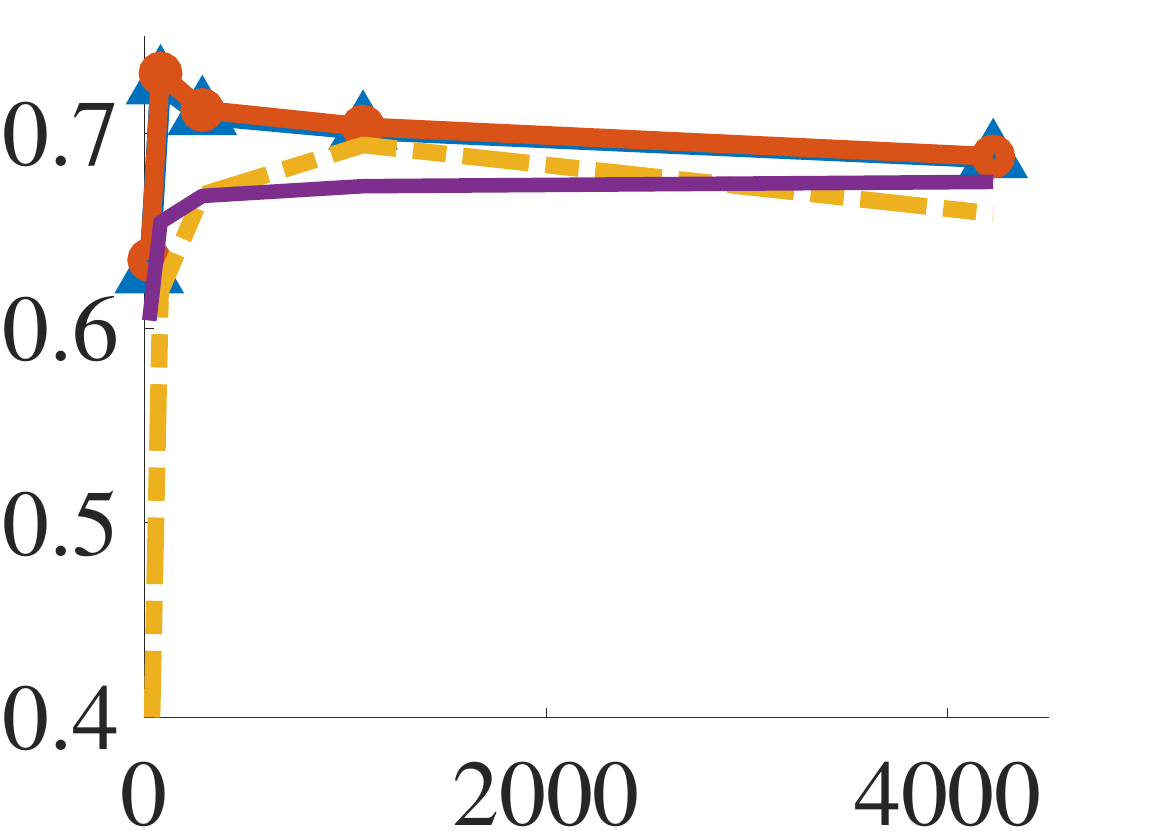}
\includegraphics[width=.225\linewidth, trim={30 190 25 200}, clip]{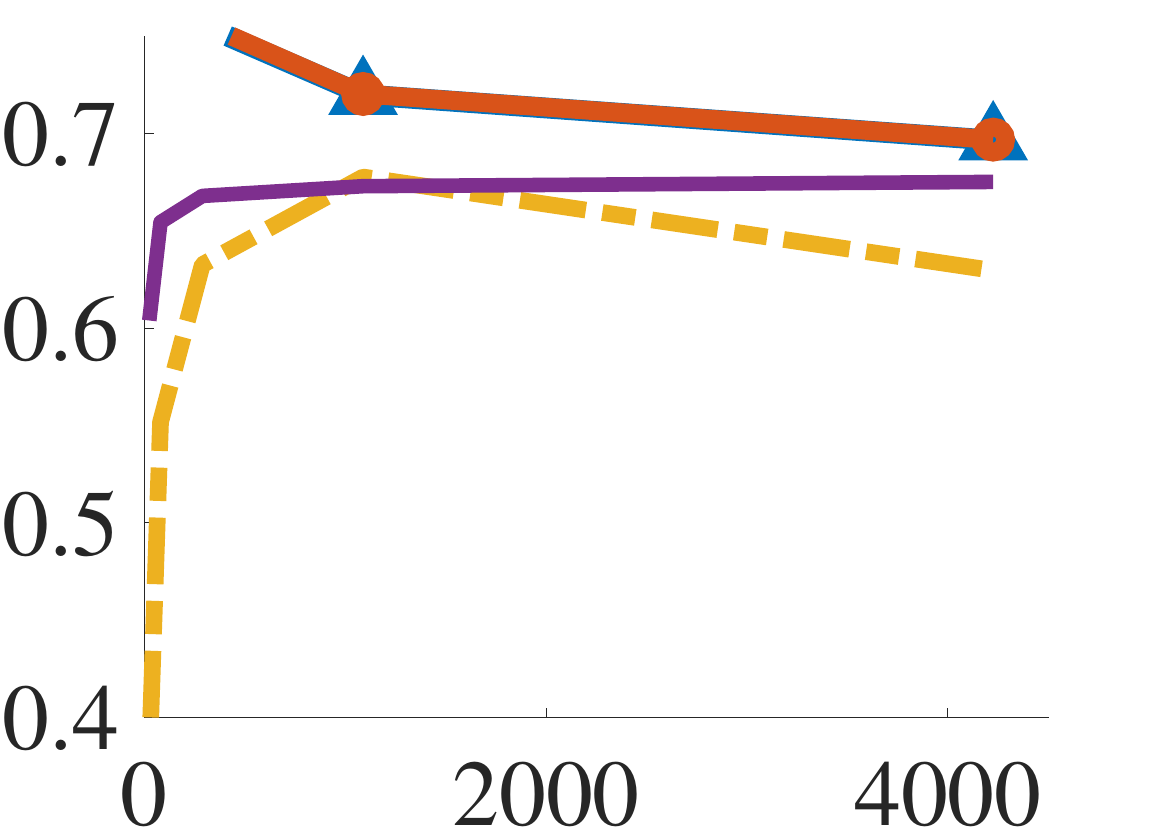}
\includegraphics[width=.225\linewidth, trim={30 190 25 200}, clip]{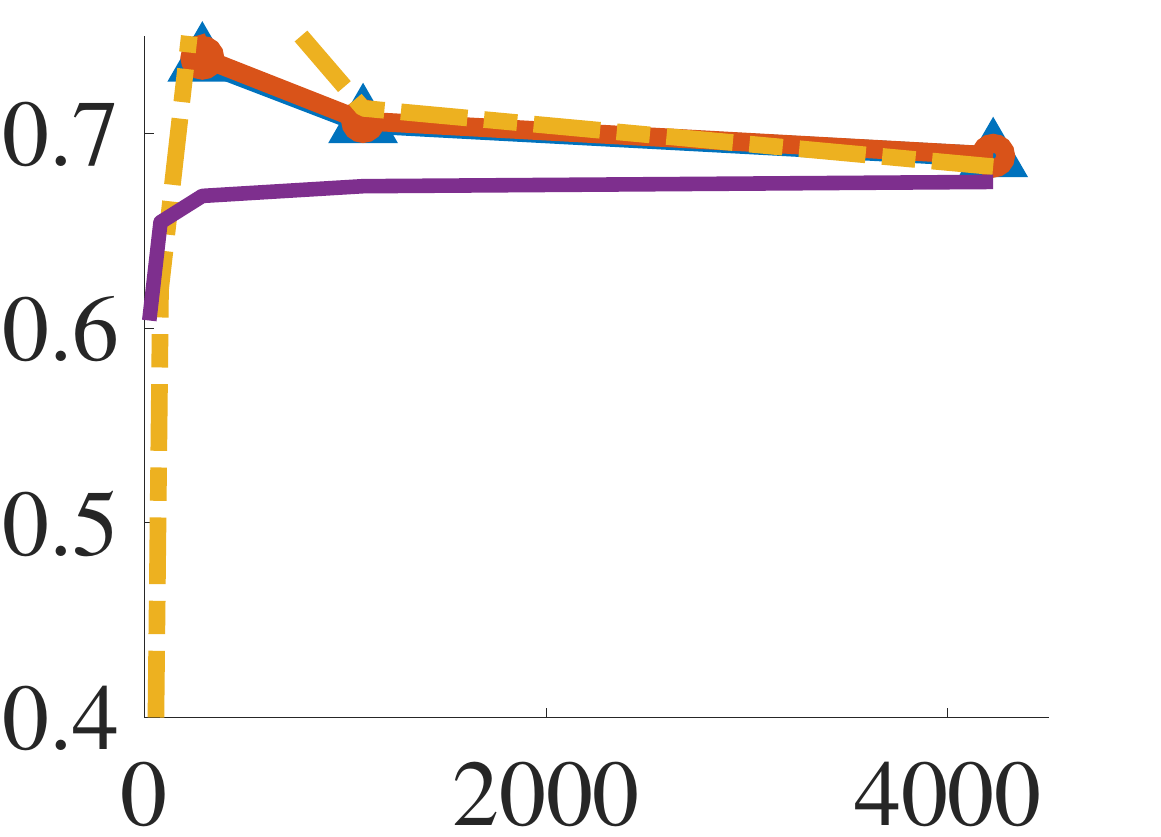}
\includegraphics[width=.225\linewidth, trim={30 190 25 200}, clip]{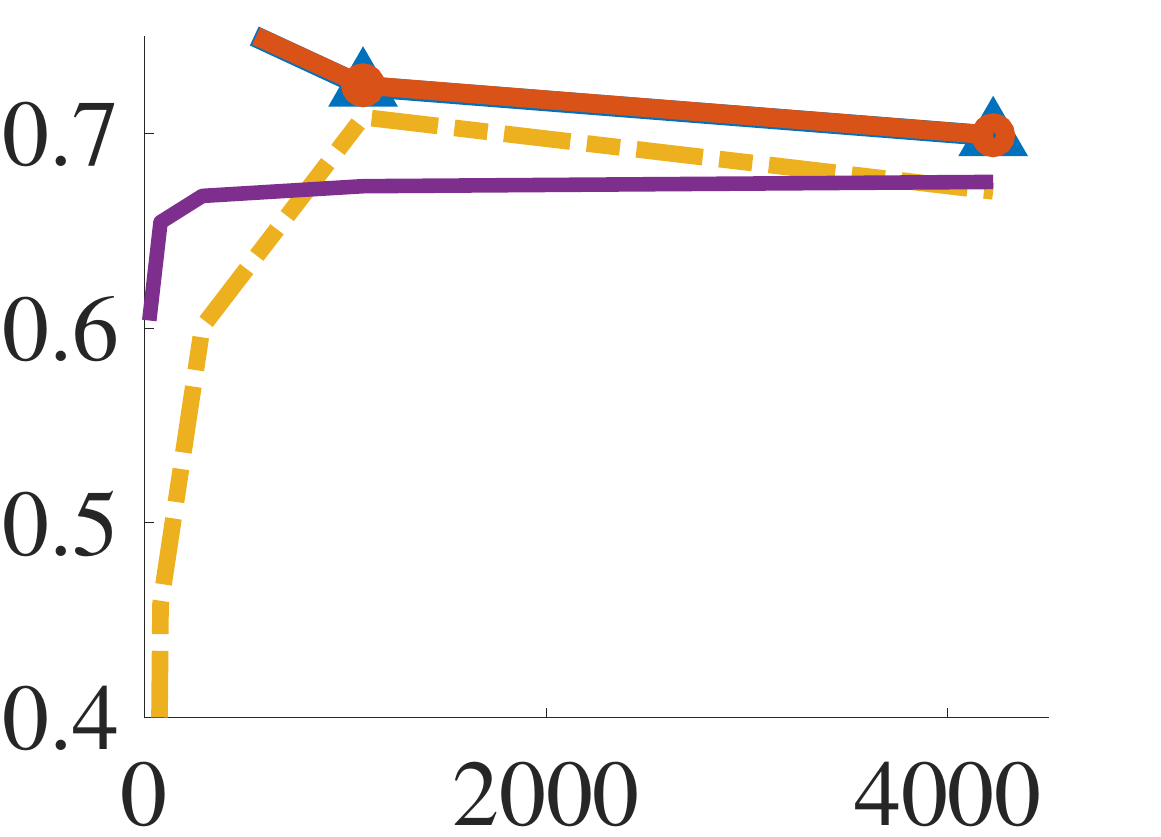}\\

$\qquad\qquad\quad$ \# Solid DOF $\qquad\qquad\quad\;$ \# Solid DOF $\qquad\qquad\quad$ \# Solid DOF $\qquad\qquad\quad\;$ \# Solid DOF
\caption{Corner $y$-displacement for different number of solid DOF for the Cook's membrane benchmark (Section \ref{Cook's Membrane}) for different choices of elements and numerical Poisson ratios. The solid DOF range from $m = 25$ to $4225$. Notice that each row has the same extents. If a value of $\nus$ is close to $\frac{1}{2}$, low order elements produce volumetric locking, and higher order elements are needed for convergence at reasonable numbers of DOF.}
\label{cooks_disp}
\end{figure}

\begin{figure}
$\qquad\qquad\qquad\;\;\;\;$ \textbf{P1} $\qquad\qquad\qquad\qquad\quad$  \textbf{Q1} $\qquad\qquad\qquad\qquad\;\;\;$  \textbf{P2} $\qquad\qquad\qquad\qquad\quad\;$ \textbf{Q2}\\
\rotatebox{90}{$\quad$ \textbf{$\nus = .49995$} }
   \rotatebox{90}{$\quad$ Vol Change \% }
\includegraphics[width=.225\linewidth, trim={30 190 25 200}, clip]{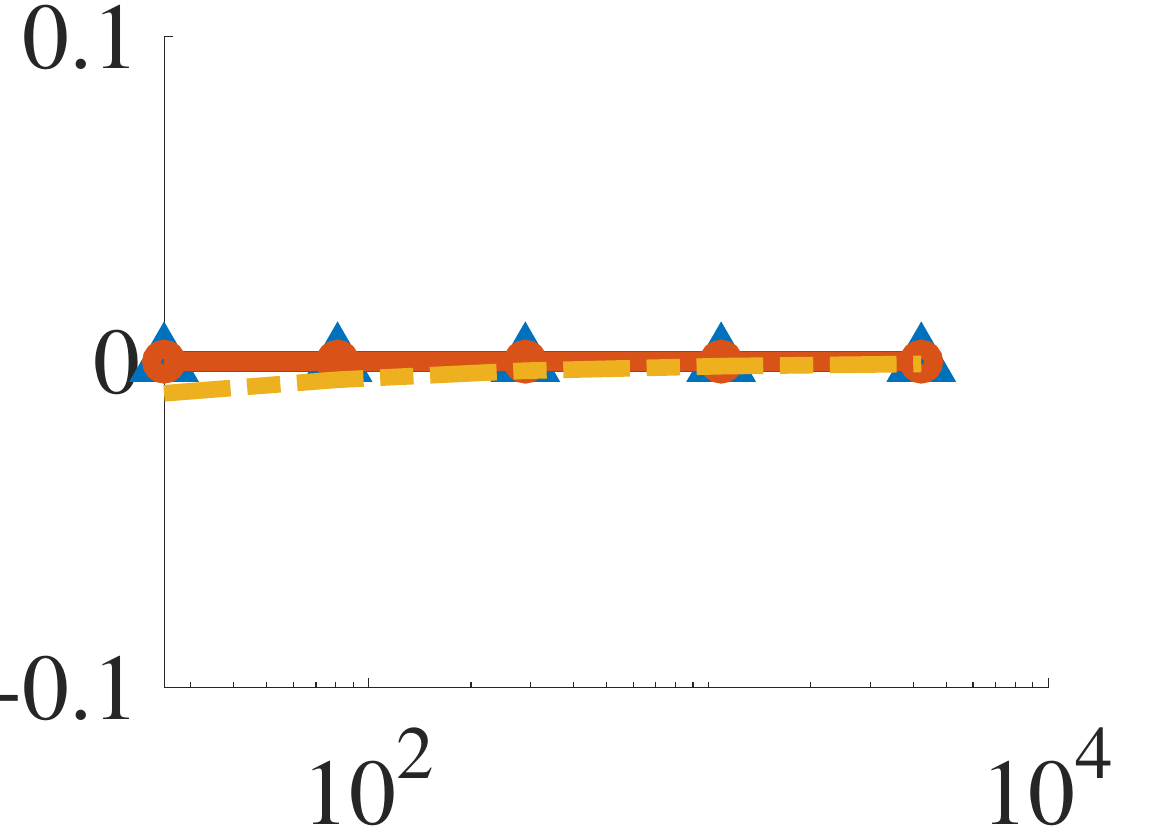} 
\includegraphics[width=.225\linewidth, trim={30 190 25 200}, clip]{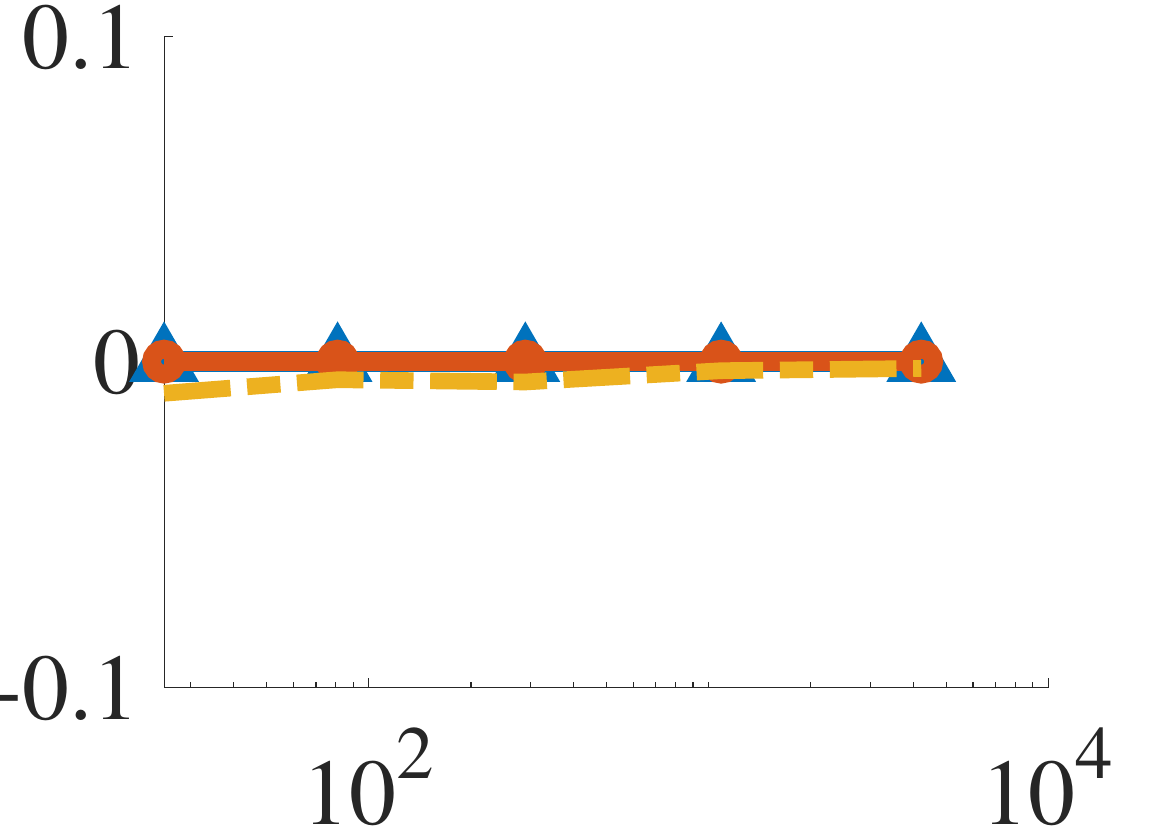} 
\includegraphics[width=.225\linewidth, trim={30 190 25 200}, clip]{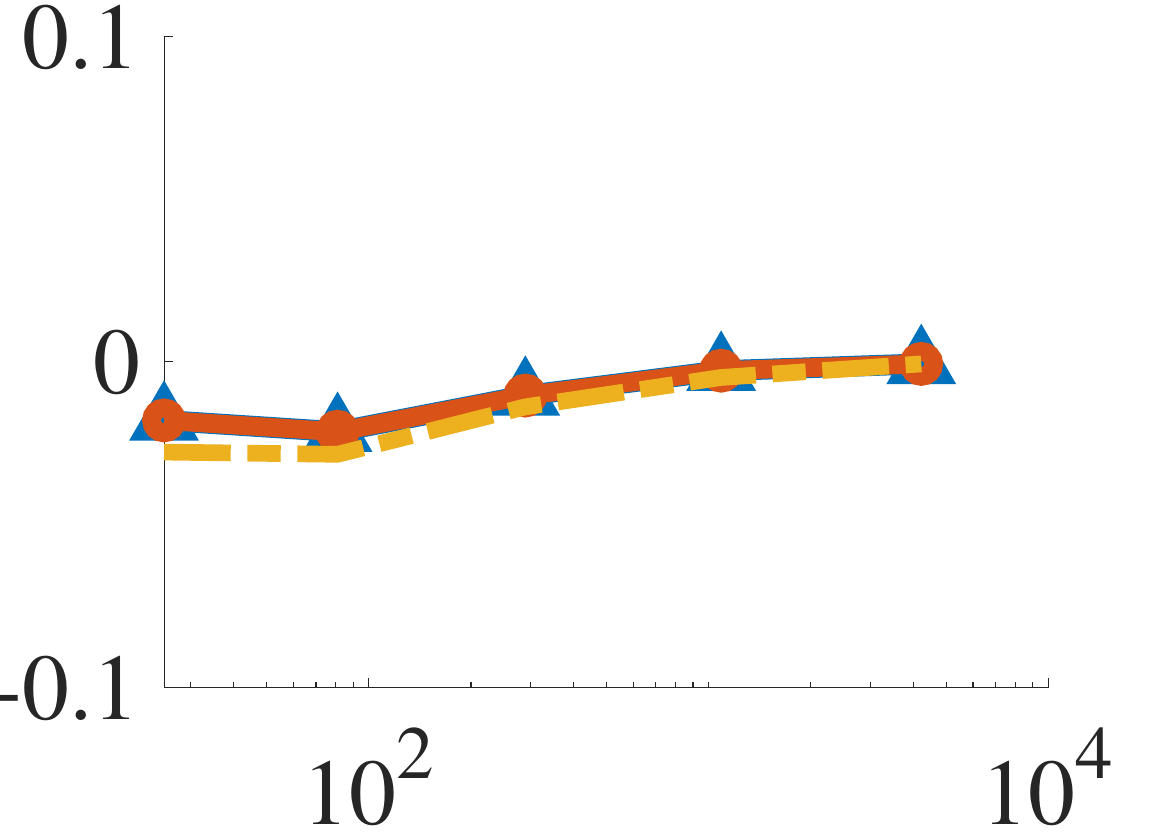} 
\includegraphics[width=.225\linewidth, trim={30 190 25 200}, clip]{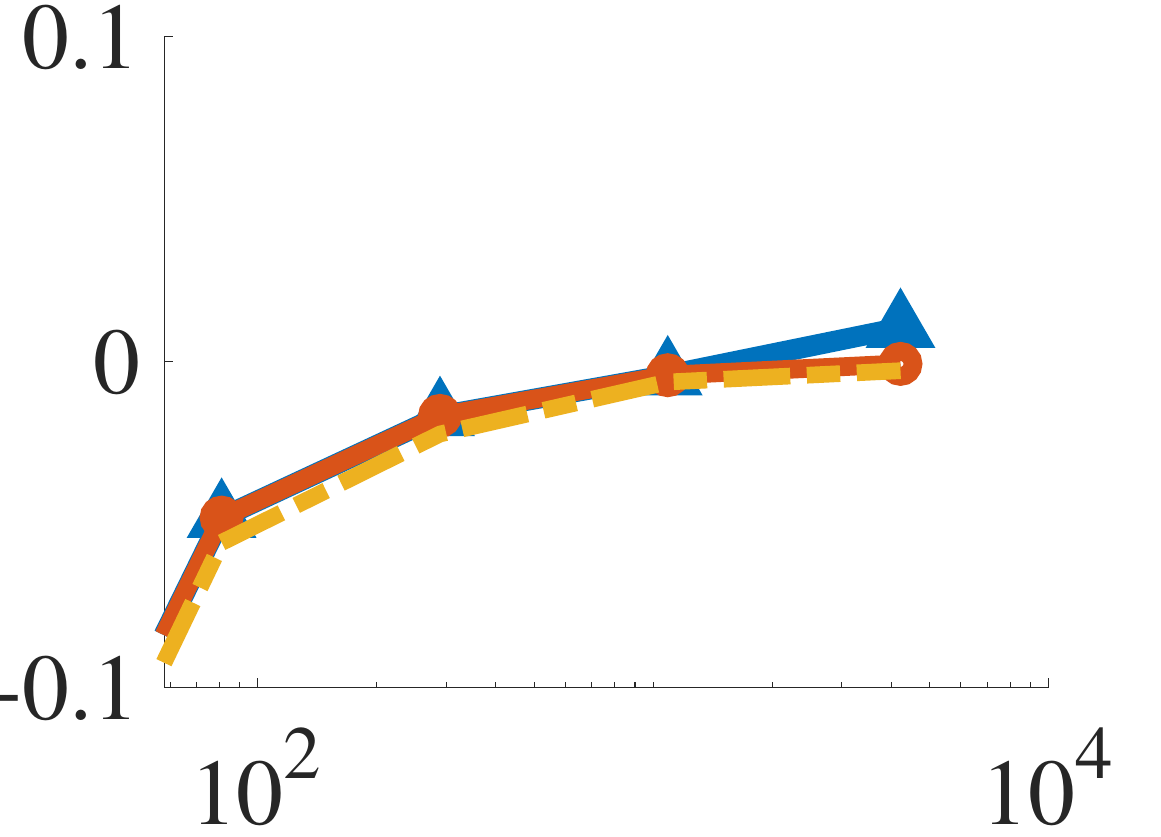} \\
\rotatebox{90}{$\qquad\;$ \textbf{$\nus = .4$} }
   \rotatebox{90}{$\quad$ Vol Change \% }
\includegraphics[width=.225\linewidth, trim={30 190 25 200}, clip]{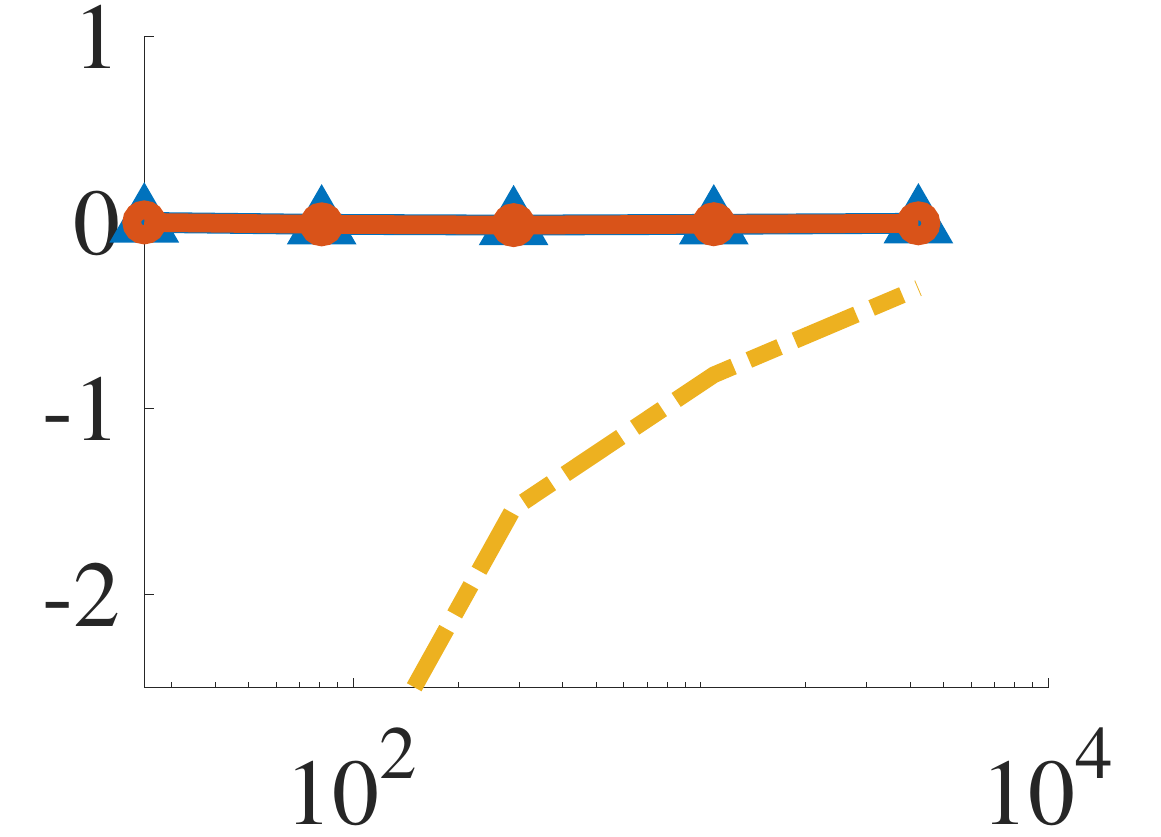} 
\includegraphics[width=.225\linewidth, trim={30 190 25 200}, clip]{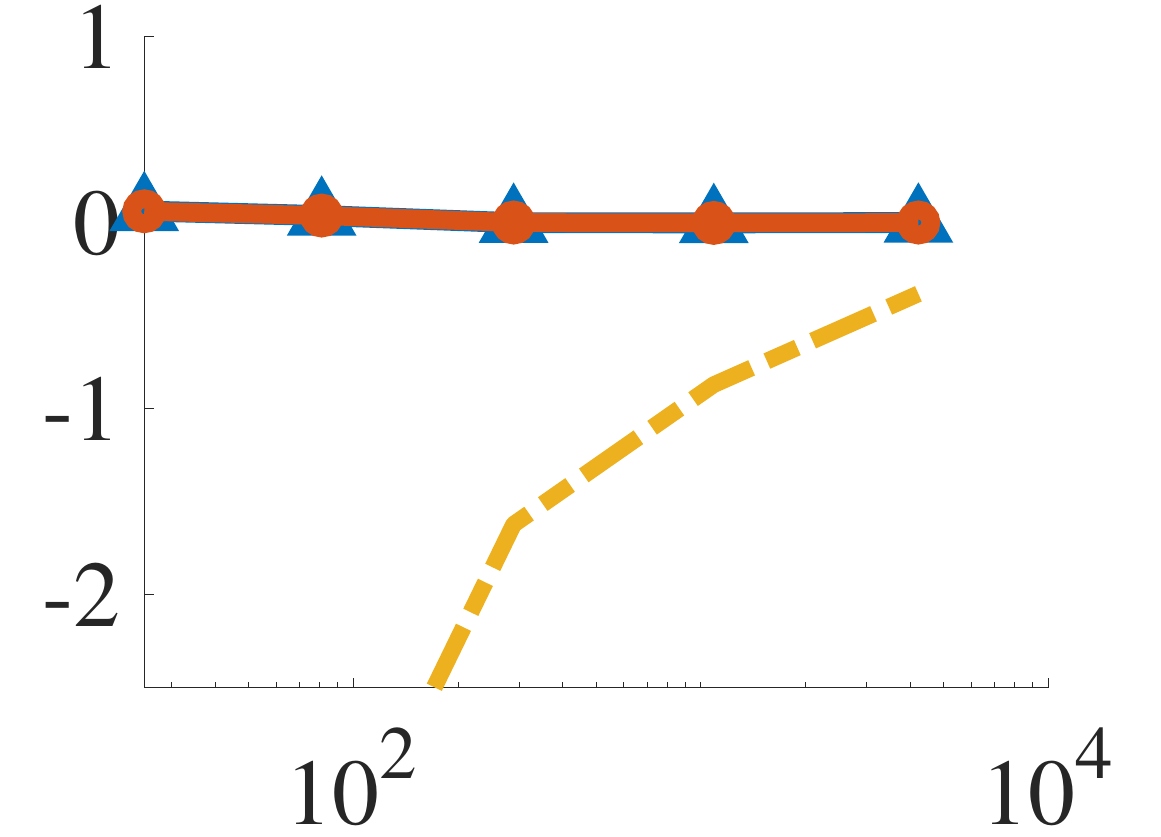}
\includegraphics[width=.225\linewidth, trim={30 190 25 200}, clip]{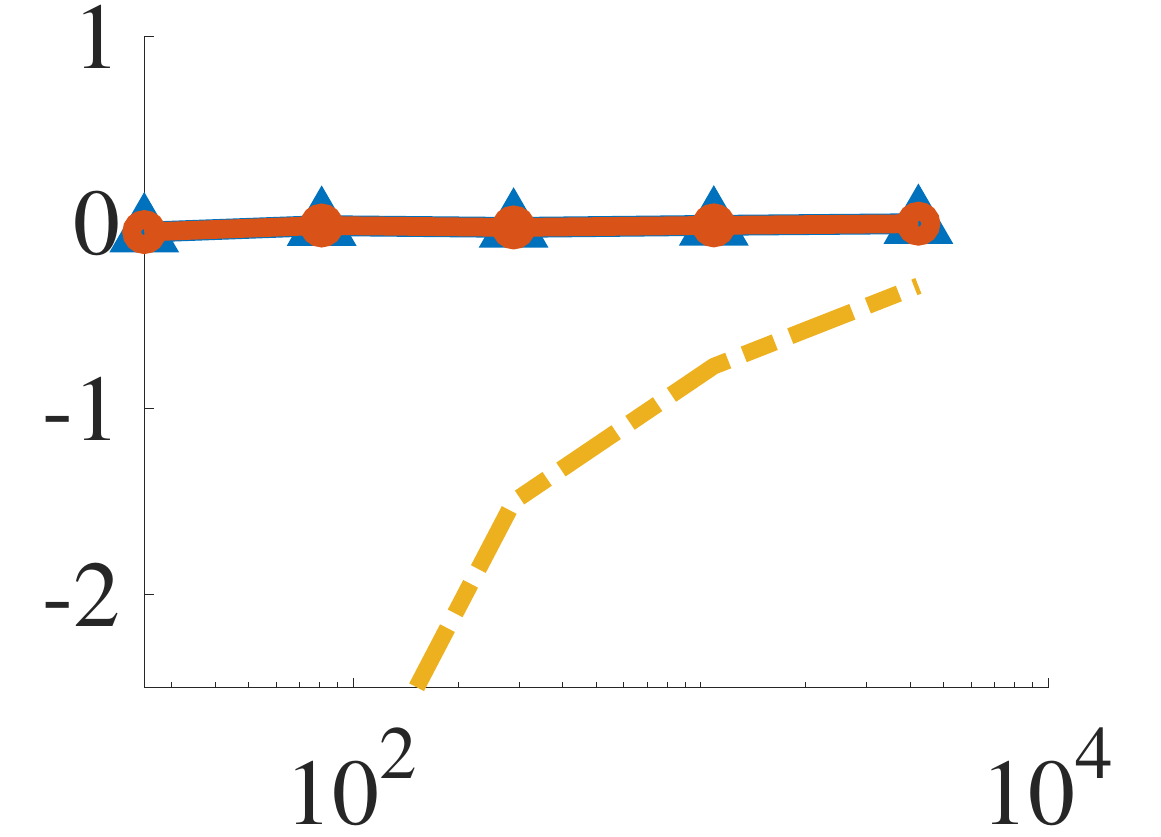}
\includegraphics[width=.225\linewidth, trim={30 190 25 200}, clip]{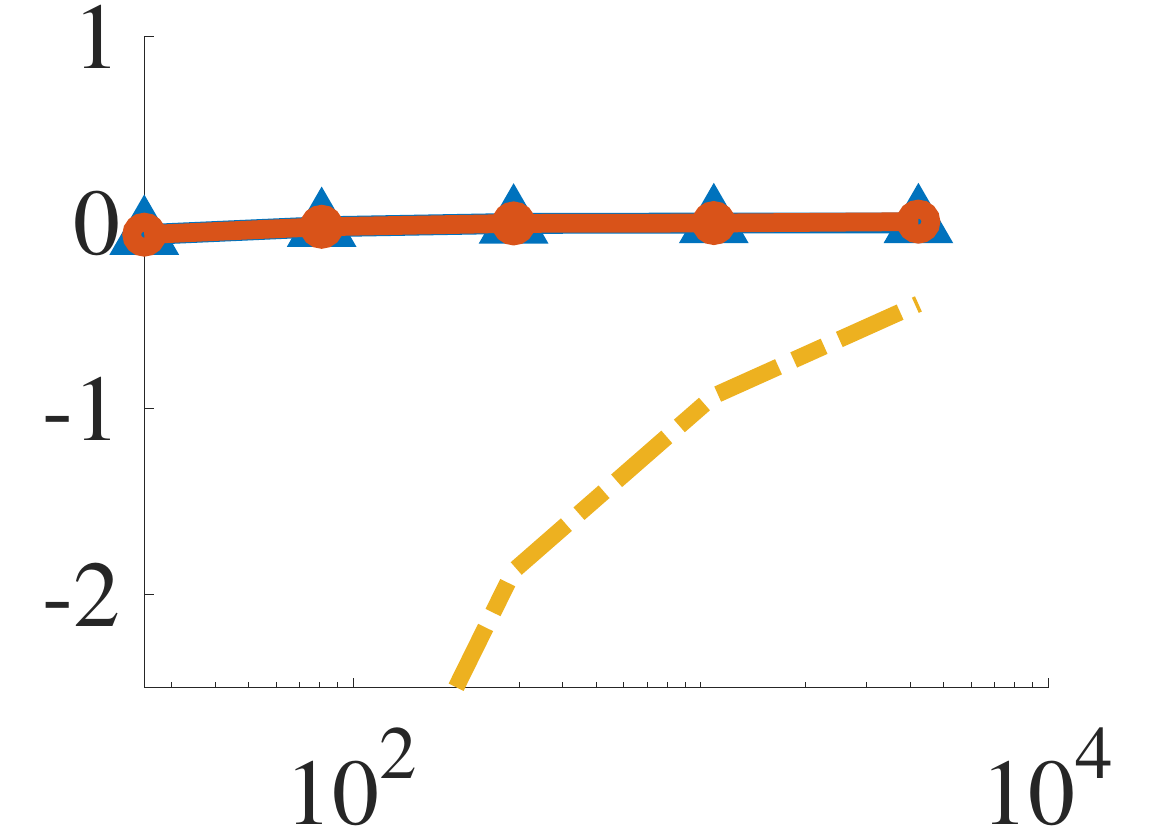}\\
\rotatebox{90}{$\qquad\;$ \textbf{$\nus = 0$} }
   \rotatebox{90}{$\quad$ Vol Change \% }
\includegraphics[width=.225\linewidth, trim={30 190 25 200}, clip]{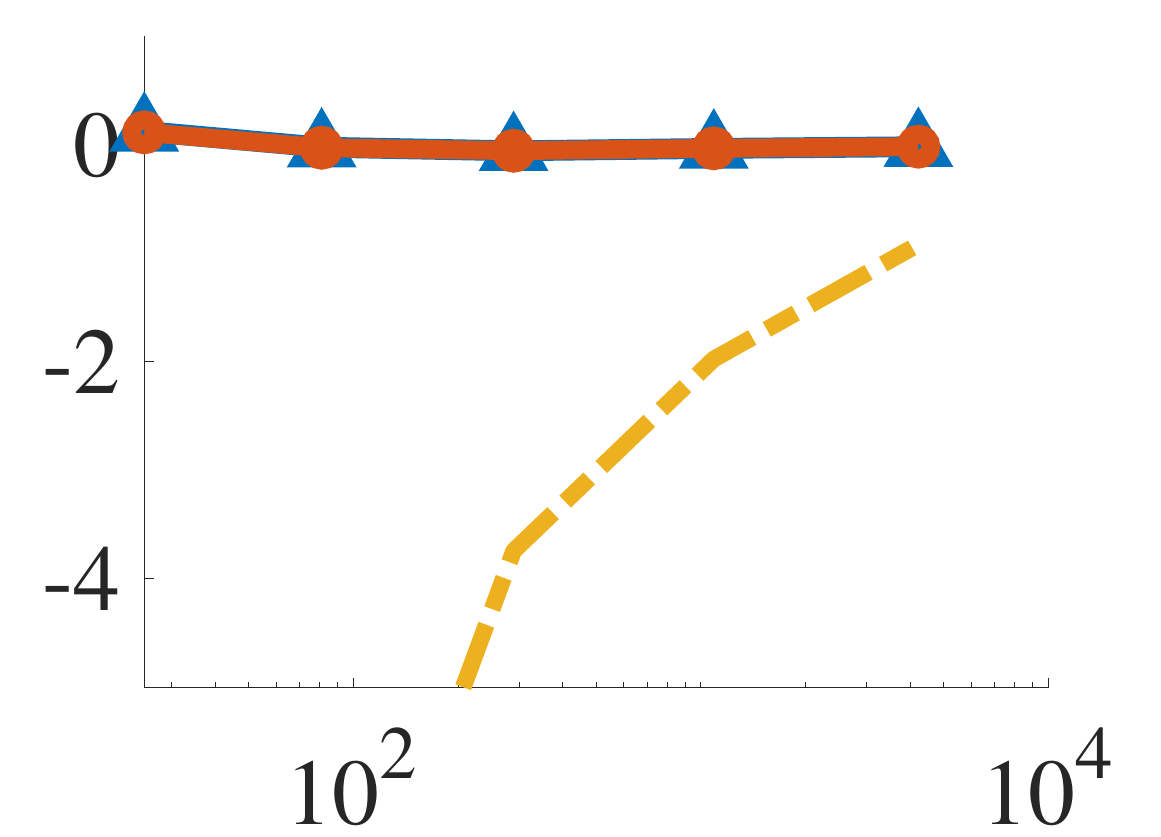} 
\includegraphics[width=.225\linewidth, trim={30 190 25 200}, clip]{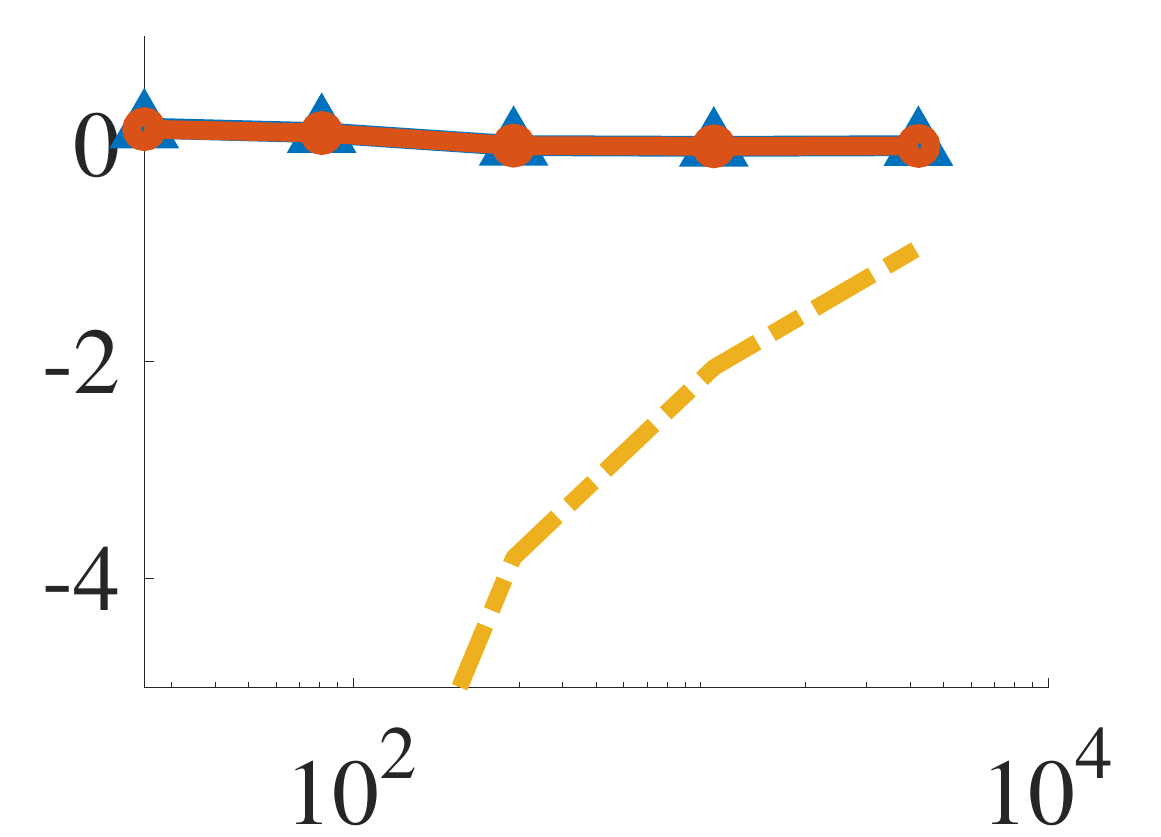} 
\includegraphics[width=.225\linewidth, trim={30 190 25 200}, clip]{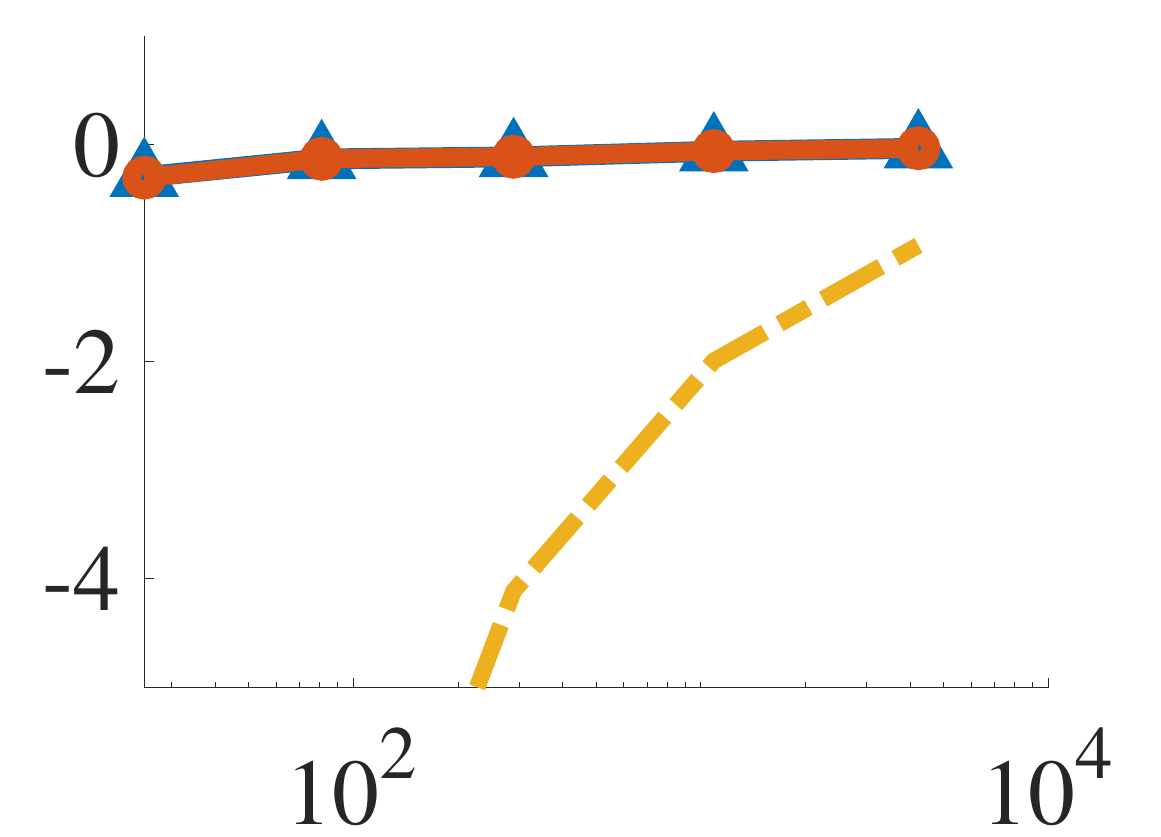} 
\includegraphics[width=.225\linewidth, trim={30 190 25 200}, clip]{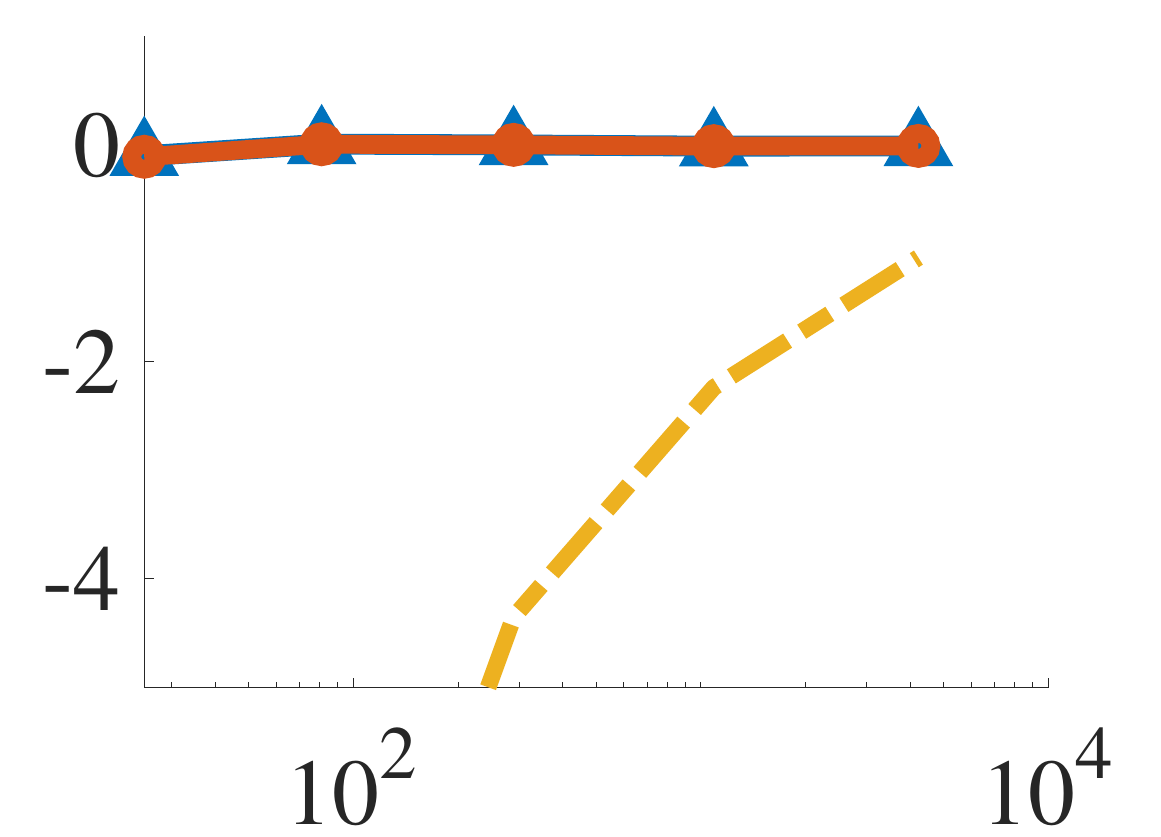}\\
\rotatebox{90}{$\qquad$ \textbf{$\nus = -1$} }
   \rotatebox{90}{$\quad$ Vol Change \% }
\includegraphics[width=.225\linewidth, trim={30 190 25 200}, clip]{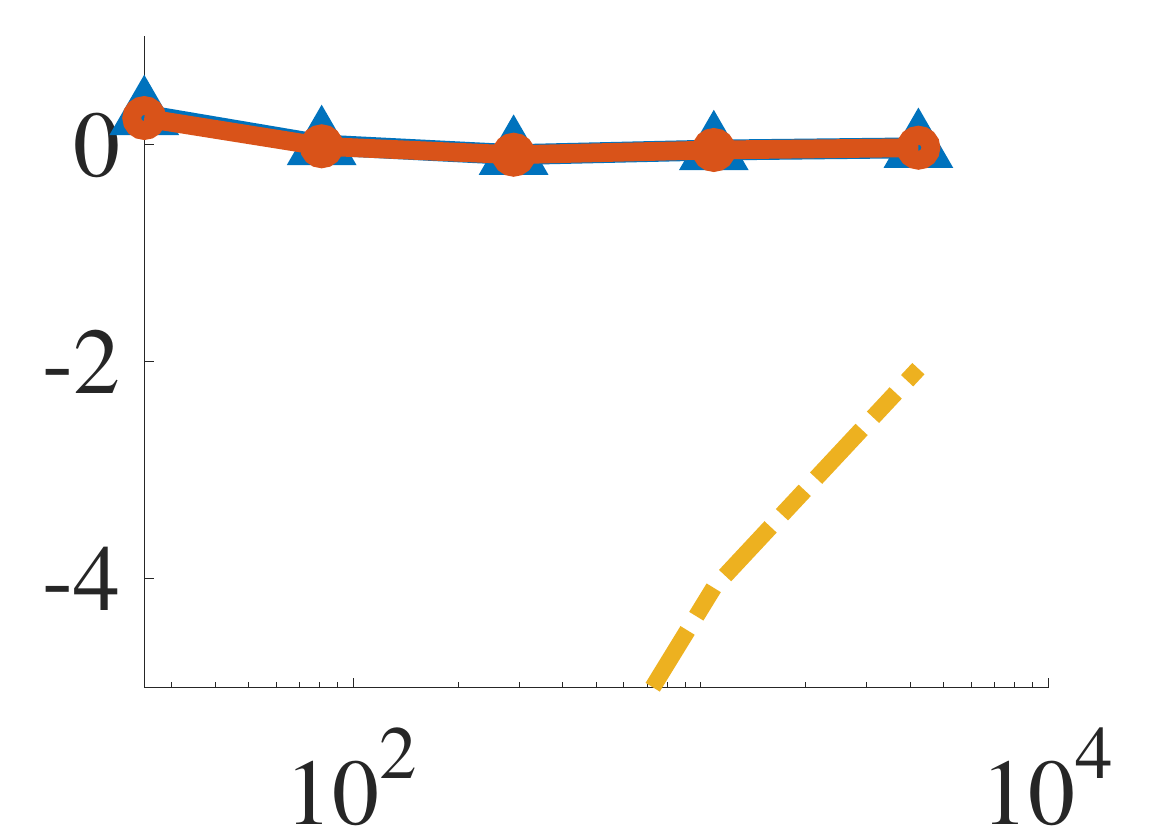} 
\includegraphics[width=.225\linewidth, trim={30 190 25 200}, clip]{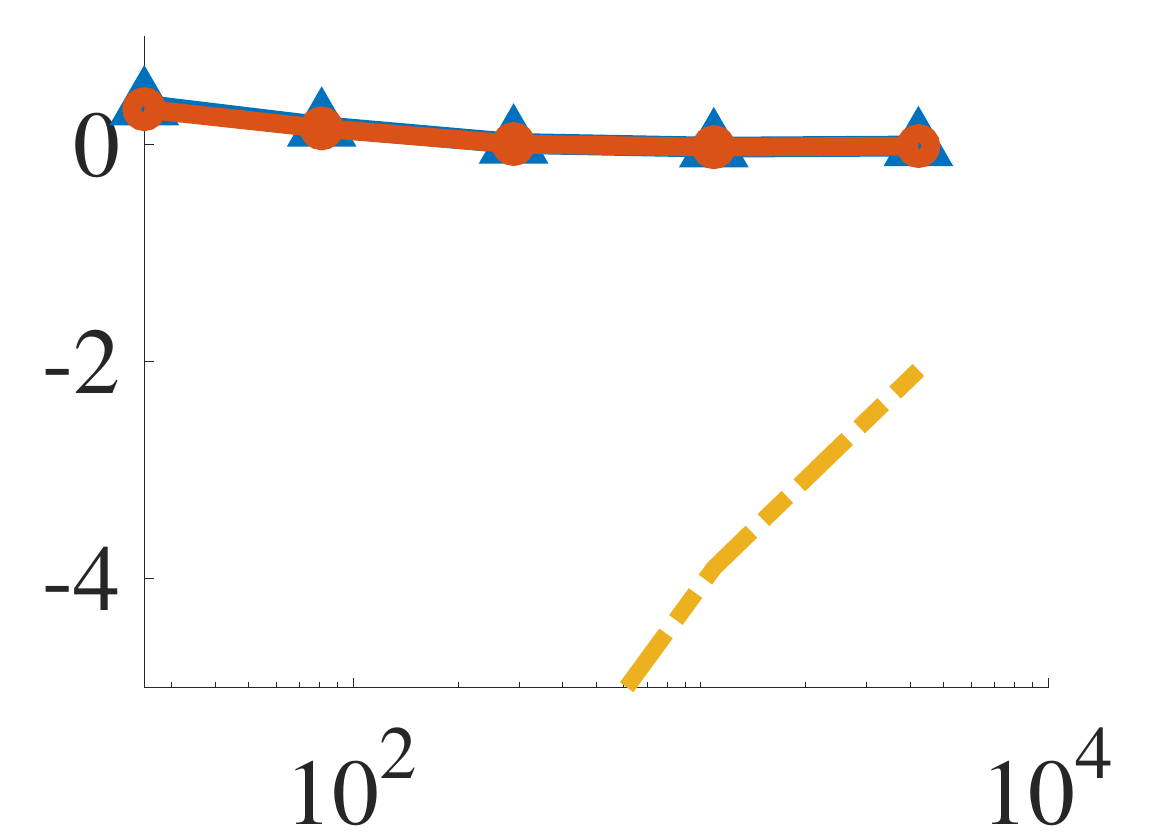} 
\includegraphics[width=.225\linewidth, trim={30 190 25 200}, clip]{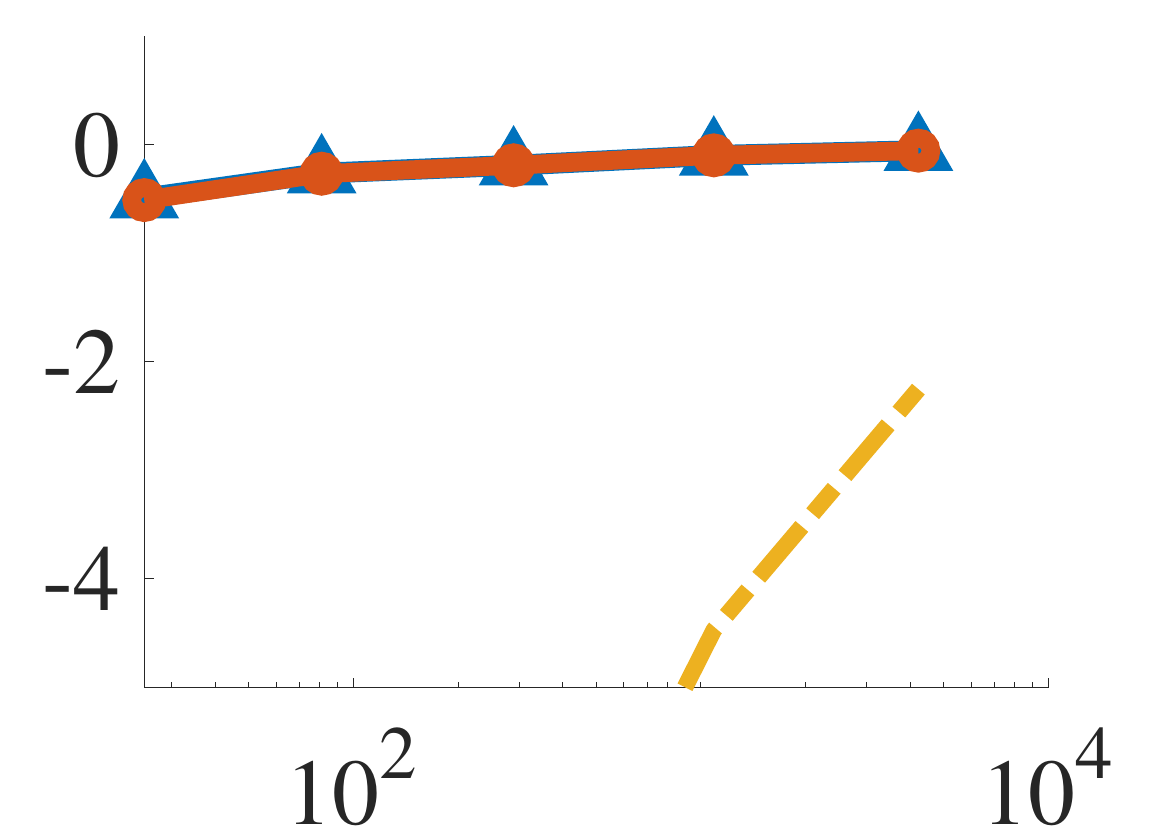} 
\includegraphics[width=.225\linewidth, trim={30 190 25 200}, clip]{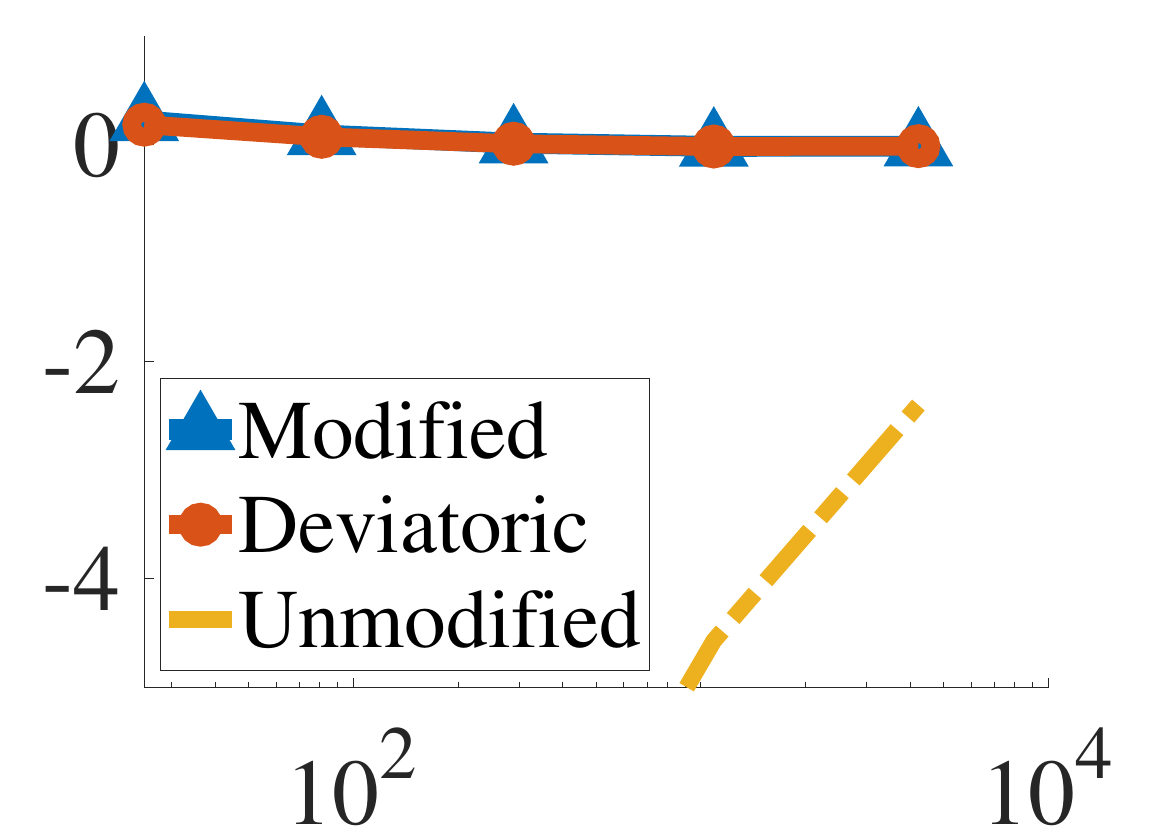}\\

$\qquad\qquad\quad$ \# Solid DOF $\qquad\qquad\quad\;$ \# Solid DOF $\qquad\qquad\quad$ \# Solid DOF $\qquad\qquad\quad\;$ \# Solid DOF
\caption{Percent change in total area for different numbers of solid DOF for the Cook's membrane benchmark (Section \ref{Cook's Membrane}) after deformation. The DOF range from $m = 25$ to $4225$, and the $x$ axis is on a log scale. Note the different scales on the $y$-axes. Omitting the coarsest discretizations ($m = 25$), the largest deviations in total area among all element types used are approximately $.10 \%$ for the modified case, $7.45 \%$ for the unmodified case, and $.10 \%$ for the deviatoric case.}
\label{cooks_area}
\end{figure}

Cook's membrane is a classical plane strain problem involving a swept and tapered quadrilateral. The dimensions of the solid domain and overall problem specification are shown in Figure (\ref{cooks}). This benchmark was first proposed by Cook \etal~\cite{RDCook1974} and is common in testing numerical methods for incompressible elasticity. An upward loading traction is applied to the right side, and the left hand is fixed in place; see Figure ($\ref{cooks}$). All other structural boundaries have stress-free boundary conditions applied. The upward traction is given as $6.25 $~$\frac{\text{dyn}}{\text{cm}^2}$. The $y$-displacement of the top right corner is measured at $T_{\text{f}} = 50$ s. The load time is $T_{\text{l}} = 20$ s. The neo-Hookean material model, equations ($\ref{nh_energy}$) -- ($\ref{nh_stress_dev}$), is used with a shear modulus of $G = 83.3333$ $\frac{\text{dyn}}{\text{cm}^2}$; this value is equivalent to using a Young's modulus of $E = 250 $~$\frac{\text{dyn}}{\text{cm}^2}$ if $\nu = \frac{1}{2}$. Damping is set to $\eta = 4.16667 \, \frac{\text{g}}{\text{s}}$ for this test. The computational domain is $\Omega = [0, L]^2$ with $L = 10 \, \text{cm}$. The numbers of solid DOF range from $m = 25$ to $m = 4225$. As was the case for the compressed block benchmark, we use a sequence of meshes that yields the same node locations for all element types considered.\\
\indent Deformations and results for this benchmark are shown in Figures ($\ref{cm}$) -- ($\ref{cooks_area}$). Figure ($\ref{cm}$) shows the deformations of the structure along with the elemental Jacobian determinant $J$. Note that, at least qualitatively, the deformation in the unmodified and unstabilized case is unphysical; see Figure ($\ref{cm}$d). We emphasize that the case of unmodified invariants and zero volumetric energy are the only cases in which this unphysical behavior is observed. Figure (\ref{cm_dev}) demonstrates that the deviatoric stresses for the IB computations are in agreement with those from the FE method when modified invariants and volumetric stabilization are used. As shown in Figure ($\ref{cooks_disp}$), most cases converge to the benchmark solution. With values of $\nus$ close to $\frac{1}{2}$, the solution exhibits volumetric locking as in a typical displacement based FE formulation. Note the unmodified case in the final row of plots in Figure (\ref{cooks_disp}) and that the displacement of the single point of interest does not seem to behave well for these cases. \\
\indent Figure (\ref{cooks_area}) shows the percent change in the total area of the mesh after deformation. It is clear that the modified invariants and deviatoric projection yield improved results in terms of global area conservation in comparison to the unmodified invariants. This effect becomes more pronounced as the numerical bulk modulus is decreased. It may appear as though the modified invariants and deviatoric projection cases have zero volume change, but this is not the case. The percent change in total volume for all elements considered rang between $.000021 \%$ and $.10 \%$ for modified invariants and between $0 \%$ and $.10 \%$ for the deviatoric projection. For the unmodified cases, this range was between $.000065 \%$ and $7.45 \%$.

%%%%%%%%%%%%%%%%

\subsection{Anisotropic Cook's Membrane}
\label{Anisotropic Cook's Membrane}
\begin{figure}
\centering
\includegraphics[width=8cm, trim={50 100 70 70}, clip]{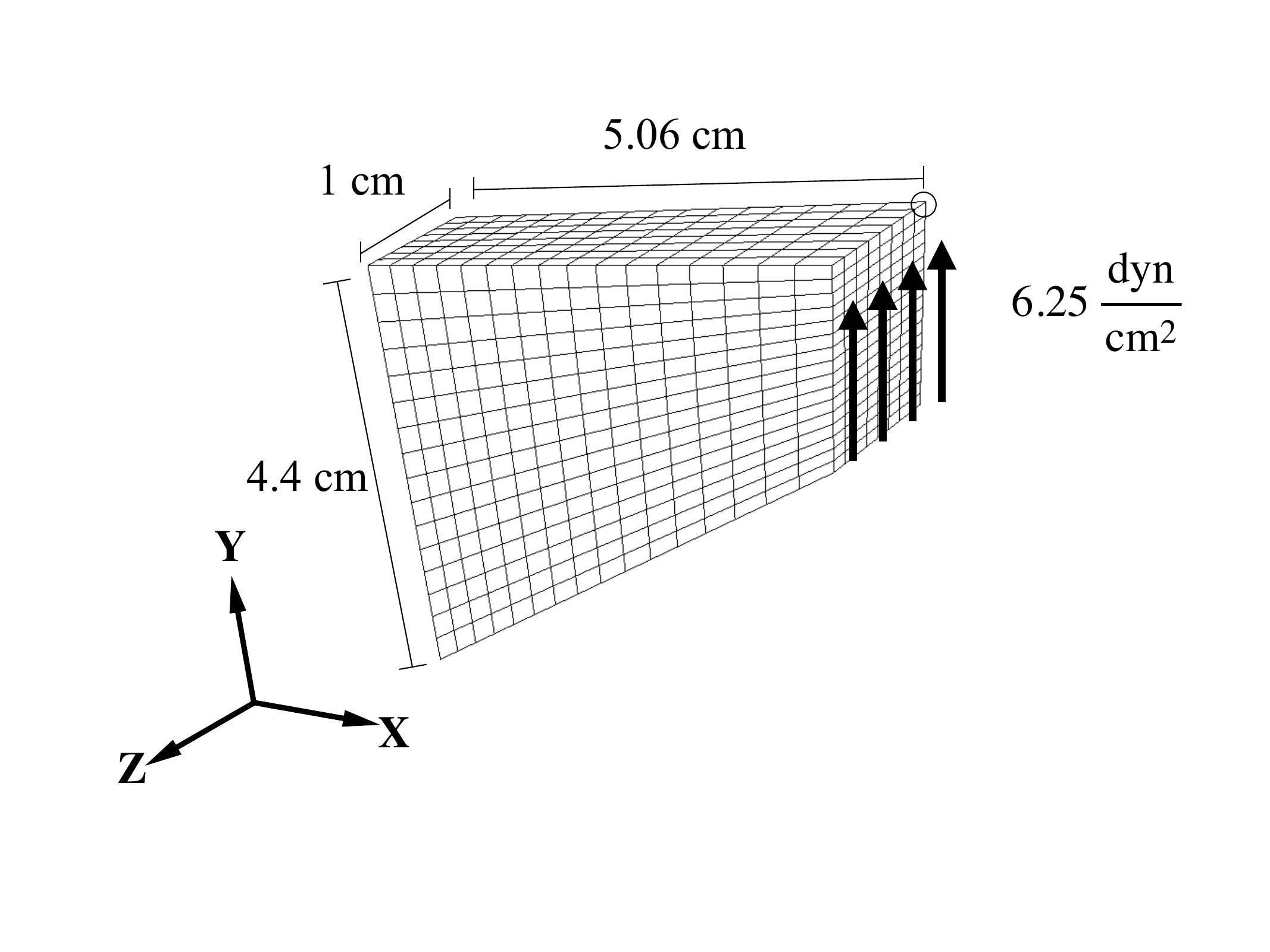}
\caption{Specifications of the anisotropic Cook's membrane benchmark (Section \ref{Anisotropic Cook's Membrane}). Traction in the $y$-direction is applied to the smallest face, and the opposite face is kept fixed. The quantity of interest is $y$-displacement as measured at the encircled point. To simplify the diagram of this three dimensional test, we omit the computational domain in this figure describing the problem setup. In the IBFE model, however, the structure is contained within a computational domain with dimensions $\Omega = [0, L]^2$ with $L = 12 \ \text{cm}$, and the solid mesh is placed in the center of this domain. Zero fluid velocity is enforced along the boundary of $\Omega$.}
\label{aniso_mesh}
\end{figure}

\begin{figure}
\begin{tabular}{l c c}
& \textbf{Modified Invariants} & \textbf{Unmodified Invariants} \\
\rotatebox{90}{\qquad\qquad\qquad\qquad \textbf{$\nus = .4$} }&
\subcaptionbox{\label{sfig:testa}} {\includegraphics[width=.45\linewidth, trim={100 50 70 70}, clip]{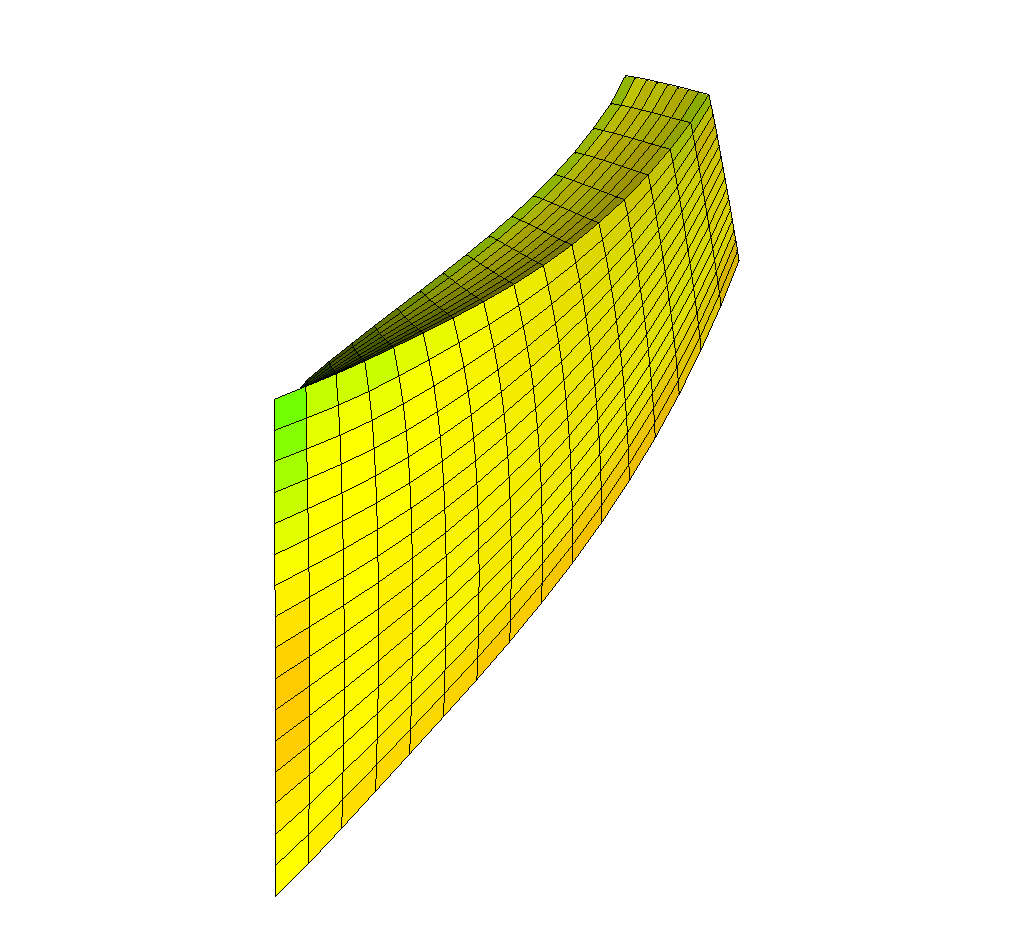}} &
\subcaptionbox{\label{sfig:testb}} {\includegraphics[width=.45\linewidth, trim={100 50 70 70}, clip]{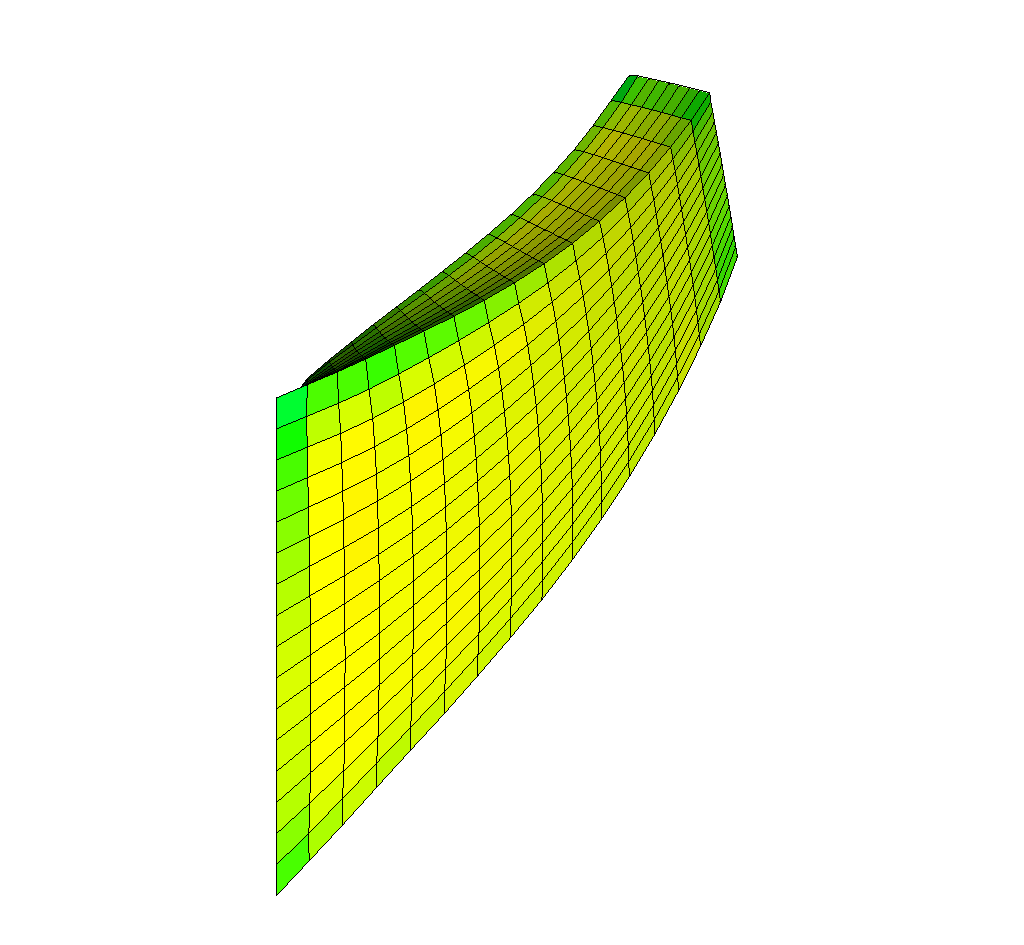}}\\

\rotatebox{90}{\qquad\qquad\qquad\qquad \textbf{$\nus = -1$}} &
\subcaptionbox{\label{sfig:testc}}{\includegraphics[width=.45\linewidth, trim={100 50 70 70}, clip]{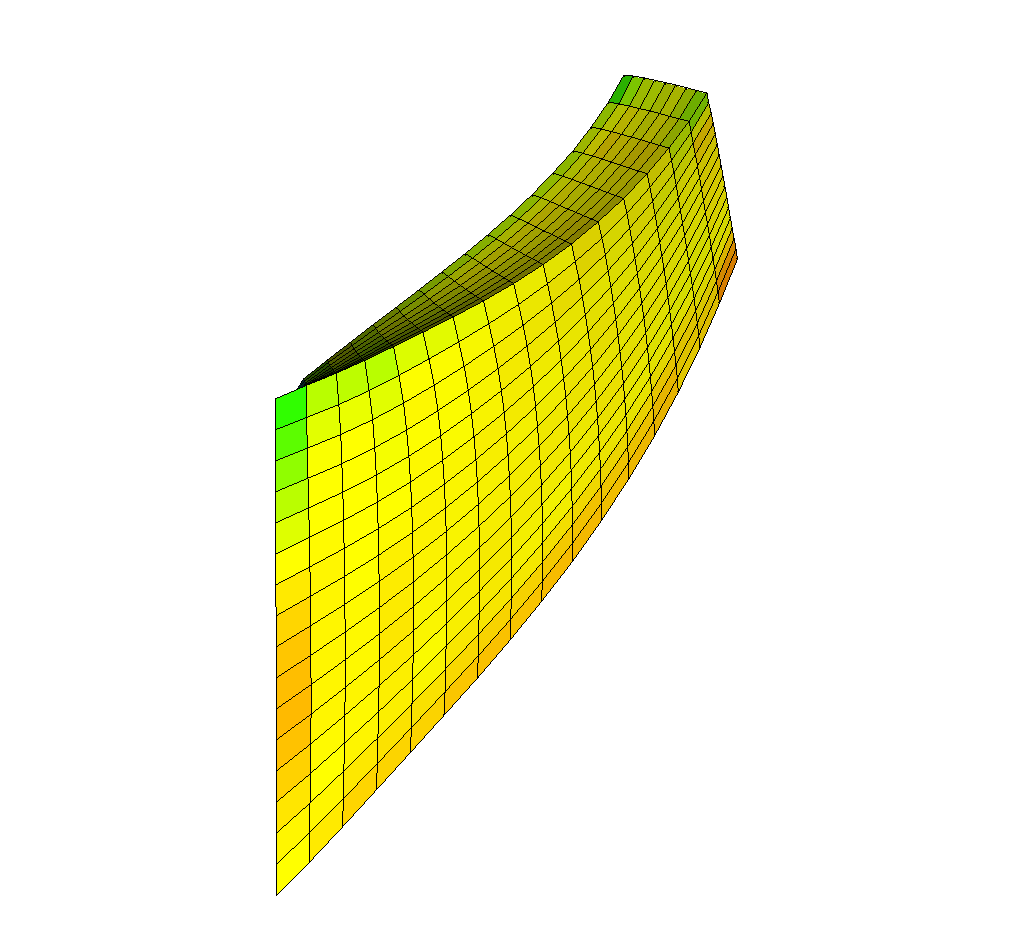}} &
\subcaptionbox{\label{sfig:testd}}{\includegraphics[width=.45\linewidth, trim={100 50 70 70}, clip]{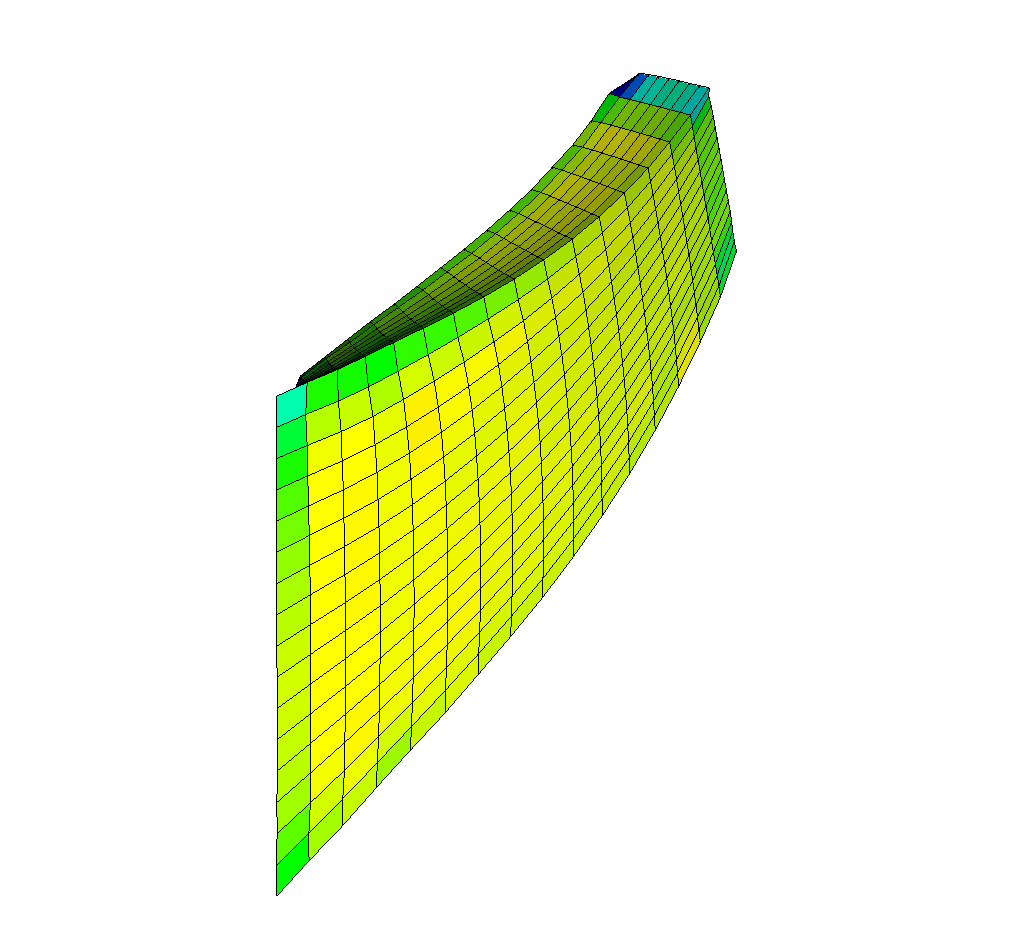}}  \\
\end{tabular}
%trim={left bottom right top}
\begin{centering}
Avg J \\
\includegraphics[width=2.5in, trim={0 5in 0 5in}, clip]{color_bar.pdf}  \\
%0.40 \ \ \ \ \ \ \ \ \ \ \ \ \ \ \ \ 1.10
0.85 $\qquad\qquad\qquad\qquad$ 1.05

\end{centering}
\caption{Deformations of the anisotropic Cook's membrane benchmark (Section \ref{Anisotropic Cook's Membrane}), along with mean values of $J$ within each element calculated via equation \eqref{avgJ}, using the modified standard reinforcing model, equations (\ref{sr_energy}) -- (\ref{sr_stress_mod}). The background Eulerian grid is not depicted. Shown here are solid meshes with \textbf{Q1} elements and $m = 2601$ solid DOF. The first row is shows cases with $\nus = .4$, and the second row shows cases with $\nus = -1$ (here equivalent to $\kappas = 0$ and no volumetric-based stabilization). The first column depicts cases with modified invariants, and the second column depicts cases with unmodified invariants. Notice that the case with unmodified invariants and zero numerical bulk modulus leads to a collapsed element on the face where the traction is applied.}

\label{aniso}
\end{figure}

\begin{figure}
\captionsetup[subfigure]{justification=centering}
\begin{tabular}{c c c c}
\begin{subfigure}{.24\linewidth}
  \centering
  \includegraphics[width=.95\linewidth, trim={350 350 350 200}, clip]{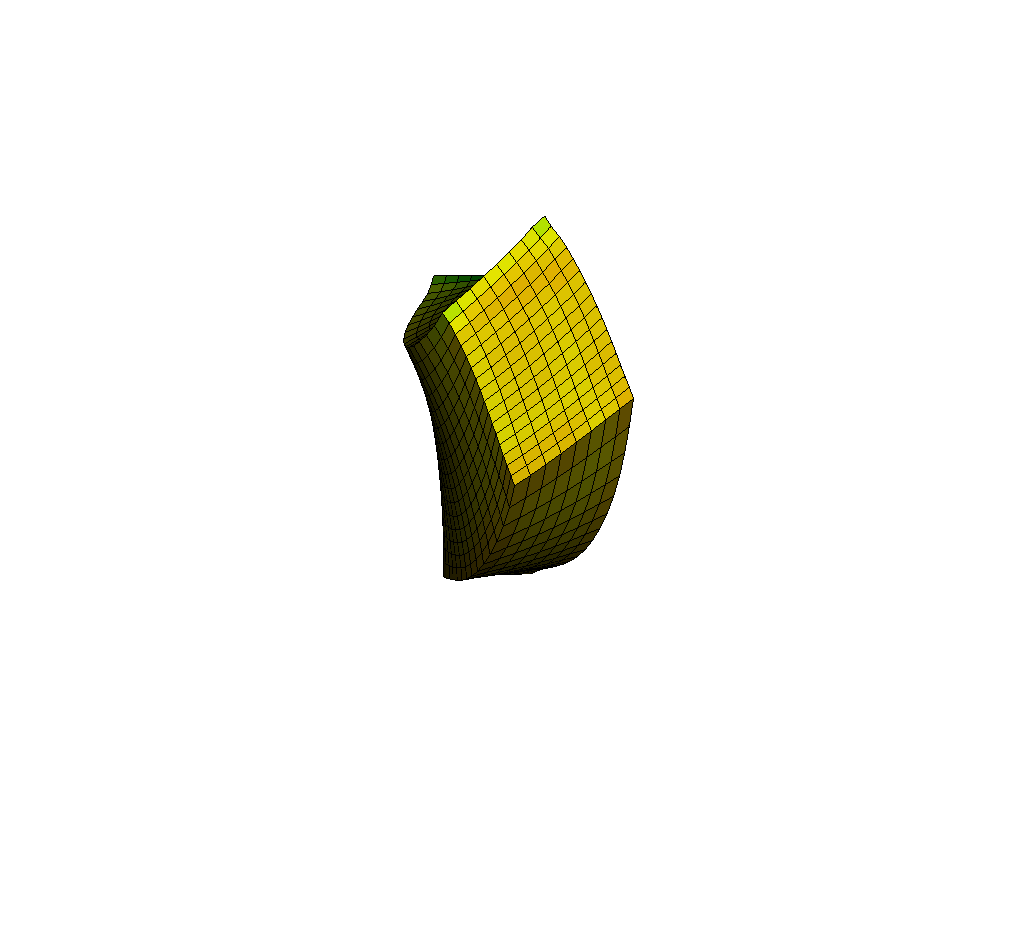}
  \caption{Modified Invariants \\ $\nus = .4$}
\end{subfigure}
\begin{subfigure}{.24\linewidth}
  \centering
\includegraphics[width=.95\linewidth, trim={350 350 350 200}, clip]{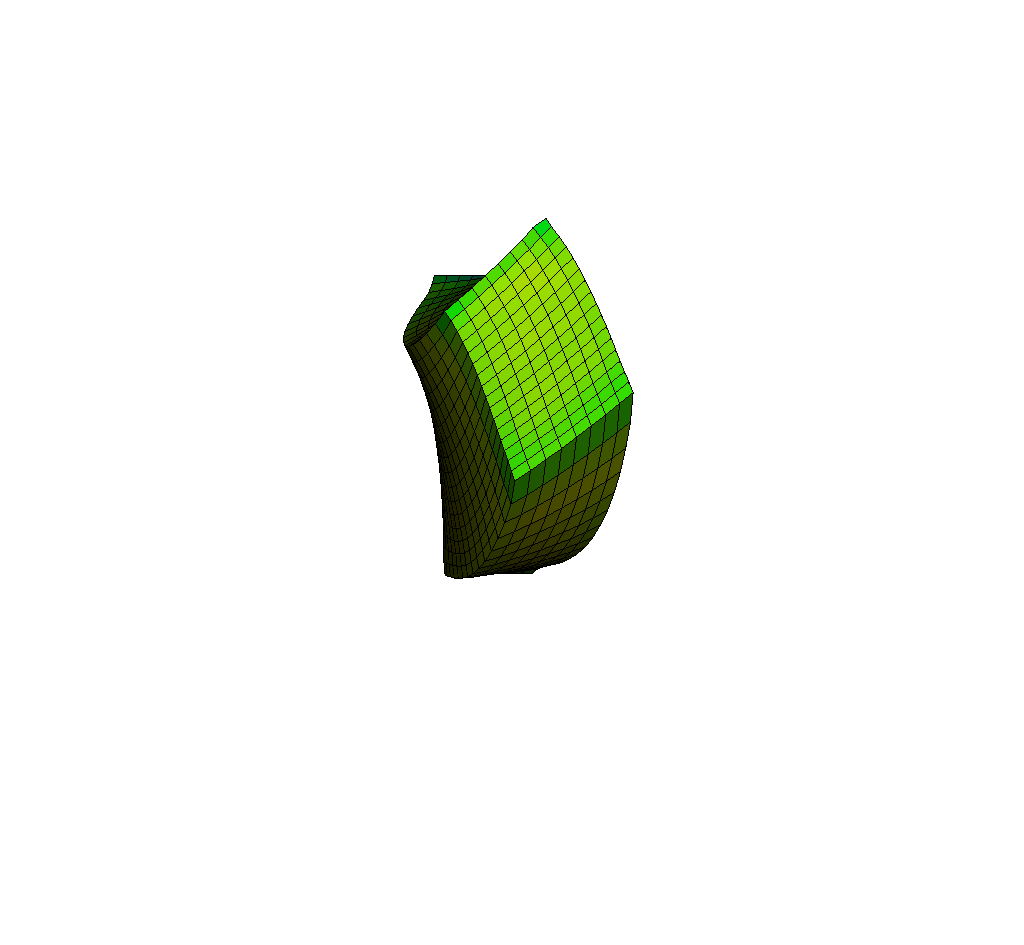}
  \caption{Unmodified Invariants \\ $\nus = .4$}
\end{subfigure}
\begin{subfigure}{.24\textwidth}
  \centering
 \includegraphics[width=.95\linewidth, trim={350 350 350 200}, clip]{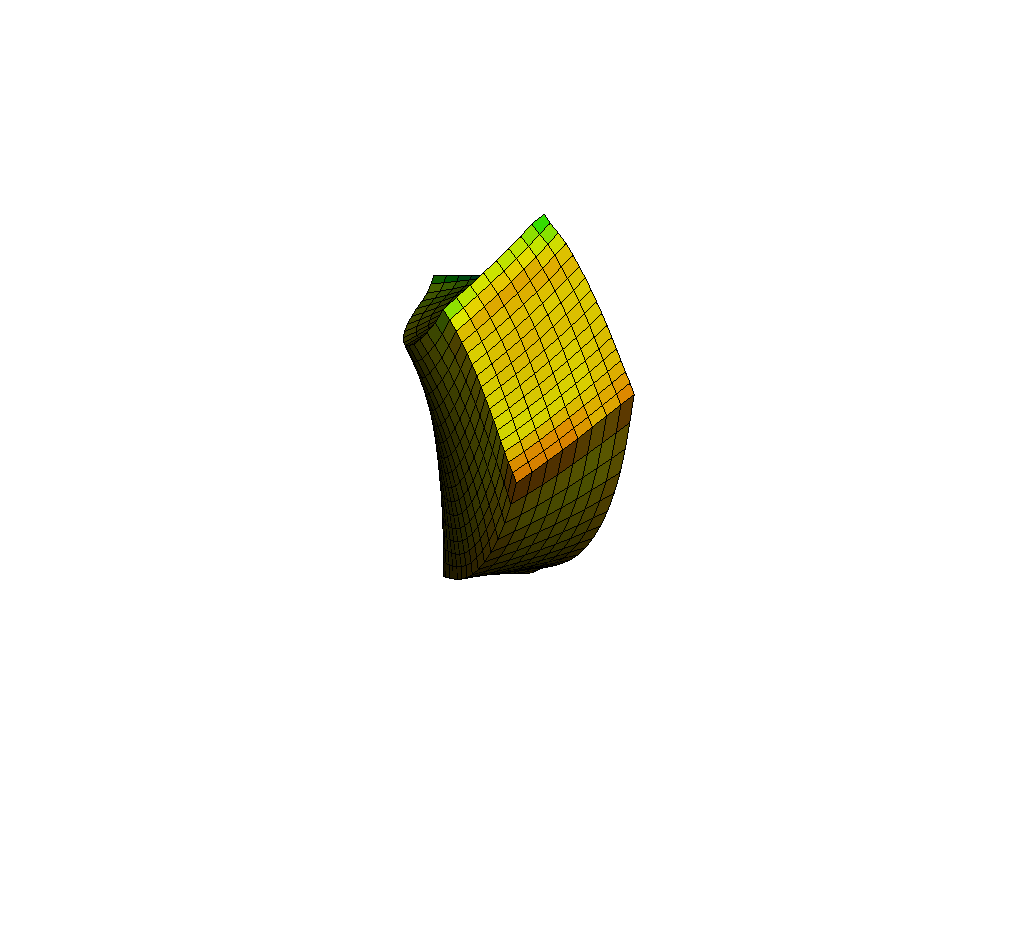}
  \caption{Modified Invariants \\ $\nus = -1$}
\end{subfigure}
\begin{subfigure}{.24\textwidth}
  \centering
\includegraphics[width=.95\linewidth, trim={350 350 350 200}, clip]{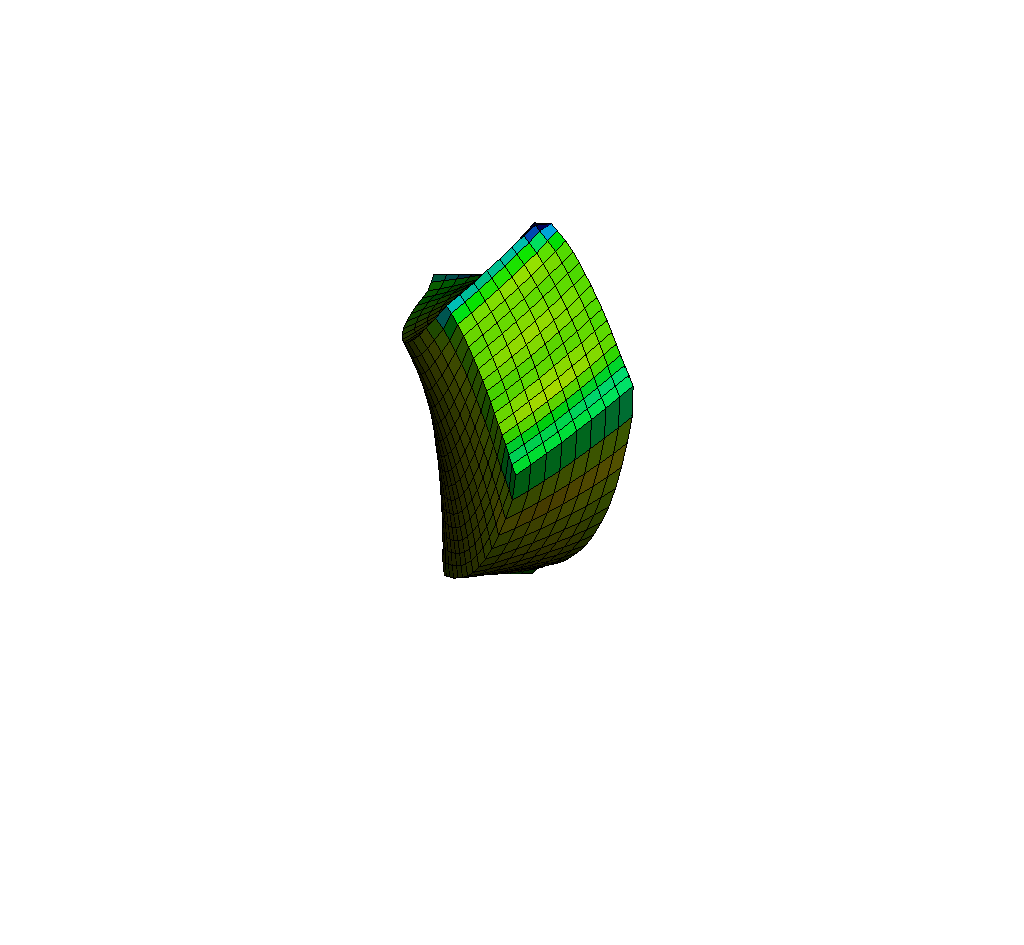}
  \caption{Unmodified Invariants \\ $\nus = -1$}
\end{subfigure}
\end{tabular}

\begin{centering}
Avg $J$ \\
\includegraphics[width=2.5in, trim={0 5in 0 5in}, clip]{color_bar.pdf}  \\
%0.40 \ \ \ \ \ \ \ \ \ \ \ \ \ \ \ \ 1.10
0.85 $\qquad\qquad\qquad\qquad$ 1.05

\end{centering}
\caption{Deformations and mean values of $J$ of the anisotropic Cook's membrane benchmark (Section \ref{Anisotropic Cook's Membrane}) with zero and finite volumetric energy from a different view; see also Figure (\ref{aniso}). The collapsed element in panel (d) is clearly visible.}
\label{aniso_alt}
\end{figure}

\begin{figure}
\begin{tabular}{l c c c}
&$\boldsymbol{\lambda_1}$& $\boldsymbol{\lambda_2} $& $\boldsymbol{\lambda_3}$ \\
\rotatebox{90}{$\qquad\quad\quad\;$ \textbf{m = 1557} }&
\subcaptionbox{\label{sfig:testa}} {\includegraphics[width=.15\linewidth]{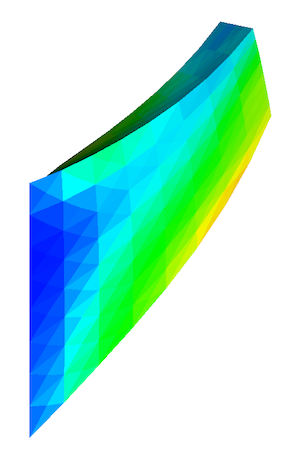}}&
\subcaptionbox{\label{sfig:testb}} {\includegraphics[width=.15\linewidth]{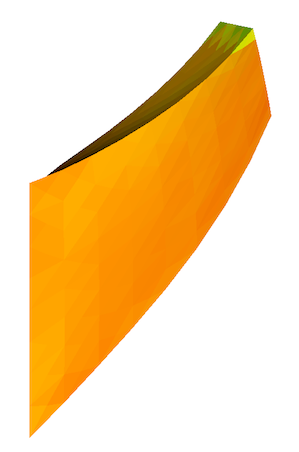}}&
\subcaptionbox{\label{sfig:testc}} {\includegraphics[width=.15\linewidth]{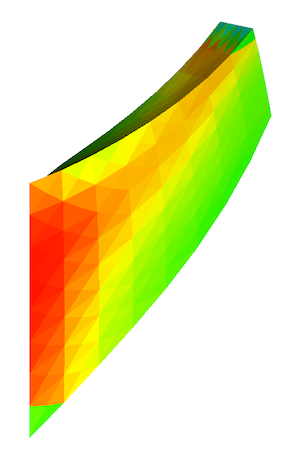}}\\ 

\rotatebox{90}{$\qquad\quad\quad\;$ \textbf{m = 11,300} }&
\subcaptionbox{\label{sfig:testa}} {\includegraphics[width=.15\linewidth]{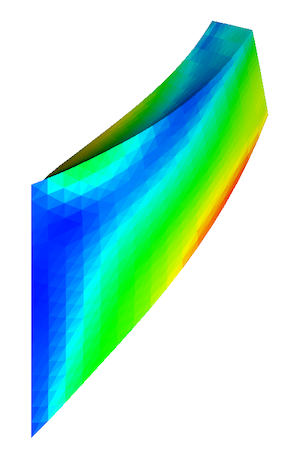}}&
\subcaptionbox{\label{sfig:testb}} {\includegraphics[width=.15\linewidth]{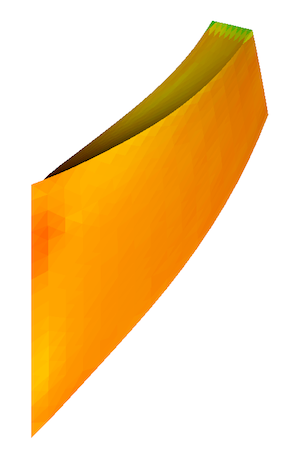}}&
\subcaptionbox{\label{sfig:testc}} {\includegraphics[width=.15\linewidth]{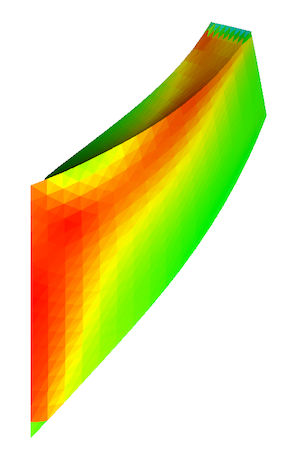}}\\ 

\rotatebox{90}{$\qquad\quad\quad\;$ \textbf{m = 36,920} }&
\subcaptionbox{\label{sfig:testa}} {\includegraphics[width=.15\linewidth]{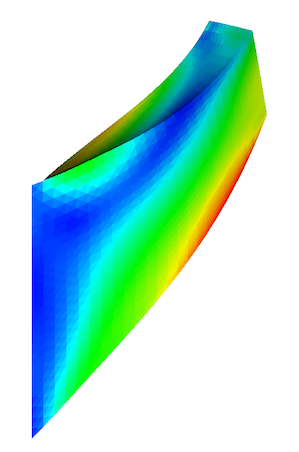}}&
\subcaptionbox{\label{sfig:testb}} {\includegraphics[width=.15\linewidth]{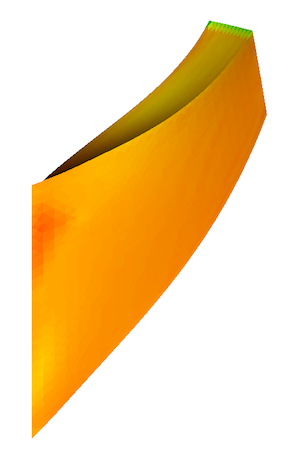}}&
\subcaptionbox{\label{sfig:testc}} {\includegraphics[width=.15\linewidth]{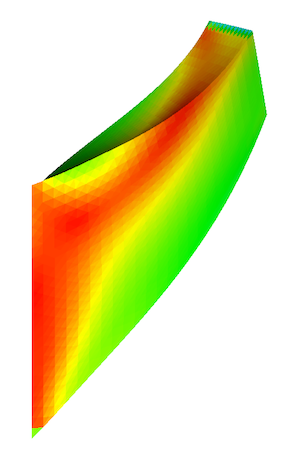}}\\ 

\\
\hline
\\
\rotatebox{90}{$\qquad\quad\quad$ \textbf{FE (P1/P1)} }&
\subcaptionbox{\label{sfig:testa}} {\includegraphics[width=.15\linewidth]{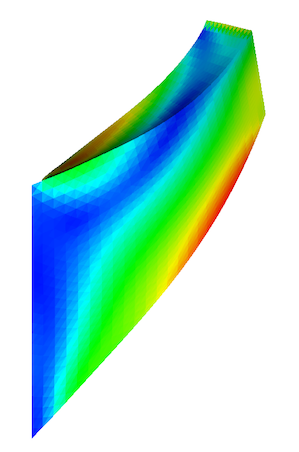}}&
\subcaptionbox{\label{sfig:testb}} {\includegraphics[width=.15\linewidth]{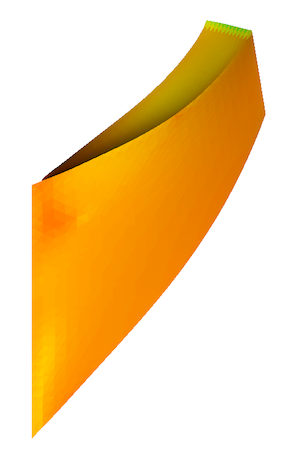}}&
\subcaptionbox{\label{sfig:testc}} {\includegraphics[width=.15\linewidth]{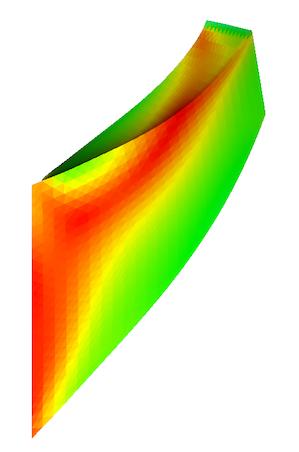}}\\ 

&\includegraphics[width=.3\linewidth, trim={0 5in 0 5in}, clip]{color_bar.pdf}&
\includegraphics[width=.3\linewidth, trim={0 5in 0 5in}, clip]{color_bar.pdf}&
\includegraphics[width=.3\linewidth, trim={0 5in 0 5in}, clip]{color_bar.pdf}  \\
&1.0 $\qquad\qquad\qquad$ 1.8& 0.75 $\qquad\qquad\qquad$ 1.05& 0.2 $\qquad\qquad\qquad$ 1.0 \\
\end{tabular}
\caption{Principal stretches (eigenvalues of $\FF$) of the anisotropic Cook's membrane benchmark (Section \ref{Anisotropic Cook's Membrane}) for the IBFE method using \textbf{P1} elements with modified invariants and volumetric stabilization ($\nus = 0.4$) and the principal stretches for the FE (\textbf{P1}/\textbf{P1}) method. The solid DOF for the IBFE method are listed in the leftmost column, and the FE (\textbf{P1}/\textbf{P1}) method uses $m = 36,920$ solid DOF.}
\label{ac_ps}
\end{figure}

\begin{figure}
$\qquad\qquad\qquad\qquad\qquad\qquad\;\;\;\;$ \textbf{P1} $\qquad\qquad\qquad\qquad\qquad\qquad\qquad\qquad\qquad\quad$  \textbf{Q1}\\
\rotatebox{90}{\qquad\qquad\qquad\quad\quad\; \textbf{$\nus = .4$}} 
   \rotatebox{90}{\qquad\qquad\qquad\qquad  Disp. (cm) }
\includegraphics[width=.45\linewidth, trim={30 190 40 200}, clip]{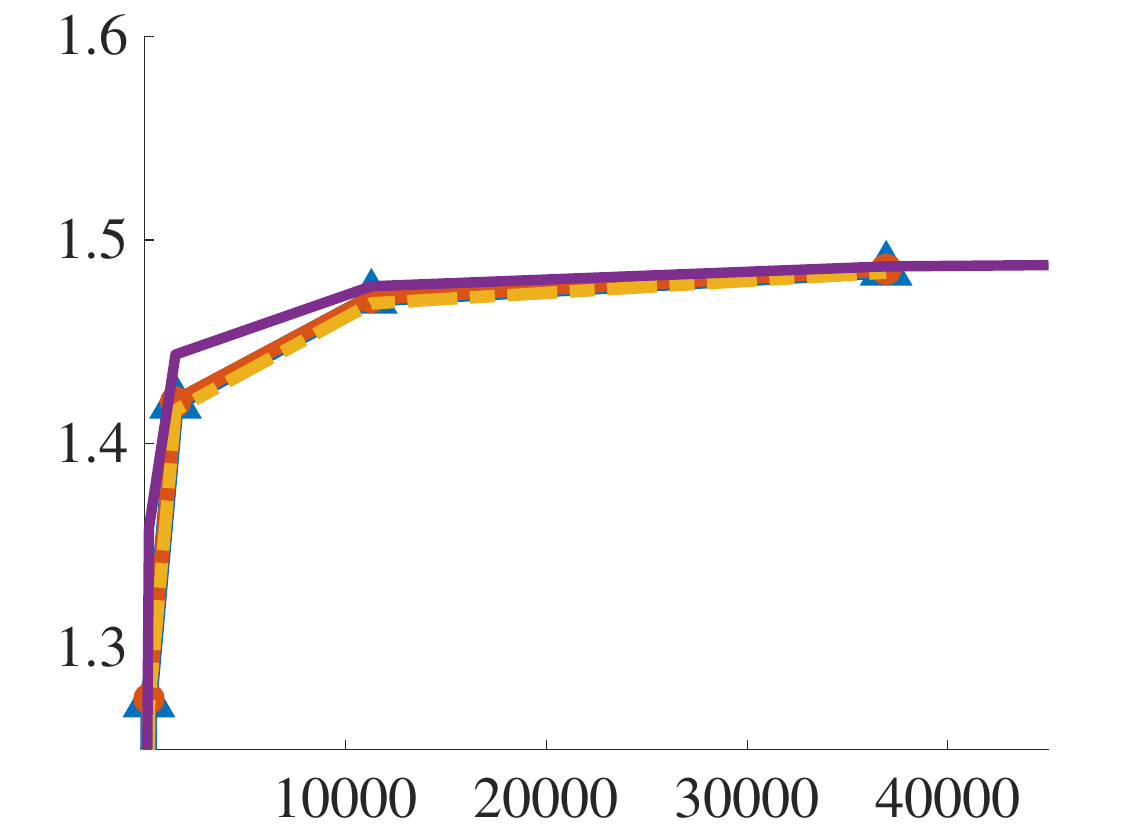} 
\includegraphics[width=.45\linewidth, trim={30 190 40 200}, clip]{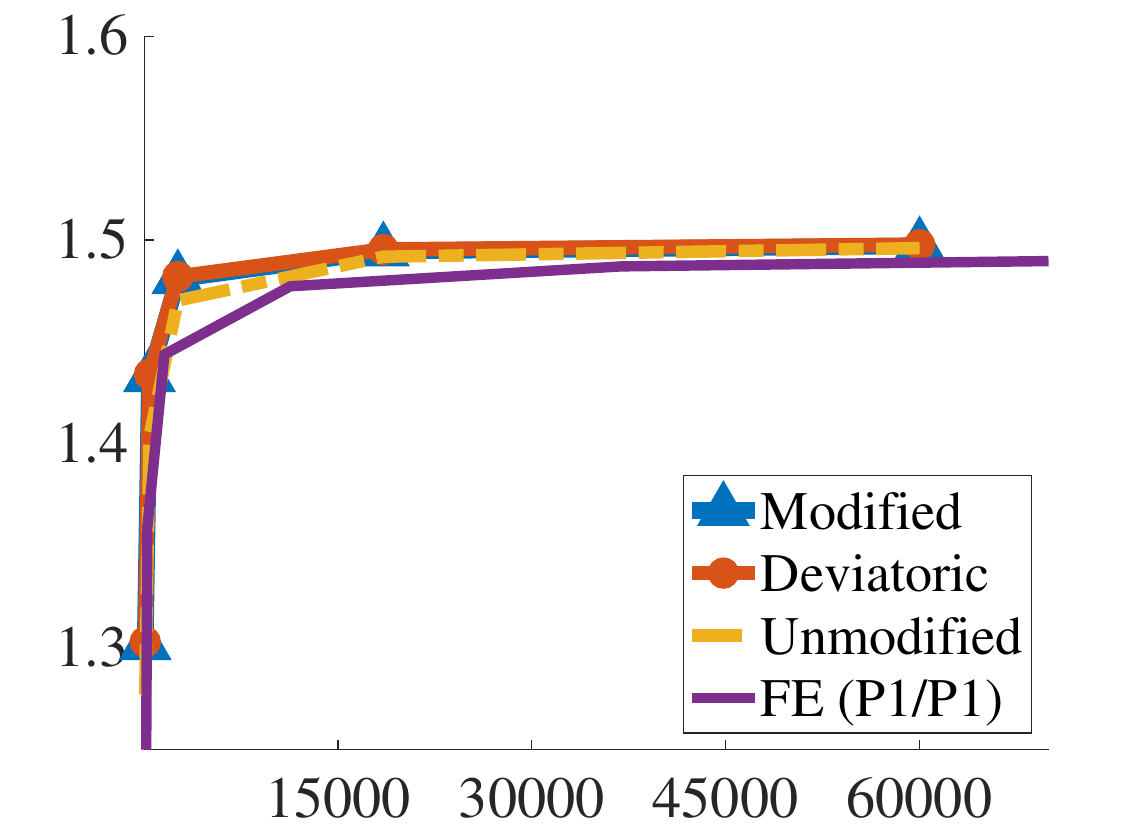}\\

\rotatebox{90}{\qquad\qquad\qquad\qquad\quad \textbf{$\nus = -1$}} 
   \rotatebox{90}{\qquad\qquad\qquad\qquad\;  Disp. (cm) }
\includegraphics[width=.45\linewidth, trim={30 190 40 200}, clip]{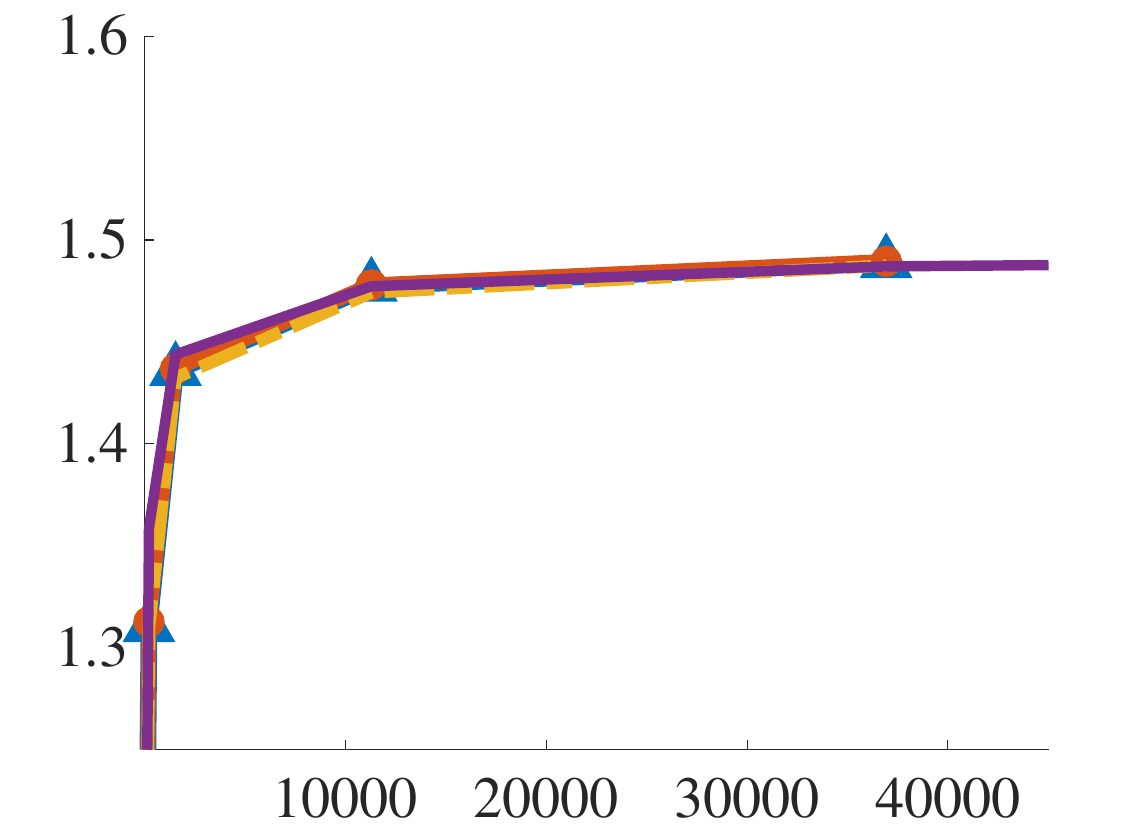} 
\includegraphics[width=.45\linewidth, trim={30 190 40 200}, clip]{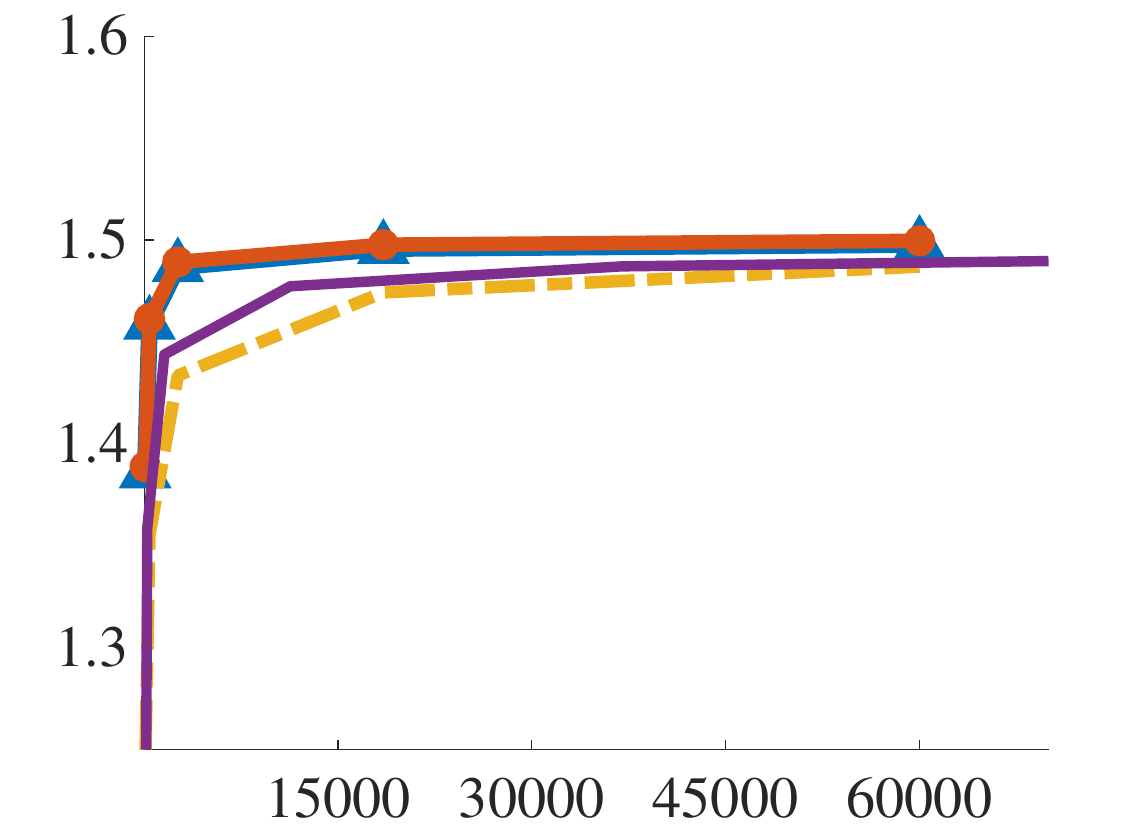}  \\

$\qquad\qquad\qquad\qquad\qquad\quad$ \# Solid DOF $\qquad\qquad\qquad\qquad\qquad\qquad\qquad\quad$  \# Solid DOF\\
\caption{Corner $y$-displacement for different numbers of DOF for the anisotropic Cook's membrane benchmark (Section \ref{Anisotropic Cook's Membrane}) for different choices of elements and numerical Poisson ratios. The solid DOF for the IBFE tests range from $m = 42$ to $60,025$. Notice that each row has the same $y$ extents, and each column has the same $x$ extents. The structural mechanics method is run with an additional discretization of $m = 86,097$ DOF. The displacement of the point of interest for the case with \textbf{P1} elements and $\nus = -1$ is in particularly good agreement with the FE (\textbf{P1}/\textbf{P1}) method here, but the deformations include numerically inaccurate artifacts, as seen in Figure (\ref{aniso_alt}).}
\label{aniso_disp}
\end{figure}

\begin{figure}
$\qquad\qquad\qquad\qquad\qquad\qquad\;\;\;\;$ \textbf{P1} $
\qquad\qquad\qquad\qquad\qquad\qquad\qquad\qquad\qquad\quad$  \textbf{Q1}\\
\rotatebox{90}{\qquad\qquad\qquad\quad \textbf{$\nus = .4$}} 
   \rotatebox{90}{\qquad\qquad\quad\;\;  Vol Change \% }
\includegraphics[width=.45\linewidth, trim={30 190 40 200}, clip]{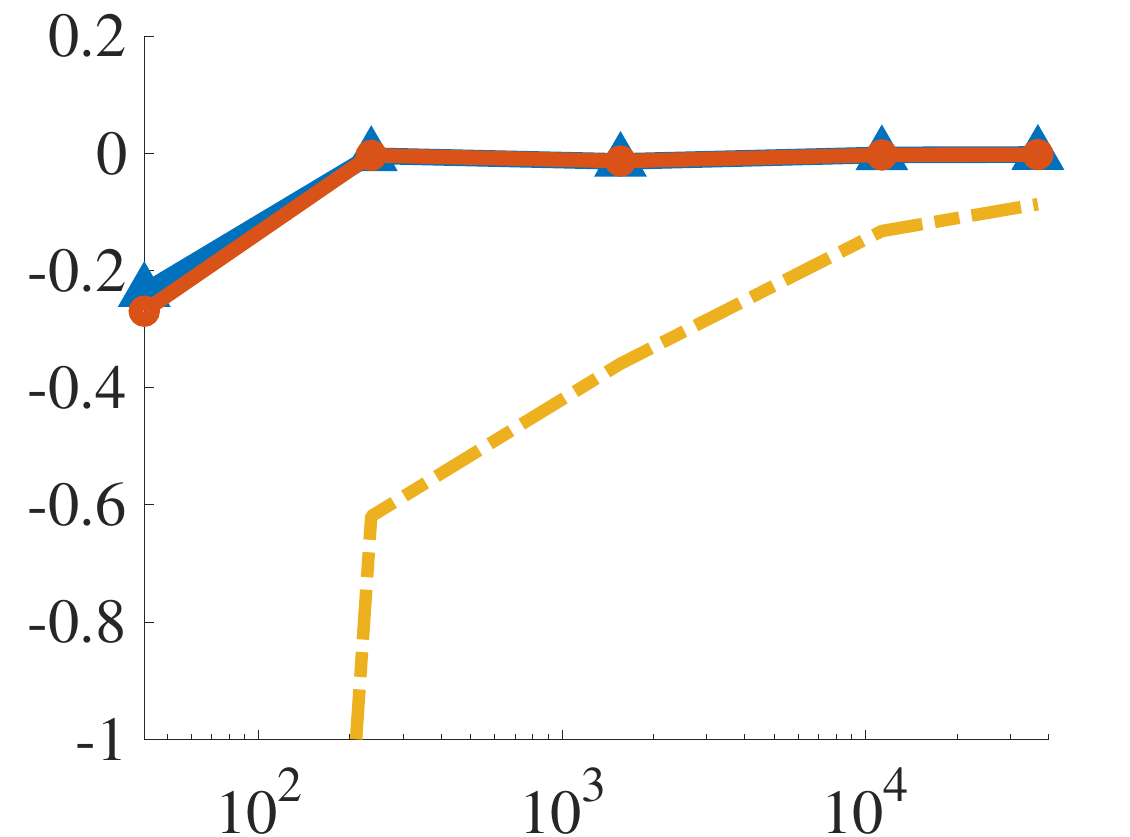}
\includegraphics[width=.45\linewidth, trim={30 190 40 200}, clip]{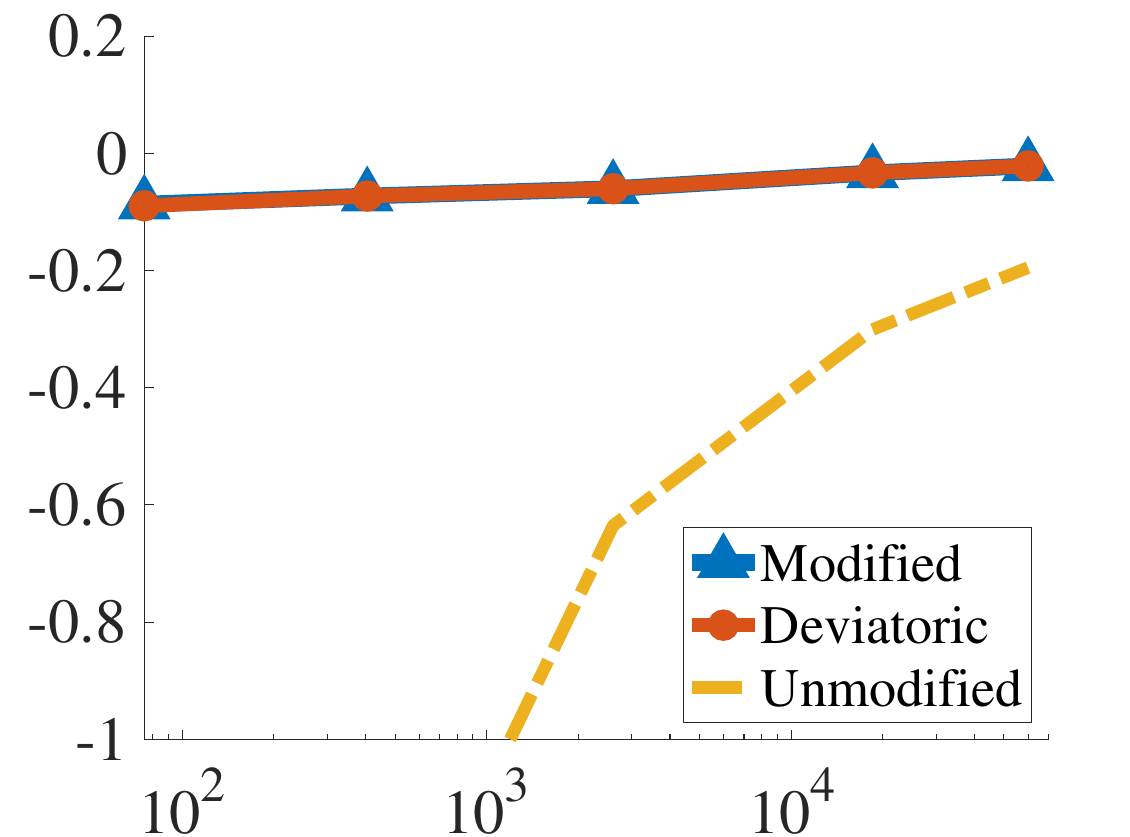}\\

\rotatebox{90}{\qquad\qquad\qquad\quad \textbf{$\nus = -1$}} 
   \rotatebox{90}{\qquad\qquad\quad\;\;  Vol Change \% }
\includegraphics[width=.45\linewidth, trim={30 190 40 200}, clip]{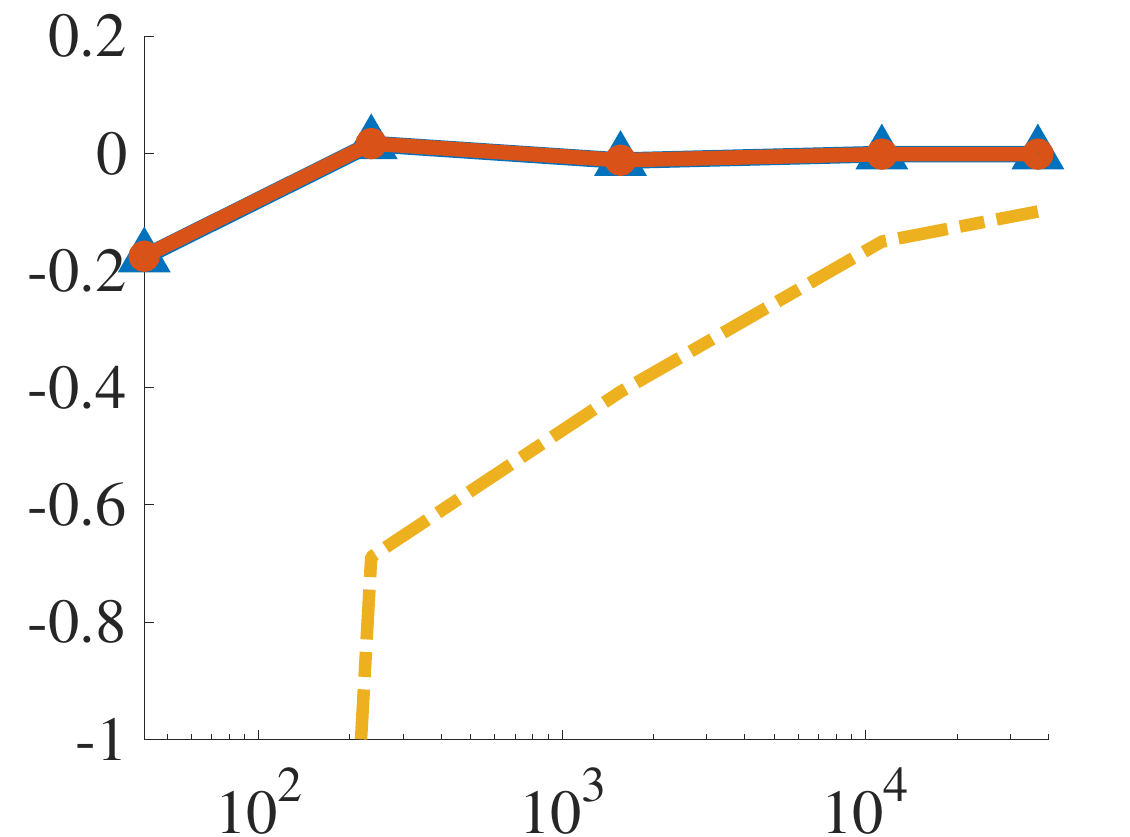}
\includegraphics[width=.45\linewidth, trim={30 190 40 200}, clip]{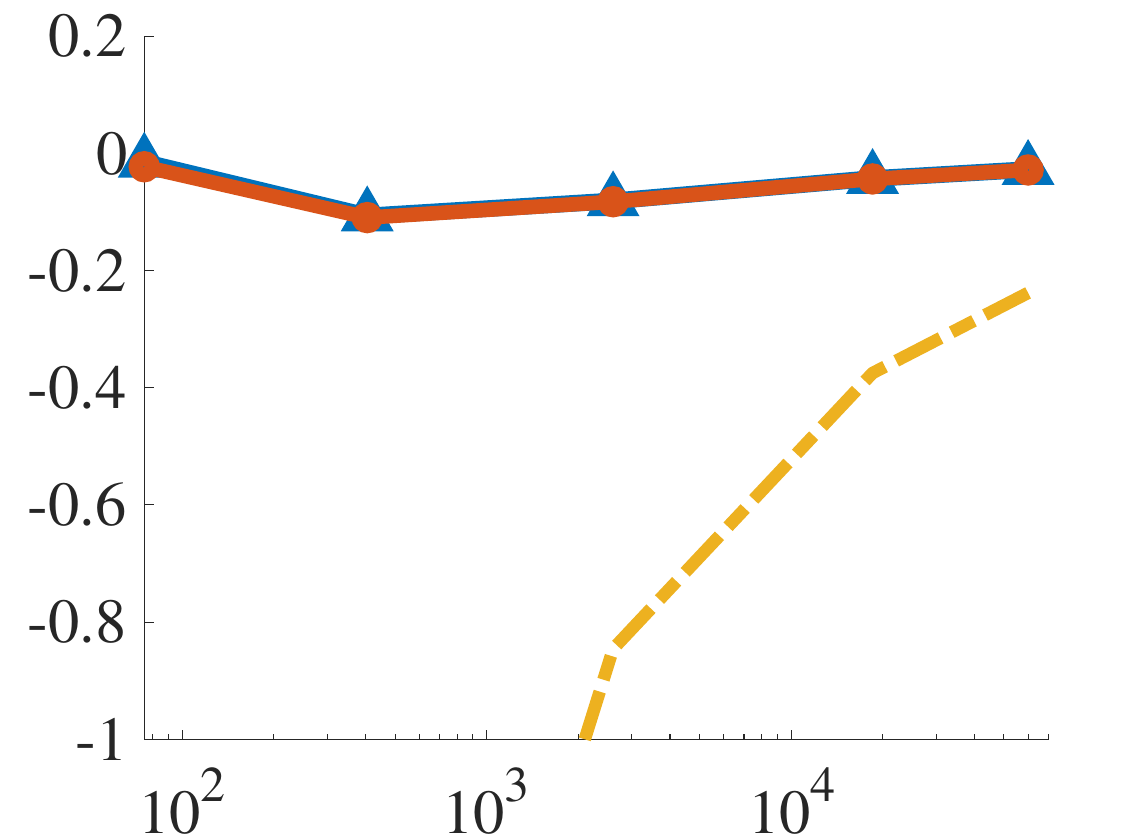}  \\

$\qquad\qquad\qquad\qquad\qquad\quad$ \# Solid DOF $\qquad\qquad\qquad\qquad\qquad\qquad\qquad\quad$  \# Solid DOF\\
\caption{Volume conservation for the anisotropic Cook's membrane benchmark (Section \ref{Anisotropic Cook's Membrane}) for different choices of elements and numerical Poisson ratio. The DOF range from $m = 42$ to $m = 60,025$, and the $x$ axis is on a log scale. Omitting the coarsest discretizations ($m=42$), the largest deviations in total volume among all element types used are approximately $.11\%$ for the modified case, $2.2\%$ for the unmodified case, and $.11\%$ for the deviatoric case.}
\label{aniso_vol}
\end{figure}

This benchmark involves a fully three-dimensional and anisotropic Cook's membrane; see Figure (\ref{aniso_mesh}). It is similar to and based upon one studied by Wriggers \etal~\cite{Wriggers2016}. The boundary conditions are the same as the two-dimensional model: an upward traction of $6.25$~$\frac{\text{dyn}}{\text{cm}^2}$ is applied to the right face, the body has zero prescribed displacement on the left face, and there is zero applied traction on all other faces. The displacement of the upper righthand corner of the right face is measured at $T_{\text{f}}=35$ s, and the load time is $T_{\text{l}} = 14$ s. This benchmark uses the standard reinforcing model, equations (\ref{sr_energy}) -- (\ref{sr_stress_dev}). Only two choices of numerical Poisson ratio are considered, $\nus = .4$ and $\nus = -1$, because of the extra computational effort required for three-dimensional simulations. Further, values of $\nus = .49995$ will exhibit locking. The fiber direction is $\Ab = \frac{1}{\sqrt{3}}(1, 1, 1)$, and we use material parameters $\Gt = 8$~$\frac{\text{dyn}}{\text{cm}^2}$, $\Gl = 160$~$\frac{\text{dyn}}{\text{cm}^2}$, and $\El = 1200$~$\frac{\text{dyn}}{\text{cm}^2}$. The density is $\rho = 1.0 \frac{\text{g}}{\text{cm}^3}$, and the fluid viscosity is $\mu = .16 \ \frac{\text{dyn} \cdot \text{s}}{\text{cm}^2}$. The larger viscosity is chosen to allow the model to more quickly reach steady state, and no damping is used in this test. The computational domain is $\Omega = [0, L]^3$ with $L = 12 \ \text{cm}$. The numbers of solid DOF range from $m = 42$ to $m = 60,025$. In our three-dimensional computations, we opt for structured tetrahedral meshes. Effectively, this means that FE nodes will have different locations for \textbf{Q1} and \textbf{P1} elements, and the sequence of meshes for each element type will have different numbers of solid DOF. IBFE computations using \textbf{P1} elements use the same meshes as those for the FE computations, whereas this is not possible for \textbf{Q1} elements.\\
%P1: 42         235        1557       11305       36925       86097, Q1: 75         405        2601       18513       60025
\indent As in the other cases considered, the behavior for the case of zero volumetric penalization with unmodified invariants yields unphysical deformations. In this case, the poor behavior is located at one of the corners on the face where the traction is applied; see Figures (\ref{aniso}) and (\ref{aniso_alt}). Specifically, the element at this location collapses; two of the FE nodes are approximately in the same location. We also show the principal stretches (eigenvalues of $\FF$) for this test when using modified invariants and volumetric stabilization in Figure (\ref{ac_ps}), which appear to be converging to results from the FE computations. Figure (\ref{aniso_disp}) shows plots of the $y$-displacement, which is measured at the encircled point in Figure (\ref{aniso_mesh}). Finally, as in the other cases considered, the case of unmodified invariants is associated with poor volume conservation. Figure (\ref{aniso_vol}) depicts the percent change in total volume for these cases. The percent change for all element types considered ranges between $.0087\%$ and $2.2\%$ for unmodified invariants. For the modified invariants and the deviatoric projection, both ranges were $.00014\%$ and $.11\%$.\\
\indent Unlike the other cases considered, the computation with zero volumetric penalization seems to perform nearly as well as or better than the case with volumetric penalization; see Figure (\ref{aniso_disp}). For the modified invariants, however, Figures (\ref{aniso}a) and (\ref{aniso}c) show that using a nonzero numerical bulk modulus produces a more uniform distribution of $J$ that is closer to $J = 1$. Overall, the differences among all cases in the results presented for this test are fairly minimal, with the exception that omitting volumetric penalization and using unmodified invariants yields unphysical deformations at the corners.

%\rotatebox{90}{\ \ \ \ \ \ \ \ \ \ \ \ \ \ \ \ \ \ \ \ \ \ \ \ \ \ \ \ \ \ \ \ \ \ \ \ \ \ \ \ Tip \ Displacement }

%%%%%%%%%%%%%%%%%%%%%

\subsection{Torsion}
\label{Torsion}
\begin{figure}
\centering
\includegraphics[width=10cm, trim={20 100 70 150}, clip]{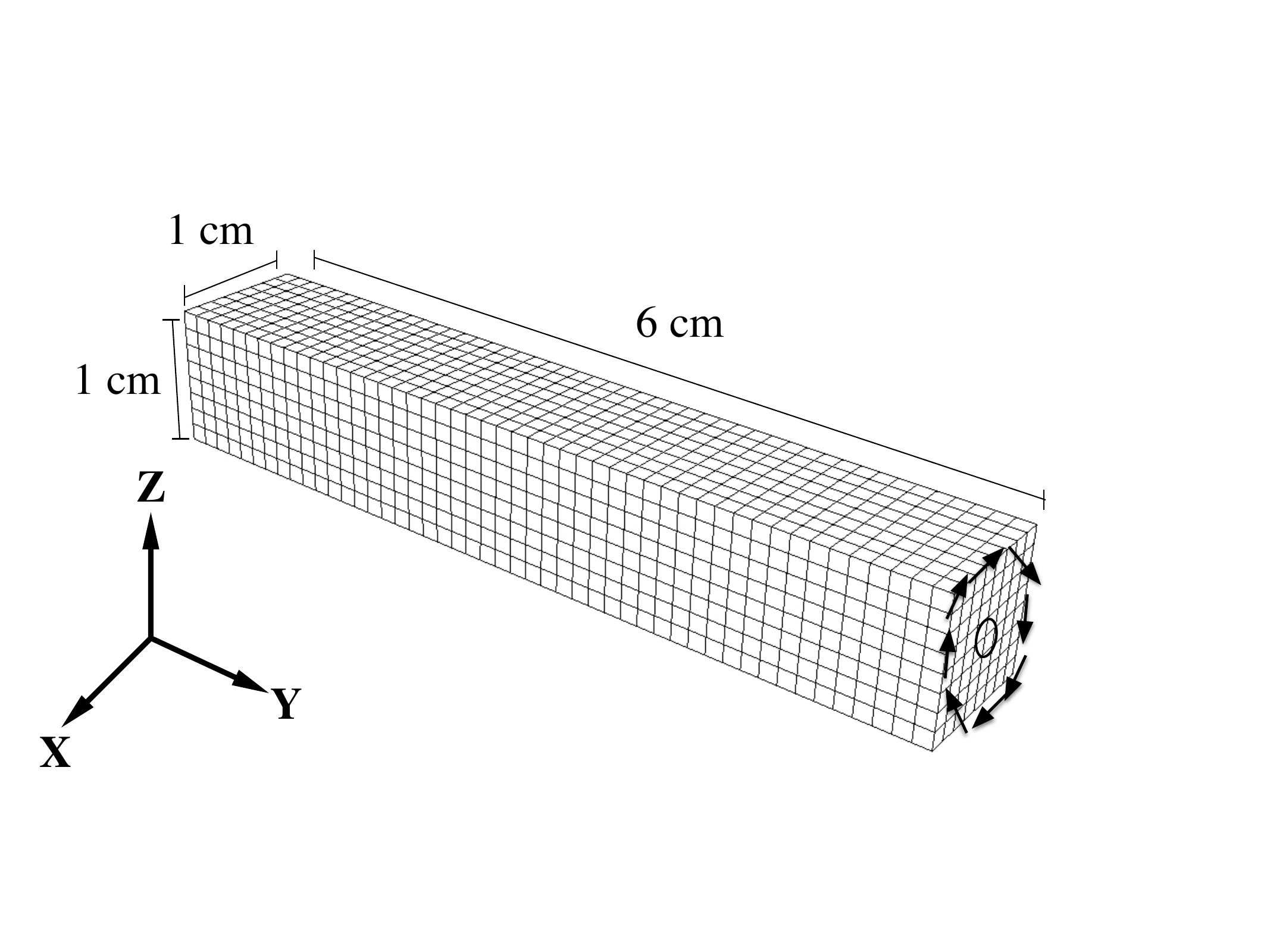}
\caption{Specifications of the torsion benchmark (Section \ref{Torsion}). The face opposite the applied torsion is kept fixed. The quantity of interest is the $y$-displacement as measured at the encircled area. To simplify the diagram of this three dimensional test, we omit the computational domain in this figure describing the problem setup. In the IBFE model, however, the structure is contained within a computational domain with dimensions $\Omega = [0,L]^3$ and $L = 9 \ \text{cm}$, and the solid mesh is placed in the center of this domain. Zero fluid velocity is enforced on the boundary of $\Omega$.}
\label{torsion_diag}
\end{figure}

\begin{figure}
\begin{tabular}{l c c}
& \textbf{Modified Invariants} & \textbf{Unmodified Invariants} \\
\rotatebox{90}{\qquad\qquad\qquad \textbf{$\nus = .4$} }&
\subcaptionbox{\label{sfig:testa}} {\includegraphics[width=.45\linewidth, trim={100 200 50 250}, clip]{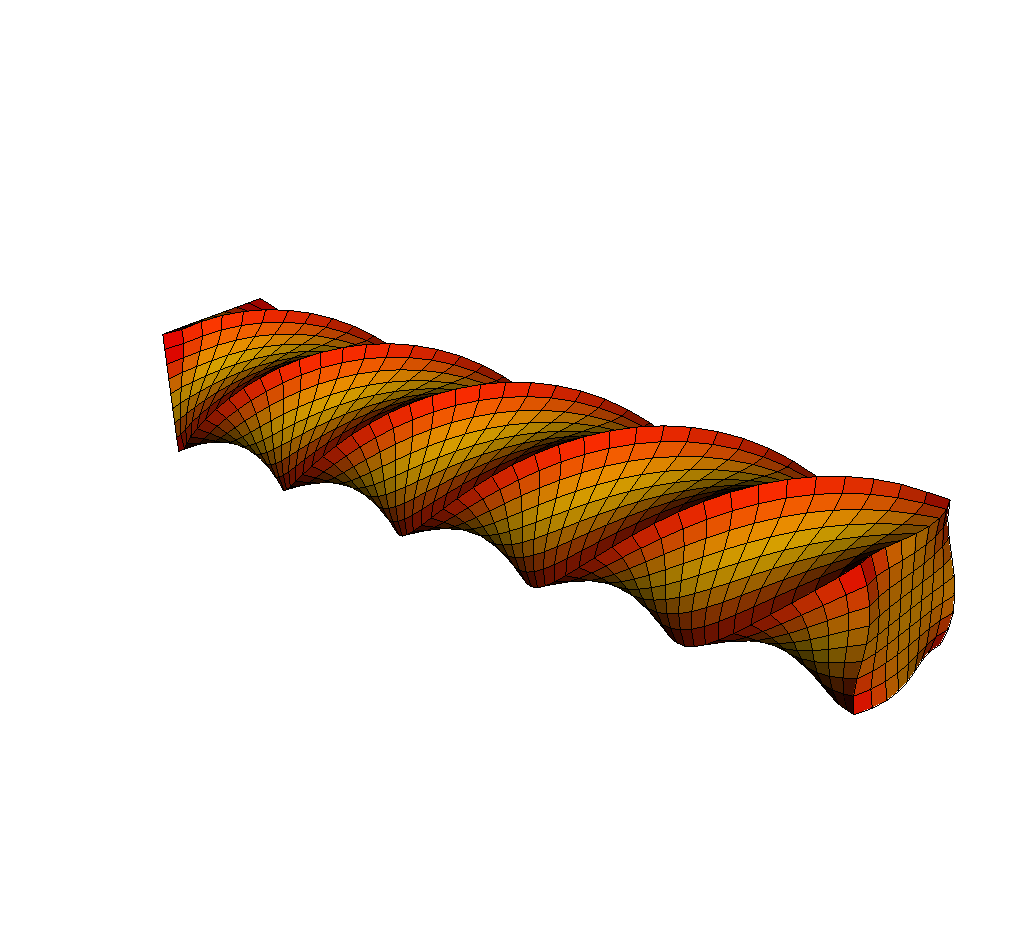}} &
\subcaptionbox{\label{sfig:testb}} {\includegraphics[width=.45\linewidth, trim={100 200 50 250}, clip]{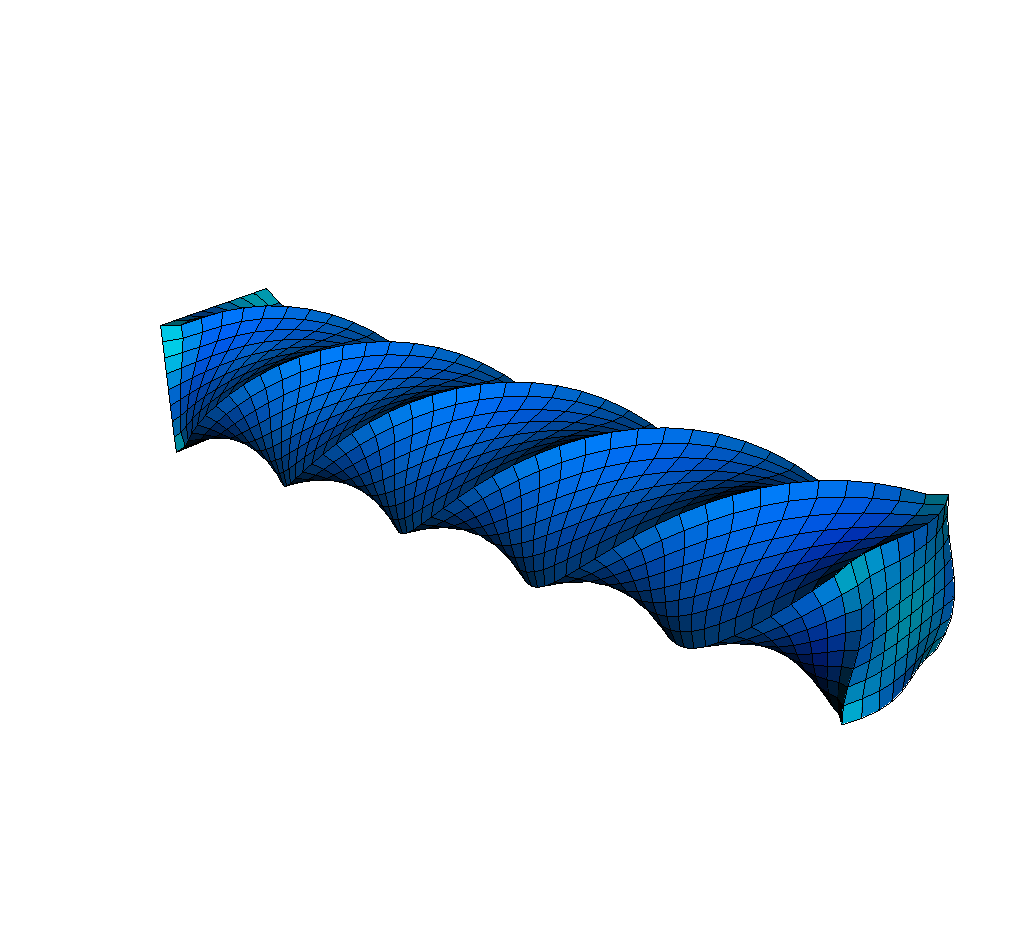}}\\

\rotatebox{90}{\qquad\qquad\qquad\textbf{$\nus = -1$}} &
\subcaptionbox{\label{sfig:testc}}{\includegraphics[width=.45\linewidth, trim={100 200 50 250}, clip]{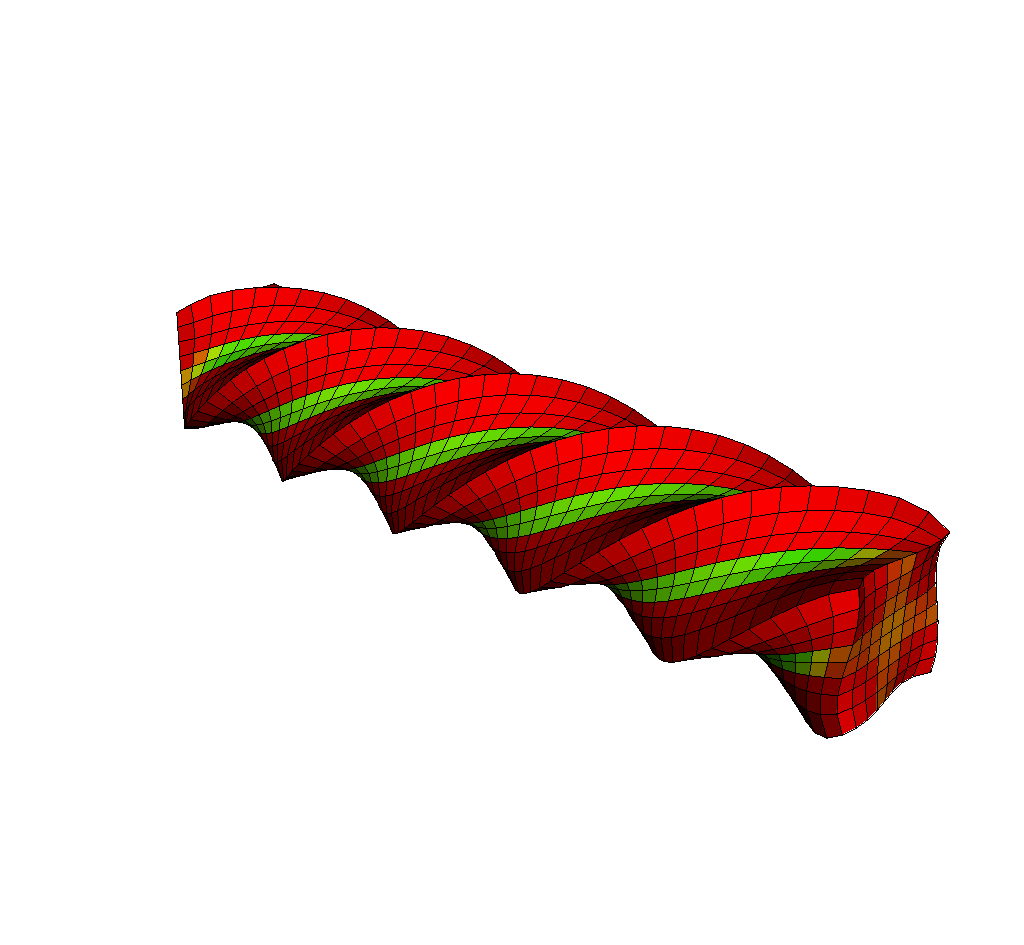}} &
\subcaptionbox{\label{sfig:testd}}{\includegraphics[width=.45\linewidth, trim={100 170 50 250}, clip]{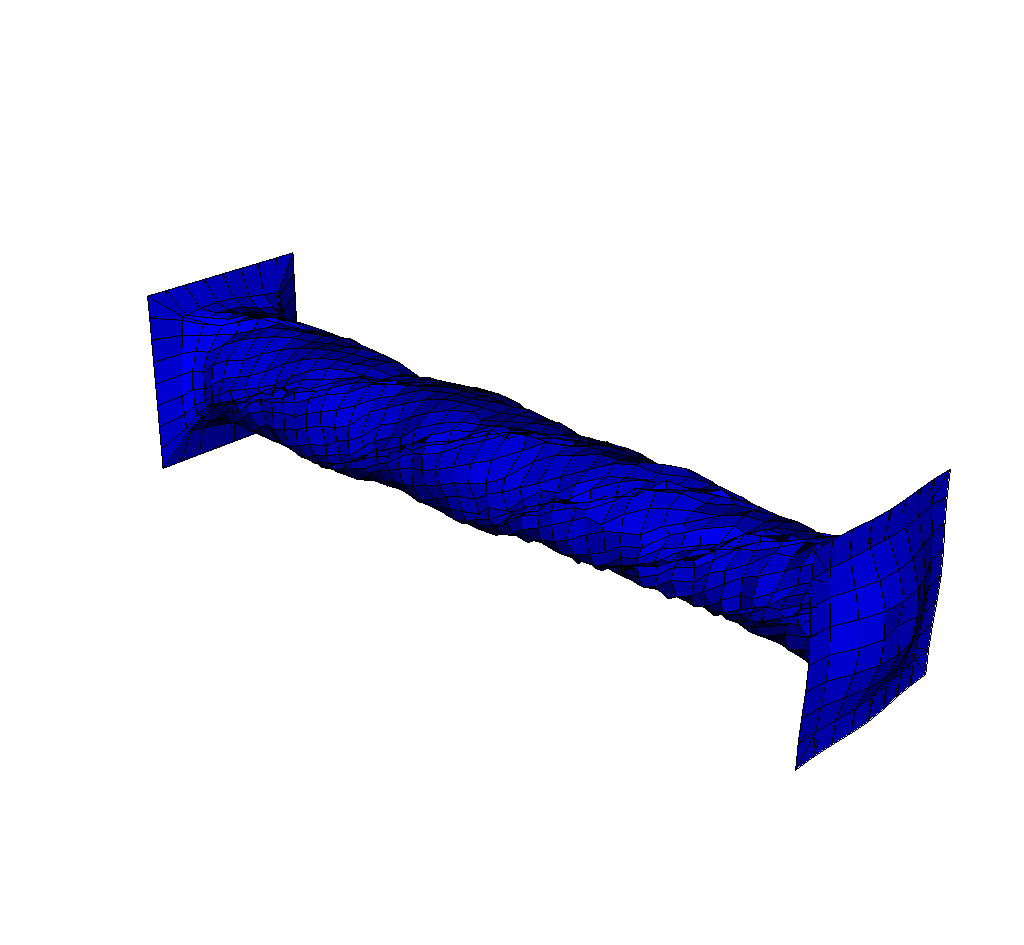}}  \\
\end{tabular}
%trim={left bottom right top}
\begin{centering}
Avg $J$ \\
\includegraphics[width=2.5in, trim={0 5in 0 5in}, clip]{color_bar.pdf}  \\
%0.00 \ \ \ \ \ \ \ \ \ \ \ \ \ \ \ \ 2.00
0.75 $\qquad\qquad\qquad\qquad$ 1.05

\end{centering}
\caption{Deformations of the torsion test (Section \ref{Torsion}), along with mean values of $J$ within each element calculated via equation \eqref{avgJ}, using the Mooney-Rivlin material model, equations (\ref{mr_energy}) -- (\ref{mr_stress_mod}), with $c_1, c_2 = 9000 \ \frac{\text{dyn}}{\text{cm}^2}$. The background Eulerian grid is not shown. Shown here are solid meshes with \textbf{Q1} elements and $m = 3969$ solid DOF. The first row shows cases with $\nus = .4$, and the second row shows cases with $\nus = -1$ (here equivalent to $\kappas = 0$ and no volumetric-based stabilization). The first column shows cases with modified invariants, and the second column depicts cases with unmodified invariants. Note the extremely unphysical deformations shown in panel (d) generated when using unmodified invariants without volumetric-based stabilization.}
\label{torsion}
\end{figure}

\begin{figure}
\begin{tabular}{l c c c}
&$\boldsymbol{\lambda_1}$& $\boldsymbol{\lambda_2} $& $\boldsymbol{\lambda_3}$ \\
\rotatebox{90}{$\qquad\quad\quad\;$ \textbf{m = 361} }&
\subcaptionbox{\label{sfig:testa}} {\includegraphics[width=.3\linewidth]{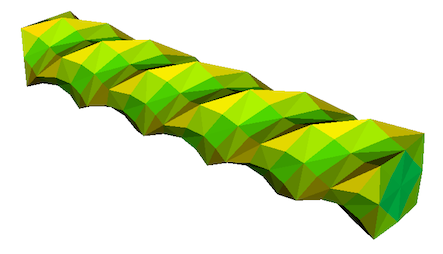}}&
\subcaptionbox{\label{sfig:testb}} {\includegraphics[width=.3\linewidth]{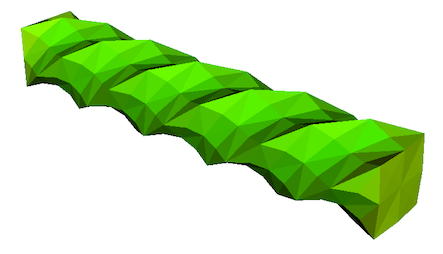}}&
\subcaptionbox{\label{sfig:testc}} {\includegraphics[width=.3\linewidth]{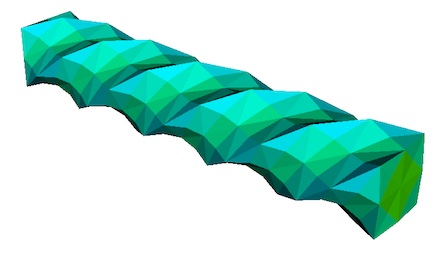}}\\ 

\rotatebox{90}{$\qquad\quad\quad\;$ \textbf{m = 2369} }&
\subcaptionbox{\label{sfig:testa}} {\includegraphics[width=.3\linewidth]{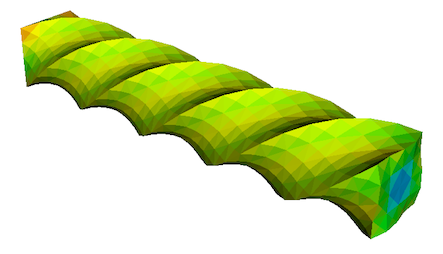}}&
\subcaptionbox{\label{sfig:testb}} {\includegraphics[width=.3\linewidth]{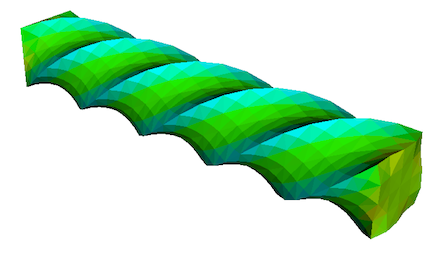}}&
\subcaptionbox{\label{sfig:testc}} {\includegraphics[width=.3\linewidth]{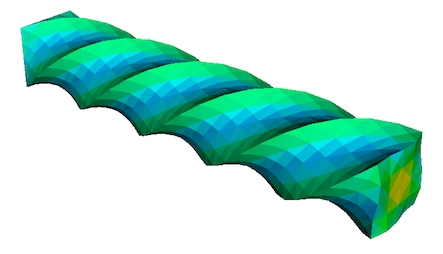}}\\ 

\rotatebox{90}{$\qquad\quad\quad\;$ \textbf{m = 7465} }&
\subcaptionbox{\label{sfig:testa}} {\includegraphics[width=.3\linewidth]{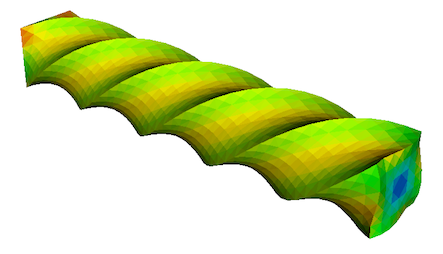}}&
\subcaptionbox{\label{sfig:testb}} {\includegraphics[width=.3\linewidth]{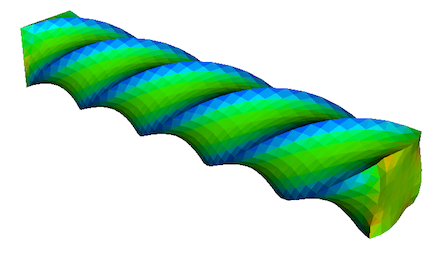}}&
\subcaptionbox{\label{sfig:testc}} {\includegraphics[width=.3\linewidth]{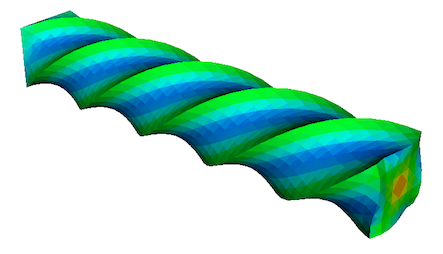}}\\ 

\\
\hline
\\
\rotatebox{90}{$\qquad\quad\quad$ \textbf{FE (P1/P1)} }&
\subcaptionbox{\label{sfig:testa}} {\includegraphics[width=.3\linewidth]{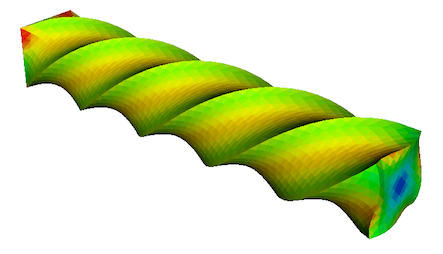}}&
\subcaptionbox{\label{sfig:testb}} {\includegraphics[width=.3\linewidth]{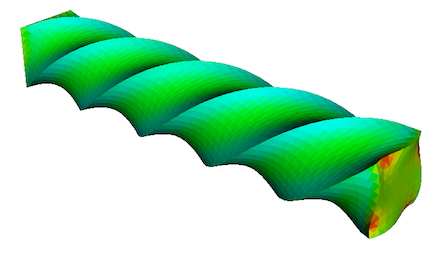}}&
\subcaptionbox{\label{sfig:testc}} {\includegraphics[width=.3\linewidth]{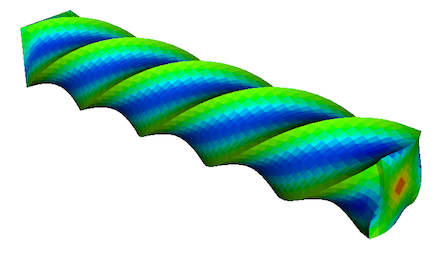}}\\ 

&\includegraphics[width=.3\linewidth, trim={0 5in 0 5in}, clip]{color_bar.pdf}&
\includegraphics[width=.3\linewidth, trim={0 5in 0 5in}, clip]{color_bar.pdf}&
\includegraphics[width=.3\linewidth, trim={0 5in 0 5in}, clip]{color_bar.pdf}  \\
&1.0 $\qquad\qquad\qquad$ 2.4& 0.8 $\qquad\qquad\qquad$ 1.1& 0.4 $\qquad\qquad\qquad$ 1.0 \\
\end{tabular}
\caption{Principal stretches (eigenvalues of $\FF$) of the torsion test (Section \ref{Torsion}) for the IBFE method using \textbf{P1} elements with modified invariants and volumetric stabilization ($\nus = 0.4$) and the principal stretches for the FE (\textbf{P1}/\textbf{P1}) method. The solid DOF for the IBFE method are listed in the leftmost column, and the FE (\textbf{P1}/\textbf{P1}) method uses $m = 17,089$ solid DOF. The IBFE results appear to be converging to those yielded by the high-resolution FE solution.}
\label{tt_ps}
\end{figure}

\begin{figure}
$\qquad\qquad\qquad\qquad\qquad\qquad\;\;\;\;$ \textbf{P1} $\qquad\qquad\qquad\qquad\qquad\qquad\qquad\qquad\qquad\quad$  \textbf{Q1}\\
\rotatebox{90}{\qquad\qquad\quad\;\; \textbf{$\nus = .4$}} 
   \rotatebox{90}{\qquad\qquad\quad\;  Disp. (cm) }
\includegraphics[width=.45\linewidth, trim={40 190 20 200}, clip]{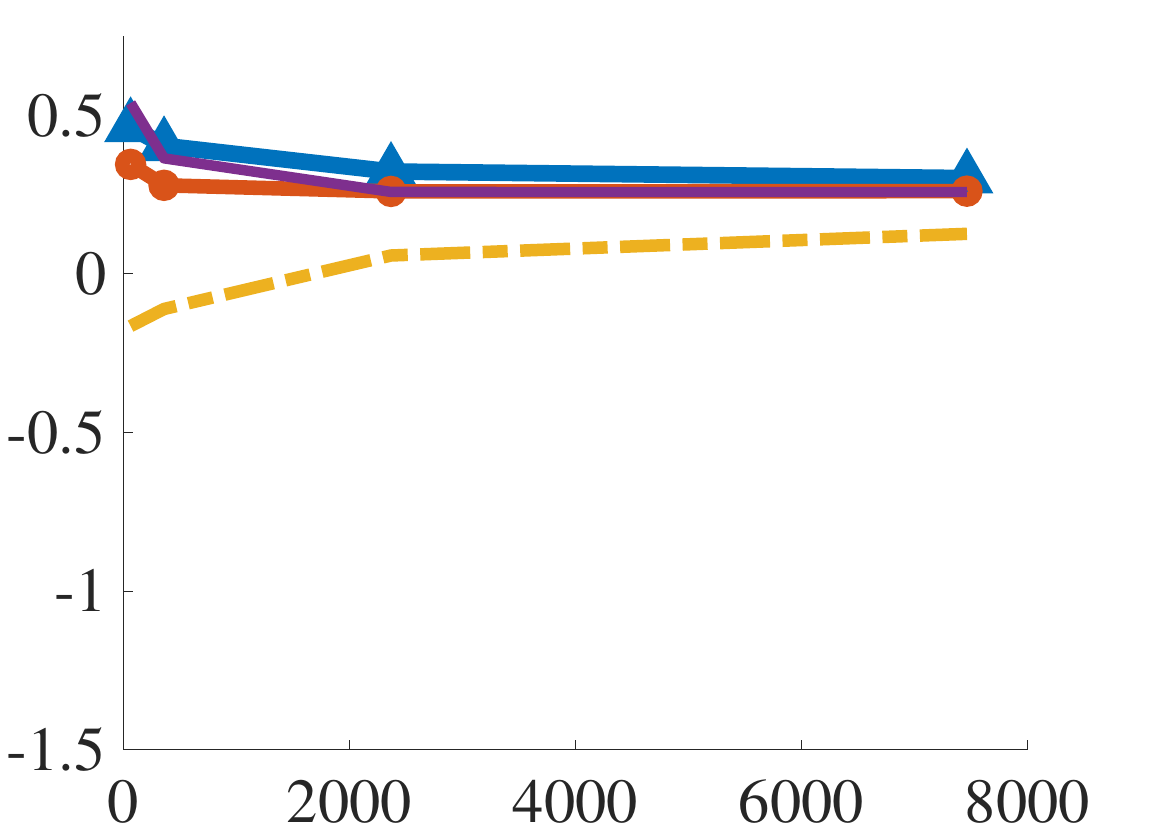}
\includegraphics[width=.45\linewidth, trim={40 190 20 200}, clip]{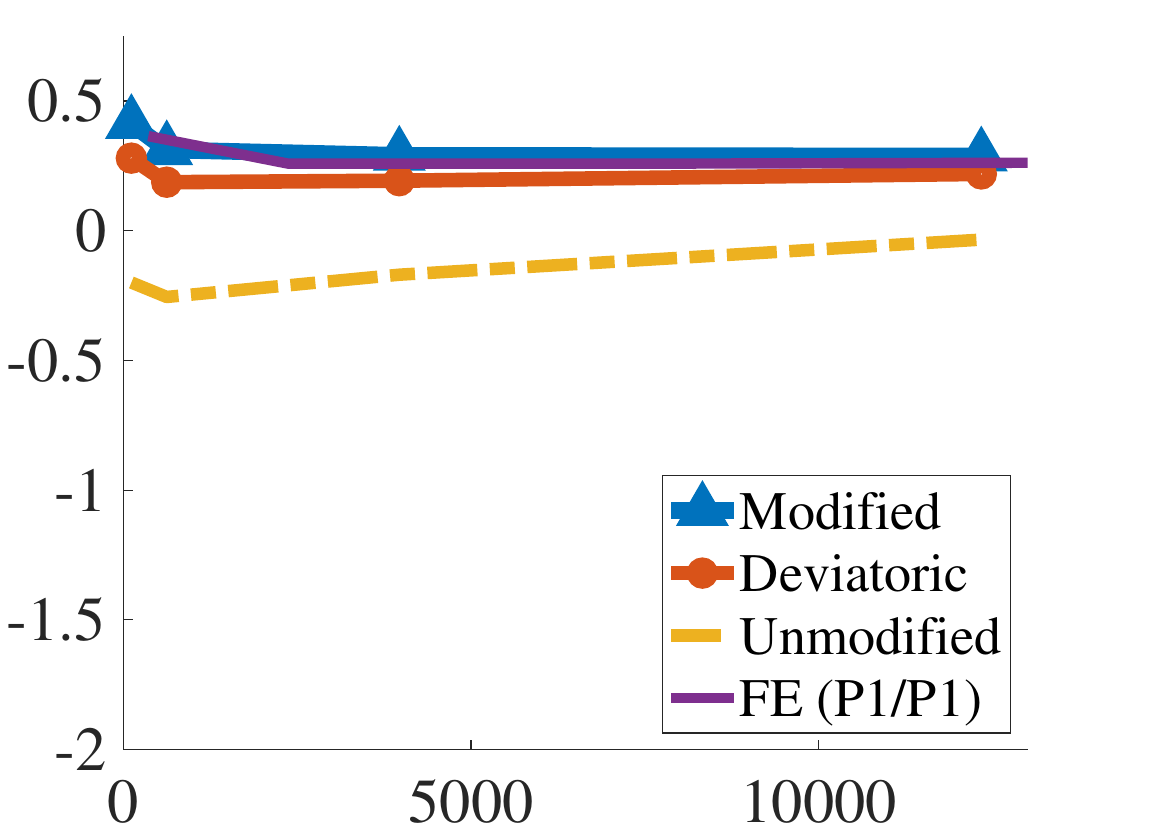}\\

\rotatebox{90}{\qquad\qquad\qquad\; \textbf{$\nus = -1$}} 
   \rotatebox{90}{\qquad\qquad\qquad\;  Disp. (cm) }
\includegraphics[width=.45\linewidth, trim={40 190 20 200}, clip]{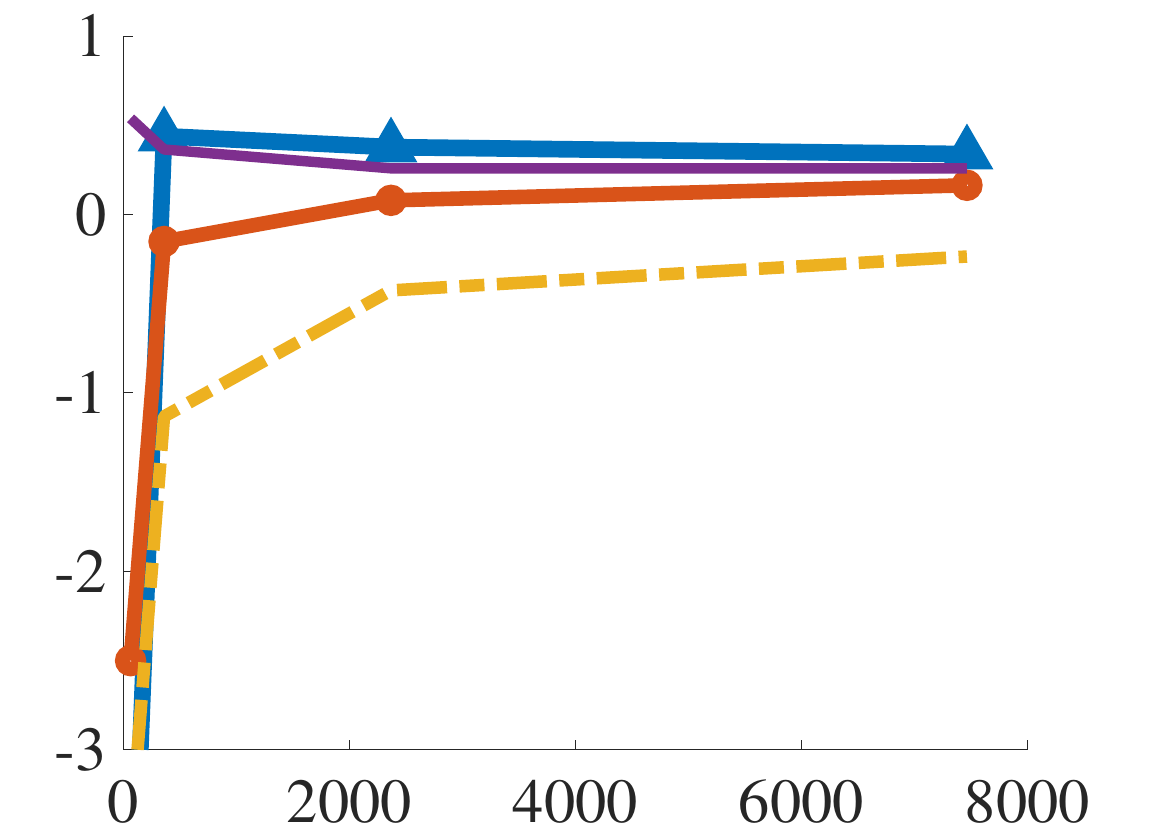}
\includegraphics[width=.45\linewidth, trim={40 190 20 200}, clip]{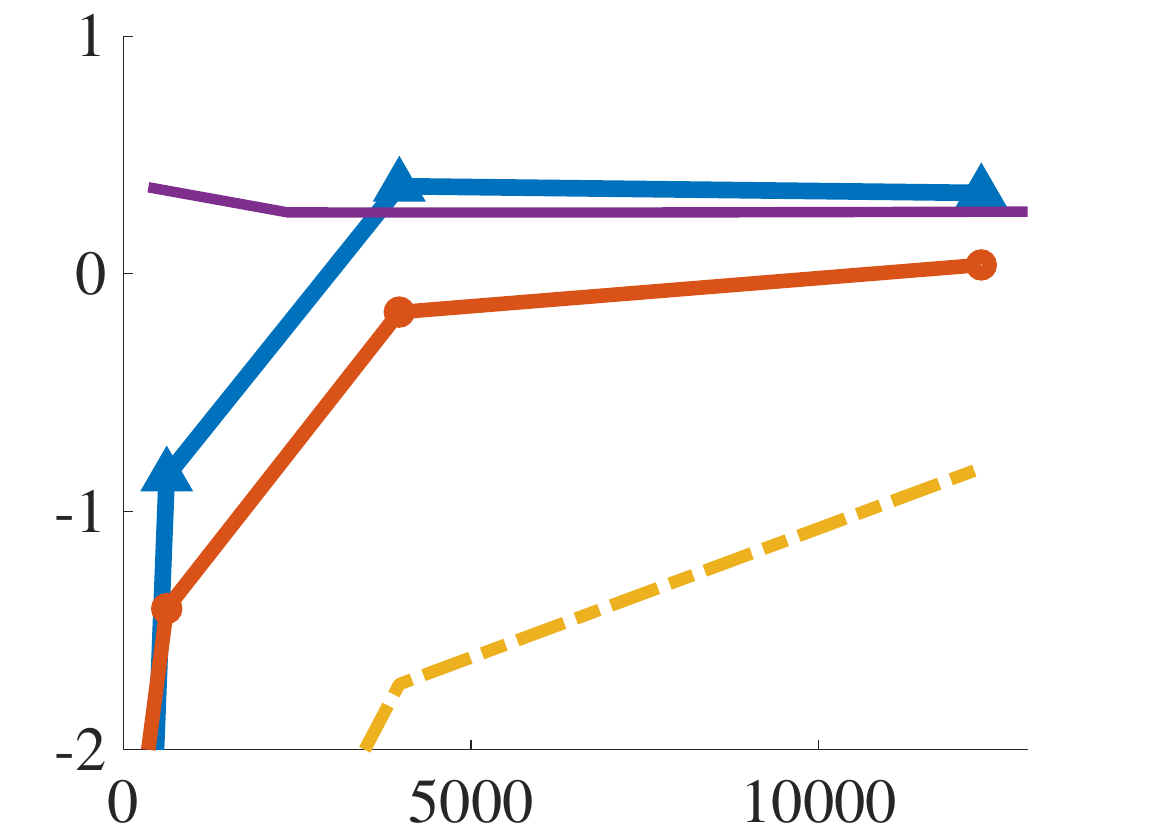}  \\

$\qquad\qquad\qquad\qquad\qquad\quad$ \# Solid DOF $\qquad\qquad\qquad\qquad\qquad\qquad\qquad\quad$  \# Solid DOF\\
\caption{Axial displacement for different solid DOF for the torsion test (Section \ref{Torsion}) for different choices of elements and numerical Poisson ratio. The solid DOF for the IB tests range from $m = 65$ to $m = 12,337$. Notice that each row has the same $y$ extents, and each column has the same $x$ extents. The structural mechanics method is run with an additional discretization of $m = 17,089$ DOF.}
\label{tor_disp}
\end{figure}

\begin{figure}
$\qquad\qquad\qquad\qquad\qquad\qquad\;\;\;\;$ \textbf{P1} $\qquad\qquad\qquad\qquad\qquad\qquad\qquad\qquad\qquad\quad$  \textbf{Q1}\\
\rotatebox{90}{\qquad\qquad\qquad\; \textbf{$\nus = .4$}} 
   \rotatebox{90}{\qquad\qquad\quad\;  Vol Change \% }
\includegraphics[width=.45\linewidth, trim={50 190 30 200}, clip]{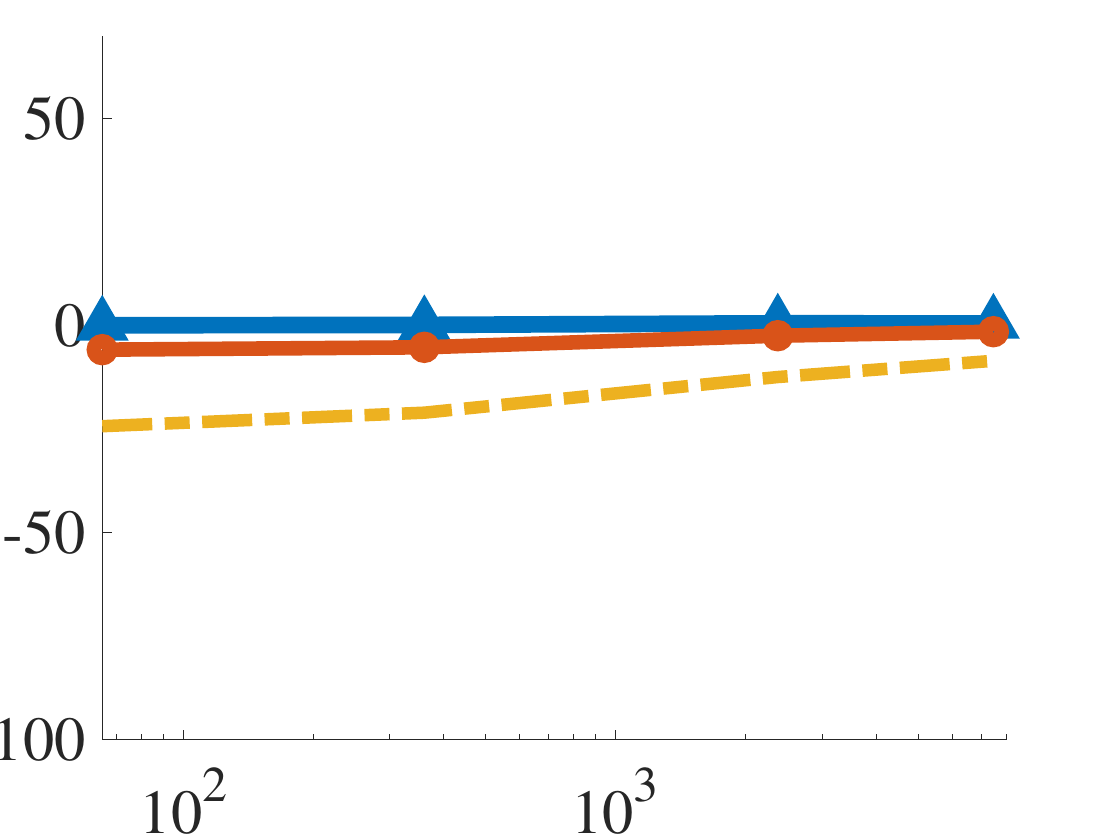}
\includegraphics[width=.45\linewidth, trim={50 190 30 200}, clip]{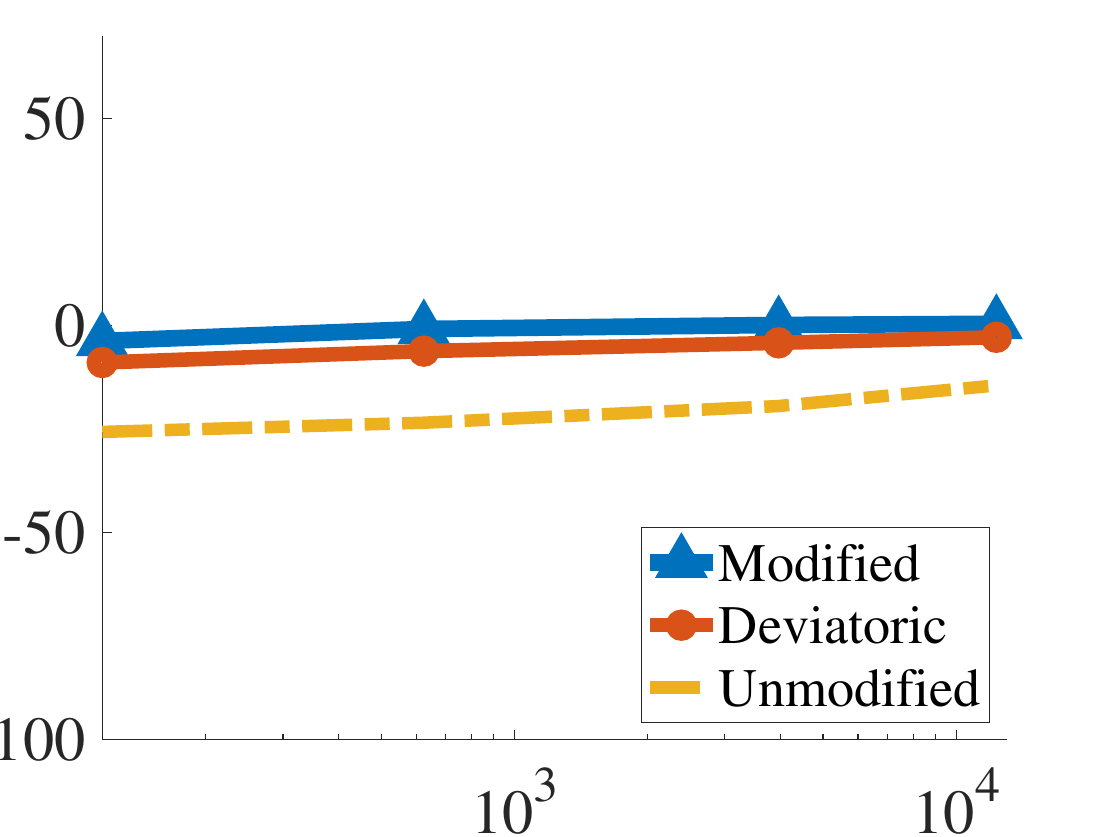}\\

\rotatebox{90}{\qquad\qquad\qquad\; \textbf{$\nus = -1$}} 
   \rotatebox{90}{\qquad\qquad\quad\;  Vol Change \% }
\includegraphics[width=.45\linewidth, trim={50 190 30 200}, clip]{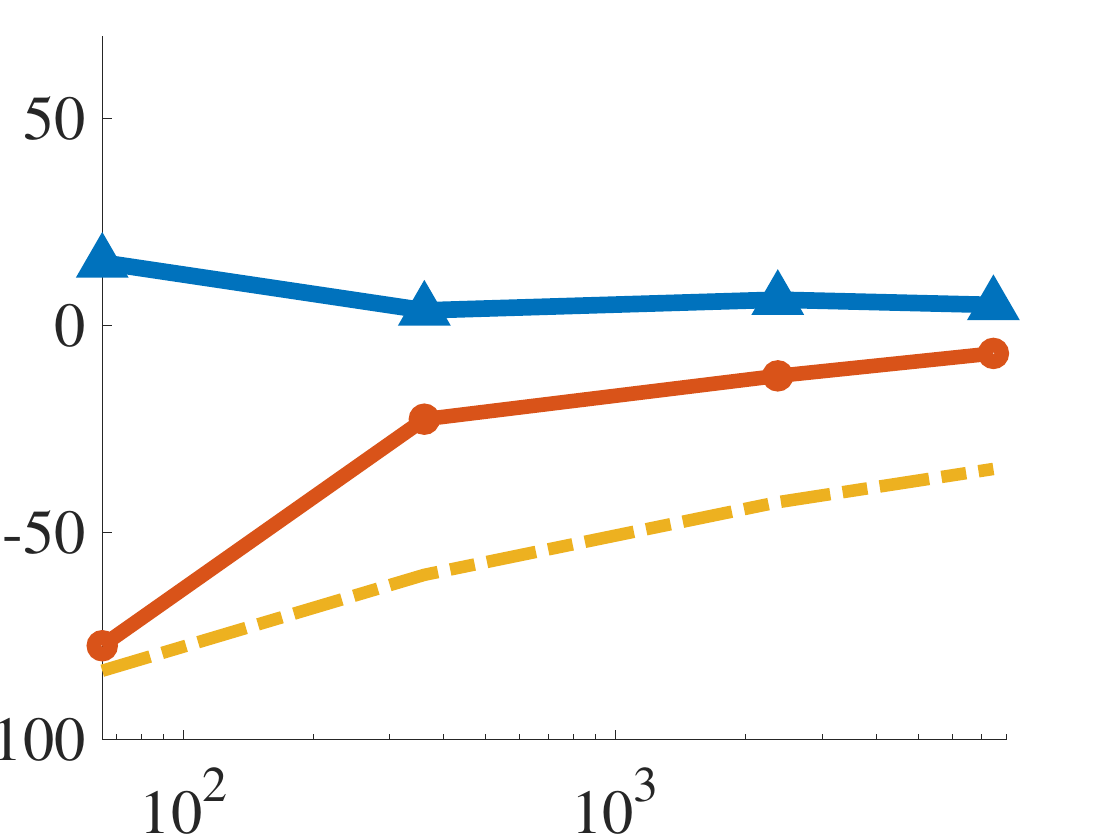}
\includegraphics[width=.45\linewidth, trim={50 190 30 200}, clip]{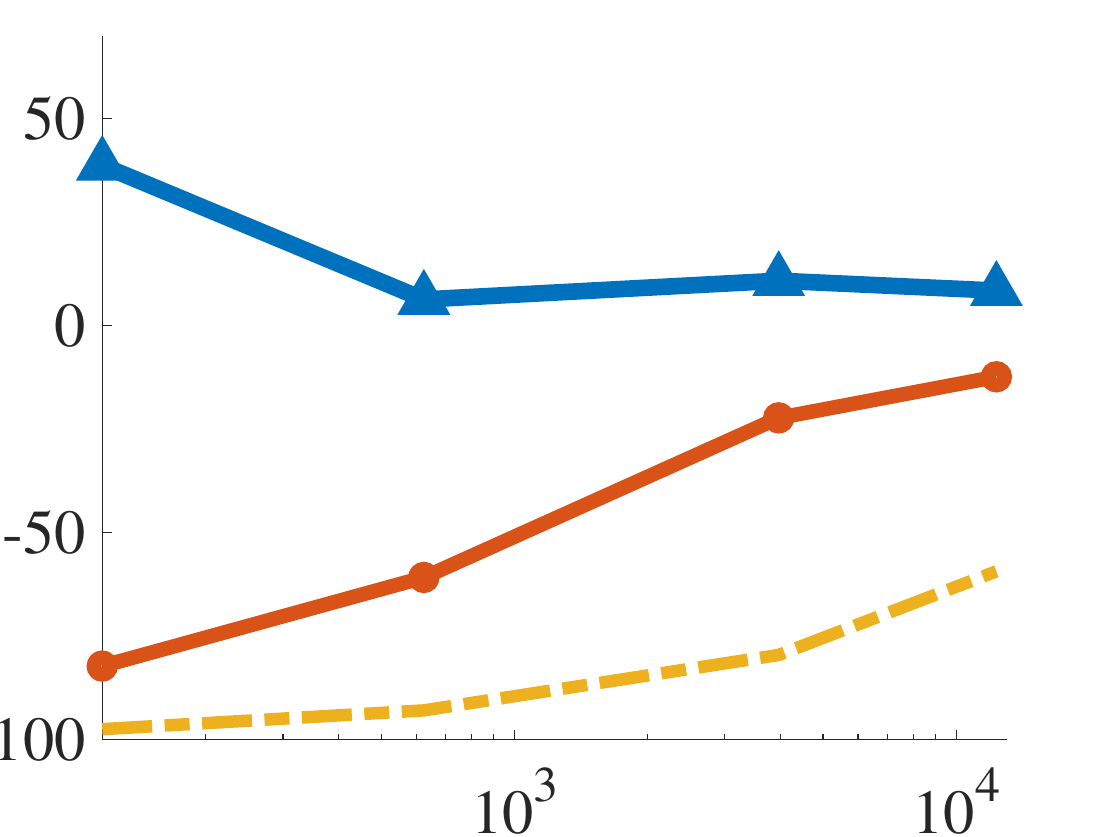}  \\

$\qquad\qquad\qquad\qquad\qquad\quad$ \# Solid DOF $\qquad\qquad\qquad\qquad\qquad\qquad\qquad\quad$  \# Solid DOF\\
\caption{Volume conservation for the torsion test (Section \ref{Torsion}) for different choices of elements and numerical Poisson ratio. The DOF range from $m = 65$ to $m = 12,337$, and the $x$ axis is on a log scale. Omitting the coarsest discretizations ($m=65$), the largest deviations in total volume among all element types used  are approximately $11\%$ for the modified case, $93\%$ for the unmodified case, and $61\%$ for the deviatoric case. In fact, in all cases where there was no volumetric stabilization, the change in volume was greater than $5\%$ for all discretizations. Further, the modified cases with volumetric stabilization were the only ones where the error dropped below $1\%$, and this was achieved for $m = 625$ DOF.}
\label{tor_vol}
\end{figure}

This benchmark is based on a similar test by Bonet \etal~\cite{Bonet2015}.  It involves applying torsion to the top face of an elastic beam, while the opposite face is fixed in place; see Figure (\ref{torsion_diag}). All other faces have zero traction applied. The torsion is applied via displacement boundary conditions, and this face is rotated by $\theta_{\text{f}} = 2.5 \pi$. The angle of rotation $\theta (t)$ increases linearly in time from $0$ to $\theta_{\text{f}}$ and reaches $\theta_{\text{f}}$ at $t = 0.4T_{\text{f}}$, with $T_{\text{f}} = 5.0$~s. We use a Mooney-Rivlin material model, equations (\ref{mr_energy}) -- (\ref{mr_stress_dev}), and with material parameters $c_1 = 9000 \, \frac{\text{dyn}}{\text{cm}^2}$ and $c_2 = 9000 \, \frac{\text{dyn}}{\text{cm}^2}$. The density is $\rho = 1.0 \frac{\text{g}}{\text{cm}^3}$, and the fluid viscosity is set to $\mu = .04 \ \frac{\text{dyn} \cdot \text{s}}{\text{cm}^2}$. The larger viscosity is chosen to allow the model to reach steady state more quickly. The choices of numerical Poisson ratio are the same as the anisotropic Cook's membrane because the computations are in three spatial dimensions. No damping is used. The computational domain is $\Omega = [0, L]^3$ with $L = 9\ \text{cm}$. The numbers of solid DOF range from $m = 65$ to $m = 12,337$. As with the anisotropic Cook's membrane test, \textbf{P1} and \textbf{Q1} meshes use different numbers of solid DOF. \\
\indent Figure (\ref{torsion}) shows the computed deformations for modified invariants, unmodified invariants, and the deviatoric projection as well as for different values of the numerical Poisson ratio. The cases of unmodified invariants and zero numerical bulk modulus lead to the most extremely unphysical deformations in all benchmarks studied; see Figure (\ref{torsion} d). As shown in Figure (\ref{tt_ps}), the principal stretches for the IB method with modified invariants and volumetric stabilization are approximately the same as those from the FE method. Figure ($\ref{tor_disp}$) shows the displacement in the $y$-direction at the center point of the twisted face. Notice that in these plots, the cases with unmodified invariants and zero volumetric energy clearly delineate themselves from other cases. Unique to this test, the convergence of the computed displacement of this case is not deceptive; the convergence is poor and the deformations are also poor. Additionally, the effect of volumetric penalization is more drastic in this benchmark: the percent change in volume is generally much larger than the previous tests; see Figure (\ref{tor_vol}). Specifically, the range of percent change for all element types considered is between $.16\%$ and $11\%$ for the modified invariants, between $8.5\%$ and $93\%$ for the unmodified invariants, and between $1.5\%$ and $61\%$ for the deviatoric projection. The choice of numerical Poisson ratio also has a large effect on the displacement of the twisted face, which can be seen in Figure (\ref{tor_disp}). The differences between displacement curves for with and without volumetric penalization is more apparent, with the case of volumetric penalization performing much better here. We contrast that with the anisotropic Cook's membrane benchmark, which has only slight differences between the displacement curves for $\nus = .4$ and $\nus = -1$. \\
\indent Finally, in this benchmark the deviatoric projection delineates itself from the modified invariants. In this test both the volume conservation and the displacement performed worse for the deviatoric projection than the modified invariants; see Figures (\ref{tor_disp}) and (\ref{tor_vol}). For the other benchmarks, there was a negligible difference between the solution produced by the modified invariants and that produced by the deviatoric projection.

%%%%%%%%%%%%%%%%%%%%%

\subsection{FSI Test: Elastic Band}
\label{Elastic Band}
\begin{figure}
\centering
\includegraphics[width=.7\linewidth, trim={10 50 10 0}, clip]{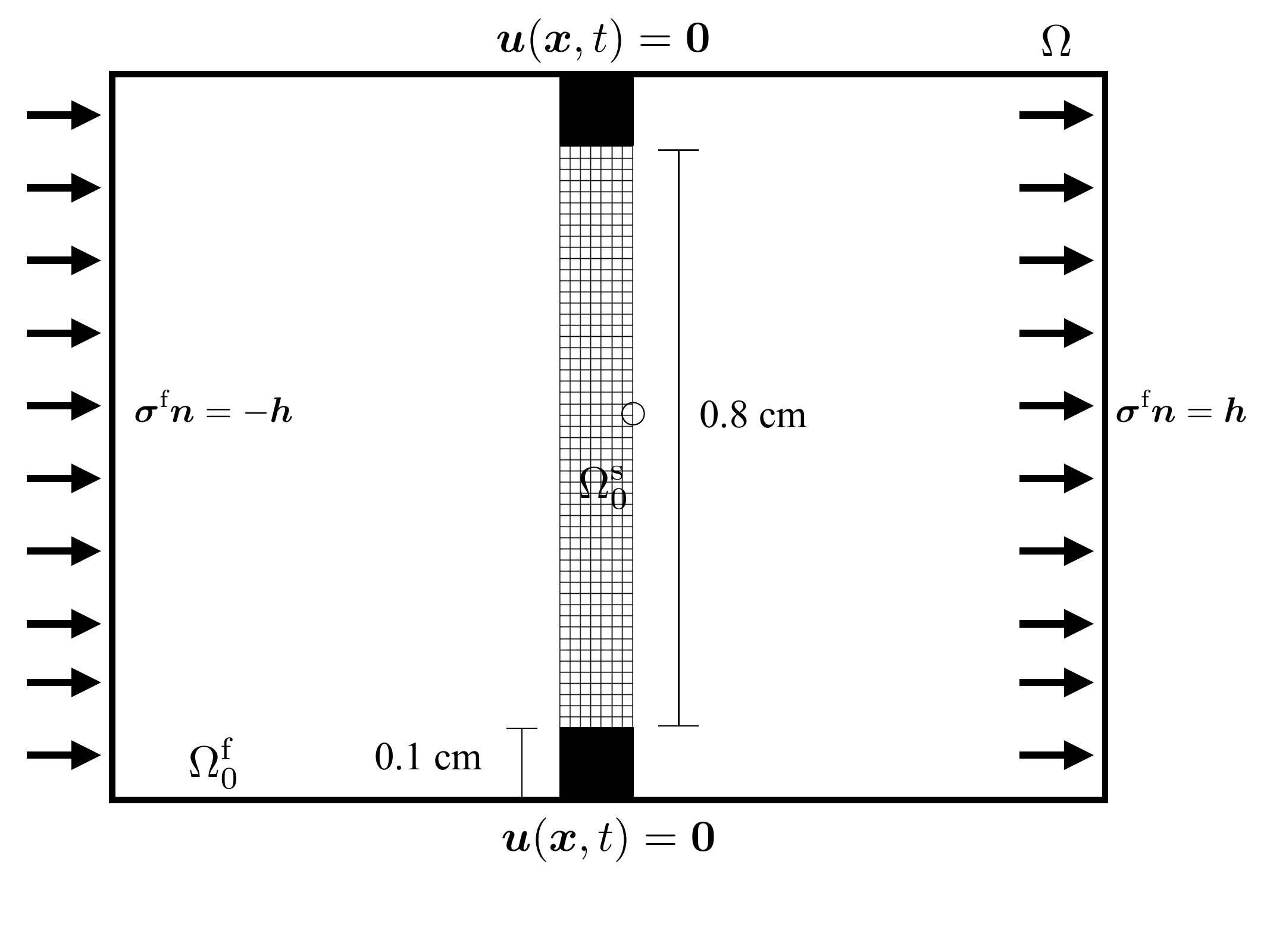}
\caption{Specifications of the fluid-loaded elastic band benchmark (Section \ref{Elastic Band}). The block boxes are rigid structures that help tether the band in place. Fluid traction boundary conditions are applied on the left and right faces. Tests described herein set $\hb(t) = \left(10 \sin \left( \frac{\pi t}{2 T_{\text{l}}}\right),0 \right)^T \frac{\text{dyn}}{\text{cm}^2}$ when $t < T_{\text{l}}$ and $\hb(t) = (10, 0)^T \frac{\text{dyn}}{\text{cm}^2}$ otherwise. The quantity of interest is the $x$-displacement as measured at the encircled point.}
\label{elastic_diag}
\end{figure}

\begin{figure}
\captionsetup[subfigure]{justification=centering}
\begin{subfigure}{.24\textwidth}
  \centering
  \includegraphics[width=.95\linewidth]{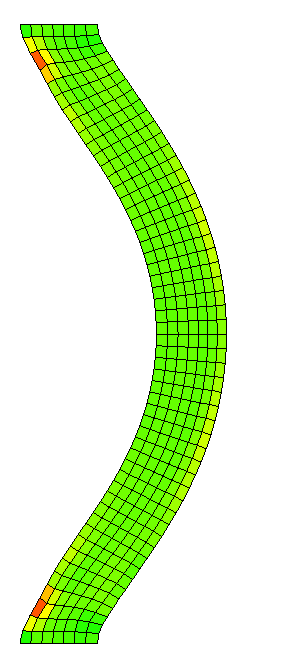}
  \caption{Modified Invariants \\ $\nus = .4$}
%  \label{fig:fixed_end_A}
\end{subfigure}
\begin{subfigure}{.24\textwidth}
  \centering
\includegraphics[width=.95\linewidth]{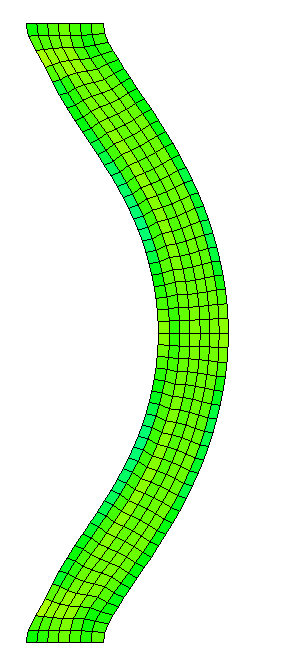}
  \caption{Unmodified Invariants \\ $\nus = .4$}
%  \label{fig:fixed_end_B}
\end{subfigure}
\begin{subfigure}{.24\textwidth}
  \centering
 \includegraphics[width=.95\linewidth,trim={0 0 0 0},clip]{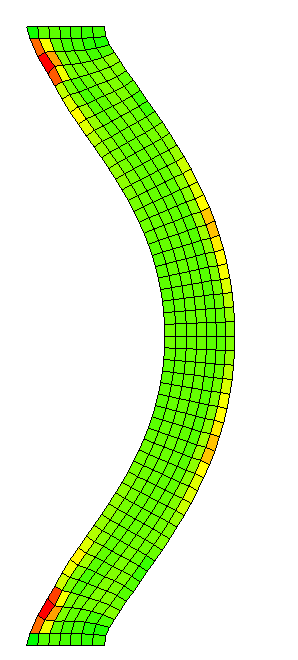}
  \caption{Modified Invariants \\ $\nus = -1$}
%  \label{fig:fixed_end_C}
\end{subfigure}
\begin{subfigure}{.24\textwidth}
  \centering
\includegraphics[width=.95\linewidth,trim={0 0 0 0},clip]{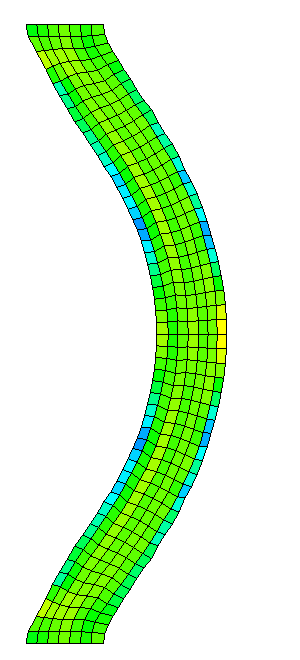}
  \caption{Unmodified Invariants \\ $\nus = -1$}
%  \label{fig:fixed_end_D}
\end{subfigure}

\begin{centering}
Avg $J$ \\
\includegraphics[width=2.5in, trim={0 5in 0 5in}, clip]{color_bar.pdf}  \\
%0.00 \ \ \ \ \ \ \ \ \ \ \ \ \ \ \ \ 2.00
0.85 $\qquad\qquad\qquad\qquad$ 1.10\\
\end{centering}
\caption{Deformations of the elastic band (Section \ref{Elastic Band}), along with mean values of $J$ within each element calculated via equation \eqref{avgJ}, for the steady state version of the test. This test uses a neo-Hookean model, equations (\ref{nh_energy}) -- (\ref{nh_stress_mod}), for the solid, equations, with $G = 200 \frac{\text{dyn}}{\text{cm}^2}$. The background Eulerian grid is not shown. Shown here are solid meshes with \textbf{Q1} elements and $m = 424$ solid DOF. Deformations are depicted at final time, $t = 15$ s. We show cases with $\nus = .4$ and $\nus = -1$ (here equivalent to $\kappas = 0$ and no volumetric-based stabilization) and cases with modified invariants and unmodified invariants.}
\label{eb}
\end{figure}

\begin{figure}
\begin{tabular}{c r}
\begin{subfigure}{.44\textwidth}
  \includegraphics[width=.85\linewidth,trim={0 0 0 0},clip]{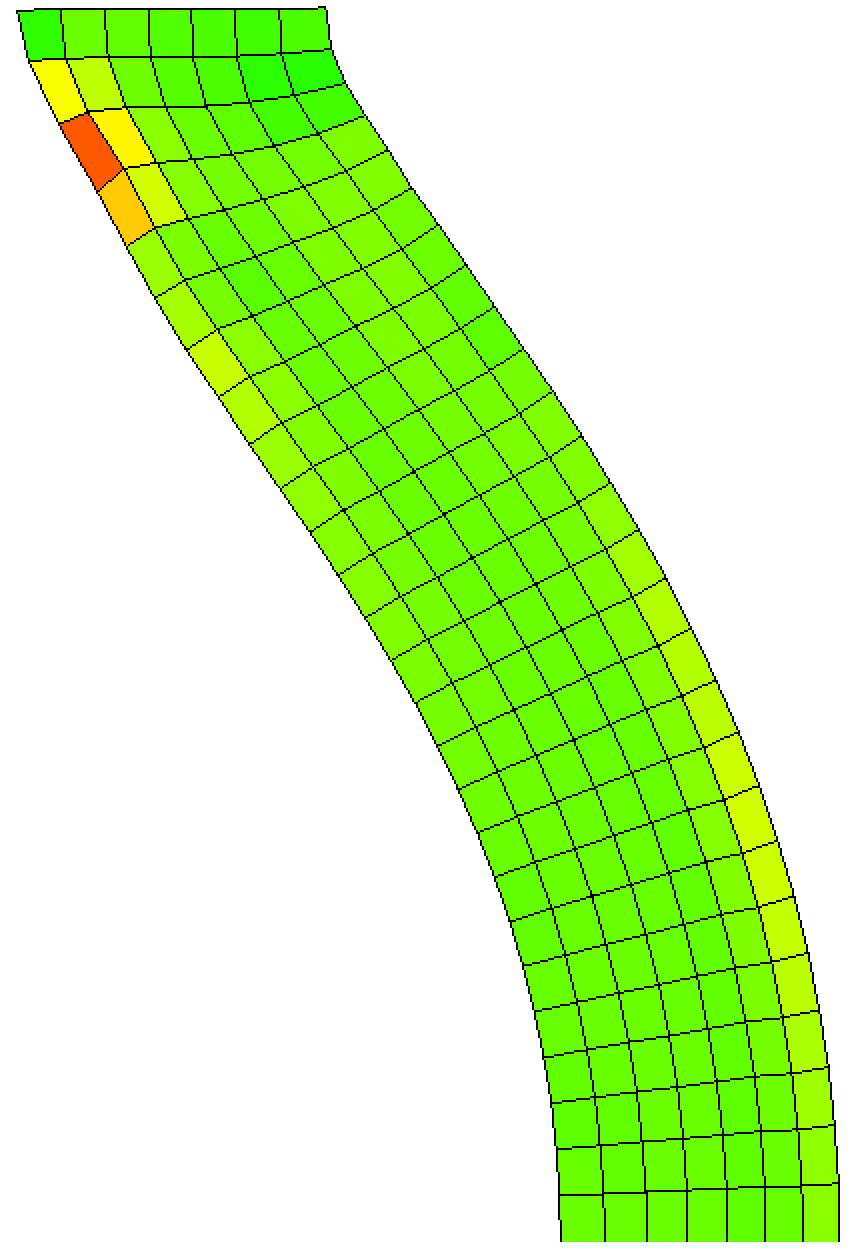}
\caption{Modified Invariants \\ $\nus = .4$}
\end{subfigure}

\begin{subfigure}{.44\textwidth}
\includegraphics[width=.85\linewidth,trim={0 0 0 0},clip]{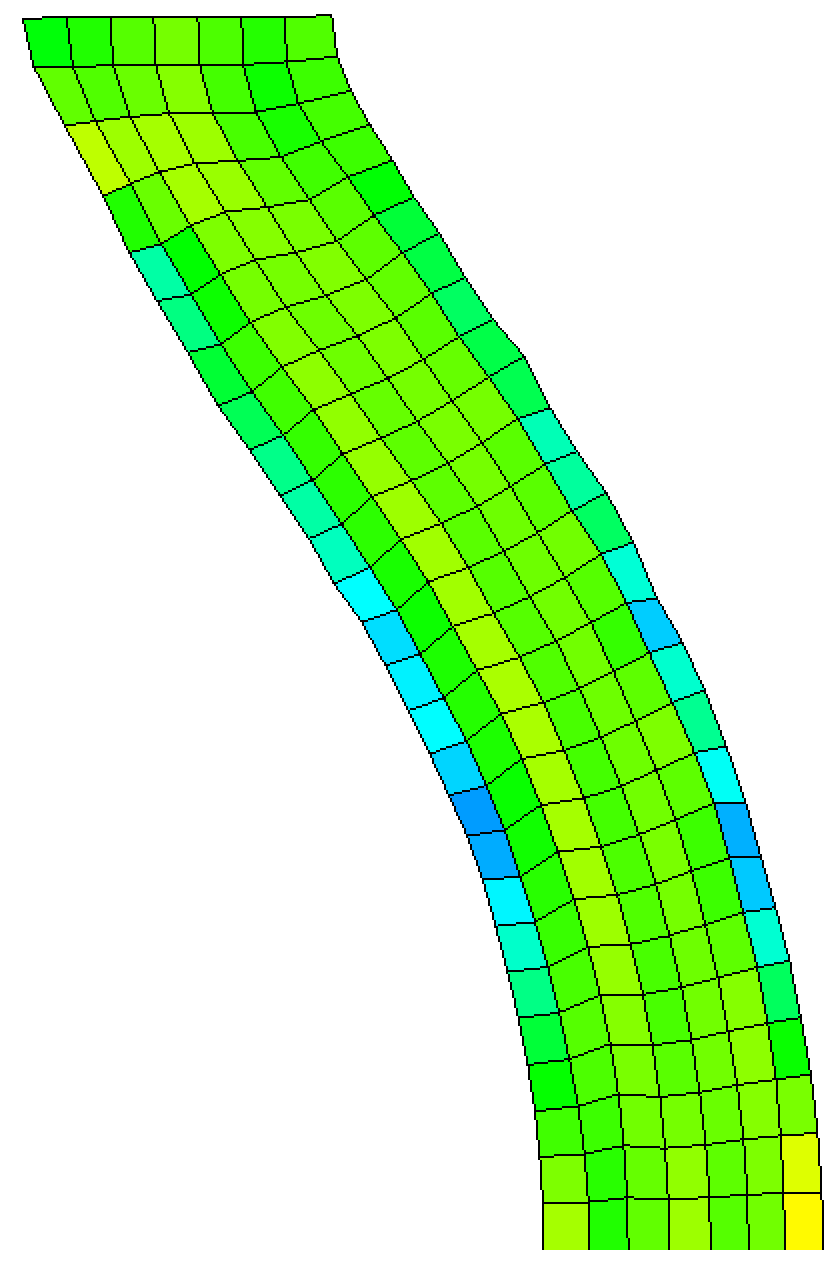}
\caption{Unmodified Invariants \\ $\nus = -1$}
\end{subfigure}

\end{tabular}
\vskip .25cm
\begin{centering}
Avg $J$ \\
\includegraphics[width=2.5in, trim={0 5in 0 5in}, clip]{color_bar.pdf}  \\
%0.00 \ \ \ \ \ \ \ \ \ \ \ \ \ \ \ \ 2.00
0.85 $\qquad\qquad\qquad\qquad$ 1.10\\
\end{centering}
\caption{Deformations and mean values of $J$ of the elastic band (Section \ref{Elastic Band}), as also shown in Figure (\ref{eb}). Note the smoother deformations in the panel (a), the modified and stabilized case, as compared to panel (b), the unmodified and unstabilized case.}
\label{eb_zoom}
\end{figure}

\begin{figure}
%trim={left bottom right top}
$\qquad\qquad\qquad\;\;\;\,$ \textbf{P1} $\qquad\qquad\qquad\qquad\quad$  \textbf{Q1} $\qquad\qquad\qquad\qquad\;\;\;$  \textbf{P2} $\qquad\qquad\qquad\qquad\quad\;$ \textbf{Q2}\\
\rotatebox{90}{$\quad$ \textbf{$\nus = .49995$} }
   \rotatebox{90}{$\quad\;\;$ Disp. (cm) }
\includegraphics[width=.225\linewidth, trim={10 190 25 200}, clip]{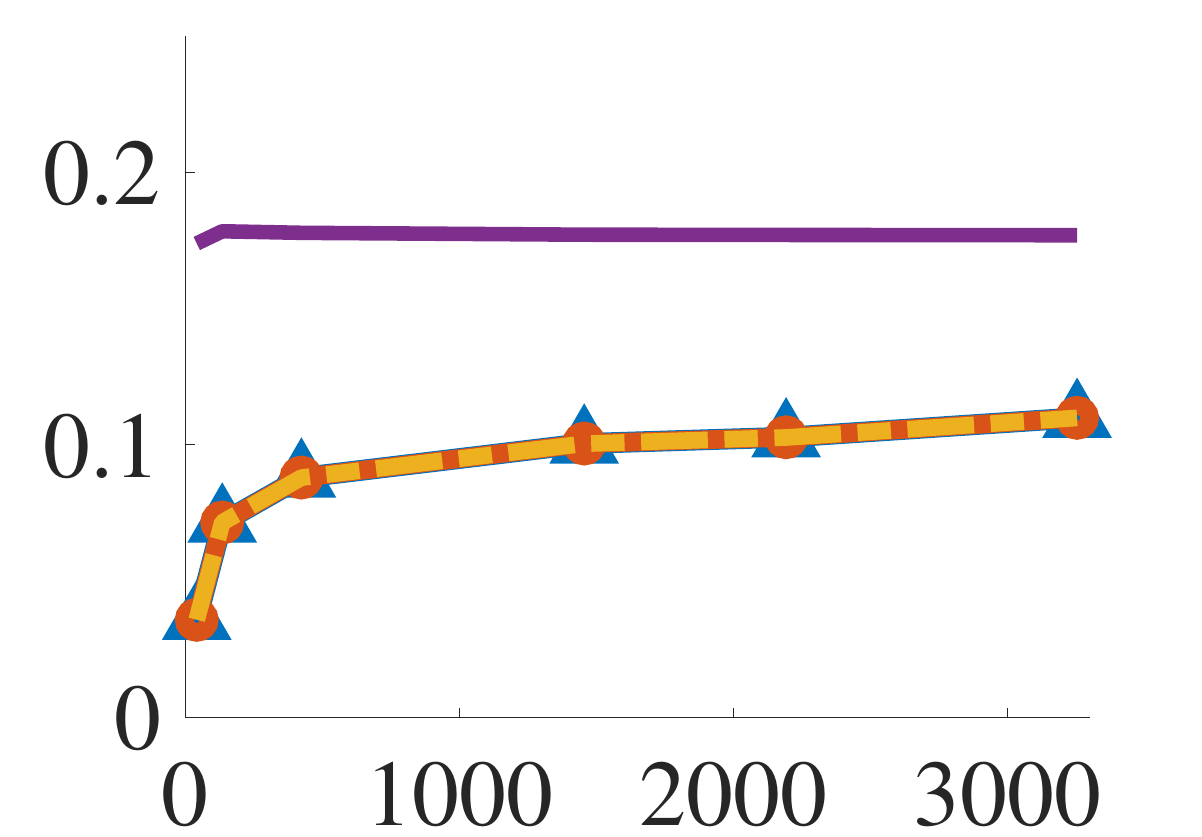} 
\includegraphics[width=.225\linewidth, trim={10 190 25 200}, clip]{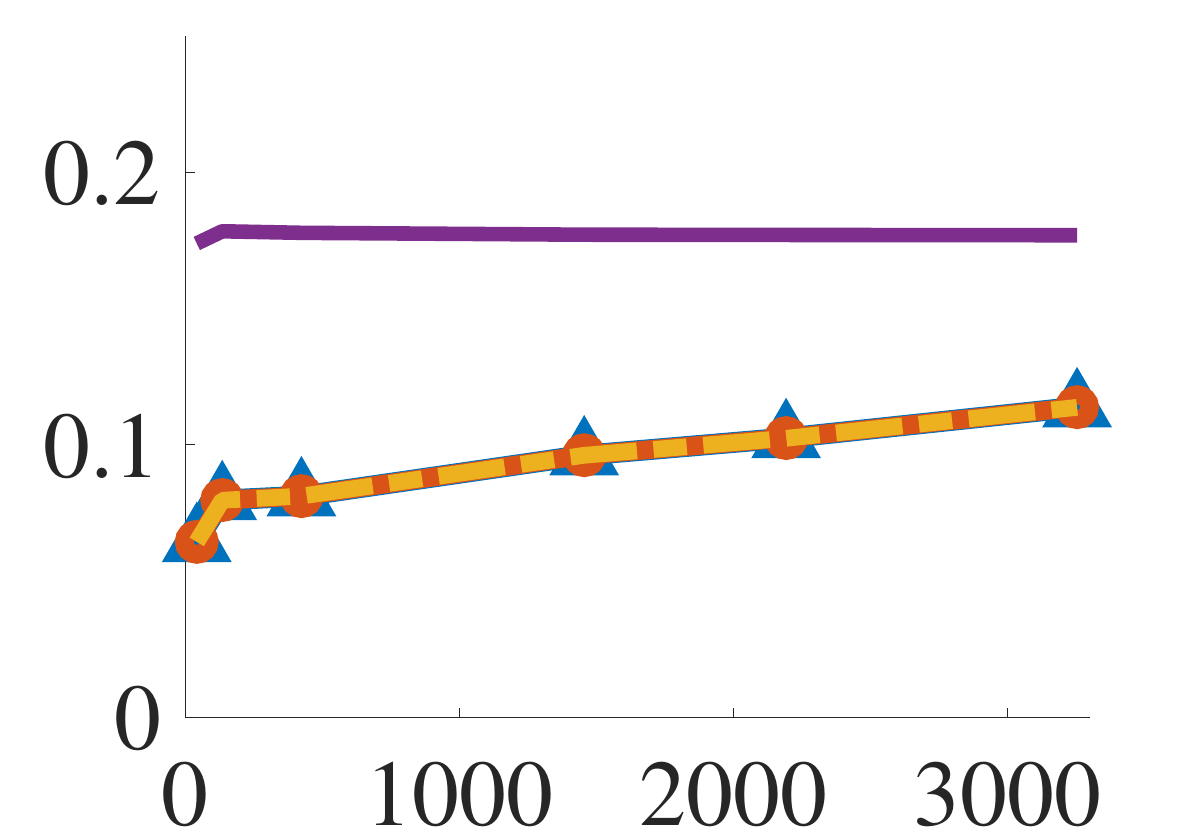}
\includegraphics[width=.225\linewidth, trim={10 190 25 200}, clip]{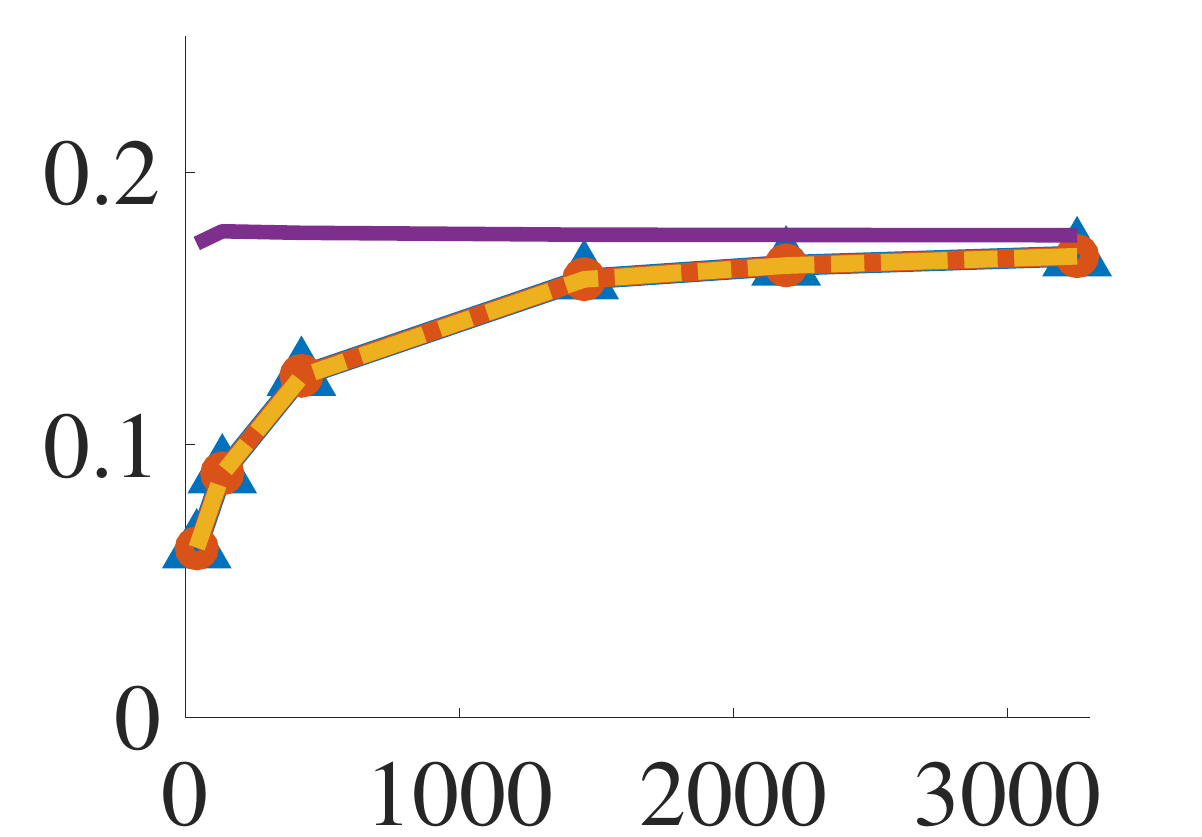}
\includegraphics[width=.225\linewidth, trim={10 190 25 200}, clip]{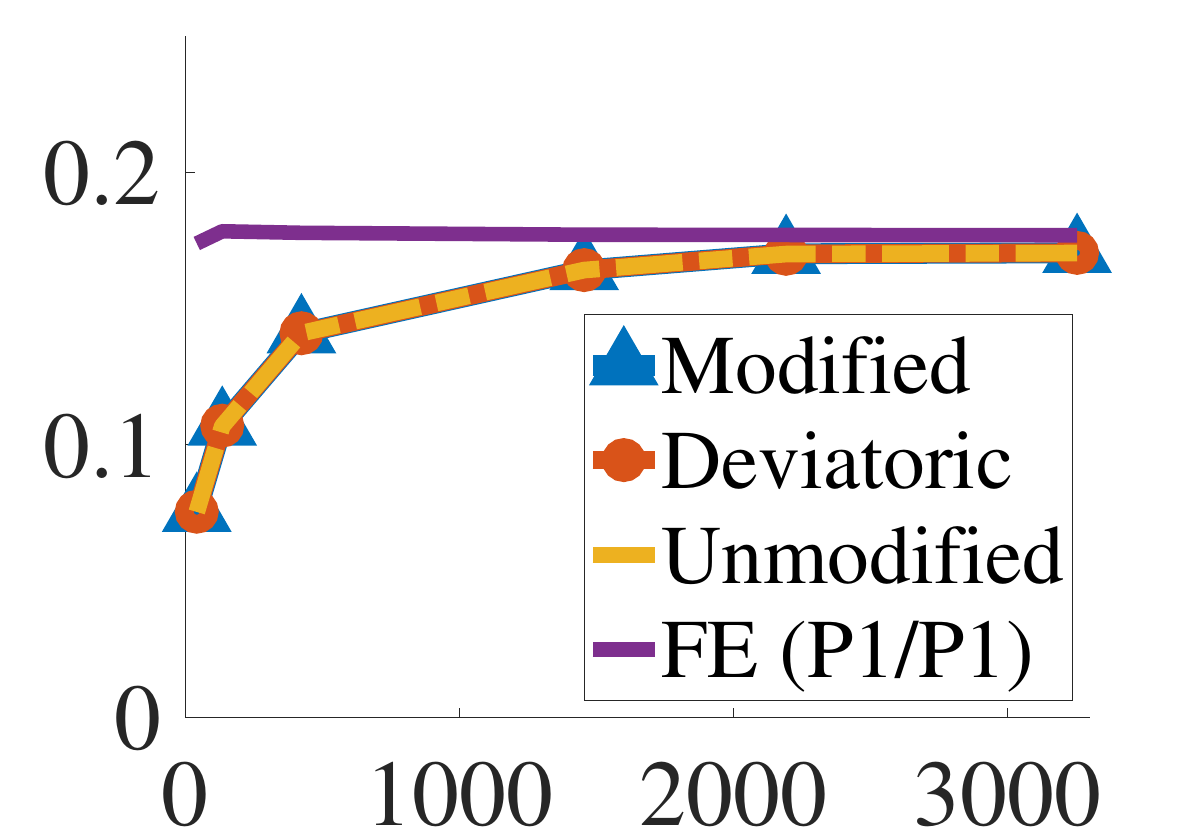} \\

\rotatebox{90}{$\quad\;\;\;$ \textbf{$\nus = .4$} }
   \rotatebox{90}{$\quad\;\;$ Disp. (cm) }
\includegraphics[width=.225\linewidth, trim={10 190 25 200}, clip]{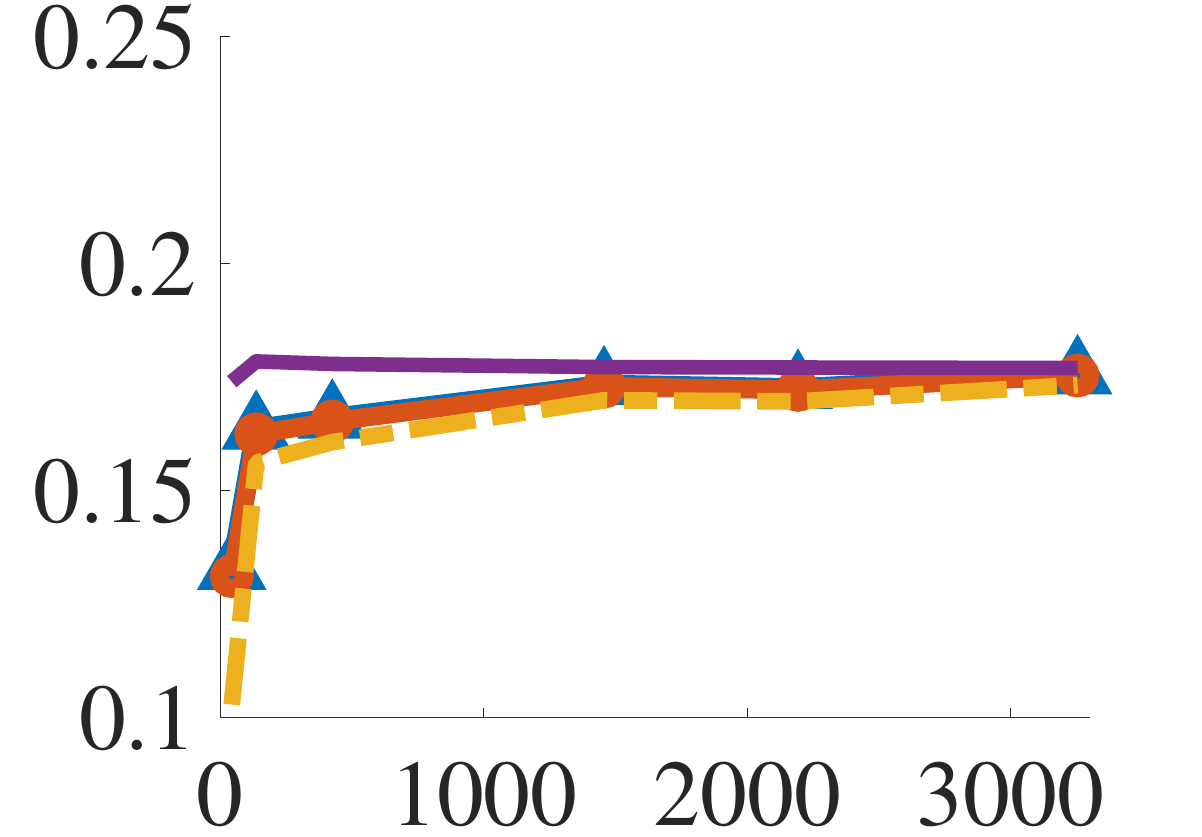} 
\includegraphics[width=.225\linewidth, trim={10 190 25 200}, clip]{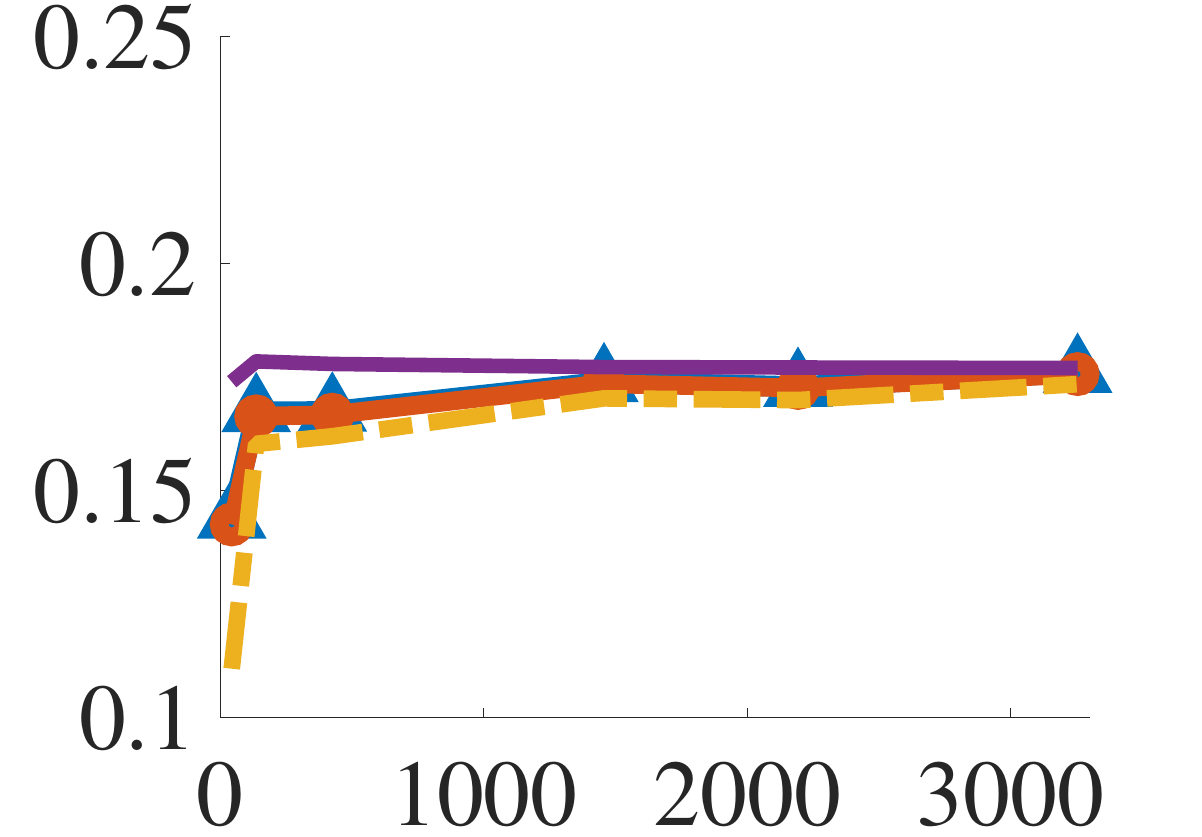}
\includegraphics[width=.225\linewidth, trim={10 190 25 200}, clip]{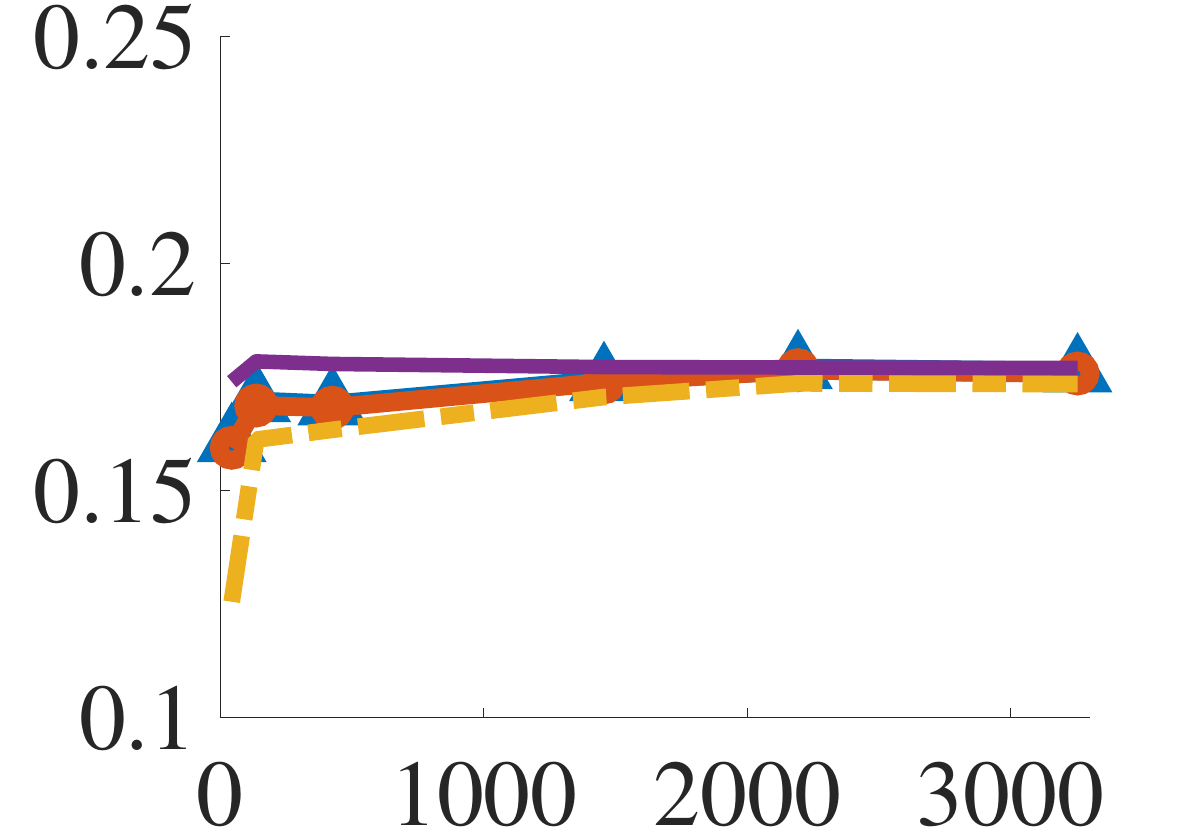}
\includegraphics[width=.225\linewidth, trim={10 190 25 200}, clip]{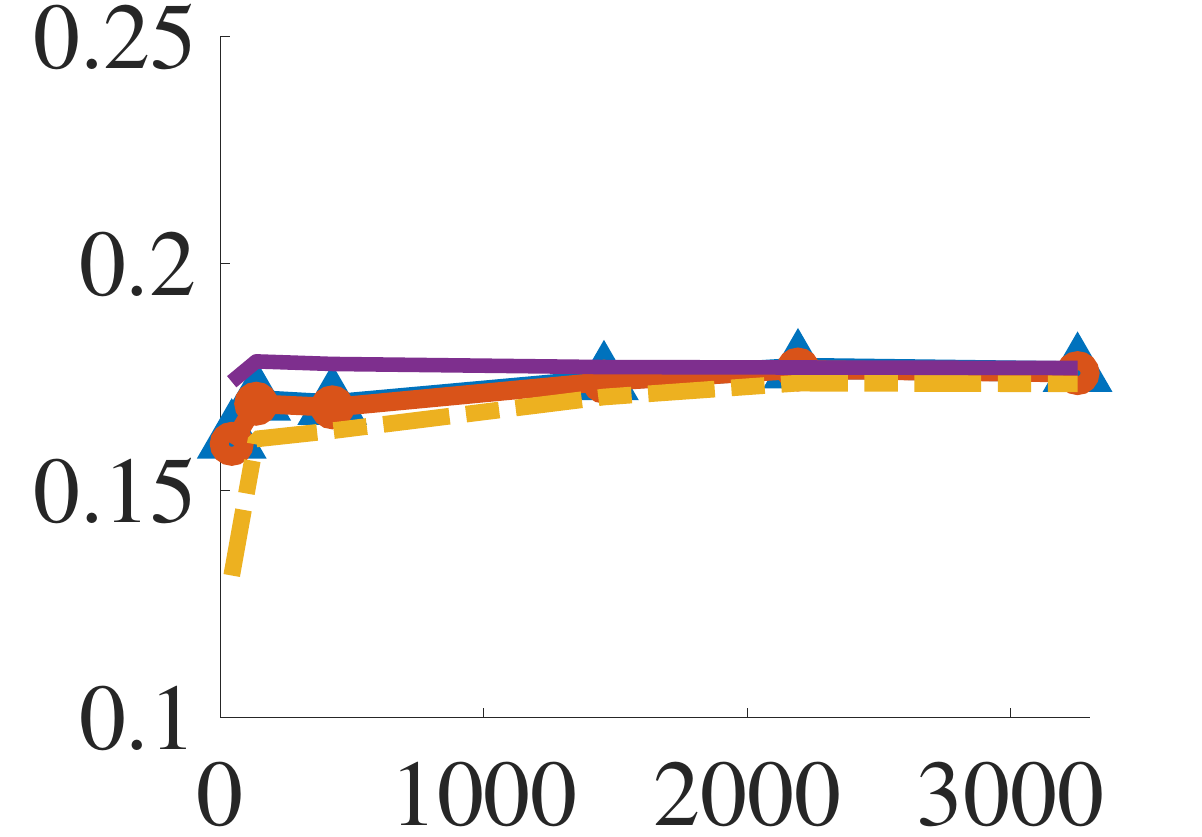}\\

\rotatebox{90}{$\qquad$ \textbf{$\nus = 0$} }
   \rotatebox{90}{$\quad\;\;$ Disp. (cm) }
\includegraphics[width=.225\linewidth, trim={10 190 25 200}, clip]{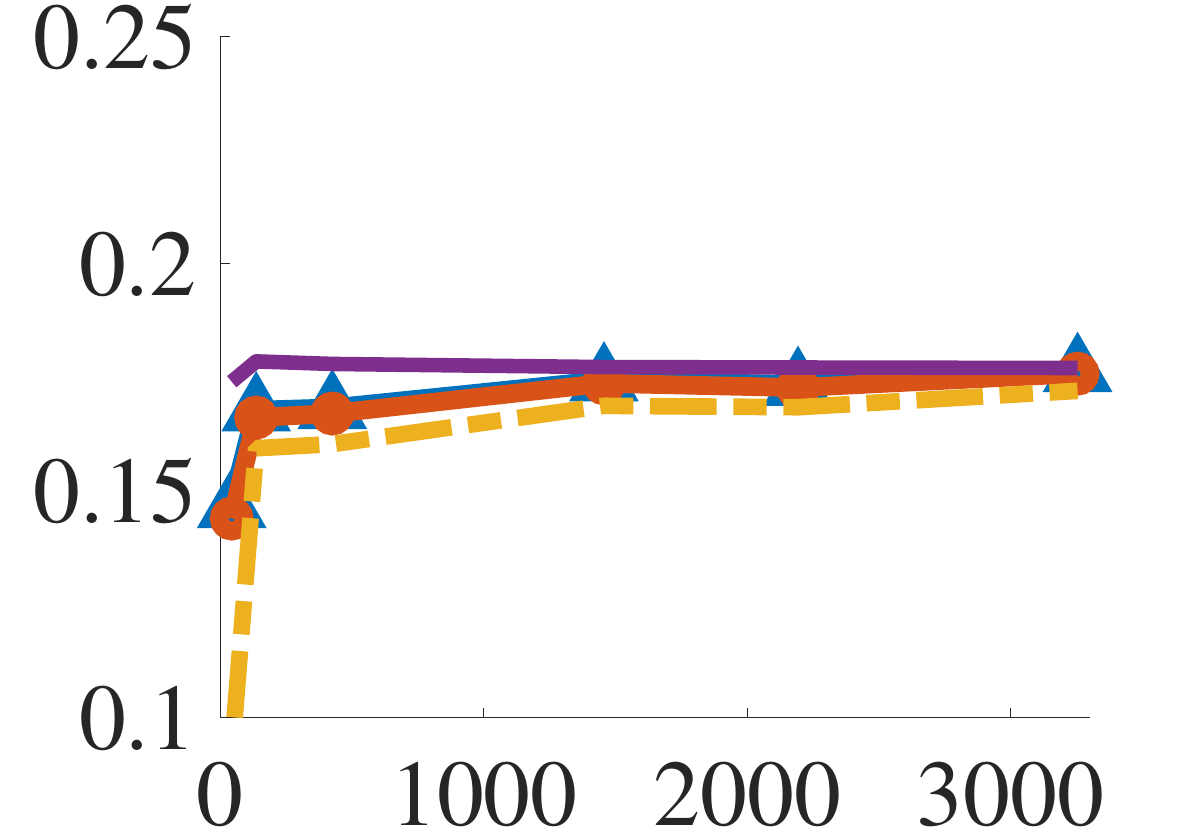}
\includegraphics[width=.225\linewidth, trim={10 190 25 200}, clip]{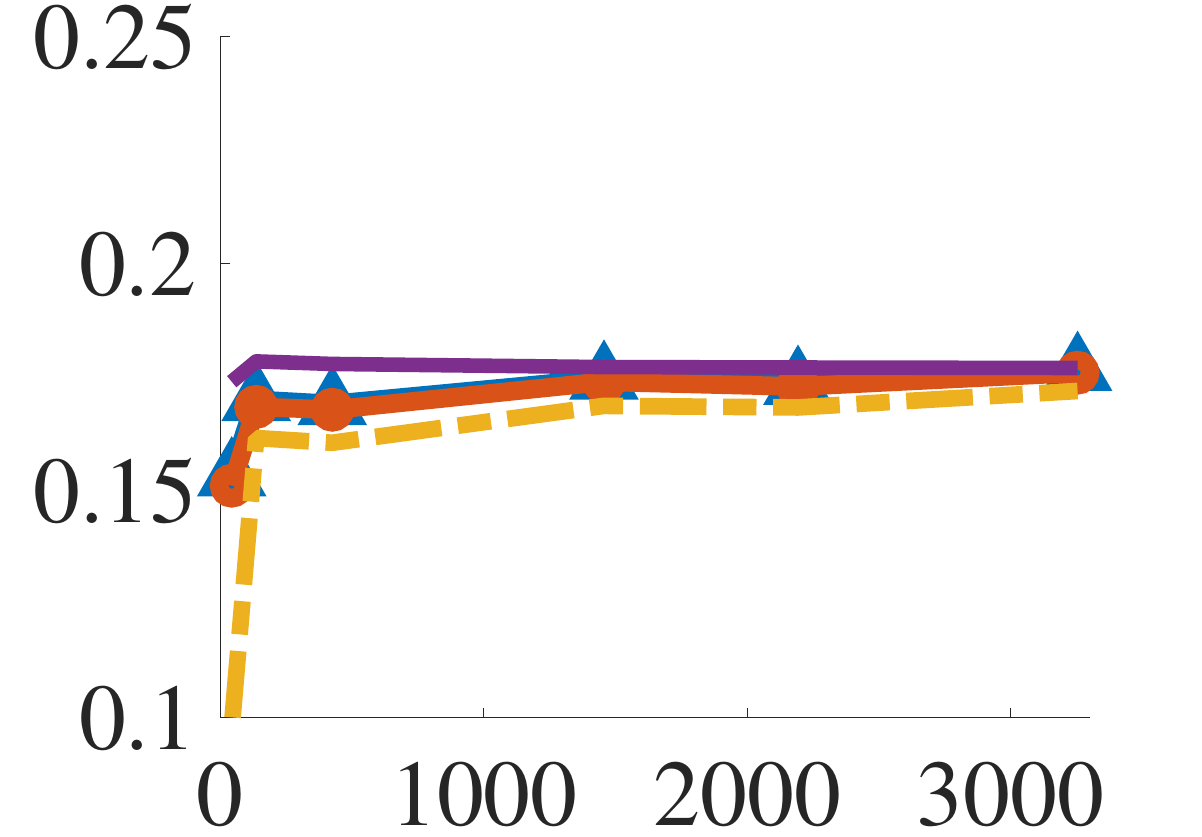}
\includegraphics[width=.225\linewidth, trim={10 190 25 200}, clip]{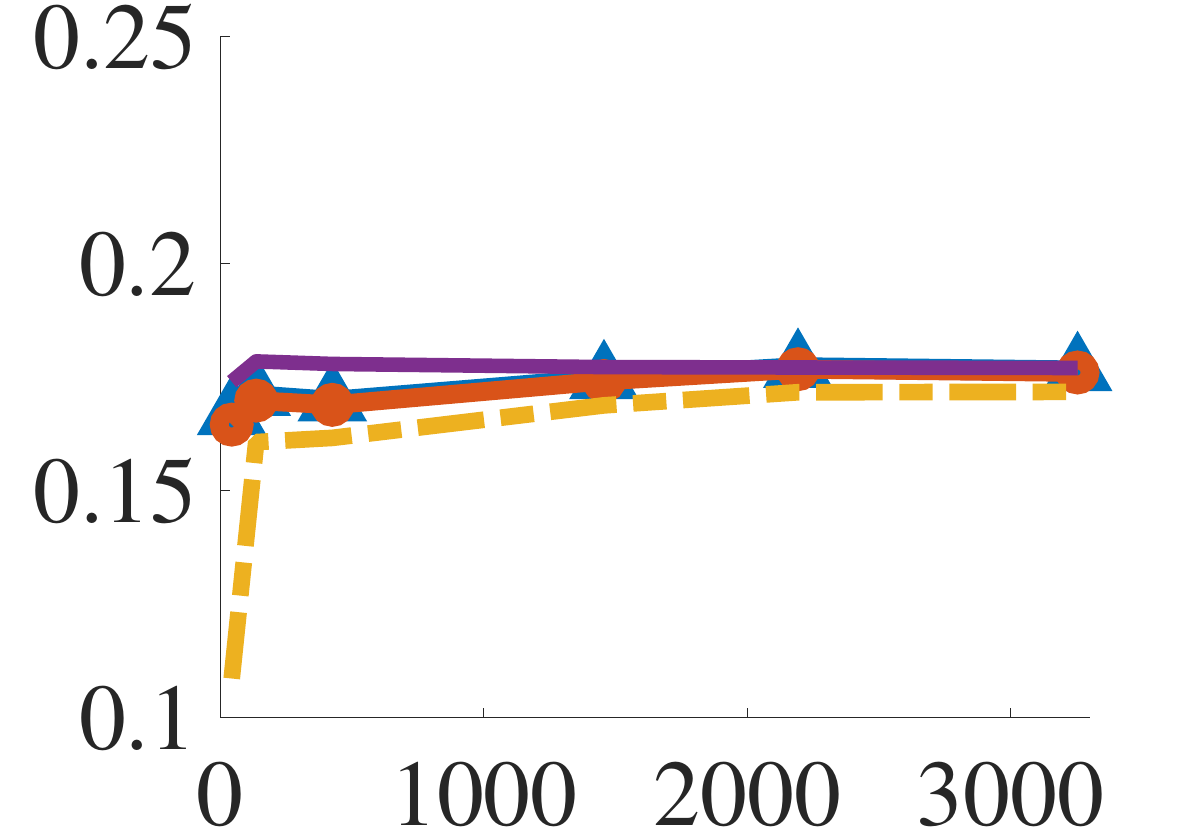}
\includegraphics[width=.225\linewidth, trim={10 190 25 200}, clip]{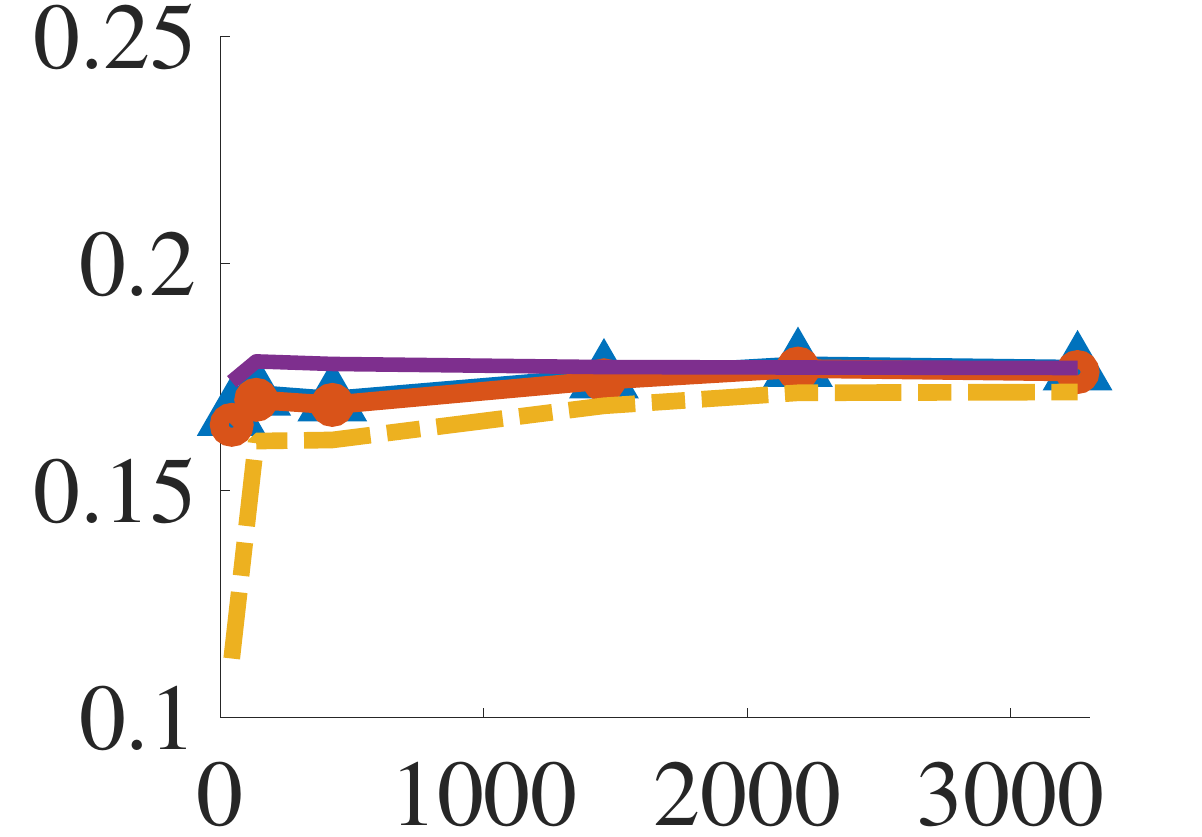}\\

\rotatebox{90}{$\quad\;\;\,$ \textbf{$\nus = -1$} }
   \rotatebox{90}{$\quad\;\;$ Disp. (cm) }
\includegraphics[width=.225\linewidth, trim={10 190 25 200}, clip]{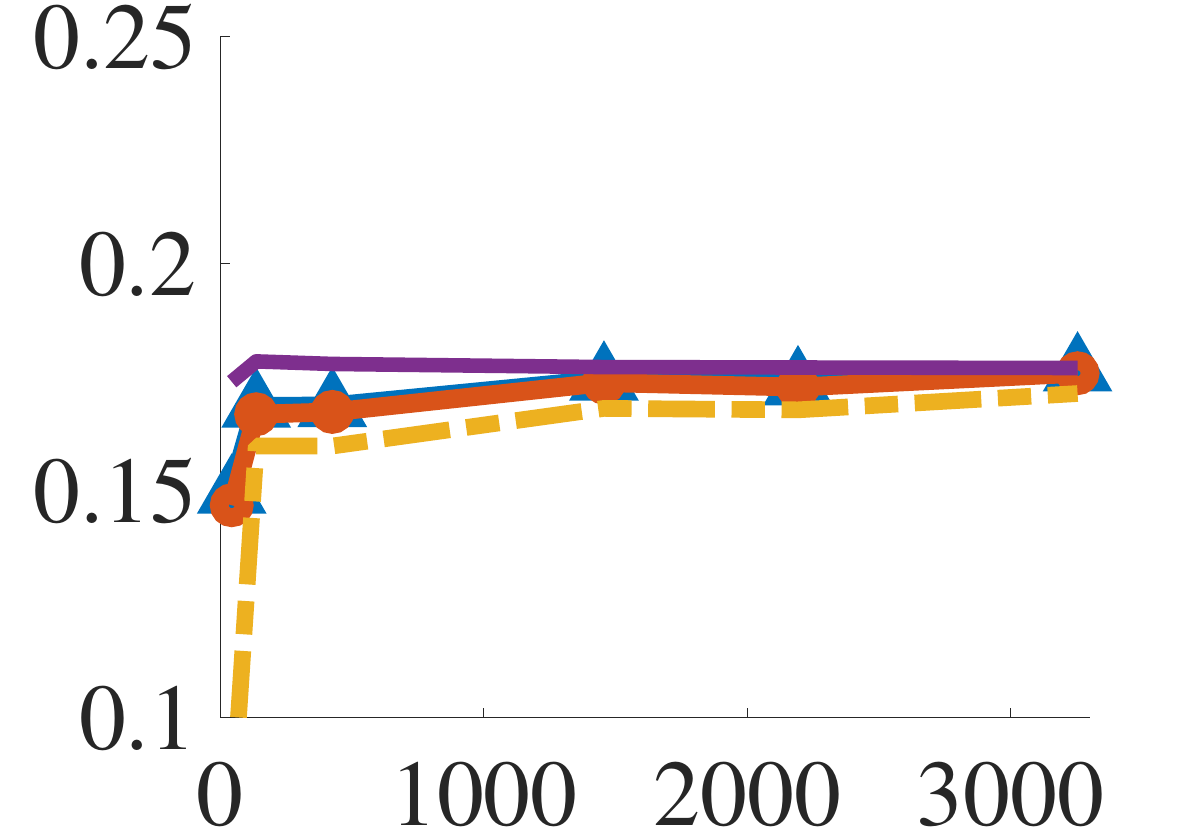}
\includegraphics[width=.225\linewidth, trim={10 190 25 200}, clip]{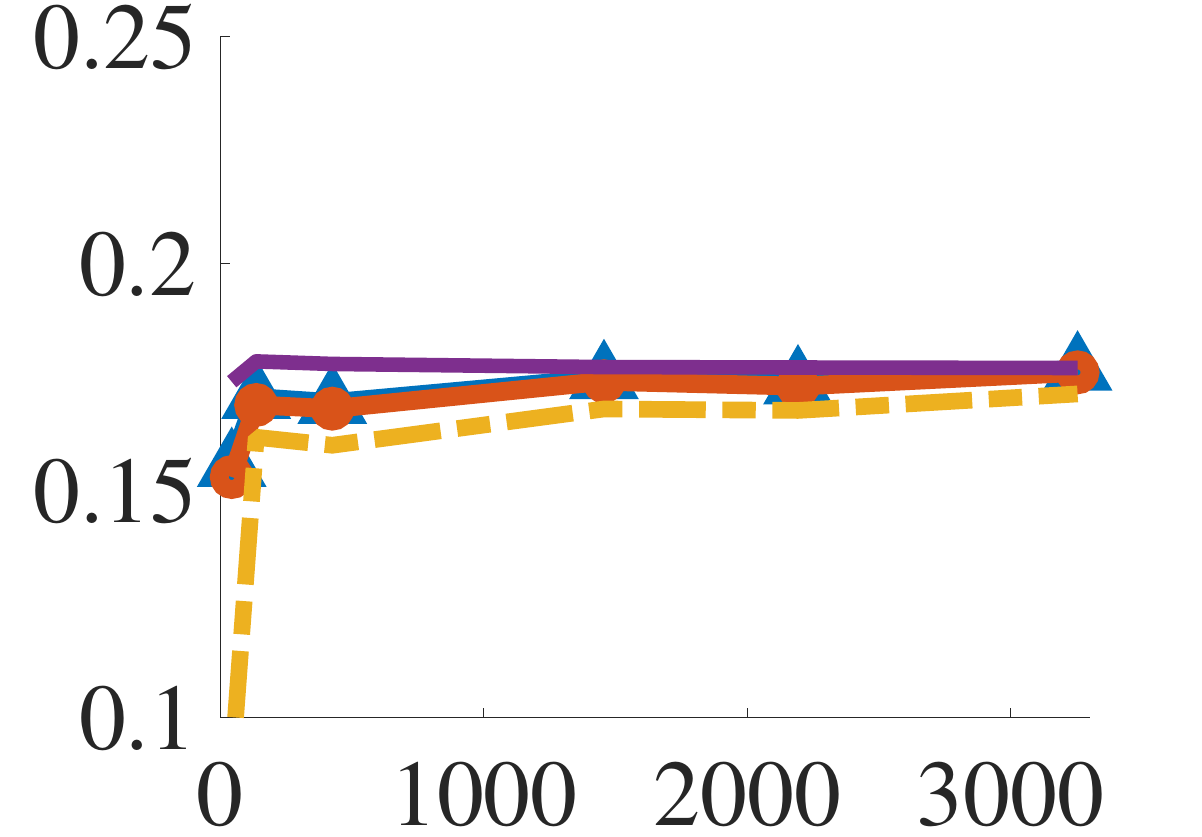}
\includegraphics[width=.225\linewidth, trim={10 190 25 200}, clip]{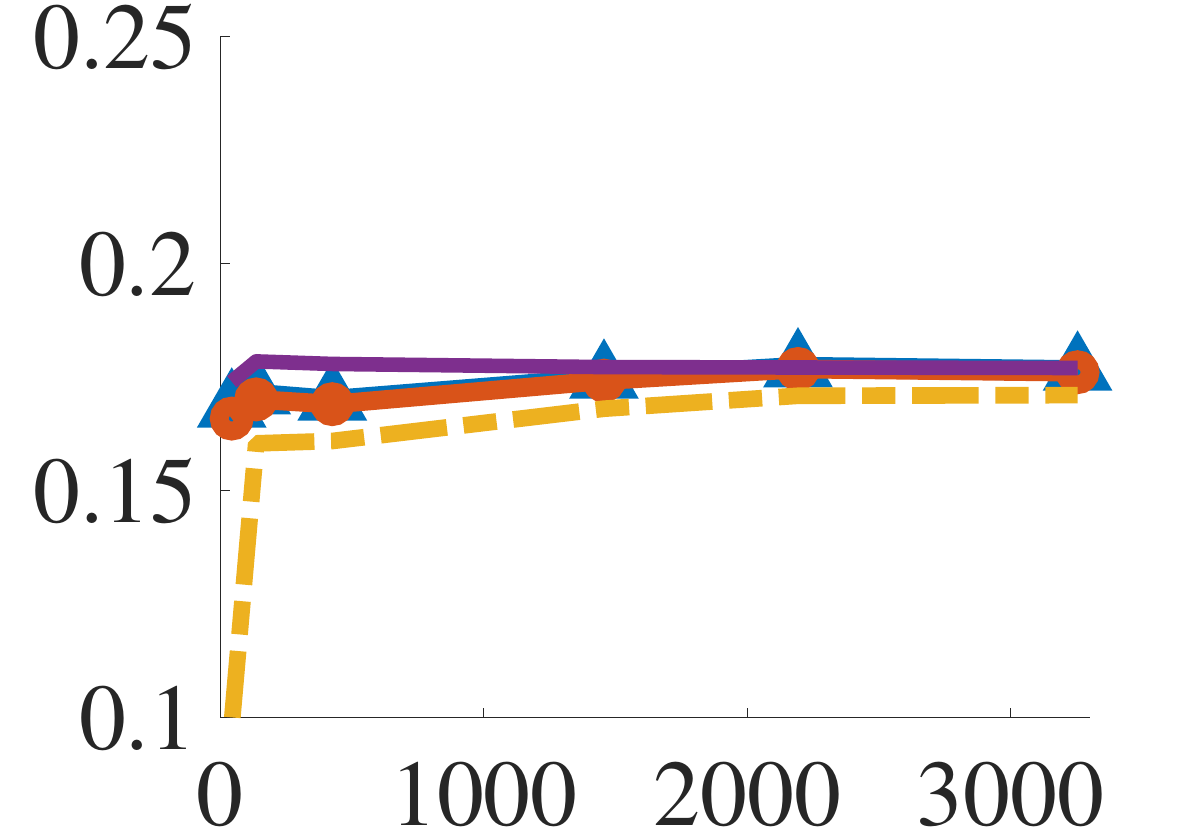}
\includegraphics[width=.225\linewidth, trim={10 190 25 200}, clip]{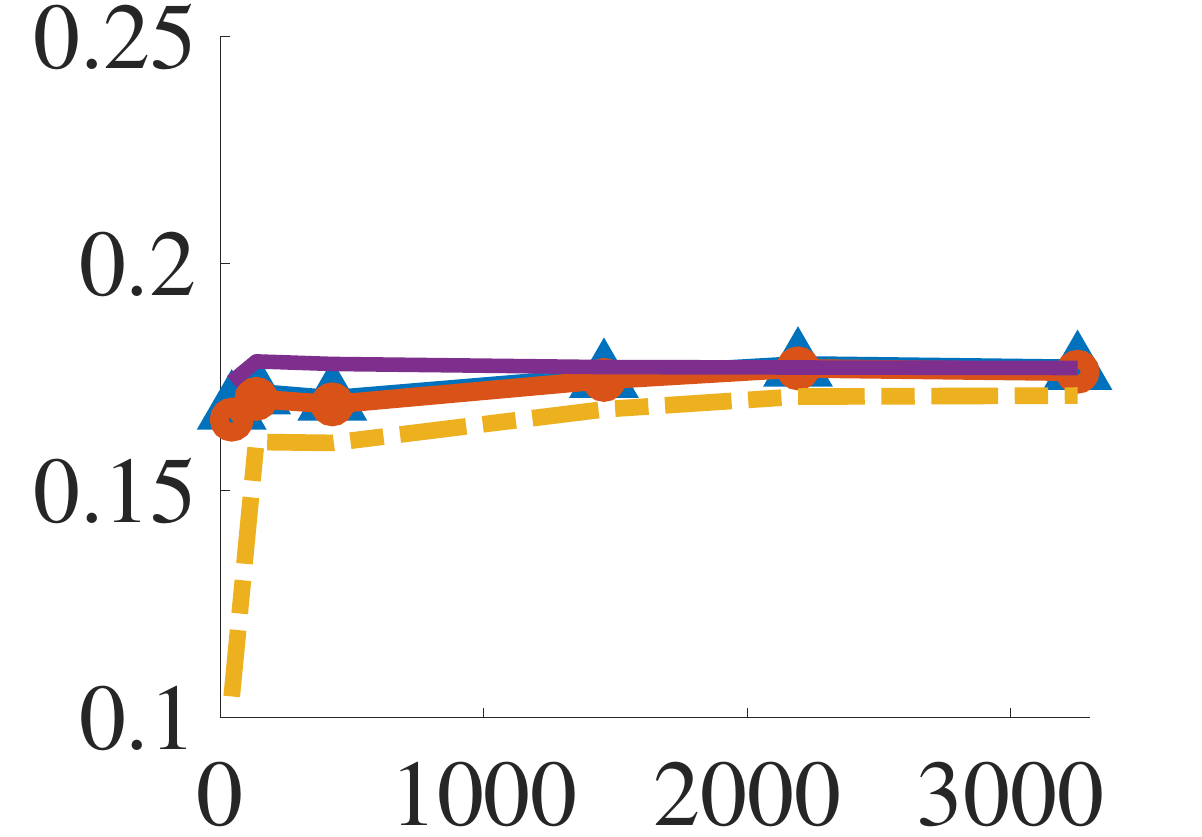}\\

$\qquad\qquad\quad$ \# Solid DOF $\qquad\qquad\quad\;$ \# Solid DOF $\qquad\qquad\quad$ \# Solid DOF $\qquad\qquad\quad\;$ \# Solid DOF
\caption{$x$-displacement of the center point in Figure (\ref{elastic_diag}) for the steady state elastic band benchmark (Section \ref{Elastic Band}) for different choices of elements and numerical Poisson ratio. The solid DOF range from $m = 42$ to $m = 3255$. Note that for $\nus = .49995$ low order elements produce volumetric locking, and higher order elements are needed for convergence at reasonable numbers of DOF. Also note that the differences between different methods are somewhat less pronounced in these tests.}
\label{elastic_disp}
\end{figure}

\begin{figure}
\label{elastic-area}
$\qquad\qquad\qquad\;\;\;\;$ \textbf{P1} $\qquad\qquad\qquad\qquad\quad$  \textbf{Q1} $\qquad\qquad\qquad\qquad\;\;\;$  \textbf{P2} $\qquad\qquad\qquad\qquad\quad\;$ \textbf{Q2}\\
\rotatebox{90}{$\quad$ \textbf{$\nus = .49995$} }
   \rotatebox{90}{$\quad$ Area Change \% }
\includegraphics[width=.225\linewidth, trim={30 190 25 200}, clip]{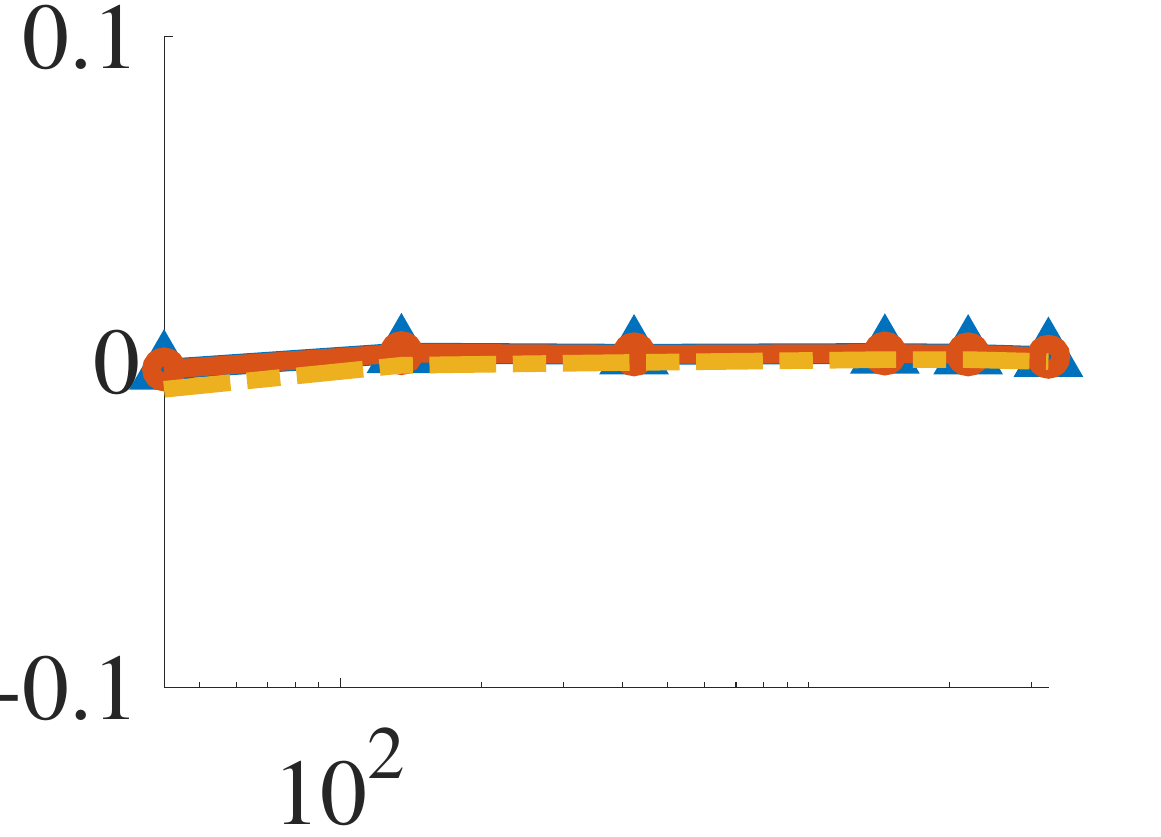} 
\includegraphics[width=.225\linewidth, trim={30 190 25 200}, clip]{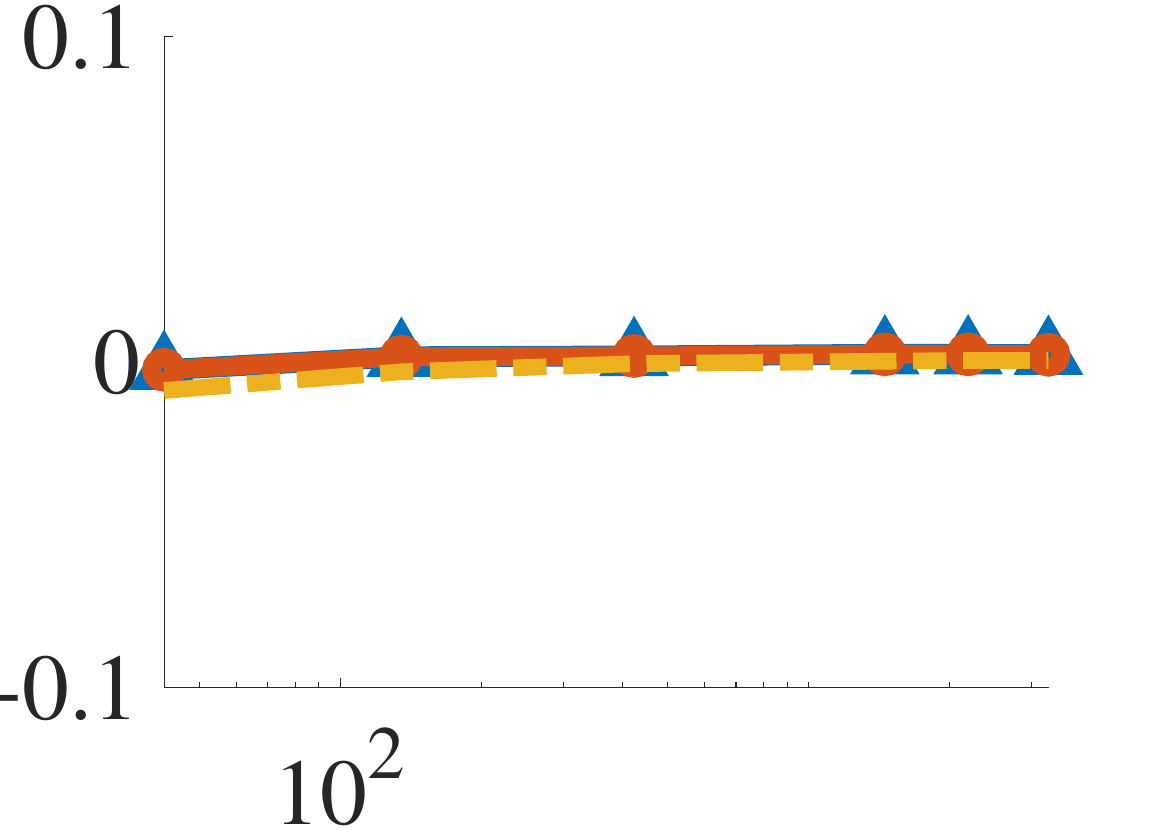} 
\includegraphics[width=.225\linewidth, trim={30 190 25 200}, clip]{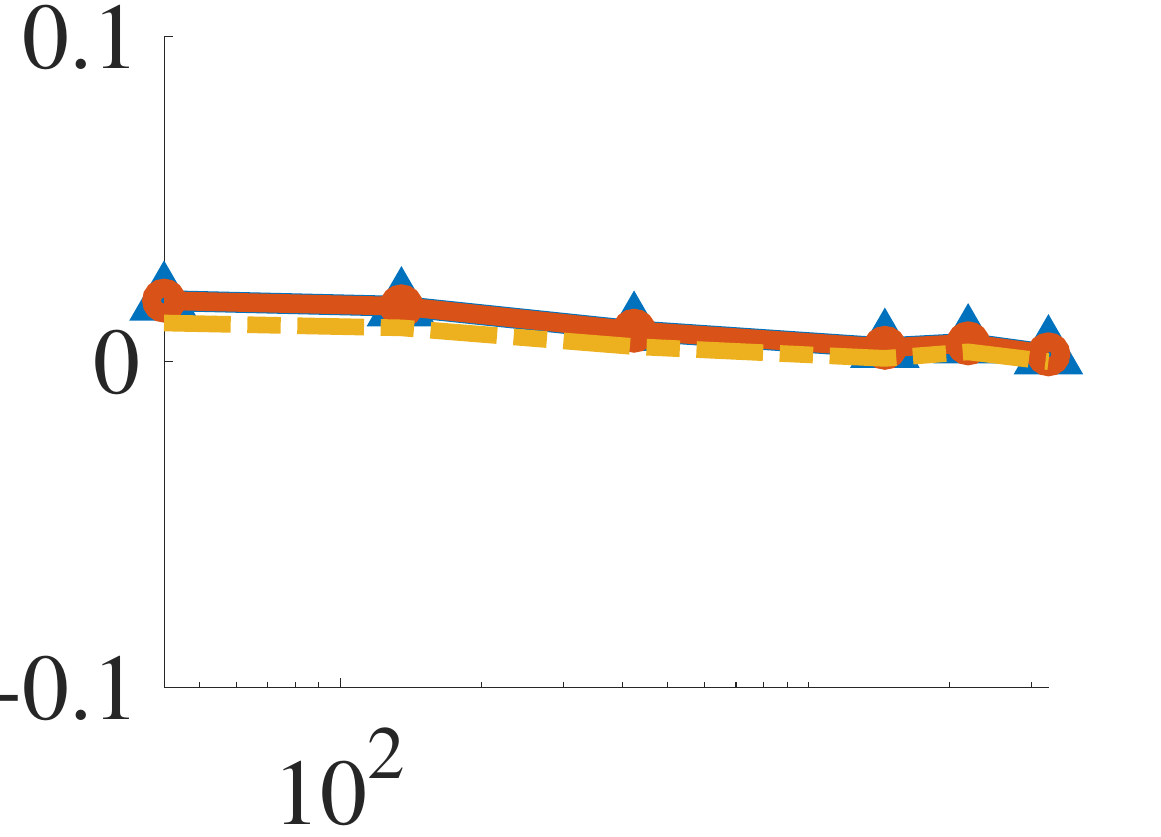} 
\includegraphics[width=.225\linewidth, trim={30 190 25 200}, clip]{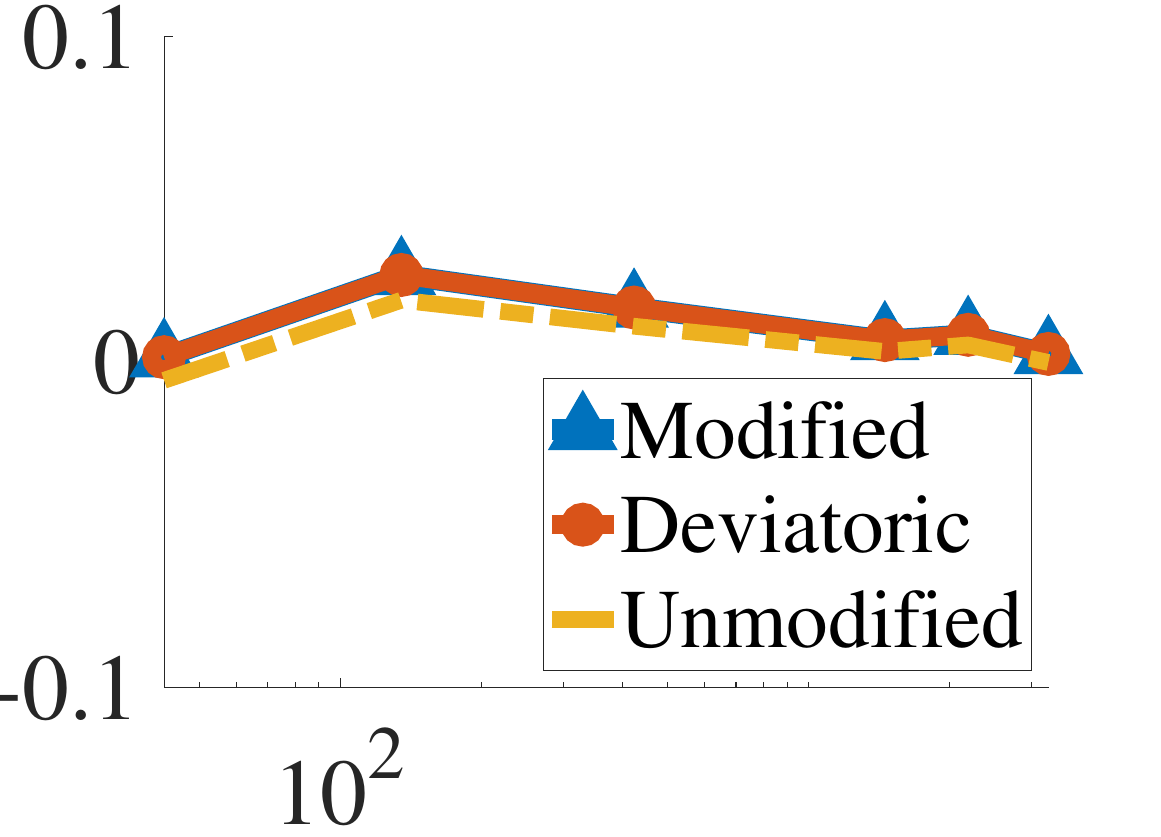} \\
\rotatebox{90}{$\qquad$ \textbf{$\nus = .4$} }
   \rotatebox{90}{$\quad$ Area Change \% }
\includegraphics[width=.225\linewidth, trim={30 190 25 200}, clip]{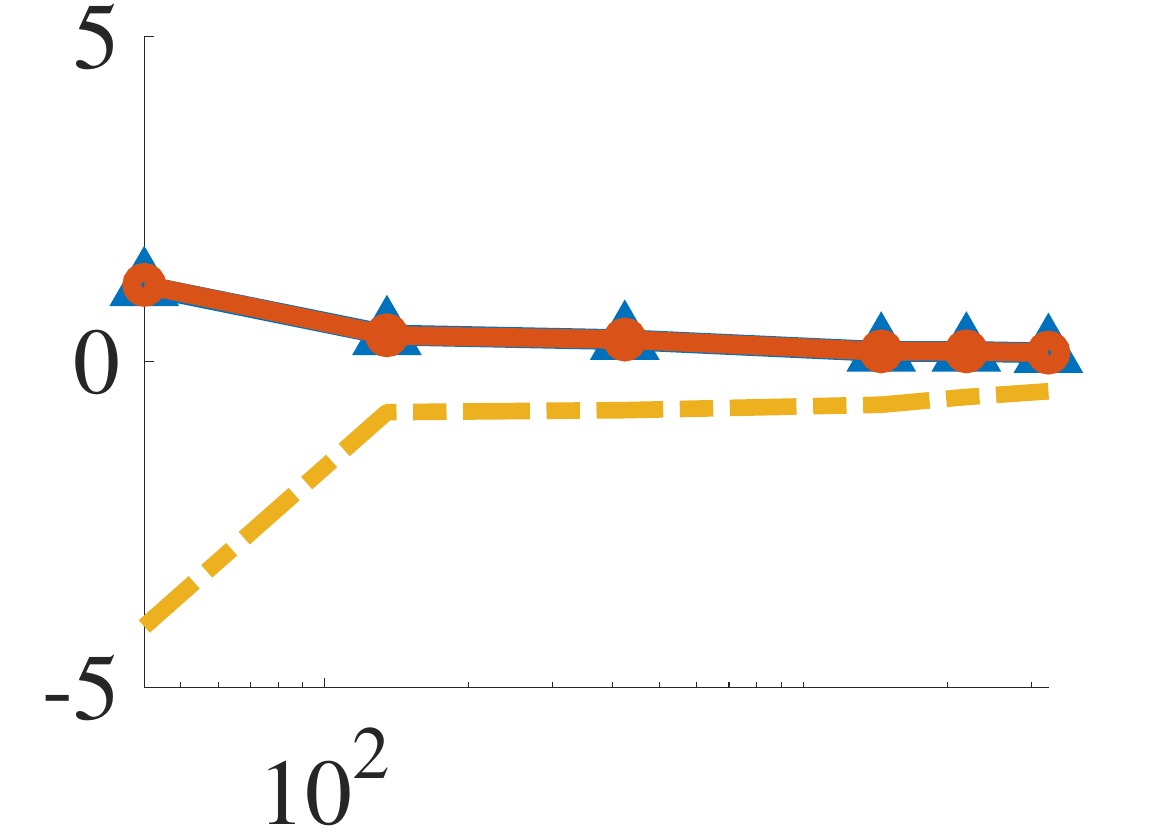} 
\includegraphics[width=.225\linewidth, trim={30 190 25 200}, clip]{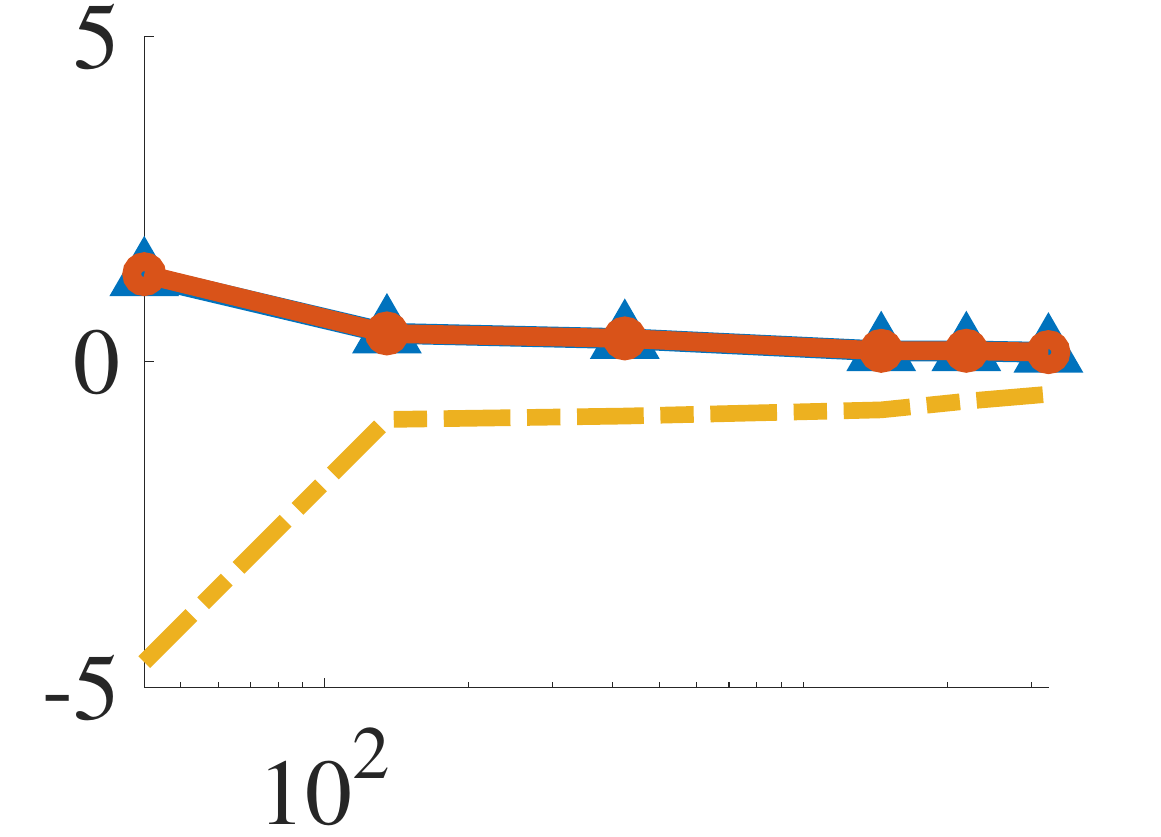}
\includegraphics[width=.225\linewidth, trim={30 190 25 200}, clip]{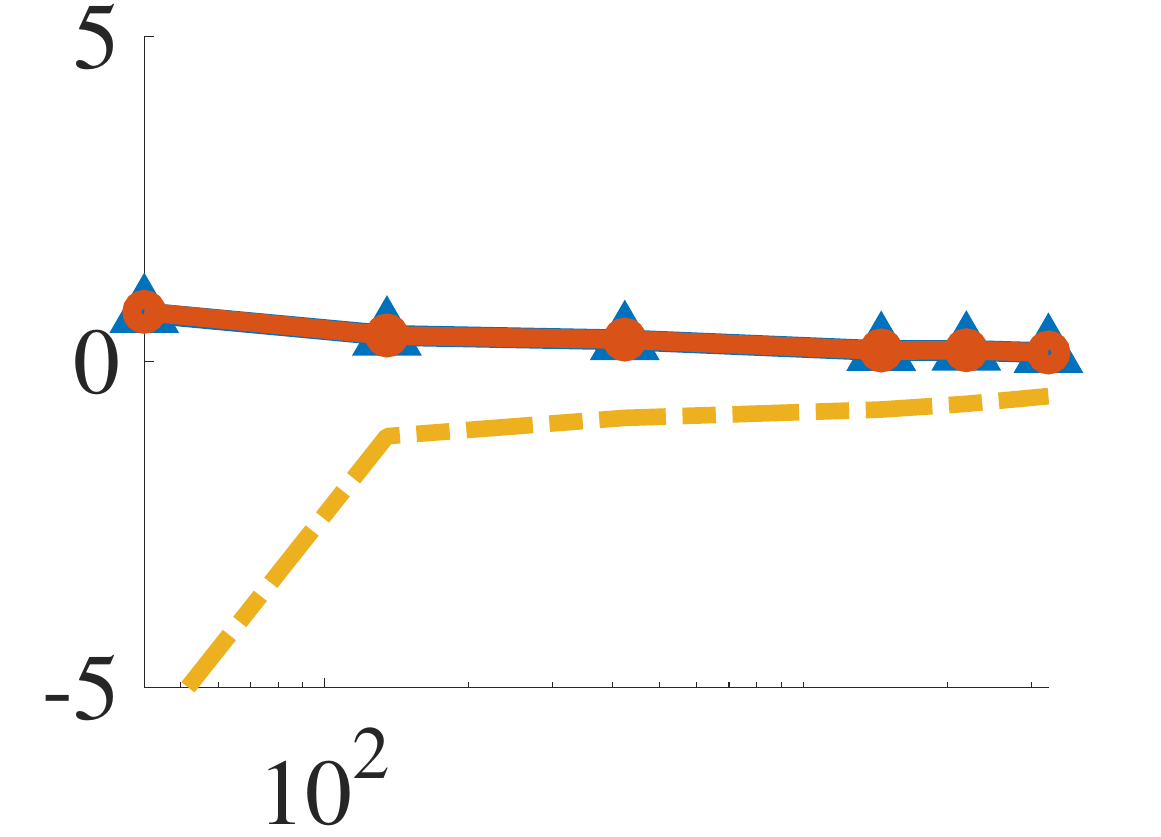} 
\includegraphics[width=.225\linewidth, trim={30 190 25 200}, clip]{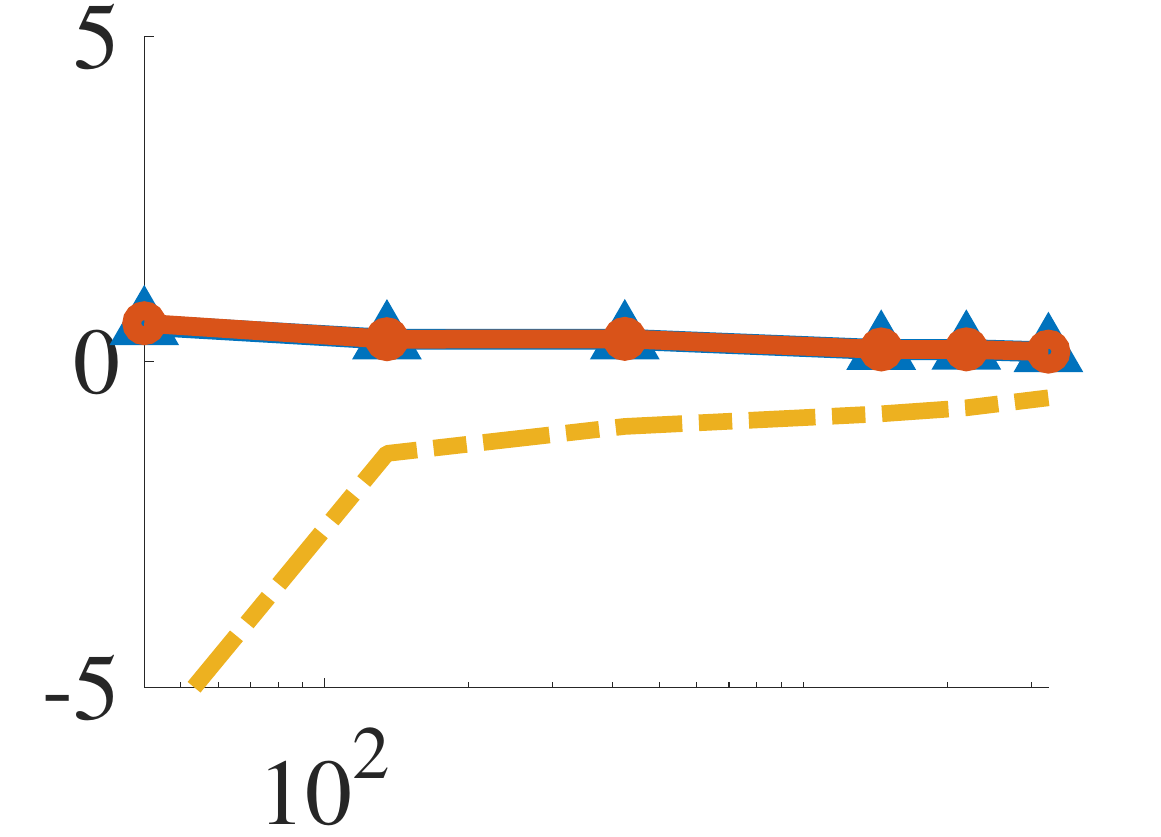}\\

\rotatebox{90}{$\qquad\;$ \textbf{$\nus = 0$} }
   \rotatebox{90}{$\quad$ Area Change \% }
\includegraphics[width=.225\linewidth, trim={30 190 25 200}, clip]{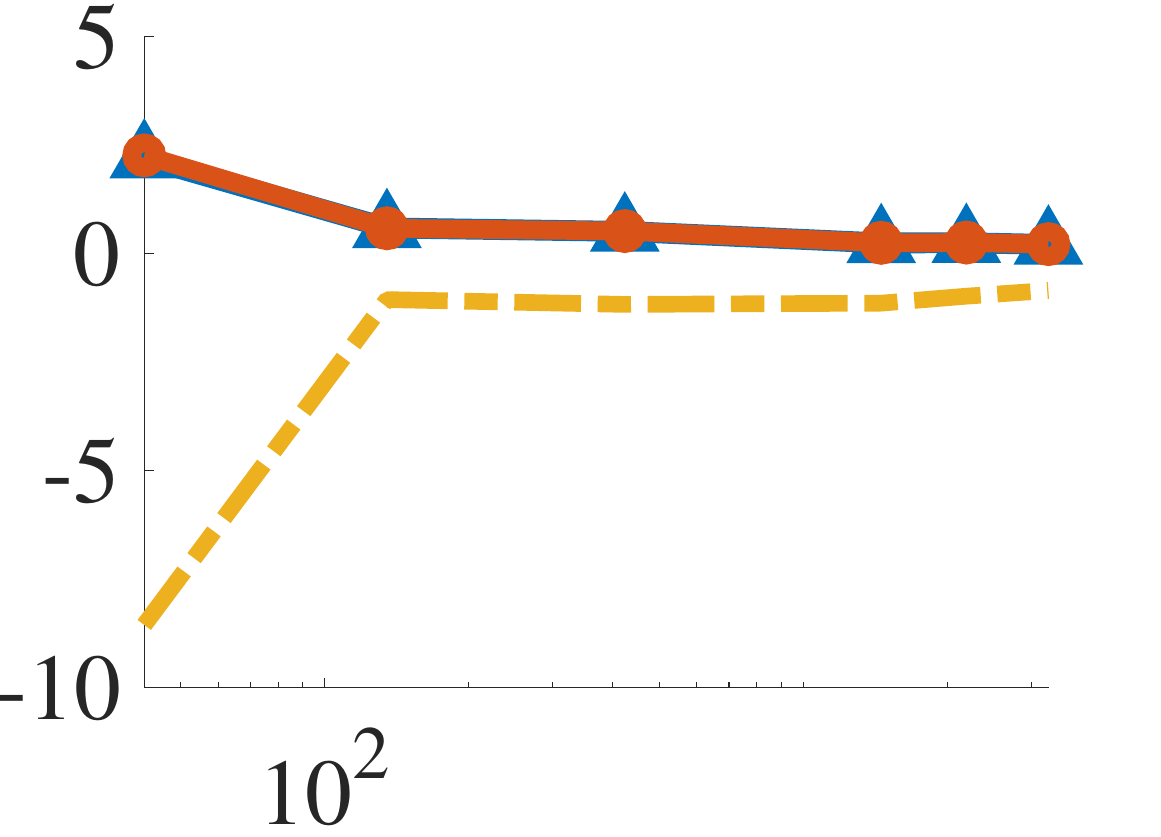} 
\includegraphics[width=.225\linewidth, trim={30 190 25 200}, clip]{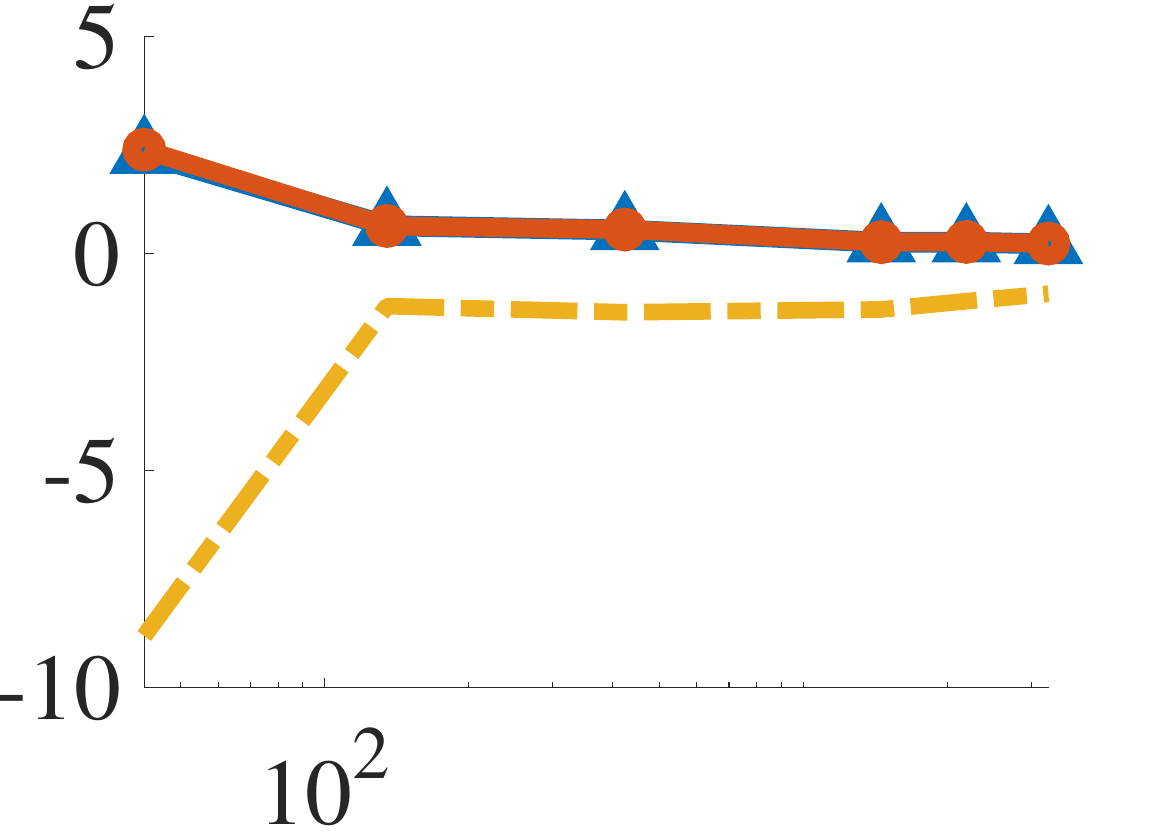} 
\includegraphics[width=.225\linewidth, trim={30 190 25 200}, clip]{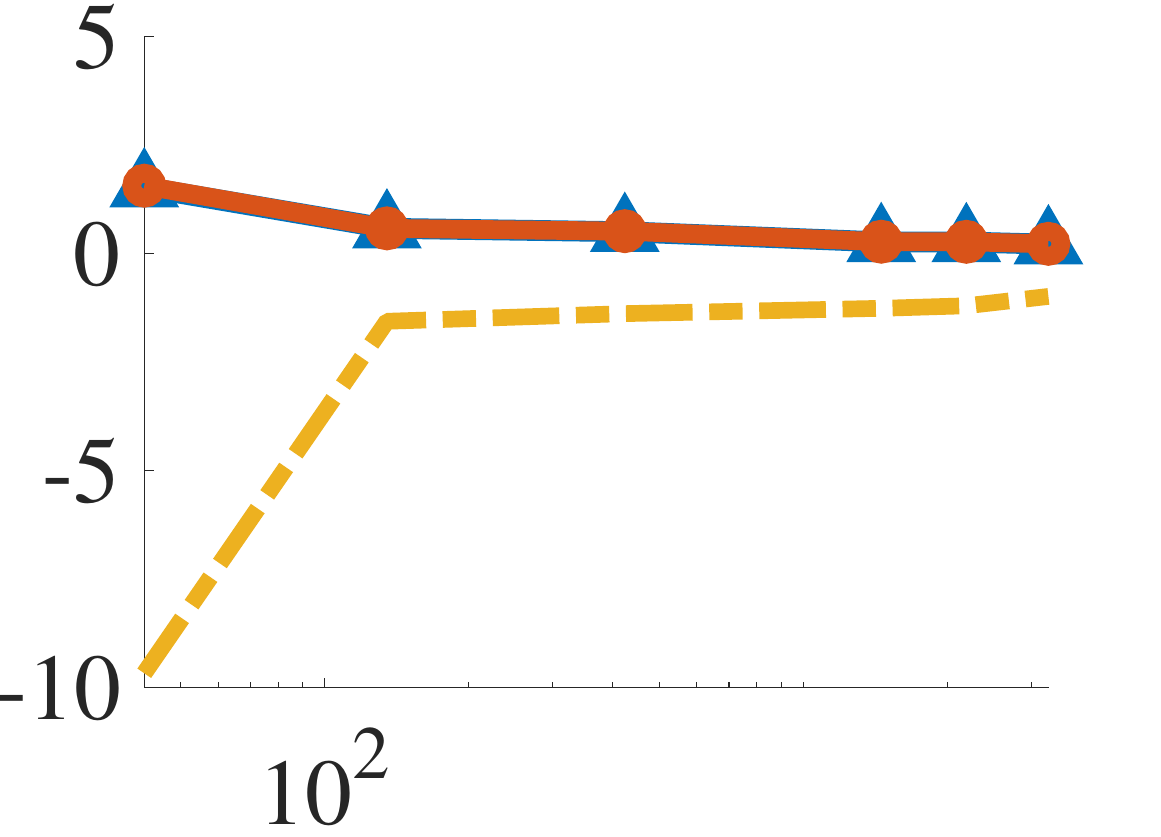} 
\includegraphics[width=.225\linewidth, trim={30 190 25 200}, clip]{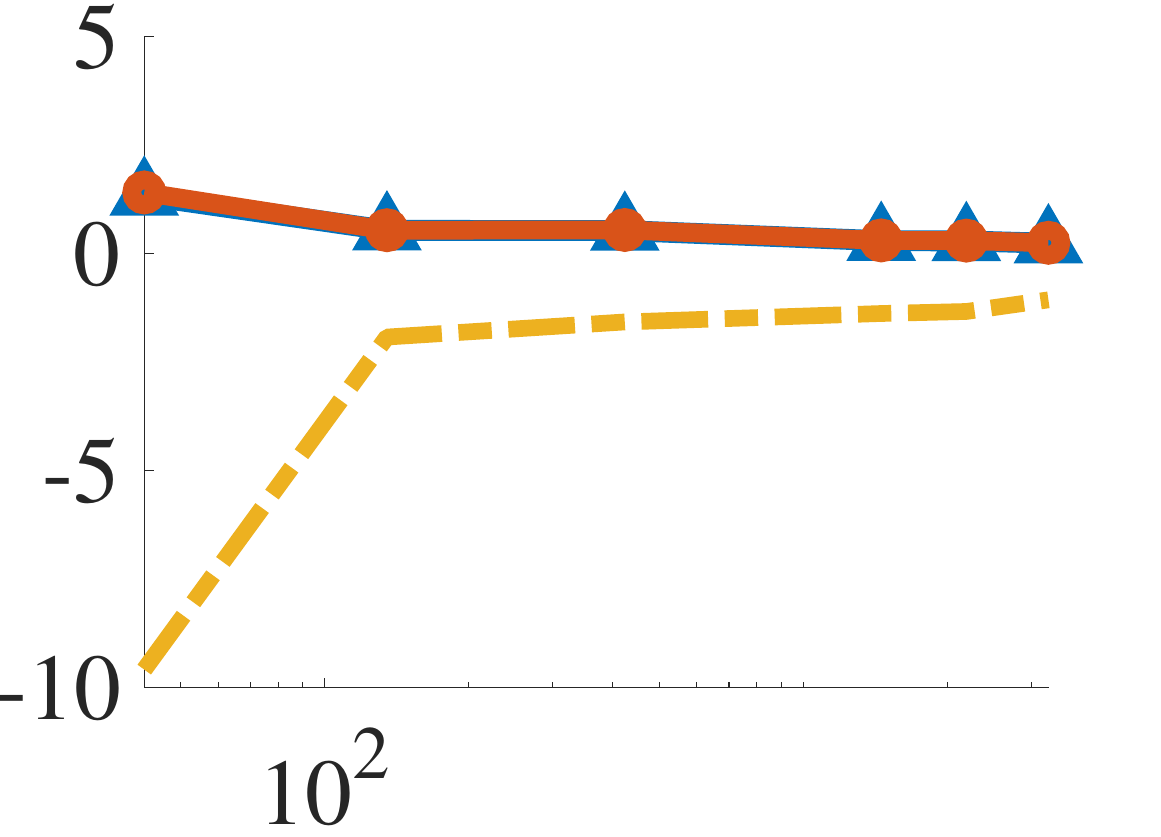}\\

\rotatebox{90}{$\qquad$ \textbf{$\nus = -1$} }
   \rotatebox{90}{$\quad$ Area Change \% }
\includegraphics[width=.225\linewidth, trim={30 190 25 200}, clip]{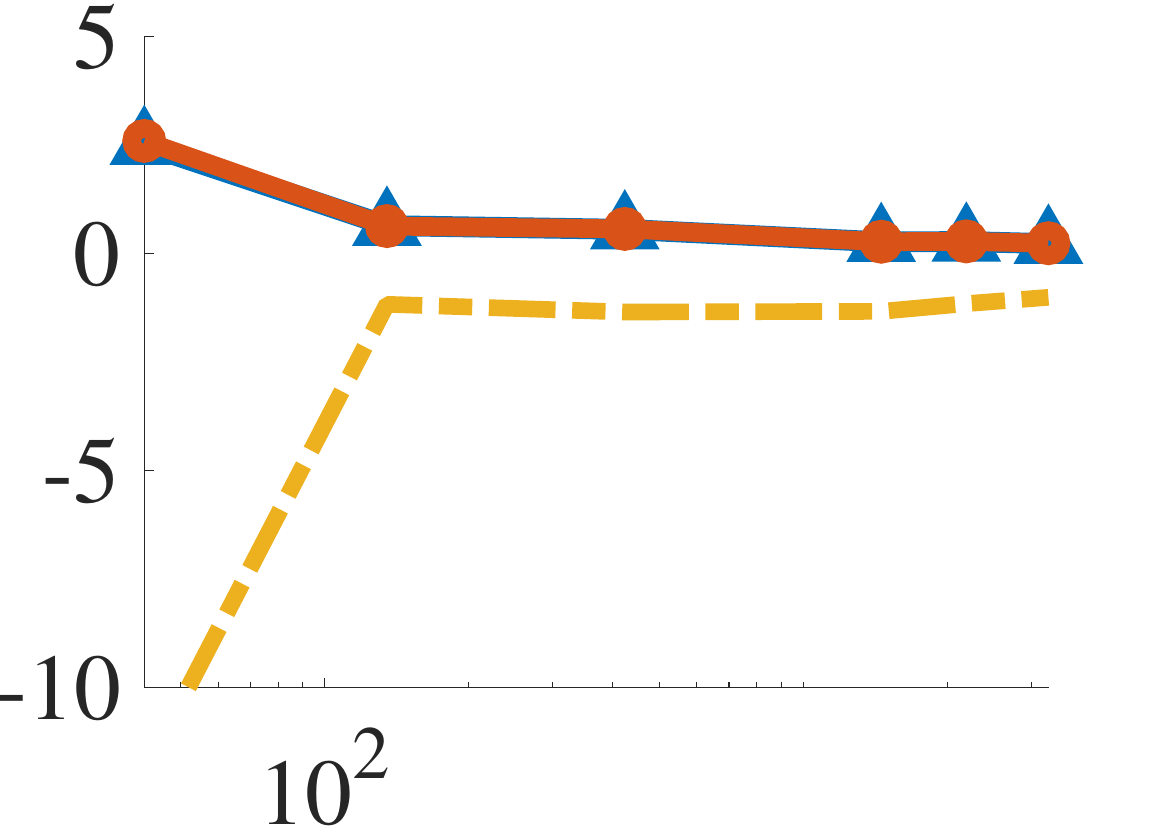} 
\includegraphics[width=.225\linewidth, trim={30 190 25 200}, clip]{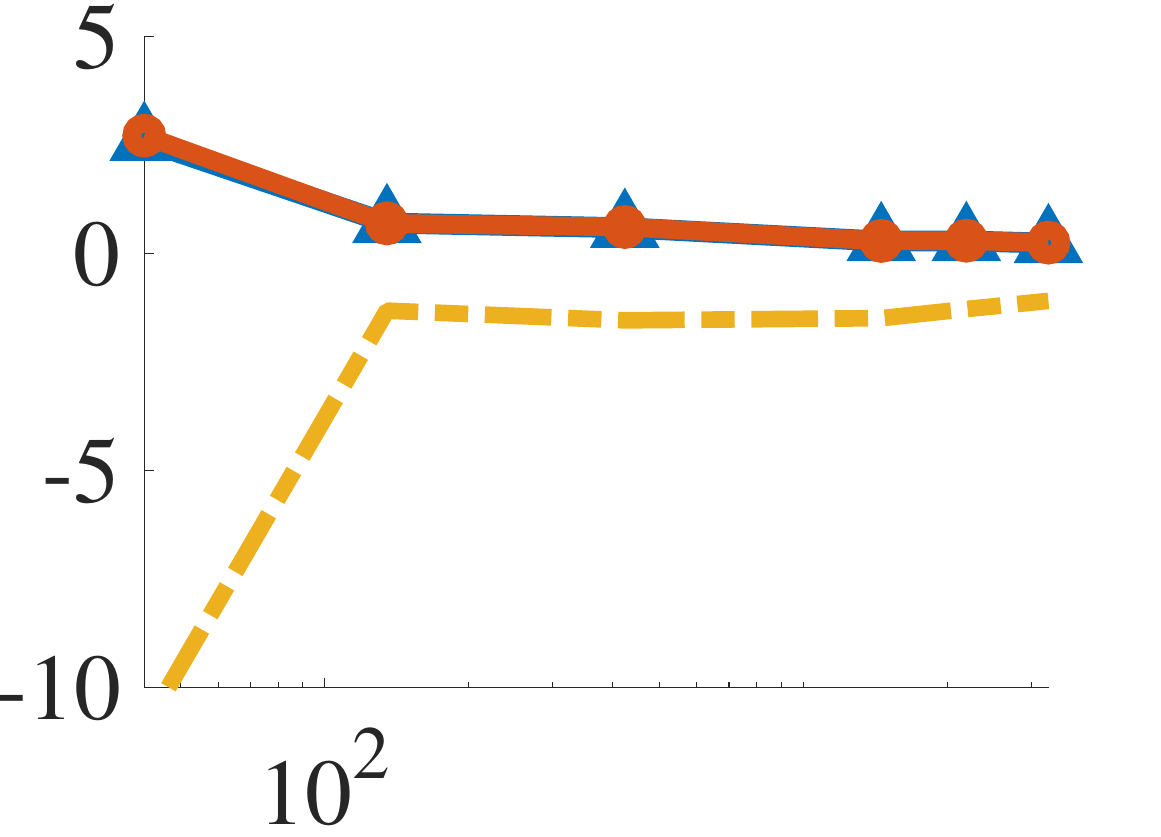} 
\includegraphics[width=.225\linewidth, trim={30 190 25 200}, clip]{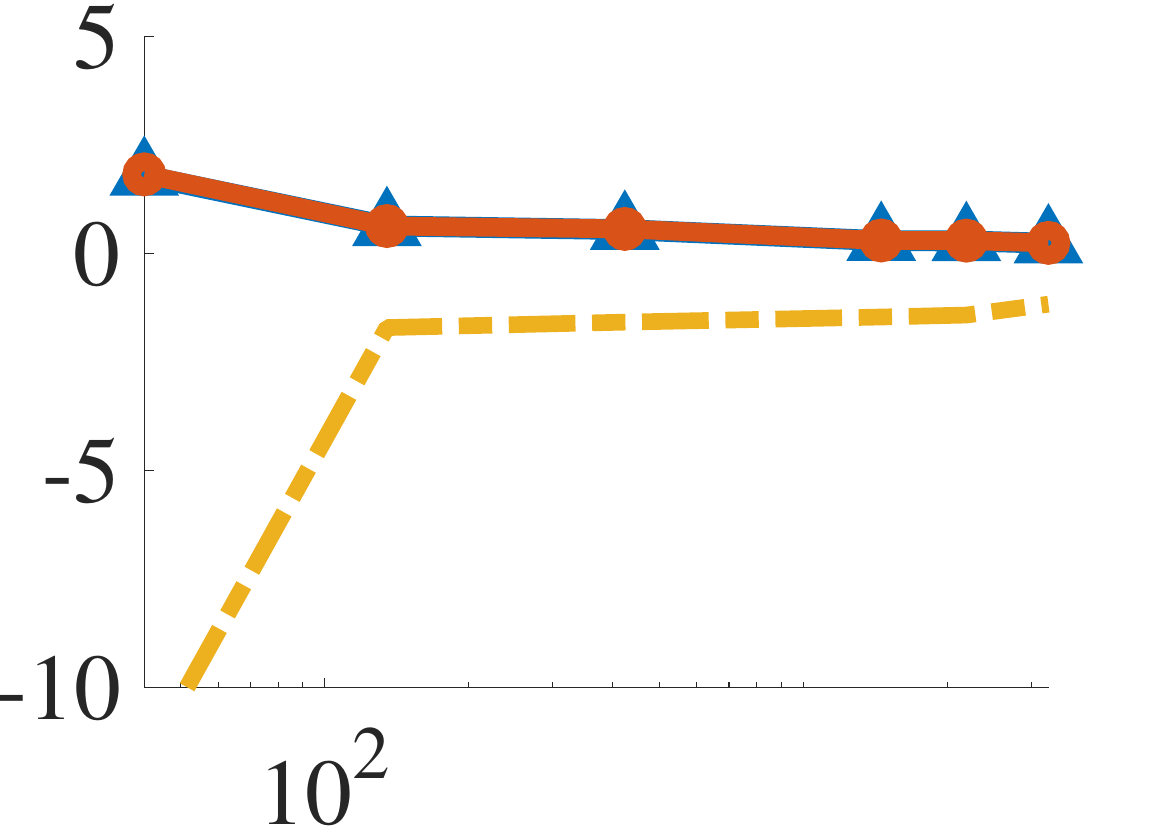} 
\includegraphics[width=.225\linewidth, trim={30 190 25 200}, clip]{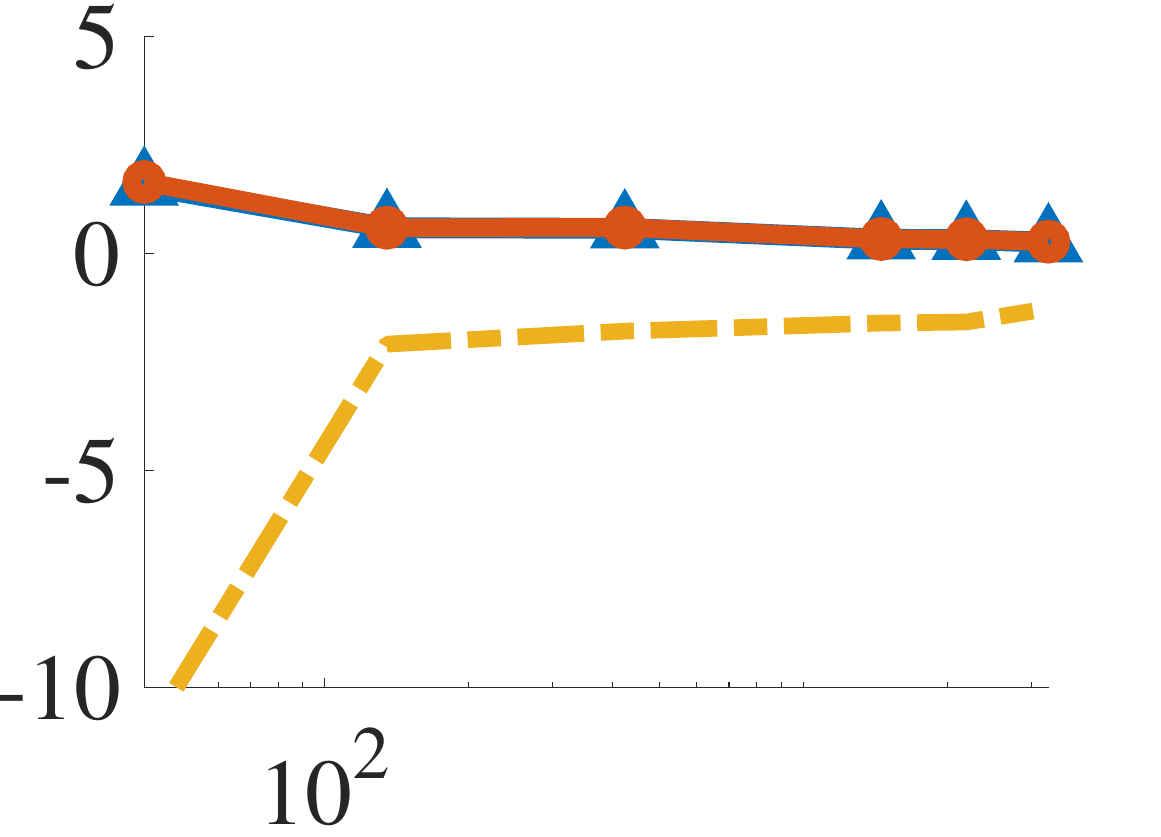}\\

$\qquad\qquad\quad$ \# Solid DOF $\qquad\qquad\quad\;$ \# Solid DOF $\qquad\qquad\quad$ \# Solid DOF $\qquad\qquad\quad\;$ \# Solid DOF
\caption{Percent change in total area for different numbers of solid DOF for the steady state elastic band benchmark (Section \ref{Elastic Band}) after deformation. The DOF range from $m = 42$ to $3255$, and the $x$ axis is on a log scale. Omitting the coarsest discretizations ($m=42$), the largest deviations in total volume among all element types used are approximately $0.69\%$ for the modified case, $2.1\%$ for the unmodified case, and $0.69\%$ for the deviatoric case.}
\label{elastic_area}
\end{figure}

\begin{figure}
\centering
\includegraphics[width=.7\linewidth, trim={10 120 10 120}, clip]{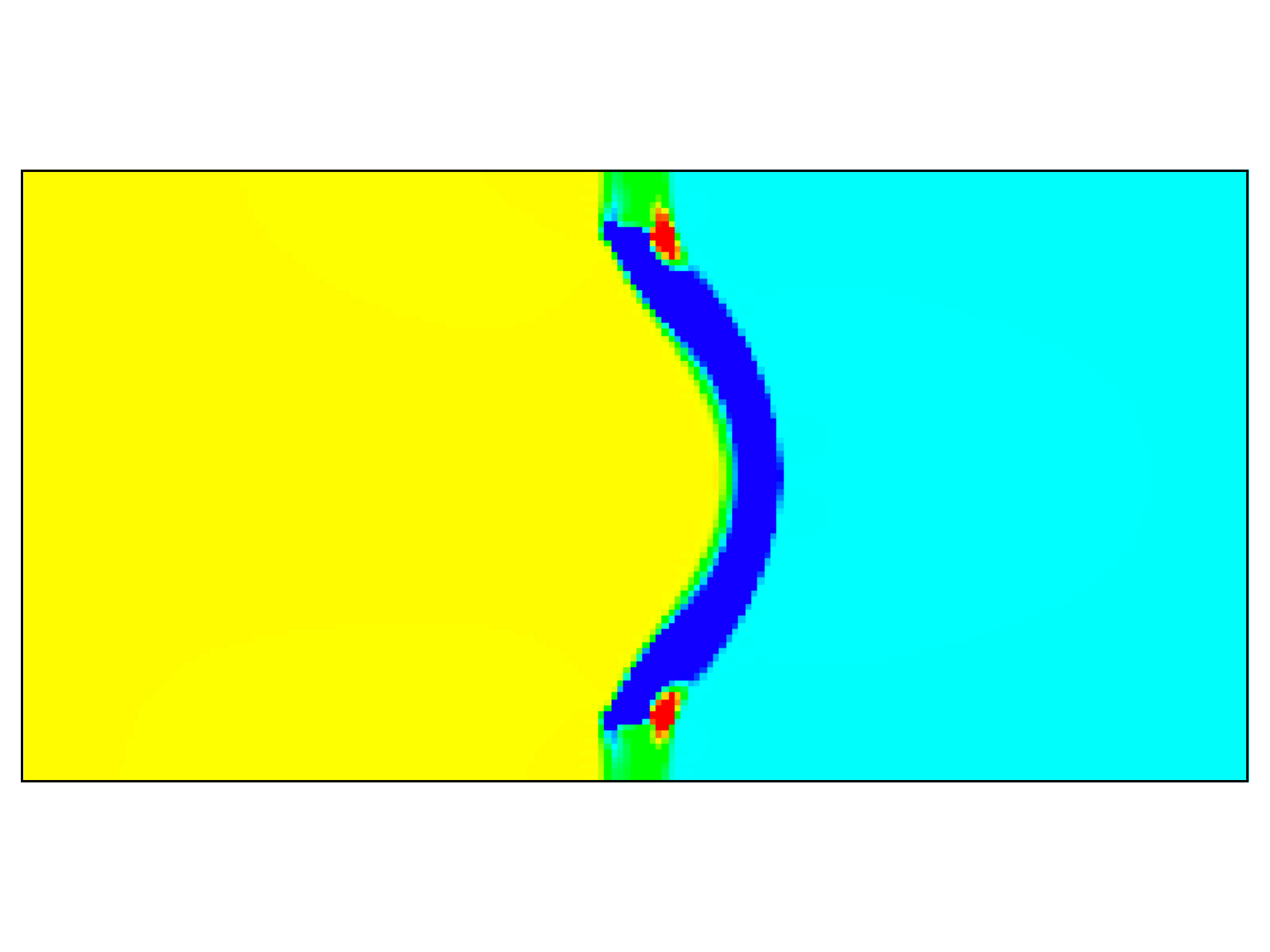}
\vskip .1cm
\begin{centering}
$\peul(\xb,t) \ \left(\frac{\text{dyn}}{\text{cm}^2}\right)$ \\
\includegraphics[width=2.5in, trim={0 5in 0 5in}, clip]{color_bar.pdf}  \\
%0.00 \ \ \ \ \ \ \ \ \ \ \ \ \ \ \ \ 2.00
-20 $\qquad\qquad\qquad\qquad$ 20\\
\end{centering}
\caption{Eulerian Lagrange multiplier field $\peul(\xb,t)$ of the steady state elastic band benchmark (Section \ref{Elastic Band}). The case shown here uses \textbf{Q1} elements, an Eulerian grid of $192 \times 96$, modified invariants, and $\nus = 0.4$. The Lagrangian mesh is not shown. At steady state, the Lagrange multiplier field is constant in the two regions separated by the elastic body. The region on the left has pressure $\peul(\xb,t) = 10 \  \frac{\text{dyn}}{\text{cm}^2}$ and the region on the right has pressure $\peul(\xb,t) = -10 \ \frac{\text{dyn}}{\text{cm}^2}$ at the end of the computation ($t = 15$ s). Within the region occupied by the solid, $\peul$ is non-constant, as shown in the figure. Note that $\peul$ corresponds to the physical pressure only in $\fluiddom$. In the solid region $\soliddom$, the physical pressure includes an additional contribution from the solid model, $p = \peul + \pstab$. The contribution of $\pstab$ is not depicted in this figure.}
\label{eb_pressure}
\end{figure}

\begin{figure}
%\begin{tabular}{c c c}
%\subcaptionbox{\label{sfig:testa}} {\includegraphics[width=.3\linewidth, trim={100 150 50 150}, clip]{eb_flow_t=02.pdf}}&
%\subcaptionbox{\label{sfig:testb}} {\includegraphics[width=.3\linewidth, trim={100 150 50 150}, clip]{eb_flow_t=1.pdf}} &
%\subcaptionbox{\label{sfig:testb}} {\includegraphics[width=.3\linewidth, trim={100 150 50 150}, clip]{eb_flow_t=10.pdf}}\\
%
%\end{tabular}
%\begin{centering}
%\subcaptionbox{\label{sfig:testa}} {\includegraphics[width=.6\linewidth, trim={10 50 10 20}, clip]{eb_flow_t=02.pdf}}\\
%\subcaptionbox{\label{sfig:testb}} {\includegraphics[width=.6\linewidth, trim={10 50 10 0}, clip]{eb_flow_t=1.pdf}}\\
%\subcaptionbox{\label{sfig:testb}} {\includegraphics[width=.6\linewidth, trim={10 50 10 20}, clip]{eb_flow_t=10.pdf}}\\
%\end{centering}
\begin{tabular}{l r}
\subcaptionbox{\label{sfig:testa}} {\includegraphics[width=.45\linewidth, trim={100 100 100 100}, clip]{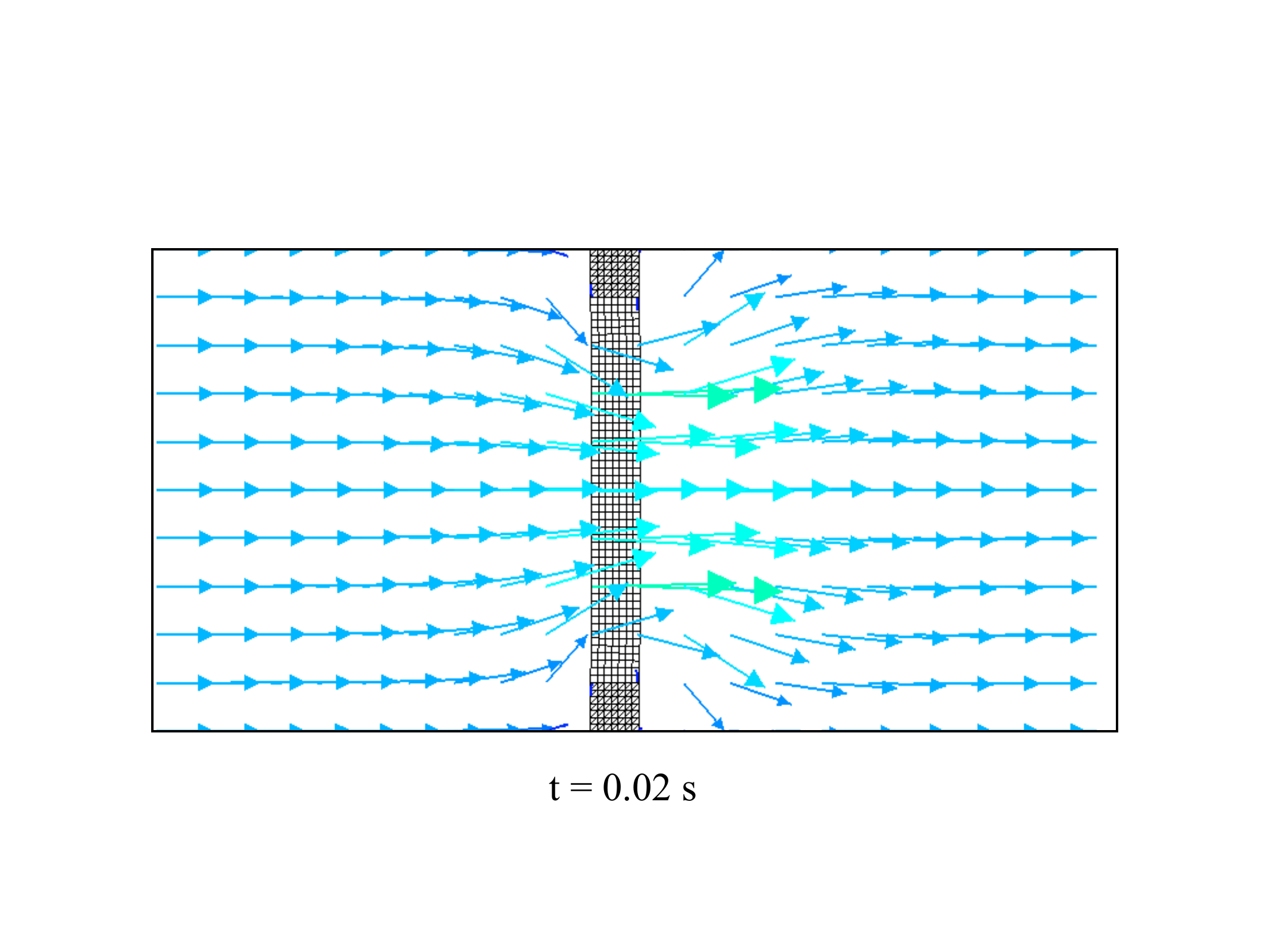}}&
\subcaptionbox{\label{sfig:testb}} {\includegraphics[width=.45\linewidth, trim={100 100 100 100}, clip]{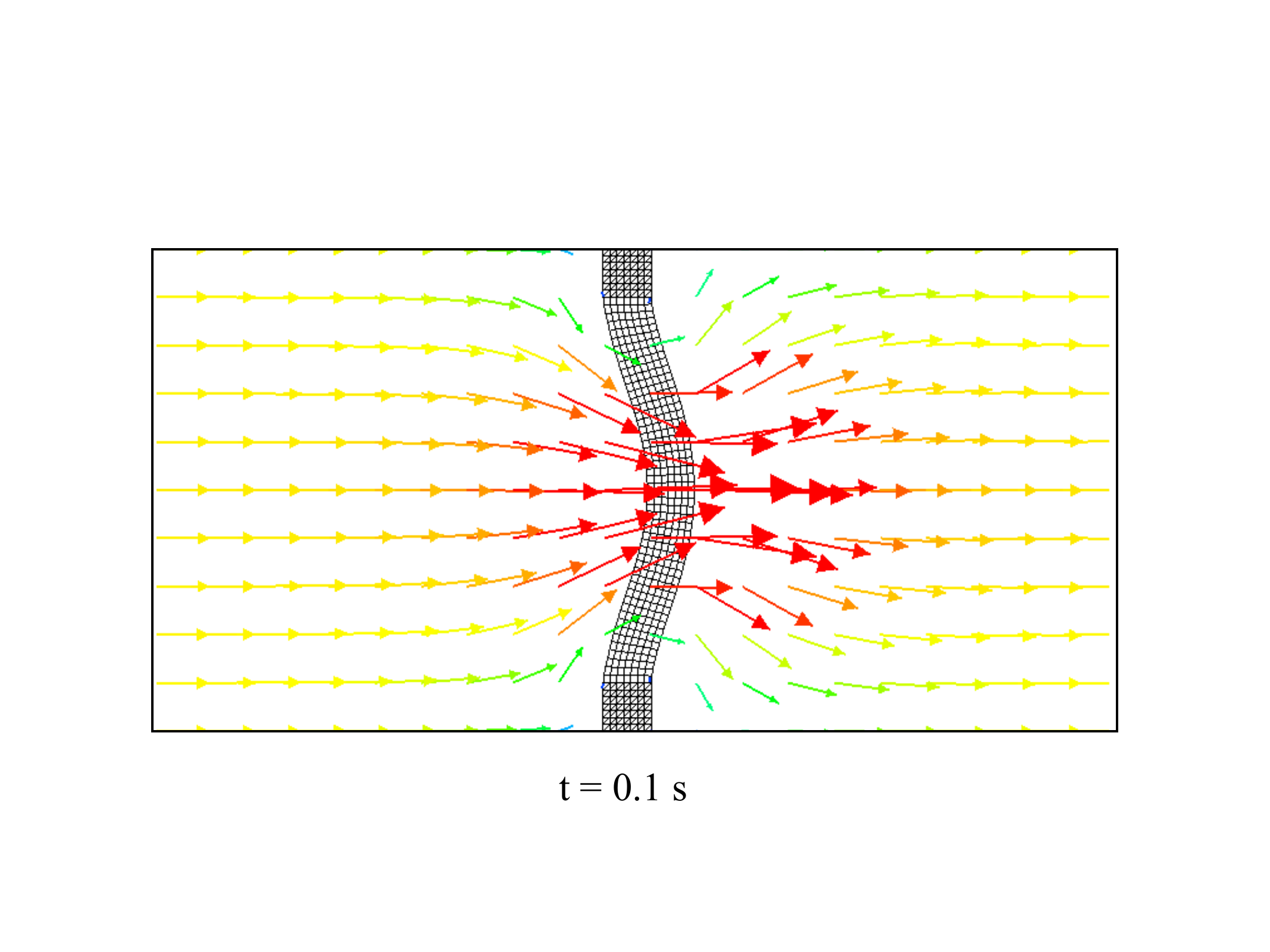}}\\

\subcaptionbox{\label{sfig:testb}} {\includegraphics[width=.45\linewidth, trim={100 100 100 100}, clip]{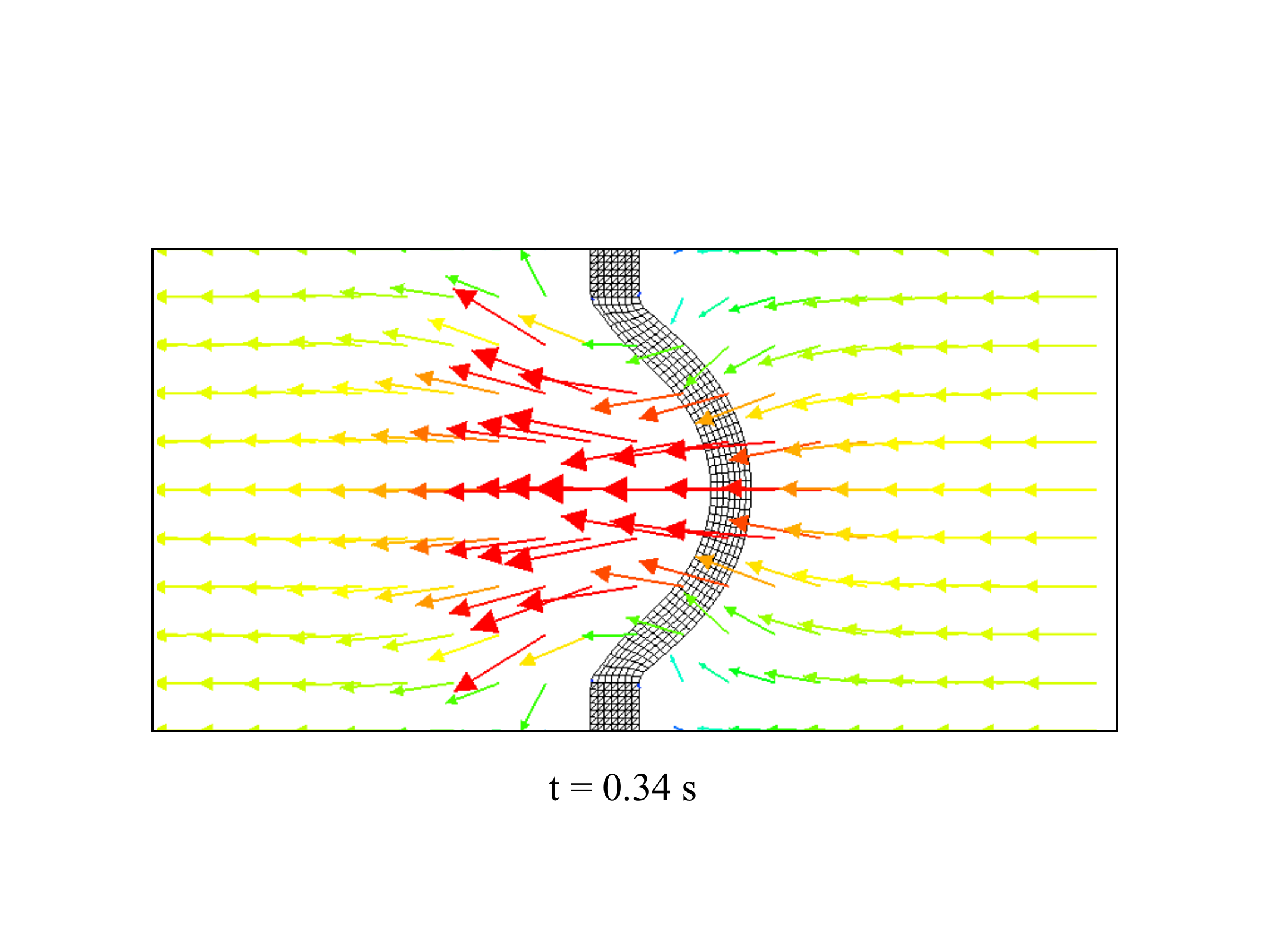}}&
\subcaptionbox{\label{sfig:testb}} {\includegraphics[width=.45\linewidth, trim={100 100 100 100}, clip]{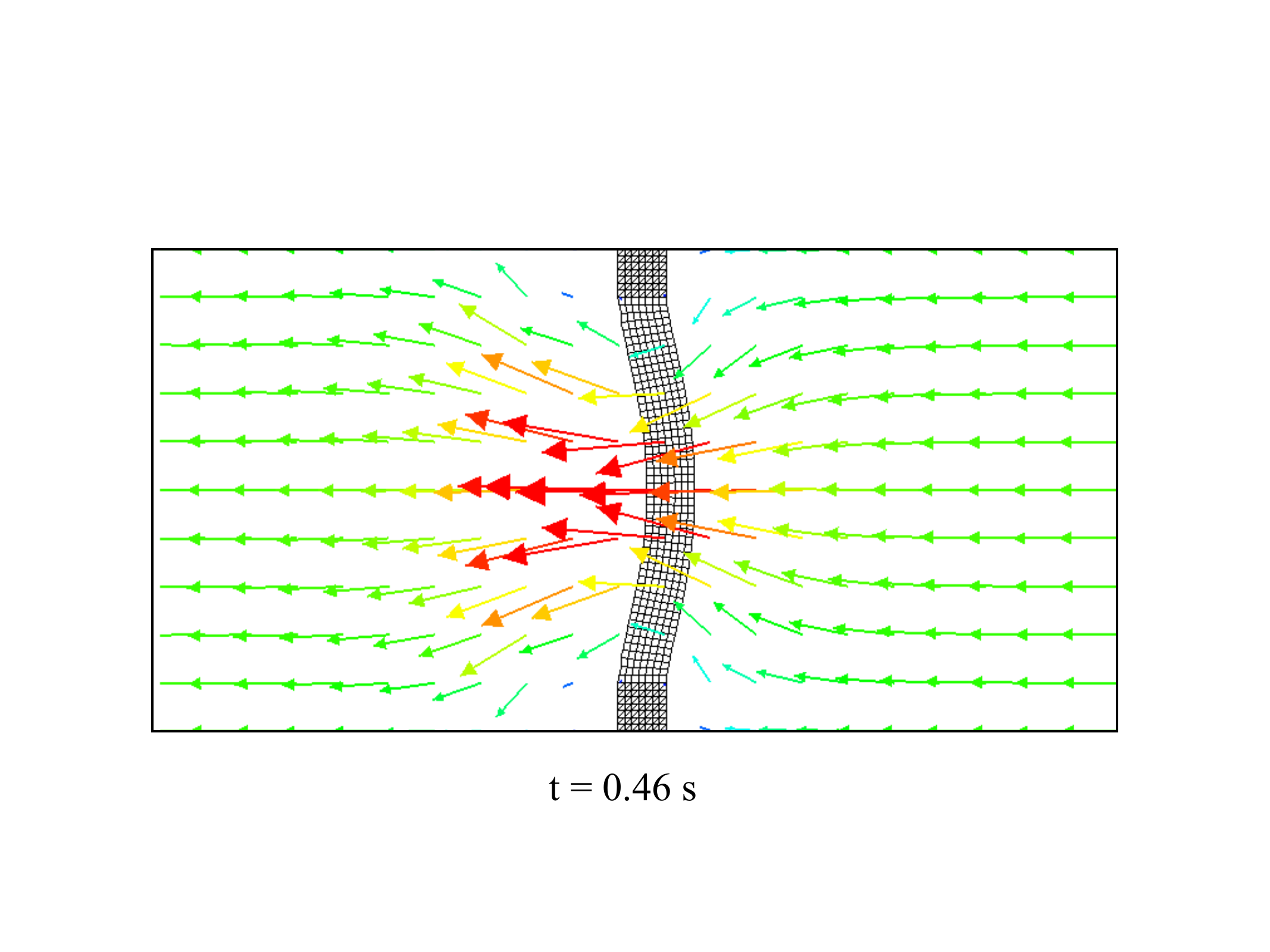}}\\
\end{tabular}

%trim={left bottom right top}
\begin{centering}
\vskip .1cm
$||\ub(\xb,t)|| \ \left(\frac{\text{cm}}{\text{s}}\right)$ \\
\includegraphics[width=2.5in, trim={0 5in 0 5in}, clip]{color_bar.pdf}  \\
%0.00 \ \ \ \ \ \ \ \ \ \ \ \ \ \ \ \ 2.00
0.0 $\qquad\qquad\qquad\qquad$ 1.0

\end{centering}
\caption{Eulerian velocity field and deformations of the fully dynamic version of the elastic band benchmark (Section \ref{Elastic Band}) at $t = 0.02$ s (a), $t = 0.1$ s (b), $t = 0.34$ s (c), and $t = 0.46$ s (d). Each of these time slices corresponds to a different characteristic deformation: (a) near initial configuration; (b) early deformation; (c) largest deformation; and (d) right before the structure enters another period of oscillation. Note that the Eulerian velocity field $\ub(\xb,t)$ corresponds to the velocity of whichever material is located at position $\xb$; it describes the velocity of the fluid as well as the structure. The color corresponds to the magnitude of the velocity at each spatial point $\xb$.}
\label{eb_flow}
\end{figure}

\begin{figure}
%%$\qquad\qquad\qquad\;\;\;\;$ \textbf{N = 16} $\qquad\qquad\qquad\qquad\quad$  \textbf{N = 32} $\qquad\qquad\qquad\qquad\;\;\;$  \textbf{N = 64} \\

$\qquad\qquad\quad\quad\quad\quad$ \textbf{N = 16} $\qquad\qquad\quad\quad\quad\quad\quad\quad\quad$ \textbf{N = 64}$\qquad\qquad\quad\quad\quad\quad\quad\quad\quad$ \textbf{N = 96}\\

\rotatebox{90}{$\qquad\quad\quad\;$ \textbf{$\nus = .4$} }
   \rotatebox{90}{$\qquad\quad\quad$ Disp. (cm) }
\includegraphics[width=.3\linewidth, trim={30 190 25 200}, clip]{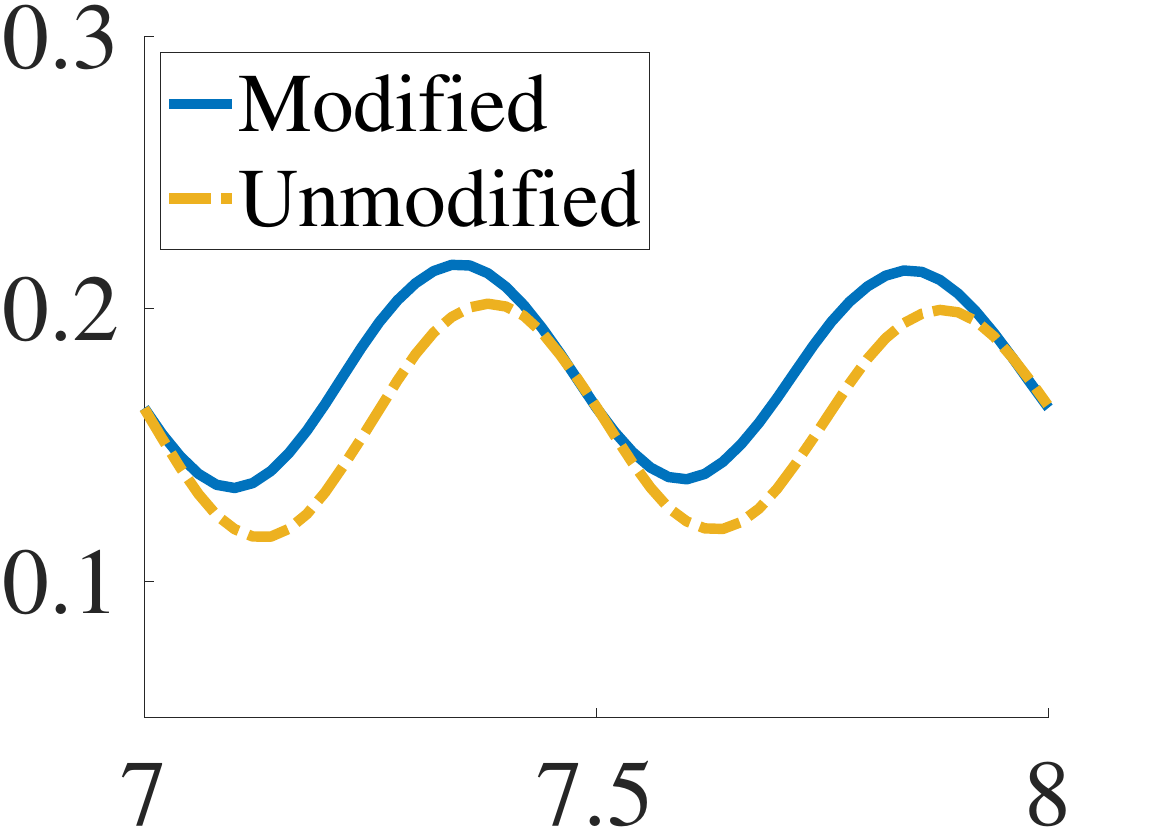}
\includegraphics[width=.3\linewidth, trim={30 190 25 200}, clip]{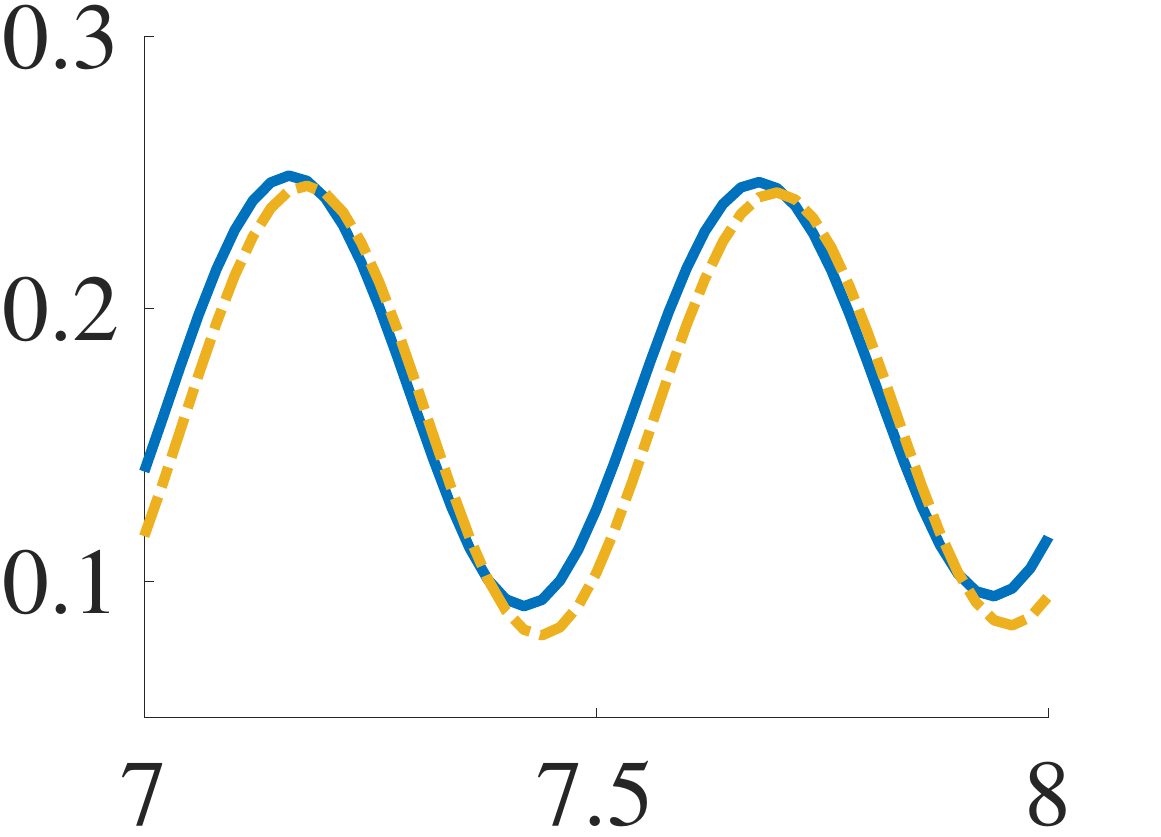}
\includegraphics[width=.3\linewidth, trim={30 190 25 200}, clip]{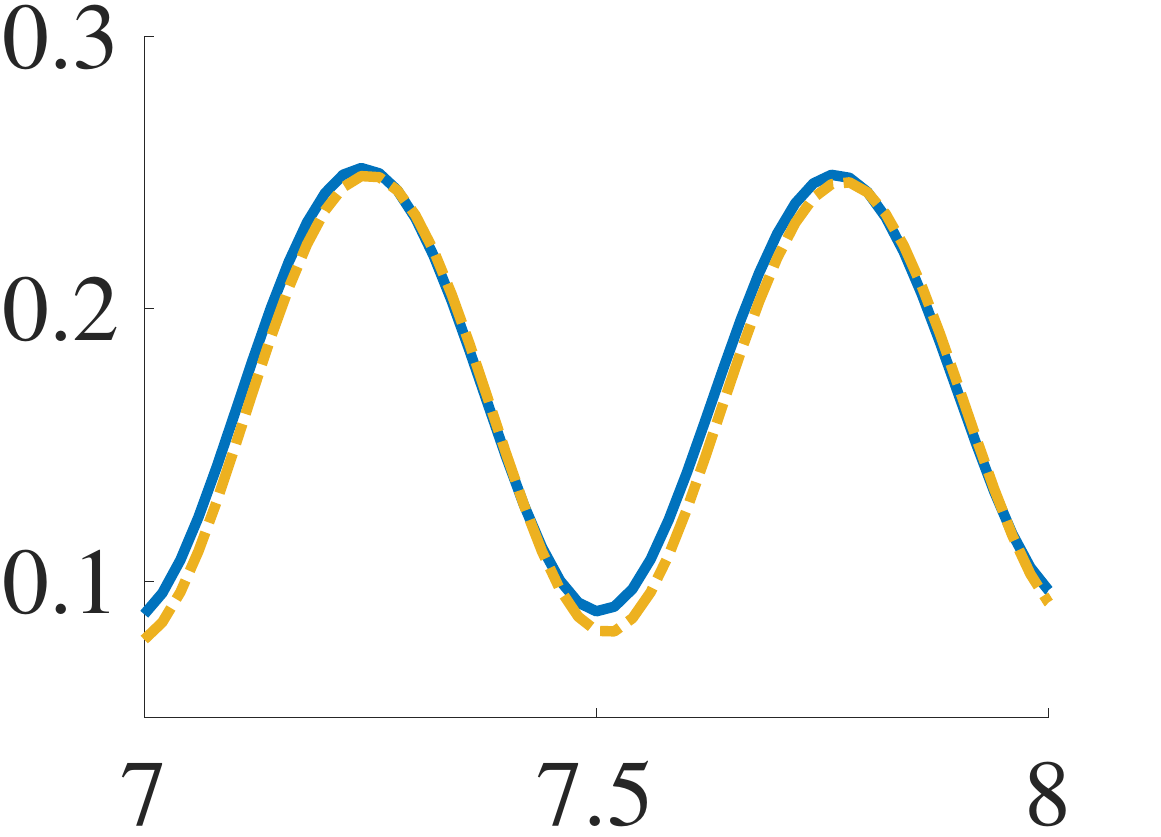}\\ 

\rotatebox{90}{$\qquad\quad\quad$ \textbf{$\nus = -1$} }
   \rotatebox{90}{$\qquad\quad\quad$ Disp. (cm)  }
\includegraphics[width=.3\linewidth, trim={30 190 25 200}, clip]{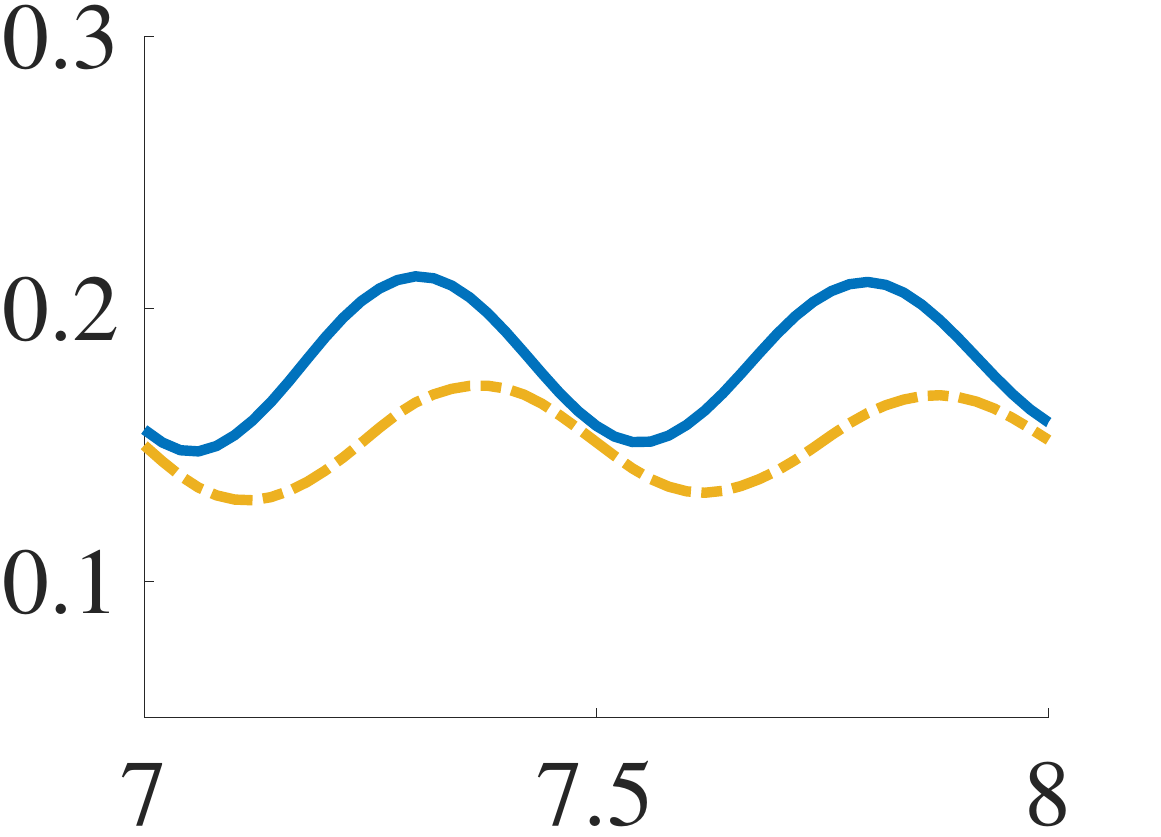} 
\includegraphics[width=.3\linewidth, trim={30 190 25 200}, clip]{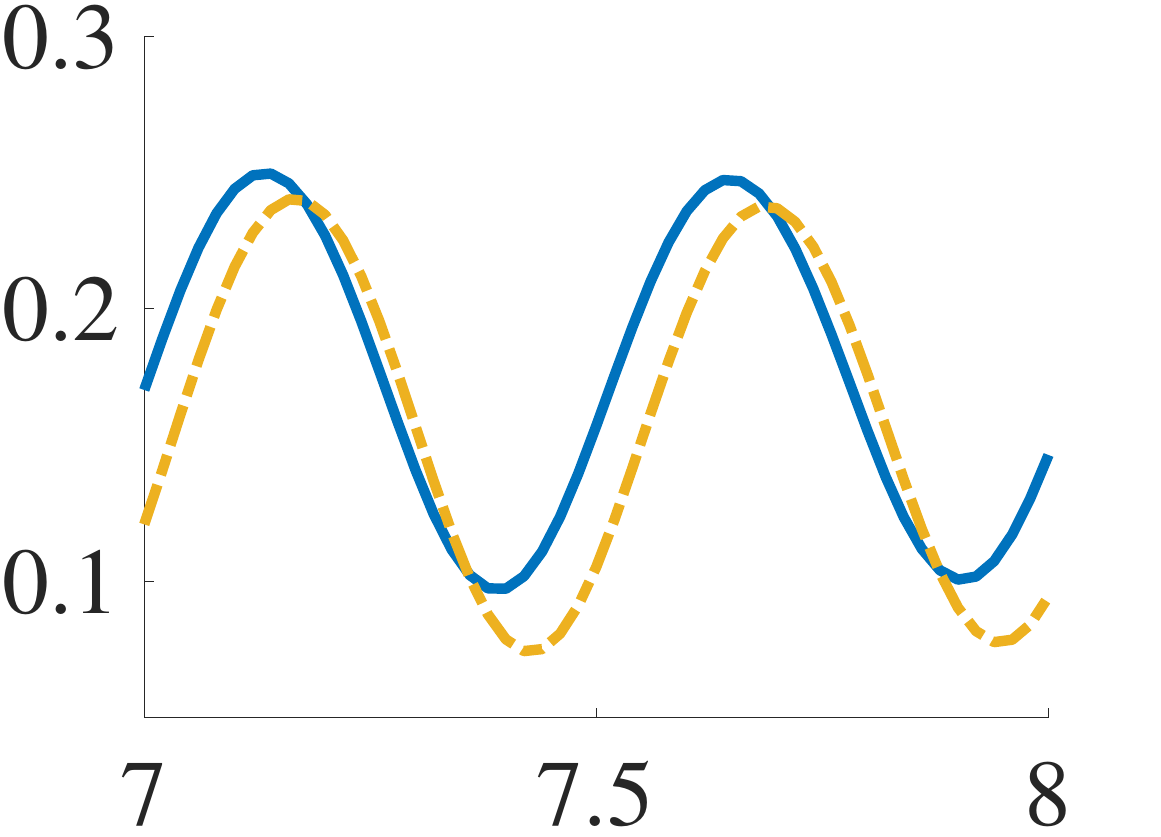}
\includegraphics[width=.3\linewidth, trim={30 190 25 200}, clip]{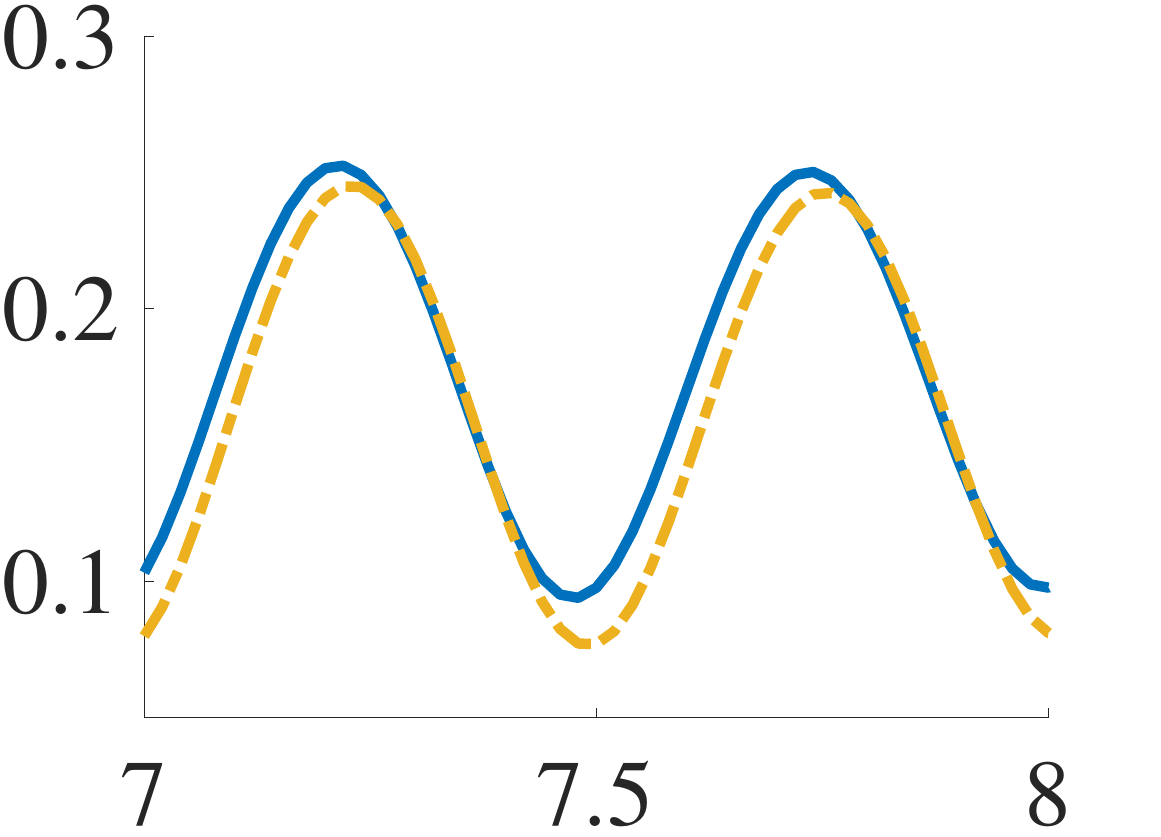}\\

$\qquad\qquad\quad\quad\quad\quad$ Time (s) $\qquad\qquad\quad\quad\quad\quad\quad\quad$ Time (s) $\qquad\qquad\quad\quad\quad\quad\quad\quad$ Time (s)
\caption{Transient behavior of the elastic band (Section \ref{Elastic Band}) for the dynamic version of the test for different discretizations between times $t = 7.0$ s and $t = 8.0$ s. Plot shows $x$-displacement of point of interest against time. $N$ describes an $2N$ by $N$ Eulerian grid. $N = 16$ corresponds to $m = 135$ solid DOF, $N = 64$ corresponds to $m = 1456$ solid DOF, and $N = 96$ corresponds to $m = 3255$ solid DOF. Results shown here are for \textbf{Q2} elements.}
\label{eb_v_time}
\end{figure}

\begin{figure}

$\qquad\qquad\quad\quad\quad\quad$ \textbf{N = 16} $\qquad\qquad\quad\quad\quad\quad\quad\quad\quad$ \textbf{N = 64}$\qquad\qquad\quad\quad\quad\quad\quad\quad\quad$ \textbf{N = 96}\\

\rotatebox{90}{$\qquad\quad\quad$ \textbf{$\nus = .4$} }
   \rotatebox{90}{$\qquad\quad$ Area Change \% }
\includegraphics[width=.3\linewidth, trim={30 190 25 200}, clip]{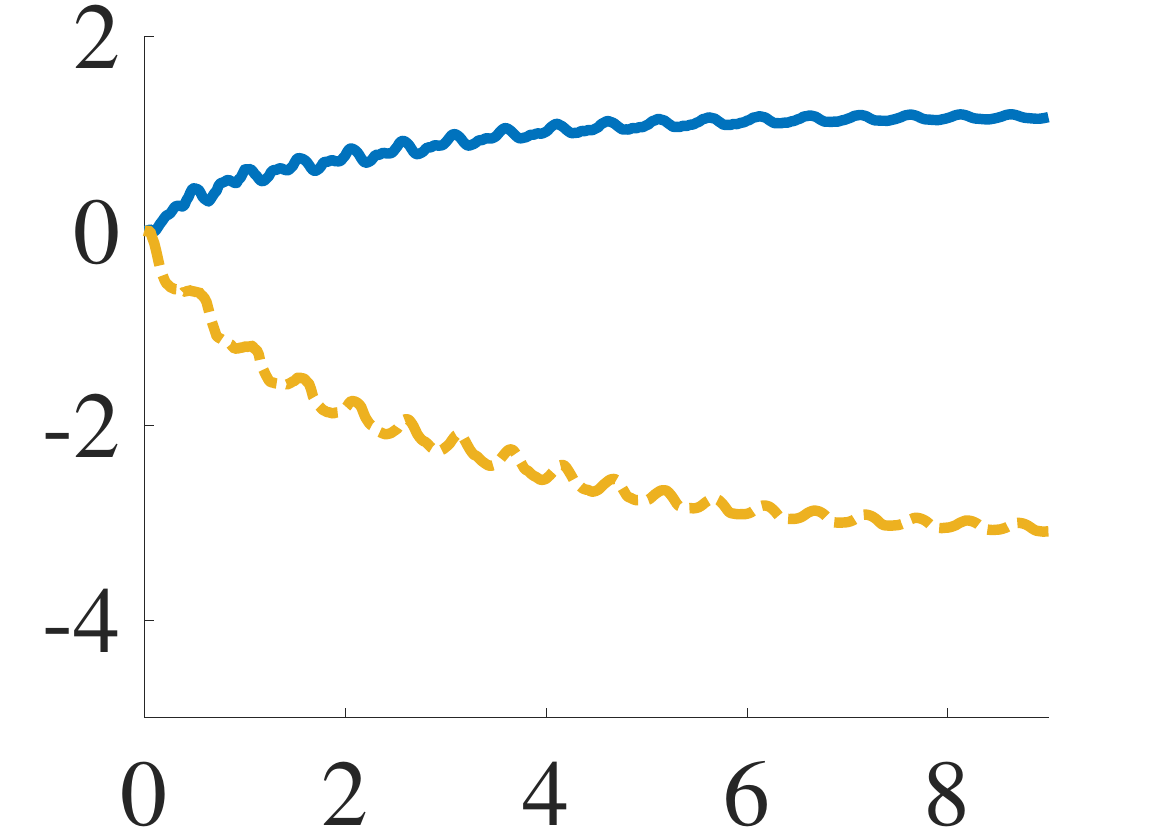} 
\includegraphics[width=.3\linewidth, trim={30 190 25 200}, clip]{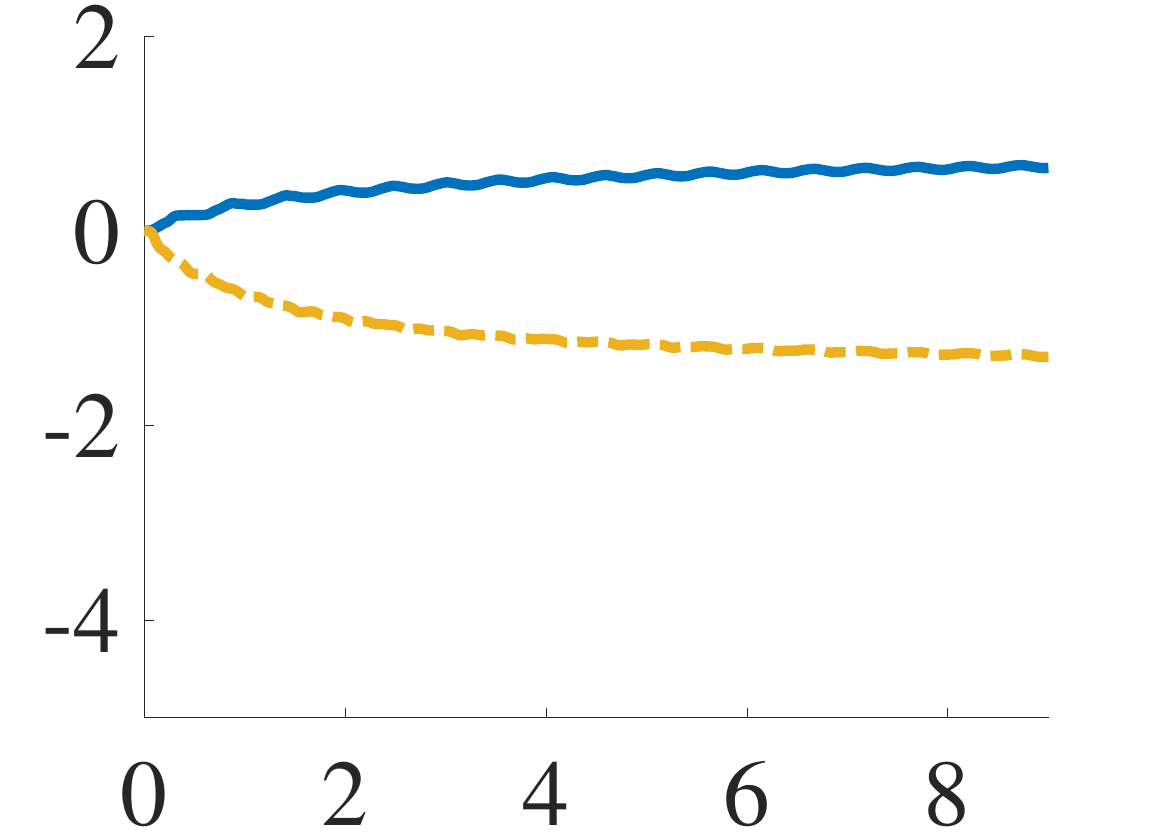}
\includegraphics[width=.3\linewidth, trim={30 190 25 200}, clip]{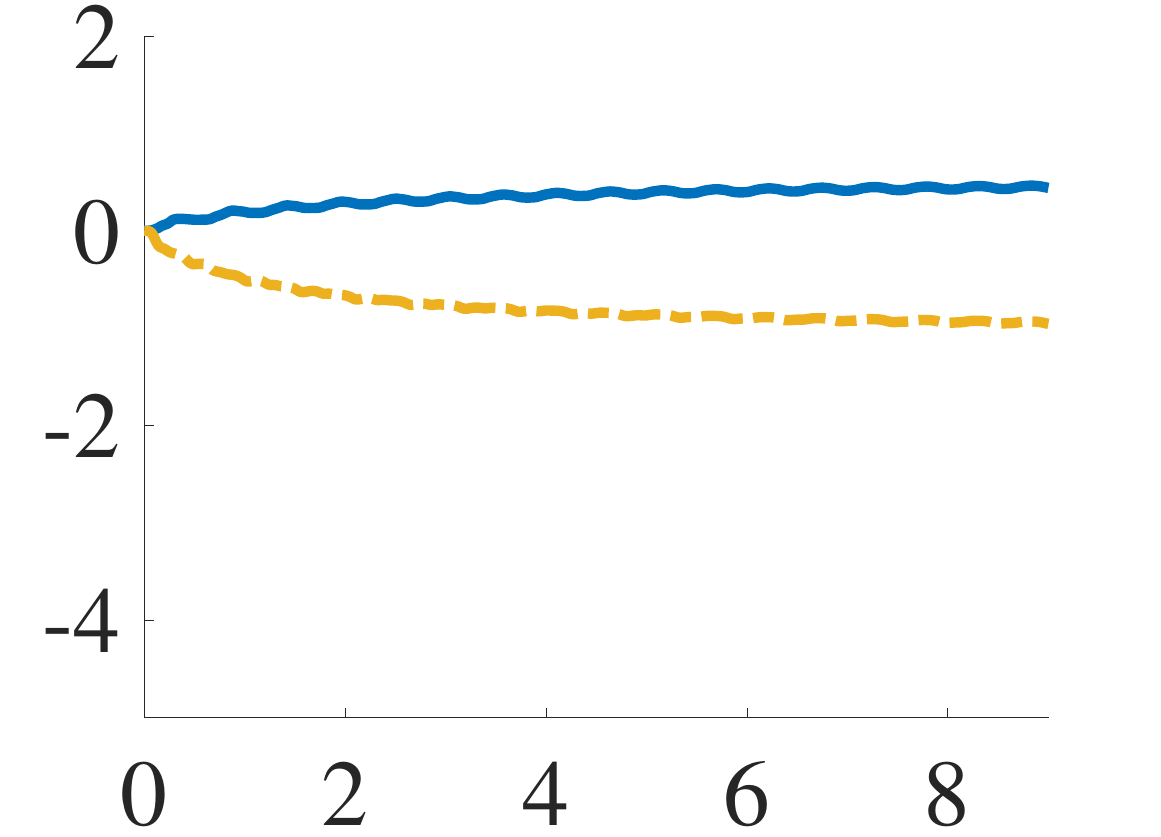}\\

\rotatebox{90}{$\qquad\quad\quad$ \textbf{$\nus = -1$} }
   \rotatebox{90}{$\qquad\quad$ Area Change \%  }
\includegraphics[width=.3\linewidth, trim={30 190 25 200}, clip]{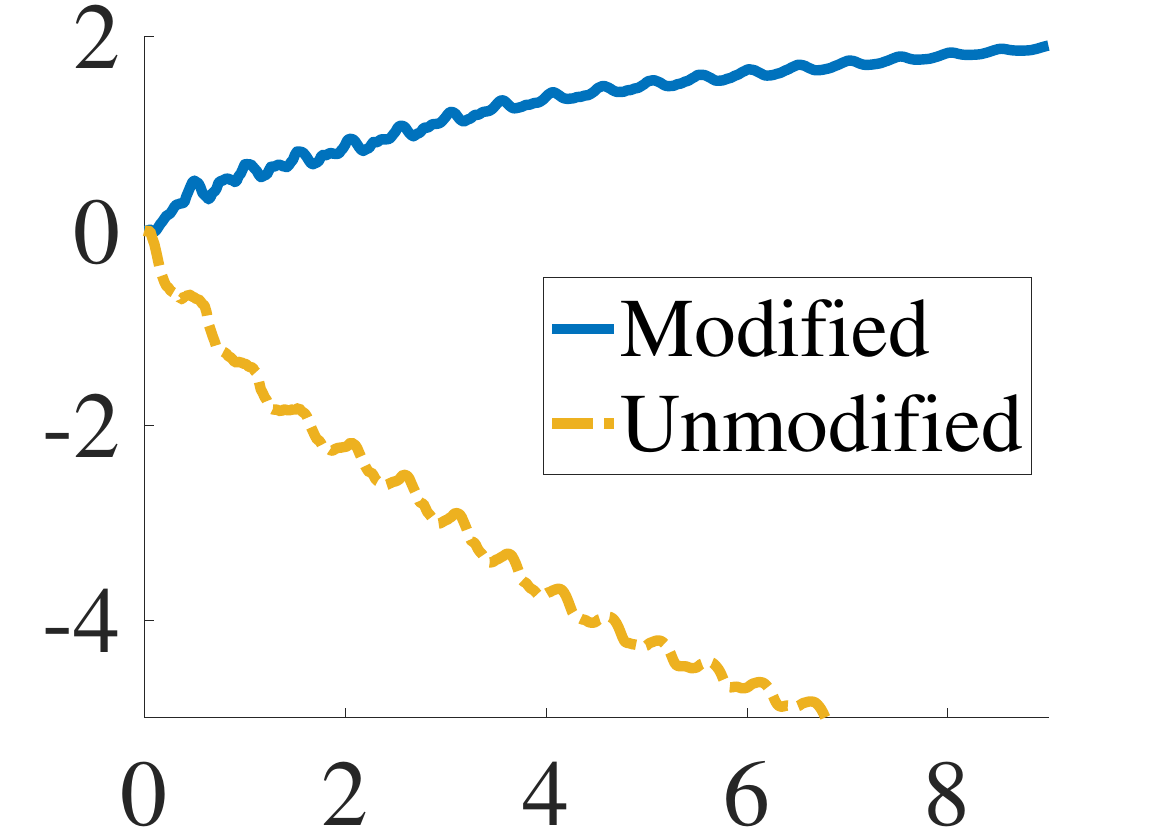} 
\includegraphics[width=.3\linewidth, trim={30 190 25 200}, clip]{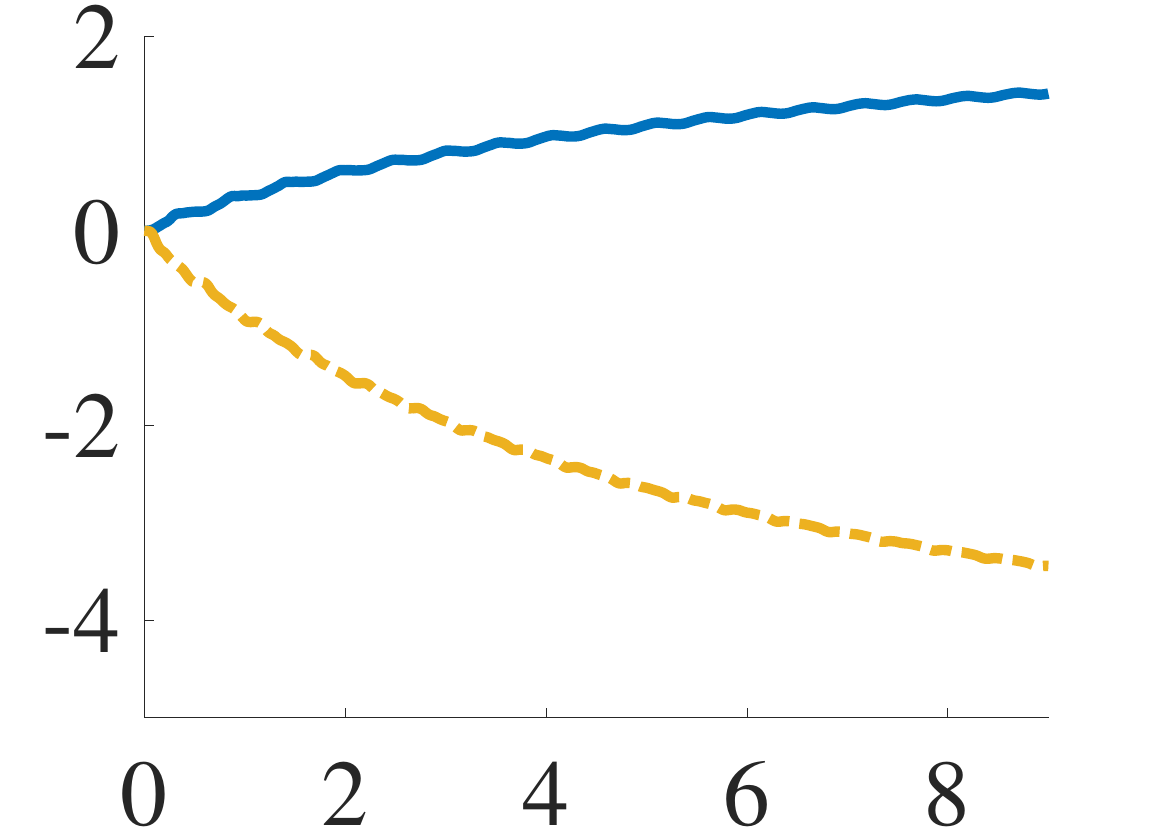}
\includegraphics[width=.3\linewidth, trim={30 190 25 200}, clip]{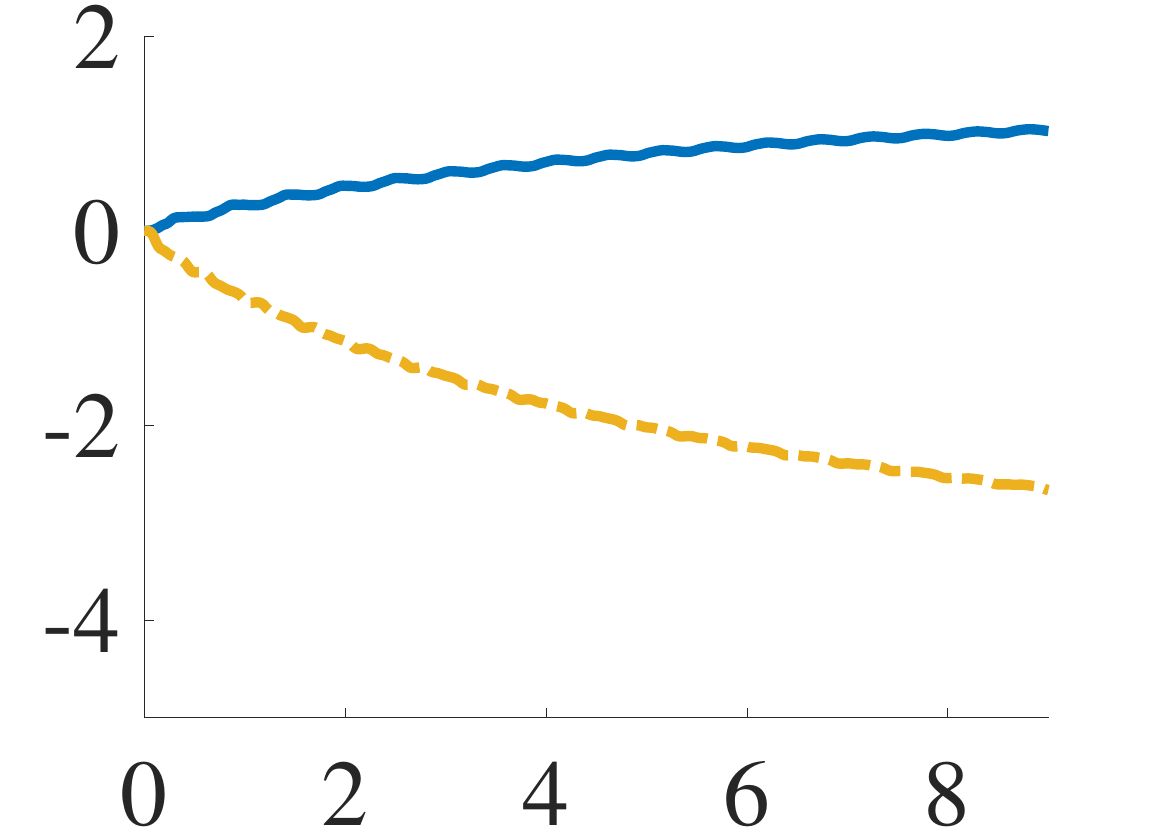}\\

$\qquad\qquad\quad\quad\quad\quad$ Time (s) $\qquad\qquad\quad\quad\quad\quad\quad\quad$ Time (s) $\qquad\qquad\quad\quad\quad\quad\quad\quad$ Time (s)
\caption{Transient behavior of the elastic band (Section \ref{Elastic Band}) for the dynamic version of the test for different discretizations for the entire simulation. Plot shows percent change in area against time.  $N = 16$ corresponds to $m = 135$ solid DOF, $N = 64$ corresponds to $m = 1456$ solid DOF, and $N = 96$ corresponds to $m = 3255$ solid DOF. Results shown here are for \textbf{Q2} elements.}
\label{eb_a_v_time}
\end{figure}

\begin{figure}
\begin{tabular}{l c c}
& \textbf{$\nus = .4$} & \textbf{$\nus = -1$} \\
\rotatebox{90}{$\qquad\quad\qquad$ Mean Disp. (cm) }&
\includegraphics[width=.45\linewidth, trim={50 190 25 200}, clip]{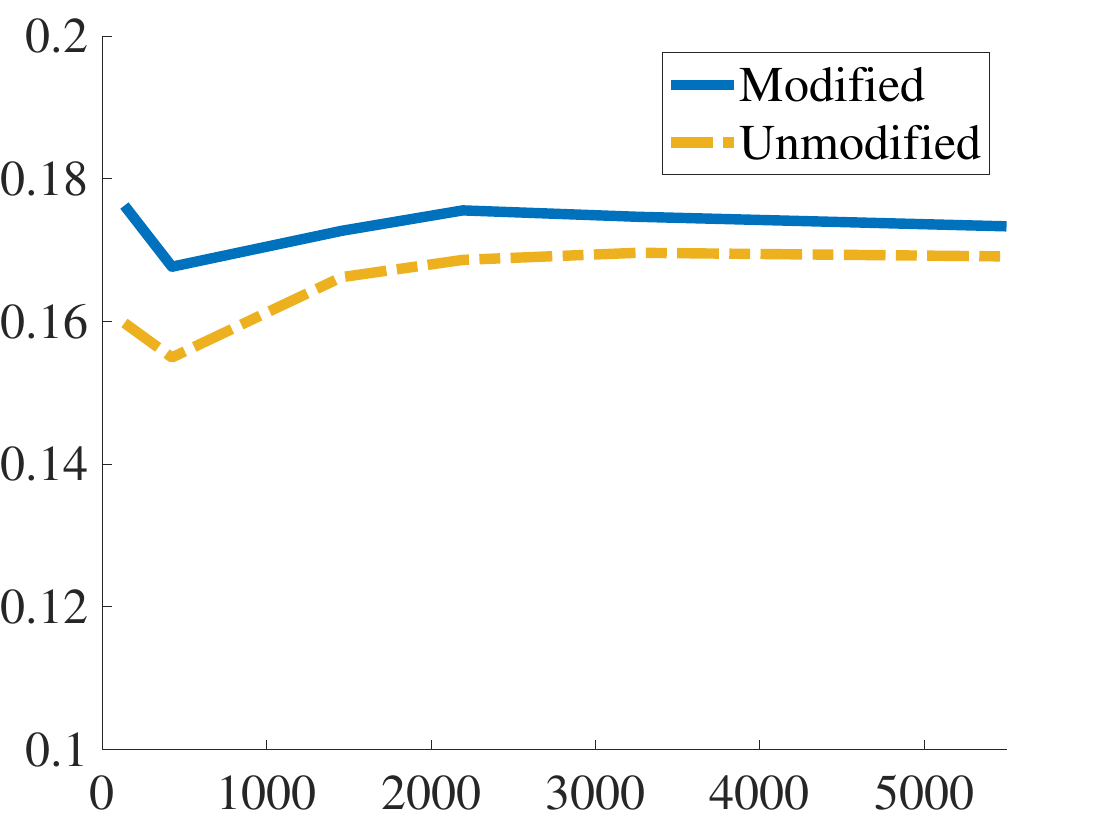} &
\includegraphics[width=.45\linewidth, trim={50 190 25 200}, clip]{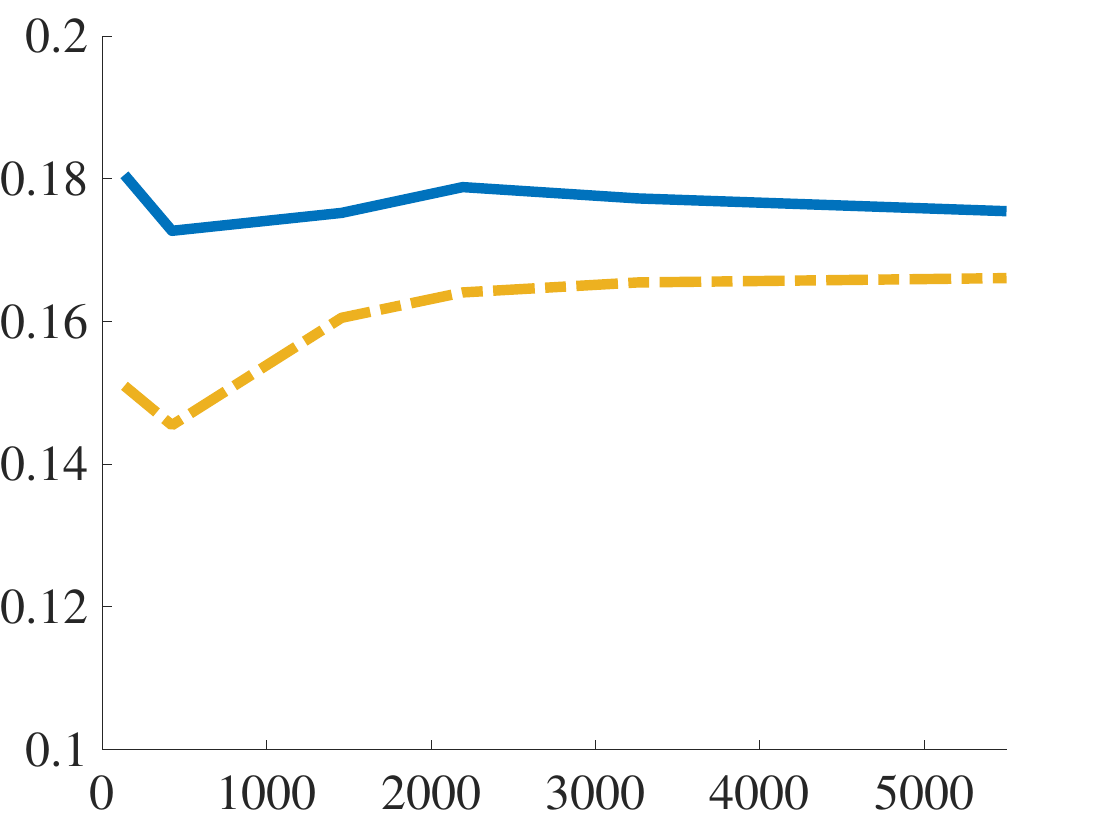}\\

\rotatebox{90}{$\qquad\quad\qquad\quad$ Max Disp. (cm) }&
\includegraphics[width=.45\linewidth, trim={50 190 25 200}, clip]{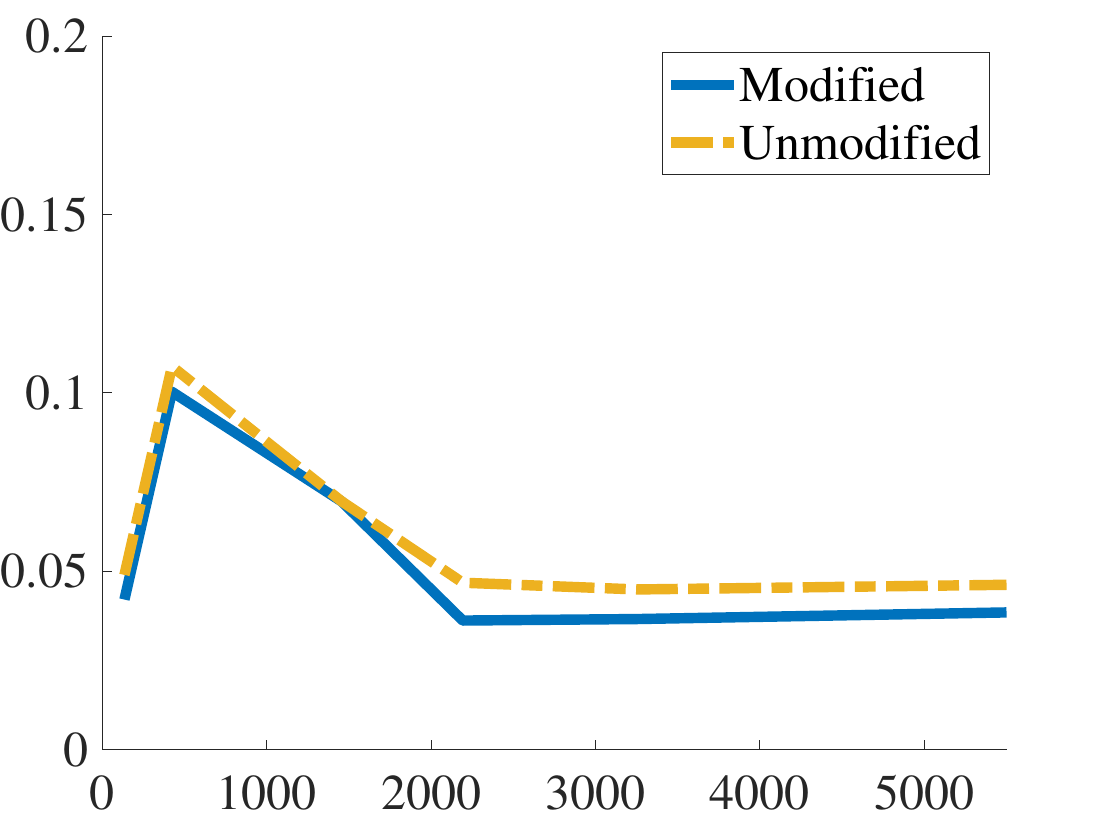} &
\includegraphics[width=.45\linewidth, trim={50 190 25 200}, clip]{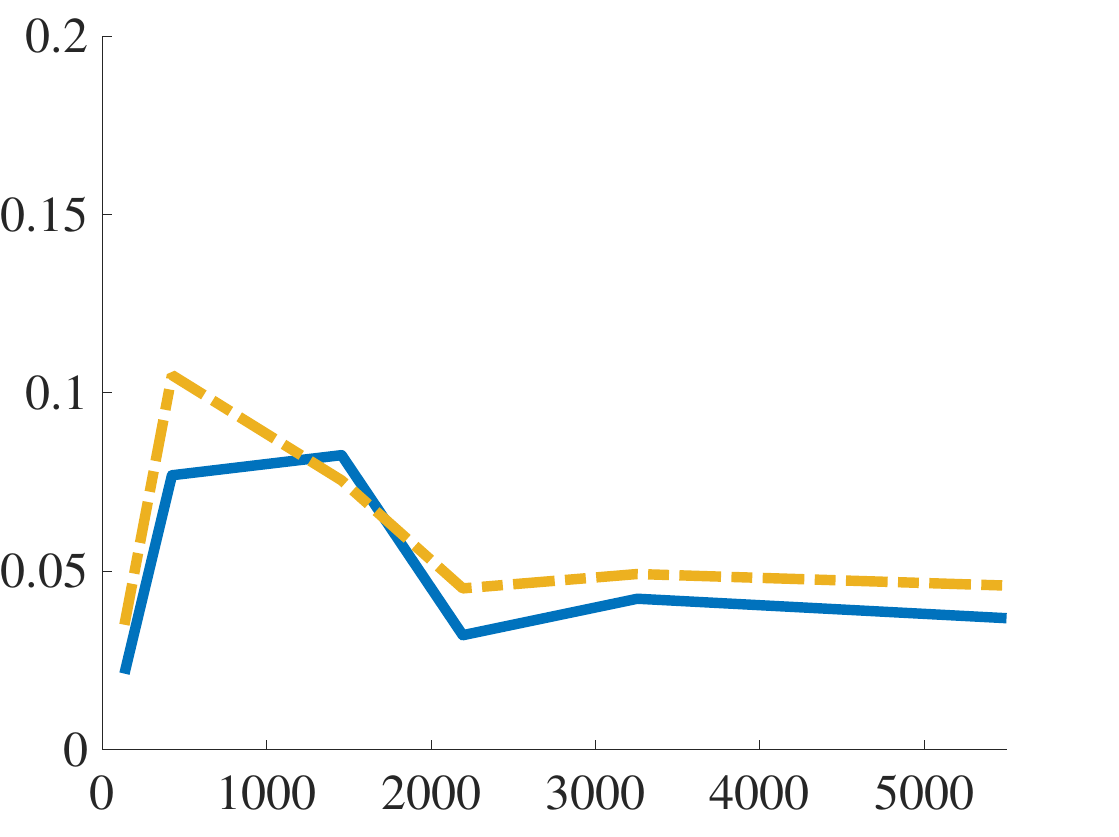}\\

& \# Solid DOF & \# Solid DOF
\end{tabular}
\caption{Mean and max displacement of the transient elastic band benchmark (Section \ref{Elastic Band}) for different solid DOF. The displacement mean and max amplitude were calculated from the data presented in Figure (\ref{eb_v_time}), which is over the time window of $t = 7$ s to $t = 8$ s. We also include another discretization with $m = 5562$ solid DOF.}
%trim={left bottom right top}
\label{eb_mean_amp}
\end{figure}

The elastic band test is a plane strain benchmark in which the loading on the structure is completely driven by fluid forces. The top and bottom of the structure are fixed in place by two rigid blocks as shown in Figure ($\ref{elastic_diag}$). Fluid forces act on the left and right side and so the traction supplied as a boundary condition on the structure is $\boldsymbol{T} = 0$. Unlike previous benchmarks, the computational domain is rectangular $\Omega = [0,L_1] \times [0,L_2]$, with $L_1 = 2 \ \text{cm}$ and $L_2 = 1\ \text{cm}$. Fluid traction boundary conditions of the form $\cauchyf(\xb,t) \nb(\xb) = -\hb(t)$ and $\cauchyf(\xb,t) \nb(\xb)  = \hb(t)$ are imposed on the left and right boundaries of the computational domain, respectively. Here, $\hb(t) = \left(10 \sin \left( \frac{\pi t}{2 T_{\text{l}}}\right),0 \right) \frac{\text{dyn}}{\text{cm}^2}$ when $t < T_{\text{l}}$ and $\hb(t) = (10, 0) \frac{\text{dyn}}{\text{cm}^2}$ otherwise, and the loading time is $T_{\text{l}} = 5$ s. Zero fluid velocity is enforced on the top and bottom boundaries of the computational domain. The $x$-displacement of the middle of the right face is measured at $T_{\text{f}} = 15~\text{s}$; see the encircled point in Figure (\ref{elastic_diag}). The neo-Hookean material model, equations ($\ref{nh_energy}$) -- ($\ref{nh_stress_dev}$), is once again used with a shear modulus of $G = 200$ $\frac{\text{dyn}}{\text{cm}^2}$. The solid viscous damping parameter is set to $\eta = 60 \ \frac{\text{g}}{\text{s}}$. The solid DOF range from $m = 42$ to $m = 3255$. As with the two-dimensional tests, the meshes all share the same node positions. \\
\indent As shown in Figure (\ref{eb}), and in more detail in Figure (\ref{eb_zoom}), the unmodified and unstabililized regime causes the structure to crinkle, whereas the modified and stabilized regime produces smoother deformations. Figures ($\ref{elastic_disp}$) and ($\ref{elastic_area}$) show steady state results for this benchmark. The range of percent change in area for all element types considered is between $0.0015\%$ and $0.69\%$ for the modified invariants, between $0.0\%$ and $2.1\%$ for the unmodified invariants, and between $0.0015\%$ and $0.69\%$ for the deviatoric projection. At steady state, the Lagrange multiplier $\peul$ will be constant within the regions separated by the structure, as shown in Figure (\ref{eb_pressure}). Consequently, we are able to compare the FSI results against results from a quasi-static solid mechanics version of the problem with pressure boundary conditions. The reported results do indicate convergence under grid refinement. Cases with $\nus = .49995$ experience volumetric locking, although it is less pronounced in this example. \\
\indent Because this benchmark includes nontrivial fluid dynamics, we also present results for the transient behavior of a dynamic version of the test. Instead of gradually applying the fluid traction in this dynamic version, it is set to $\hb(t) = (10, 0) \frac{\text{dyn}}{\text{cm}^2}$ for the entire duration of the simulation. The test parameters are the same as the steady state version with the exception of using a final time of $T_{\text{f}} = 9$ s and using no solid damping. This alternate problem specification will result in oscillations of the elastic band, and so we show the Eulerian velocity field of the entire computational domain at a collection of time slices; see Figure ($\ref{eb_flow}$). For the plots of the transient behavior of the dynamic case, see Figure (\ref{eb_v_time}) for the displacement plotted against time and Figure (\ref{eb_a_v_time}) for the total area plotted against time. It can be seen that the dynamics of the point of interest are similar between modified and unmodified invariants for $\nus = 0.4$, and there is a noticeable difference for $\nus = -1.0$. For modified and unmodified cases, Figure (\ref{eb_a_v_time}) shows the area change noticeably decreases when $\nus$ increases from $\nus = -1$ to $\nus = 0.4$. Additionally, the total area change decreases under grid refinement for the dynamic case. Figure (\ref{eb_mean_amp}) shows that the mean and amplitude of the oscillations in displacement converge under grid refinement for the stabilized case. 

%{elastic_diag, eb, eb_zoom, elastic_disp, elastic_area, eb_pressure, eb_flow, eb_v_time, eb_a_v_time, eb_mean_amp

%%%%%%%%%%%%%%%%%%%%%%%
\section{Effect of volumetric stabilization on a model of esophageal transport}
\label{eso-transport}
The phenomenon of esophageal transport involves the passage of a fluid
bolus that is accompanied by a significant amount of solid deformation
in the esophageal walls. The amount of deformation is closely related
to the velocity and pressure of the fluid, and the
wall stiffness also has an important effect on the rate of bolus
emptying through peristalsis. Thus esophageal transport is a strongly coupled fluid-structure
interaction problem. This Section demonstrates that introducing volumetric stabilization, along with using modified
invariants to describe the isotropic material response, substantially
reduces unphysical deformations during simulations
of esophageal bolus transport and emptying.

\begin{figure}
\captionsetup[subfigure]{justification=centering}
\begin{subfigure}{.7\textwidth}
  \centering
  \includegraphics[width=.95\linewidth, trim={0 0 0 0}, clip]{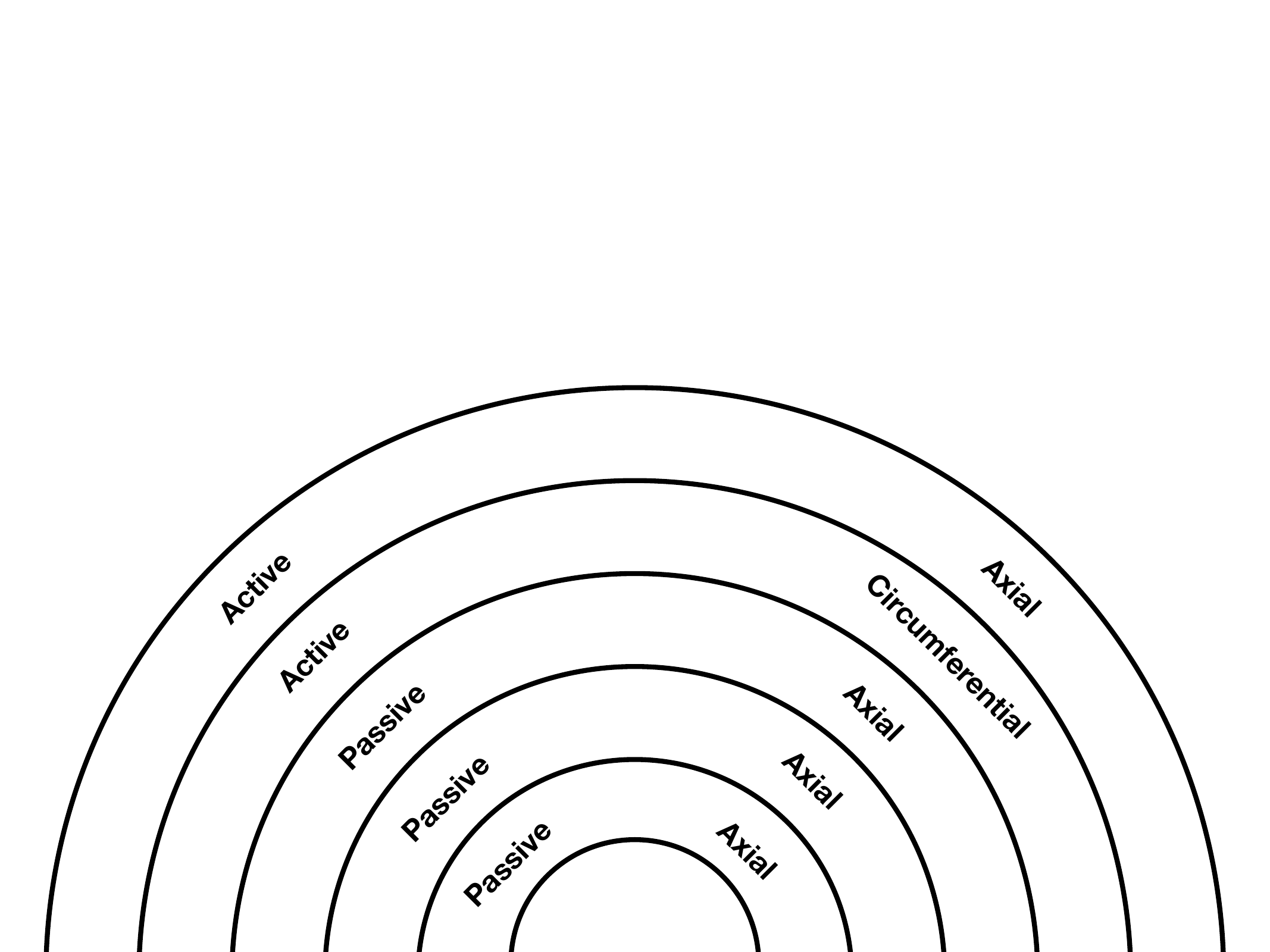}
  \caption{Tissue layer properties}
%  \label{fig:fixed_end_A}
\end{subfigure}
\begin{subfigure}{.23\textwidth}
  \centering
  \includegraphics[width=.95\linewidth, trim={0 0 0 0}, clip]{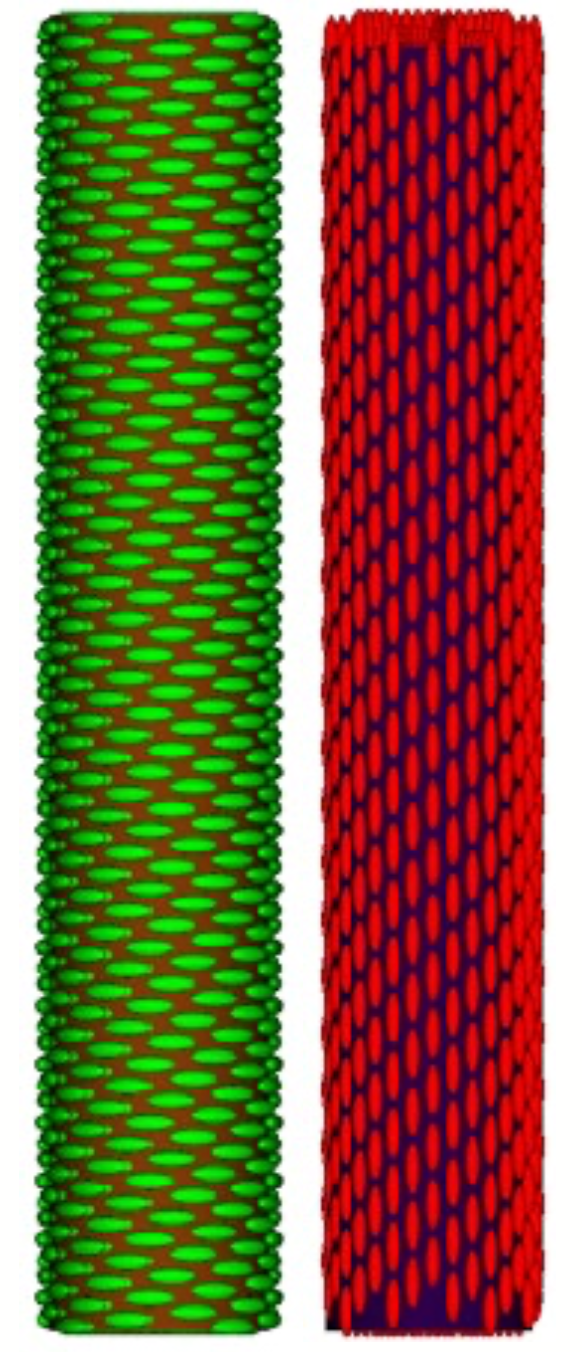}
  \caption{Fiber Directions}
%  \label{fig:fixed_end_B}
\end{subfigure}
\caption{Diagram depicting properties of each layer of the esophageal model and the fiber directions used in that model (not drawn to scale). The esophageal tube is shown from a top down view in panel (a). A label of ``passive" means that the passive model, equations (\ref{passive-energy} -- \ref{passive-stress-mod}), is used in that region, whereas a label of ``active" means that the active model, equations (\ref{active-energy} -- \ref{active-stress-mod}), is used in that region. ``Axial" and ``circumferential" mean the fibers run in those directions, respectively. Panel (b) shows fibers running in the circumferential direction, which is used for the fourth layer, and fibers running in the axial direction, which is used for all other layers. The diagrams in panel (b) were previously reported in a work by Kou \etal~\cite{Kou2017b}}
\label{eso_diag}
\end{figure}

We construct a model that follows previously published work by Kou \etal~\cite{Kou2017, Kou2017b}.
The esophageal model used for this demonstration consists of a thick-walled cylindrical
tube. The walls are divided into five layers, and each of these layers has different material properties and fiber orientations
that resemble observed tissue mechanical properties; see Figure (\ref{eso_diag}).
The outermost layer and the three inner layers are described as an anisotropic
material with fibers oriented along the axis of the tube, and the fibers of the fourth
layer have circumferential orientation. The three innermost
layers do not have any active component resembling muscular contractions.
The remaining two layers contract in a manner resembling peristalsis. This effect is achieved by changing the rest length of the fibers that are embedded in these layers.\\
\indent The three passive layers have a neo-Hookean isotropic component and
an anisotropic fiber component that is described by
\begin{align}
\Psi &=\frac{G_{1}}{2}\left(I_{1}-3\right)+\frac{G_{2}}{2}\left(I_{4}-1\right)^{2} + \frac{\kappas}{2}\left(\ln J\right)^2, \ \text{and} \label{passive-energy} \\
\PPs &= G_1\FF + 2G_2(I_4 - 1)\FF\MM + \kappas\left(\ln J\right)\FF^{-T}.
\end{align}
The modified versions are
\begin{align}
	\Psi &=\frac{G_{1}}{2}\left(\bar{I}_{1}-3\right)+\frac{G_{2}}{2}\left(I_{4}-1\right)^2 + \frac{\kappas}{2}\left(\ln J\right)^2, \ \text{and} \\
	\PPs &= G_1J^{-2/3}\left(\FF -\frac{I_1}{3}\FF^{-T}\right) + 2G_2(I_4 - 1)\FF\MM + \kappas\left(\ln J\right)\FF^{-T}. \label{passive-stress-mod}
\end{align}
For the active layers, the strain energy and stress are, respectively,
\begin{align}
	\Psi &=\frac{G_{3}}{2}\left(I_{1}-3\right)+\frac{G_{4}}{2}\left(\frac{\sqrt{I_{4}}}{\lambda}-1\right)^{2} + \frac{\kappas}{2}\left(\ln J\right)^2, \ \text{and} \label{active-energy} \\
	\PPs &= G_3\FF + G_4\left(\frac{I_4 - \lambda \sqrt{I_{4}}}{\lambda^2 I_4}\right)\FF\MM+ \kappas\left(\ln J\right)\FF^{-T},
\end{align}
in which $\lambda = \lambda(\Xb,t)$ is the resting length. The modified versions are 
\begin{align}
	\Psi &=\frac{G_{3}}{2}\left(\bar{I}_{1}-3\right)+\frac{G_{4}}{2}\left(\frac{\sqrt{I_{4}}}{\lambda}-1\right)^2 + \frac{\kappas}{2}\left(\ln J\right)^2, \ \text{and}  \\
	\PPs &= G_3J^{-2/3}\left(\FF -\frac{I_1}{3}\FF^{-T}\right) + G_4\left(\frac{I_4 - \lambda \sqrt{I_{4}}}{\lambda^2 I_4}\right)\FF\MM+ \kappas\left(\ln J\right)\FF^{-T}. \label{active-stress-mod}
\end{align}
The effect of a peristaltic contraction on the system is achieved by
changing the value of $\lambda$ in a sinusoidal, traveling wave-like
manner via
\begin{equation}
	\lambda(\Xb,t) = \lambda_0 e^{-(X_3 - ct)^2/(2\omega^2)}.
\end{equation}
This wave travels at a speed of $c = 10 \, \frac{\text{cm}}{\text{s}}$, has width $\omega = 1.5$ cm, and max amplitude $\lambda_0 = 0.4$. The motivation for these fiber strain energy functions and model of peristalsis is discussed in previous work \cite{Merodio2003, Yang2006}.\\
\indent The innermost
layers model the mucosal part of the esophagus and have
stiffnesses of $G_1 = 40.0~\frac{\text{dyn}}{\text{cm}^2}$ for the isotropic matrix and $G_2 = 400.0~\frac{\text{dyn}}{\text{cm}^2}$ for the anisotropic part. The active layers of the esophagus
are considerably stiffer, and the assigned material properties for
these two layers are $G_3 = 4000~\frac{\text{dyn}}{\text{cm}^2}$ for the isotropic part and $G_4 = 40,000~\frac{\text{dyn}}{\text{cm}^2}$ for
the anisotropic part. Additional physiological details, implementation
of peristaltic contraction, insights on esophageal transport, and the
effect of different material properties and strain energy functions
on bolus transport are discussed in prior work by Kou \etal~\cite{Kou2017}. \\
\indent As in Section \ref{Anisotropic Cook's Membrane}, we only modify $I_1$ for the modified cases. Because the layers have different stiffnesses,
different values of $\kappas$ are computed for each layer using
its respective stiffness, and we use either $\nus = 0.4$ for the stabilized case or $\nus = -1.0$ for the unstabilized case, using equation (\ref{kappa-nu}). For each layer, we select the largest elastic modulus to use in equation (\ref{kappa-nu}).\\
\indent The computational domain in this problem is $\Omega = [-0.7,0.7] \times [-0.7,0.7] \times [-1.6, 18.5]$, in which each linear dimension has units of cm. The inner and outer radii of the tube are $r_{\text{i}} = 0.03$ cm and $r_{\text{o}} = 0.53$ cm, respectively, and the tube extends from $z=0$ cm to $z=18.0$ cm in the reference configuration. The density of the fluid is $\rho = 1.0 \frac{\text{g}}{\text{cm}^3}$ and the viscosity is $\mu = 0.1~\frac{\text{dyn}\cdot s}{\text{cm}^2}$. These values of density and viscosity reflect the properties of ingested food after
it has been considerably mixed with saliva from chewing and forms
a bolus when entering the esophagus. The velocity at the top face of the computational domain is zero
and zero traction boundary conditions are enforced for the remaining
fluid boundary faces. The top end of the esophageal tube is fixed
using a penalty force and the bottom end of the tube is free
to move.\\
\indent The computational domain is discretized by $N_x = 70,\, N_y = 70$ and $N_z = 335$ grid cells in the $x,\, y$, and $z$ directions respectively. The Lagrangian domain is described using Q1 elements with $m_{\text{r}} = 7$, $m_{\text{c}} = 32$ and $m_{\text{a}} = 180$ elements in the radial, circumferential, and axial directions, respectively.
 
 \begin{figure}
\begin{centering}
\includegraphics[width=0.9\linewidth, trim={70 20 10 30}, clip]{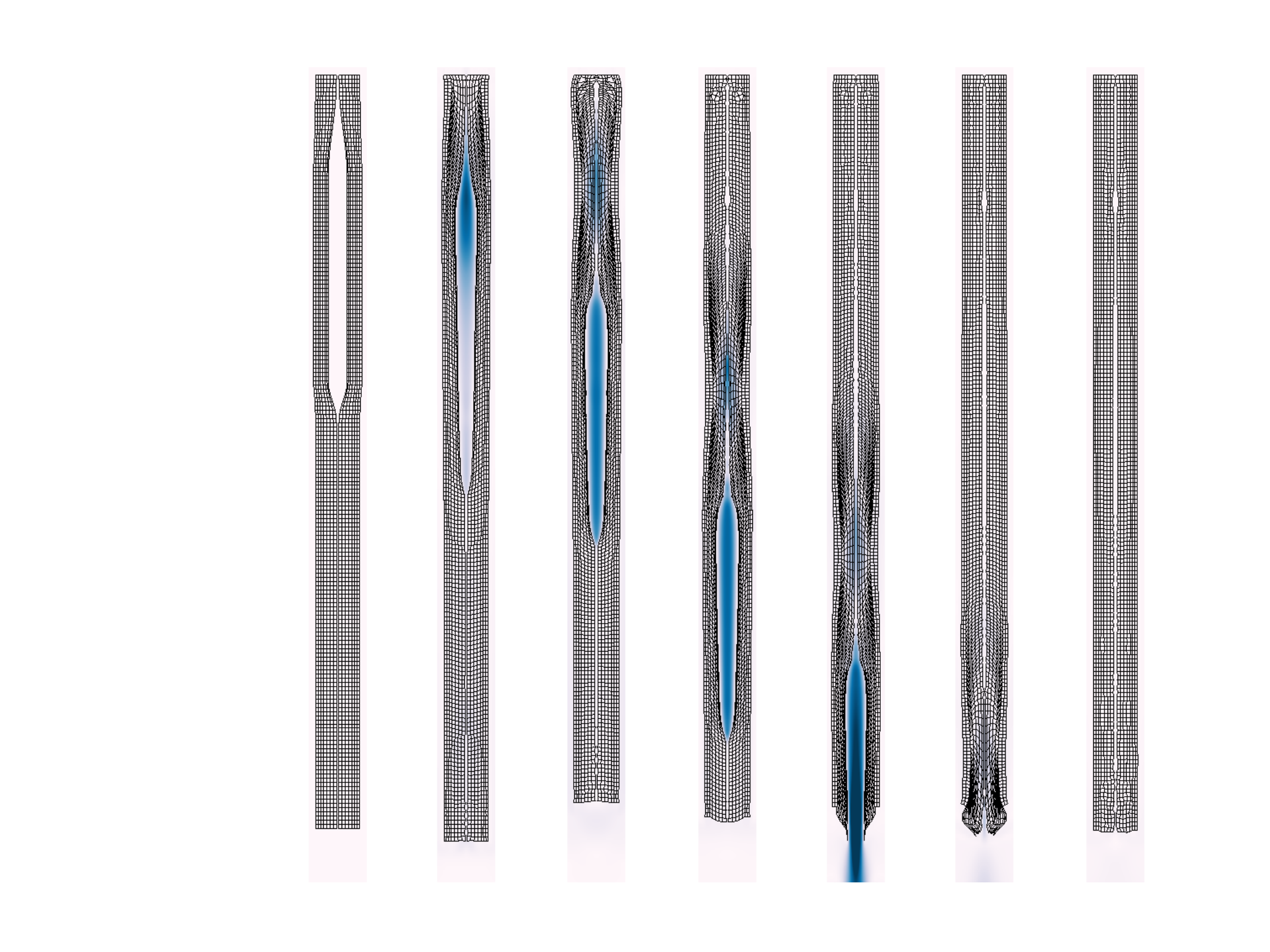}
\end{centering}
\vskip .1cm
\begin{centering}
$||\ub(\xb,t)|| \ \left(\frac{\text{cm}}{\text{s}}\right)$ \\
\includegraphics[width=2.5in, trim={0 5in 0 5in}, clip]{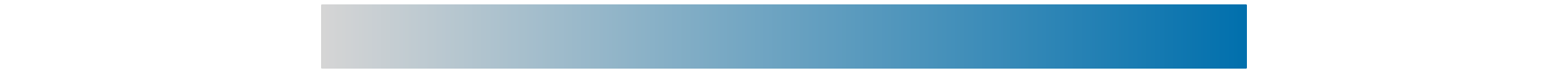}  \\
%0.00 \ \ \ \ \ \ \ \ \ \ \ \ \ \ \ \ 2.00
0.0 $\qquad\qquad\qquad\qquad$ 25.0\\
\end{centering}
%\begin{centering}
%$||\ub(\xb,t)||$ \\
%\includegraphics[width=2.5in, trim={0 5in 0 5in}, clip]{blue_color_bar.pdf}
%\vskip .1cm
%0.0 $\qquad\qquad\qquad\qquad$ 2.0
%\end{centering}
\caption{Two dimensional cross-section of a bolus traveling down a three dimensional model of the esophagus resulting from peristalsis (Section \ref{eso-transport}). The panels show the structural deformations along with contours of velocity magnitude. Shown here is the stabilized case with modified invariants. Note that after the bolus has traveled through the vessel, the geometry returns to its initial configuration with minimal artifacts.}
\label{fig:Bolus-traveling-down-sequence}
\end{figure}

\begin{figure}
\captionsetup[subfigure]{justification=centering}
\begin{subfigure}{.24\textwidth}
  \centering
  \includegraphics[width=.95\linewidth]{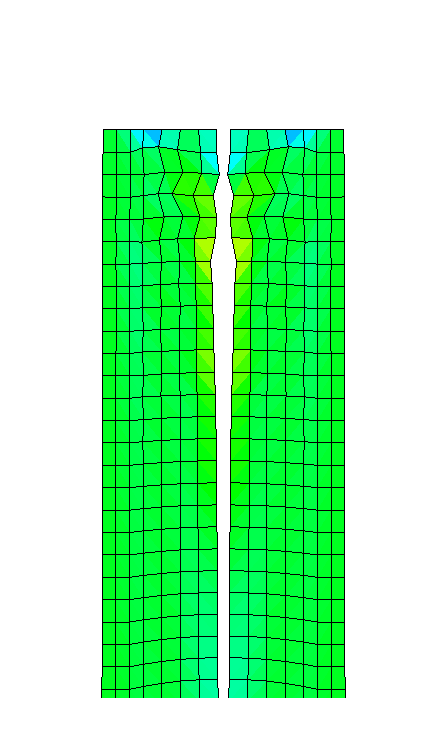}
  \caption{Modified Invariants \\ $\nus = .4$}
%  \label{fig:fixed_end_A}
\end{subfigure}
\begin{subfigure}{.24\textwidth}
  \centering
\includegraphics[width=.95\linewidth]{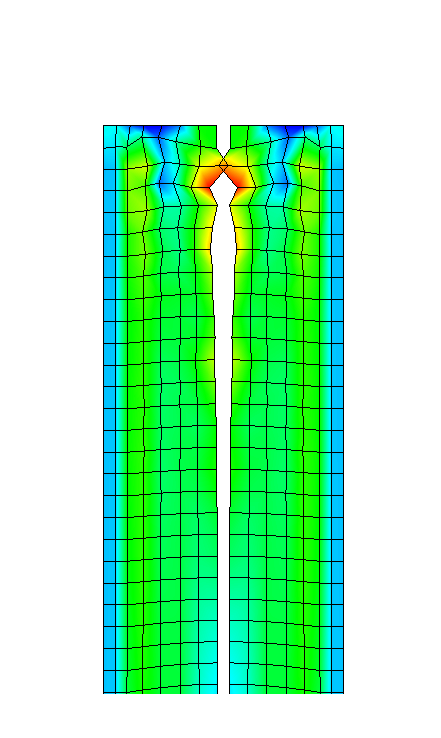}
  \caption{Unmodified Invariants\\ $\nus = .4$}
%  \label{fig:fixed_end_B}
\end{subfigure}
\begin{subfigure}{.24\textwidth}
  \centering
 \includegraphics[width=.95\linewidth,trim={0 0 0 0},clip]{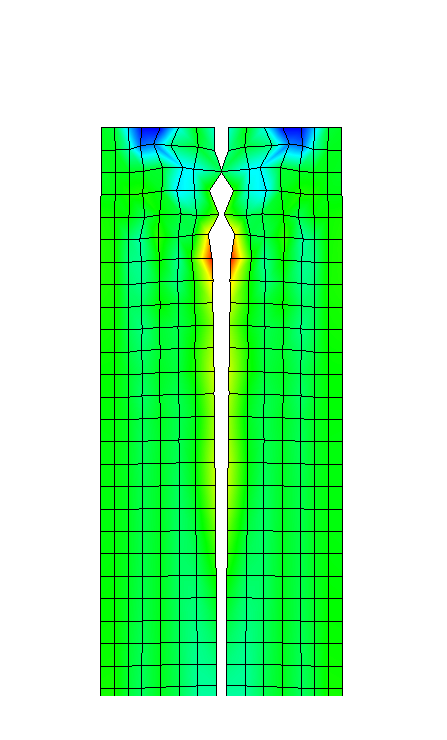}
  \caption{Modified Invariants \\ $\nus = -1$}
%  \label{fig:fixed_end_C}
\end{subfigure}
\begin{subfigure}{.24\textwidth}
  \centering
\includegraphics[width=.95\linewidth,trim={0 5 0 0},clip]{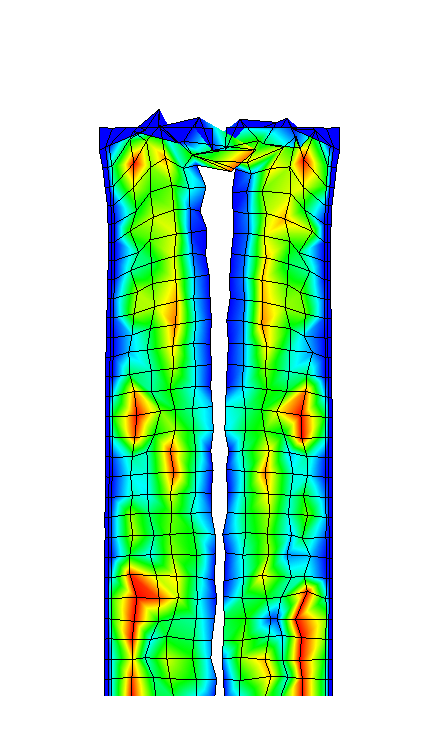}
  \caption{Unmodified Invariants \\ $\nus = -1$}
%  \label{fig:fixed_end_D}
\end{subfigure}

\begin{centering}
$\text{Avg} \ J$ \\
\includegraphics[width=2.5in, trim={0 5in 0 5in}, clip]{color_bar.pdf}  \\
%0.00 \ \ \ \ \ \ \ \ \ \ \ \ \ \ \ \ 2.00
0.75 $\qquad\qquad\qquad\qquad$ 1.30\\
\end{centering}
\caption{Two dimensional cross-sections of the deformation fields and mean values of $J$ for the esophageal transport model (Section \ref{eso-transport}), along with values of $J$, at the fixed end of the tube at time $t = 0.8$ s. Note the that panel (d) shows extreme deformations for the case of unmodified invariants and no stabilization.}
\label{deformations_fixed_end}
\end{figure}

Figure (\ref{fig:Bolus-traveling-down-sequence}) shows the transport of the bolus down the length of the esophagus for the case with modified invariants and volumetric stabilization. Once the bolus is vacated, the esophagus returns to its undeformed state with minimal deformation artifacts. For each of the four cases, we inspect cross-sections of the deformations of the esophageal structure at time $t = 0.8$ s. This occurs when the bolus has traveled slightly past halfway along its length, i.e.~the middle panel shown in Figure (\ref{fig:Bolus-traveling-down-sequence}).
During this moment, we consider the section of the esophagus behind the bolus and the deformations observed at the free and fixed ends. \\
\indent During the passage of the bolus, the tube walls undergo large deformations and eventually return to their undeformed shape. We consider the deformation
patterns at the fixed end of the tube because it is the first section that is allowed to return to its original configuration. Figure (\ref{deformations_fixed_end}) shows the Lagrangian mesh when the contraction wave has propagated partway down the tube. It is clear that the most unphysical deformations occur in the case with unmodified invariants and no volumetric stabilization. We also observe that that the recovery of the tube is imperfect in all cases, with case (a) most closely resembling the undeformed configuration.

\begin{figure}
\captionsetup[subfigure]{justification=centering}
\begin{subfigure}{.24\textwidth}
  \centering
  \includegraphics[width=.95\linewidth,trim={5 0 5 0},clip]{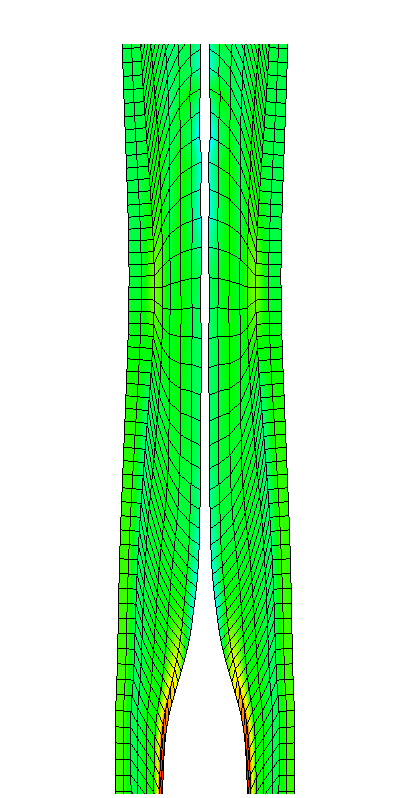}
  \caption{Modified Invariants \\ $\nus = .4$}
  \label{fig:behind_bolus_A}
\end{subfigure}
\begin{subfigure}{.24\textwidth}
  \centering
\includegraphics[width=.95\linewidth,trim={5 0 5 0},clip]{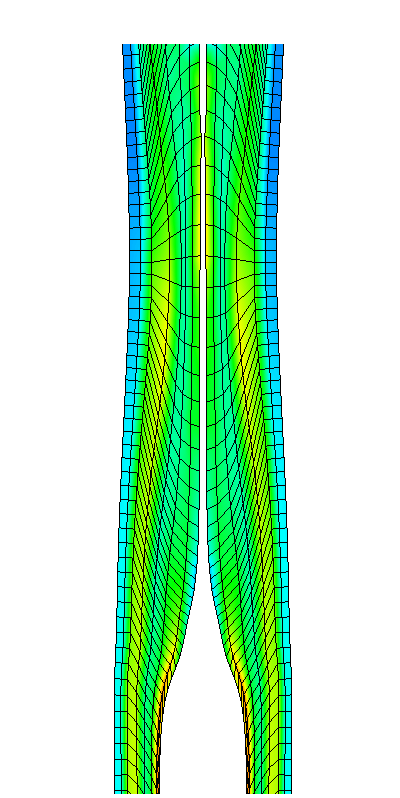}
  \caption{Unmodified Invariants \\ $\nus = .4$}
  \label{fig:behind_bolus_B}
\end{subfigure}
\begin{subfigure}{.24\textwidth}
  \centering
\includegraphics[width=.95\linewidth,trim={5 0 5 0},clip]{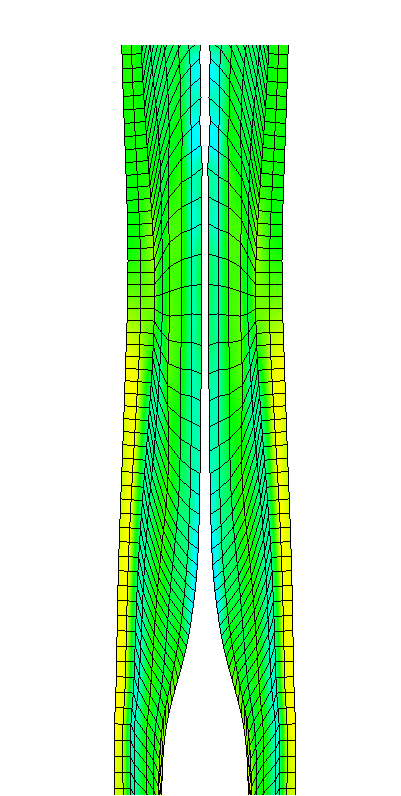}
  \caption{Modified Invariants \\ $\nus = -1$}
  \label{fig:behind_bolus_C}
\end{subfigure}
\begin{subfigure}{.24\textwidth}
  \centering
\includegraphics[width=.95\linewidth,trim={5 0 5 0},clip]{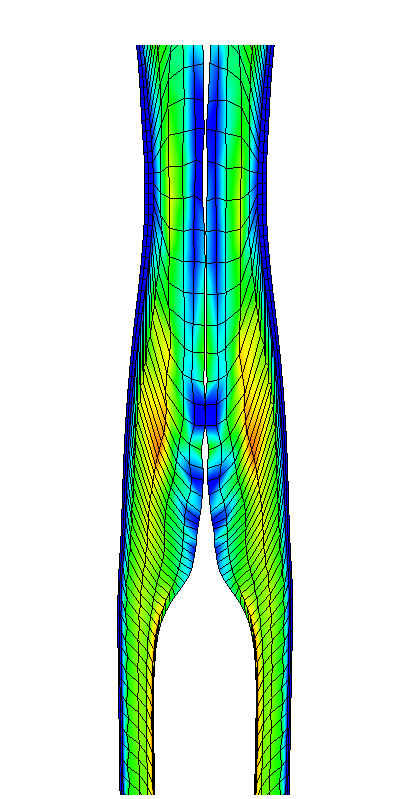}
  \caption{Unmodified Invariants \\ $\nus = -1$}
  \label{fig:behind_bolus_D}
\end{subfigure}

\begin{centering}
$\text{Avg} \ J$ \\
\includegraphics[width=2.5in, trim={0 5in 0 5in}, clip]{color_bar.pdf}  \\
%0.00 \ \ \ \ \ \ \ \ \ \ \ \ \ \ \ \ 2.00
0.75 $\qquad\qquad\qquad\qquad$ 1.30\\
\end{centering}
\caption{Deformation fields and mean values of $J$ of the esophageal transport model (Section \ref{eso-transport}) behind the bolus at $t = 0.8$ s. Here the differences between each case is less pronounced, but note the wider distribution of $J$ as found in panel (d) as compared to panel (a).}
\label{fig:deformations_behind_bolus}
\end{figure}

\begin{figure}
\captionsetup[subfigure]{justification=centering}
\begin{subfigure}{.24\textwidth}
  \centering
\includegraphics[width=.95\linewidth,trim={5 0 5 0},clip]{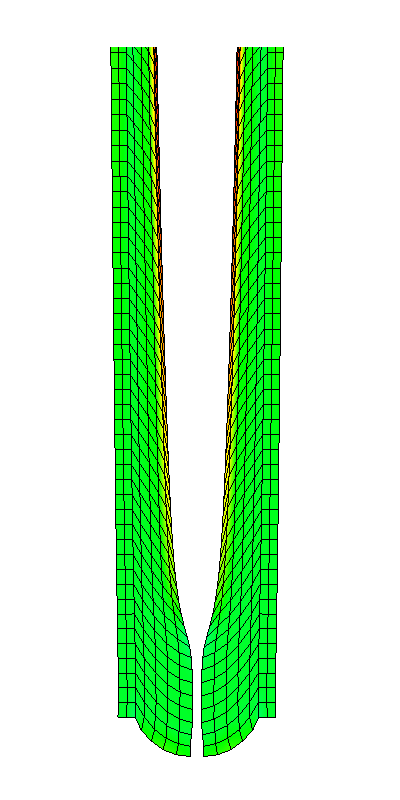}
  \caption{Modified Invariants\\ $\nus = .4$}
  \label{fig:free_end_A}
\end{subfigure}
\begin{subfigure}{.24\textwidth}
  \centering
 \includegraphics[width=.95\linewidth,trim={5 0 5 0},clip]{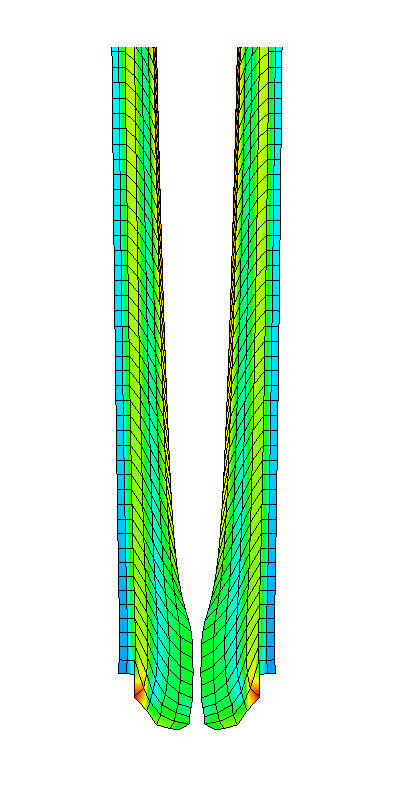}
  \caption{Unmodified Invariants\\ $\nus = .4$}
  \label{fig:free_end_B}
\end{subfigure}
\begin{subfigure}{.24\textwidth}
  \centering
\includegraphics[width=.95\linewidth,trim={5 0 5 0},clip]{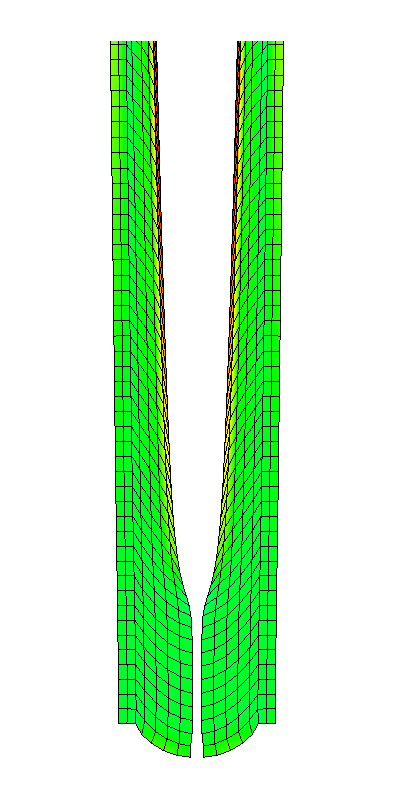}
  \caption{Modified Invariants\\ $\nus = -1$}
  \label{fig:free_end_C}
\end{subfigure}
\begin{subfigure}{.24\textwidth}
  \centering
\includegraphics[width=.95\linewidth,trim={5 0 5 0},clip]{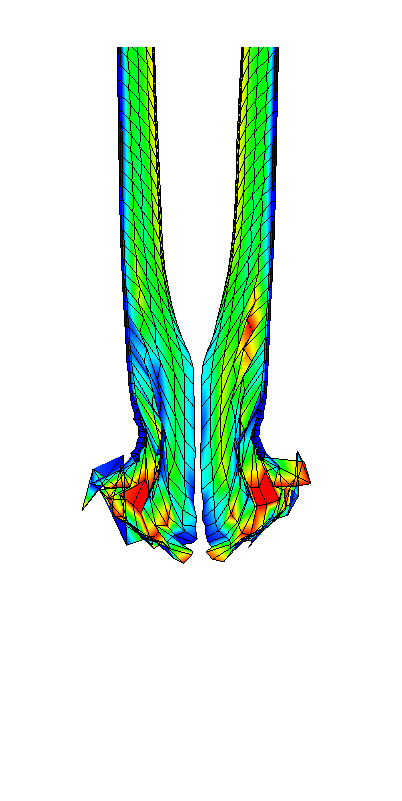}
  \caption{Unmodified Invariants\\ $\nus = -1$}
  \label{fig:free_end_D}
\end{subfigure}
\begin{centering}
\vskip .5cm
$\text{Avg} \ J$ \\
\includegraphics[width=2.5in, trim={0 5in 0 5in}, clip]{color_bar.pdf}  \\
%0.00 \ \ \ \ \ \ \ \ \ \ \ \ \ \ \ \ 2.00
0.75 $\qquad\qquad\qquad\qquad$ 1.30\\
\end{centering}
\caption{Deformations and mean values of $J$ of the esophageal transport model (Section \ref{eso-transport}) at the free end of the tube at $t = 0.8$ s. Note that the case with unmodified invariants and no stabilization once again experiences extreme deformations. Also, note how panels (a) and (c) are nearly indistinguishable.}
\label{fig:deformations_free_end}
\end{figure}

Figure (\ref{fig:deformations_behind_bolus}) shows the deformations of the tube behind the bolus after
the contraction has traveled along the tube. As before, Figure (\ref{fig:behind_bolus_D}) shows that the unstabilized and unmodified case yields unphysical deformations. The other three scenarios display reasonable deformation fields,
and the computed bolus profile is quite similar. Further differentiation between these three approaches can
be made when the tube profile at the free end is plotted as in Figure (\ref{fig:deformations_free_end}).
Observing the deformations of the tube at the free end during the bolus transport process clearly demonstrates the
effectiveness of the stabilization methods and their ability to reduce or remove unphysical deformations. When no stabilization is applied, we see completely unphysical
deformations. In our tests, the case shown in Figure (\ref{fig:free_end_D}) never reached completion because of extreme unphysical deformations. Between the cases shown in Figure \ref{fig:deformations_free_end}, panels (b) and (c), we see that, for this problem at least, it is preferable to use modified invariants with no volumetric energy as opposed to using a volumetric energy stabilization with unmodified invariants. This is because the difference between the modified stabilized case and the modified unstabilized case is less pronounced in this problem. If unmodified invariants are used, we observe irregular deformations of the mucosal layers regardless of whether or not volumetric stabilization is employed. Finally, at the final time step of each simulation the percent changes in volume were $0.7\%$ for the modified stabilized case, $5.0\%$ for the unmodified stabilized case, $1.2\%$ for the modified unstabilized case, and $262.2\%$ for the unmodified unstabilized case.

%%%%%%%%%%%%%%%%

\section{Discussion and Conclusion}
\indent This work develops simple stabilization methods for the hyperelastic IB method that correct for spurious volume changes that the immersed structure may experience in the discretized equations. Following standard nonlinear solid mechanics formulations, our strategy is to include an additional volumetric energy in the strain energy functional of the immersed structure. Although our numerical examples demonstrate that the volumetric stabilization alone may improve the accuracy of the method, the Cauchy stress can also be decomposed into deviatoric and dilatational components, offering further improvements in accuracy. We present two approaches to achieve a deviatoric Cauchy stress in the solid region: using the Flory decomposition to reformulate the strain energy with the modified strain invariants \cite{Flory1961}, and using the deviatoric projection operator (\ref{dev}). Prior benchmark studies in the solid mechanics literature have indicated that the decomposition into deviatoric and dilatational components is important for accurate simulations, especially in the nearly incompressible or fully incompressible limits. Although these changes are standard in solid mechanics, previous hyperelastic extensions of the IB method have not used them \cite{Boffi2008, Griffith2017}. To date, the effect of the formulation of the elastic stress on the accuracy of the IB method does not appear to have been systematically studied. This type of study is the main contribution of the present work. \\
\indent The effect of the volumetric stabilization is controlled by a parameter that we refer to as the numerical bulk modulus $\kappas$ that is linked to a numerical Poisson ratio $\nus$ through a standard linear elasticity relationship (\ref{kappa-nu}). Setting $\kappas = 0$ (corresponding to $\nus = -1$) describes the situation of no volumetric penalization. As described herein, the IB method is exactly incompressible in the continuum limit, and the \textit{physical} Poisson ratio of the immersed structure is automatically $\nu = \frac{1}{2}$. Thus, including a volumetric energy term makes no difference in the continuous equations. Upon discretization, however, incompressibility is exactly maintained only for the Eulerian velocity field in the present methodology. Nonetheless, even in the discrete case, it is clear in the numerical tests presented herein that the solid does inherit some incompressibility from the discretized Eulerian incompressibility constraint. Specifically, in the benchmark cases considered herein, the unstabilized formulation does appear to converge under grid refinement. In these cases, however, the stabilized method is generally much more accurate than the unstabilized method on a given Eulerian grid and Lagranginan mesh. If the solid stress is not fully deviatoric, we show that poor results can be expected if no volumetric penalization is added. The numerical Poisson ratio thereby acts as a stabilization parameter in the discrete IB equations, rather than a material parameter. Further, as is desirable for a stabilization parameter, the effect of the numerical Poisson ratio vanishes in the case that the structural deformations are exactly incompressible. \\
\indent This work documents the practical impact of the proposed stabilization method within the context of four common benchmark problems drawn from structural mechanics literature \cite{Reese1999,RDCook1974,Wriggers2016,Bonet2015}, a fully dynamic FSI benchmark, and a detailed three-dimensional model of esophageal transport. Two-dimensional plane-strain problems, the compressed block and Cook's membrane, are used to perform an in depth study of the proposed volumetric stabilization. For these problems we are able to consider many element types and choices of $\nus$. We report the components of $\cauchys$ for these two tests and the pressure fields for the compressed block, and we demonstrate that these quantities converge to the same solution offered by a high resolution FE method. We also perform extensive tests for the elastic band benchmark to investigate whether the results for quasi-static cases extend to a fluid-driven FSI test. We also use three-dimensional benchmark problems, specifically an anisotropic extension to Cook's membrane and a torsion test. For these two three-dimensional tests, we report the principal strains of the deformation and show that these quantities also converge to principal strains yielded by the FE approach. A consistent finding in all of the tests considered herein is that appropriate levels of volumetric stabilization improve volume conservation and reduce spurious, nonphysical deformations. More specifically, the results suggest that a numerical Poisson ratio of $\nus = 0.4$ yields results with accuracy that is comparable to a fully incompressible FE approach to large deformation nonlinear elasticity across a broad range of benchmarks and grid spacings. Although we stop short of completing the argument that $\nus = 0.4$ is an optimal value to determine the stabilization parameter $\kappas$, it does strike a balance between maintaining solid incompressibility and enabling accurate solutions in the displacement field for the multiple cases presented here. Investigation of a scaling law for more precisely determining $\kappas$ is left as future work.\\
\indent As $\nus \rightarrow \frac{1}{2}$, our method suffers from volumetric locking, as in many other low order displacement-based FE formulations. This phenomenon is demonstrated in our results for the compressed block and Cook's membrane benchmarks: if $\nus = .49995$, the linear and bilinear elements do not converge to the benchmark solution for the number of solid degrees of freedom considered. Reducing the numerical Poisson ratio is demonstrated to be a simple fix for these tests, and convergence behavior is improved for all element types considered, including second order elements. Additionally, this change yields only a small change in total volume conservation. As indicated by Figures (\ref{cb}), (\ref{cm}), (\ref{aniso}), and (\ref{torsion}), the reduction in pointwise volume conservation is also small. This change in numerical Poisson ratio is justified because $\nus$ is a numerical parameter, not a physical parameter. Specifically, and unlike the typical situation in nearly incompressible elasticity, changing the numerical Poisson ratio does not affect the limiting material model obtained under grid refinement. \\
\indent The effect of anisotropy is explored through study of the anisotropic Cook's membrane benchmark and model of esophageal transport. Broadly, the same conclusions which hold true for the other benchmarks also hold for these simulations. That is, unmodified isotropic invariants paired with zero volumetric penalization yield unphysical deformations. Introducing anisotropy opens the door for many more possible forms for the elastic energy functional such as modifying the anisotropic invariants or including the anisotropic invariants in the volumetric term (only for nearly incompressible materials). Though we argue that the use of $\bar{I}_4$ and $\bar{I}_5$ is not appropriate for the cases considered here \cite{Sansour2008}, the study of their use in an IB framework may still prove interesting in future work.\\
\indent Of the benchmark problems considered, the effect of the proposed changes is most pronounced for the torsion test. Introducing the volumetric term leads to extremely noticeable improvement in numerical results, especially when paired with modified invariants. More specifically, the volume conservation can be greatly improved by using this formulation, and the displacement of the point under study converges much more quickly. Even for relatively fine discretizations of this benchmark, our method demonstrates improvements in volume conservation of up to $59 \%$ when compared to the unstabilized case with unmodified invariants. \\
\indent Finally, we explored the effect of the stabilization method in the dynamic FSI regime with the elastic band benchmark and the esophageal transport model. Results from both tests are in agreement with the conclusions from the quasi-static tests. Specifically, volumetric stabilization and use of modified invariants greatly qualitatively improve the deformations for FSI problems and reduce deviations of $J$ away from 1. The esophageal test in particular demonstrates how the lack of stabilization may lead to very low quality results; see Figure (\ref{fig:free_end_D}), in which can elements are inverted, and computation can even become infeasible. \\
\indent In closing, we emphasize that using elastic energy functionals with unmodified invariants and zero volumetric penalization, given by equations (\ref{bad_energy}) and (\ref{bad_energy_aniso}), can yield unphysical deformations in our IB computations. Further, such choices of energy functional produces noticeably worse volume conservation. Even further, modified invariants generally perform better than unmodified invariants in all our results. Including a volumetric penalization term also generally produces more accurate results, with the exception of the anisotropic Cook's membrane in which the differences in performance are minimal. The present method is easy to implement and requires only a change in the way that the elastic stresses are evaluated. Further, in our experience, it has no effect on the stability of the numerical method. Thus, it results in a potentially large improvement in accuracy of computed displacements and volume conservation at negligible computational cost. Given these desirable features, this study strongly indicates using volumetric energy-based stabilization with modified invariants as the default structural formulation for IB-type methods in which the solid elastic response is hyperelastic. For other material models, the tests reported herein suggest using the stabilization technique proposed here with the deviatoric projection formulation.

%%%%%%%%%%%%%%%%%%%%%%%

\section*{Acknowledgement}

We thank Charles Puelz for a careful reading of an initial version of this manuscript and for suggestions that improved the manuscript.
We gratefully acknowledge research support through NIH Awards HL117063 and HL143336 and NSF Awards OAC 1450327 and OAC 1652541.
Computations were primarily performed using facilities provided by University of North Carolina at Chapel Hill through the Research Computing division of UNC Information Technology Services. This work also used the Extreme Science and Engineering Discovery Environment (XSEDE) resource Comet \cite{xsede} at the San Diego Supercomputer Center (SDSC) through allocation TG-ASC170023, which is supported by NSF Award OAC 1548562.

%%%%%%%%%%%%%%%%%%%%%%%%%%%%%%%%%%%%%%%%%%%%%%%%%%%%%%%%%%%%%%%%%%%%%%%%%%%%%%%%%%%%%%%%%%%%%%%%%%%%%%%%%%%%%%%%%%%%%%%%%%%

\bibliography{benchmark_paper.bib}

\end{document}